\documentclass[11pt]{amsart}
\usepackage{amscd,amssymb}
\usepackage{amsthm,amsmath,amssymb}
\usepackage[matrix,arrow]{xy}
\usepackage{fancyhdr}
\makeatletter\@addtoreset{equation}{section}
\makeatletter\@addtoreset{section}{part}

\makeatother

\newtheorem{theorem}[equation]{Theorem}
\newtheorem{theorem-definition}[equation]{Theorem-Definition}
\newtheorem{proposition}[equation]{Proposition}
\newtheorem{lemma}[equation]{Lemma}
\newtheorem{corollary}[equation]{Corollary}

\newtheorem{K3proposition}[equation]{K3-Proposition}
\theoremstyle{definition}
\newtheorem{definition}[equation]{Definition}

\newtheorem{example}{Example}

\theoremstyle{remark}
\newtheorem{remark}[equation]{Remark}

\newcommand{\recip}[1]{\frac{1}{#1}}

\newsavebox{\saveno}\def\no#1;{\sbox{\saveno}{#1}}
\newsavebox{\savereln}\def\reln#1;{\sbox{\savereln}{$\scriptstyle{#1}$}}
\newsavebox{\savegens}\def\gens#1,#2,#3,#4,#5;{\sbox{\savegens}{#1,#2,#3,#4,#5}}
\newsavebox{\savedeg}\def\degree#1;{\sbox{\savedeg}{#1}}

\newcounter{lines}

\newsavebox{\saveinv}\def\inv#1;{\sbox{\saveinv}{$#1$}\setcounter{lines}0}
\newsavebox{\saveell}\def\ell#1;{\sbox{\saveell}{$#1$}}
\newsavebox{\savekkk}\def\kkk#1;{\sbox{\savekkk}{$#1$}}
\newsavebox{\savetype}\def\type#1;{\sbox{\savetype}{$#1$}\setcounter{lines}1}
\newsavebox{\saveexcl}\def\excl#1;{\sbox{\saveexcl}{$#1$}}
\newsavebox{\saveuntw}\def\untw#1;{\sbox{\saveuntw}{#1}}
\newsavebox{\savetypea}\def\typea#1;{\sbox{\savetypea}{$#1$}\setcounter{lines}2}
\newsavebox{\saveexcla}\def\excla#1;{\sbox{\saveexcla}{$#1$}}
\newsavebox{\saveuntwa}\def\untwa#1;{\sbox{\saveuntwa}{#1}}
\newsavebox{\savetypeb}\def\typeb#1;{\sbox{\savetypeb}{$#1$}\setcounter{lines}3}
\newsavebox{\saveexclb}\def\exclb#1;{\sbox{\saveexclb}{$#1$}}
\newsavebox{\saveuntwb}\def\untwb#1;{\sbox{\saveuntwb}{#1}}
\newsavebox{\savetypec}\def\typec#1;{\sbox{\savetypec}{$#1$}\setcounter{lines}4}
\newsavebox{\saveexclc}\def\exclc#1;{\sbox{\saveexclc}{$#1$}}
\newsavebox{\saveuntwc}\def\untwc#1;{\sbox{\saveuntwc}{#1}}
\newsavebox{\savetyped}\def\typed#1;{\sbox{\savetyped}{$#1$}\setcounter{lines}5}
\newsavebox{\saveexcld}\def\excld#1;{\sbox{\saveexcld}{$#1$}}
\newsavebox{\saveuntwd}\def\untwd#1;{\sbox{\saveuntwd}{#1}}

\newenvironment{entry}%
{\medskip \no\empty; \reln\empty;
\gens\empty,\empty,\empty,\empty,\empty; \degree\empty;
\type\empty;\sbox{\saveexcl}{} \untw\empty; \typea\empty;
\sbox{\saveexcla}{} \untwa\empty; \typeb\empty;
\sbox{\saveexclb}{} \untwb\empty; \typec\empty;
\sbox{\saveexclb}{} \untwc\empty; \inv\empty; \sbox{\saveell}{}
\kkk\empty; }
{$\gimel=$\ \usebox{\saveno}:
$X_{\usebox\savereln}\subset\mathbb{P}(\usebox\savegens)$\nopagebreak
\endgraf

\ifcase\value{lines}\or
\begin{tabular}{|p{6cm}|p{6cm}|}
 \hline
$-K_X^3$=\usebox{\savedeg}& \usebox{\savekkk}\\
\hline \hline
\usebox{\savetype} & \usebox{\saveexcl}  \\
\hline
\end{tabular}\or
\begin{tabular}{|p{6cm}|p{6cm}|}
 \hline
$-K_X^3$=\usebox{\savedeg} &  \usebox{\savekkk} \\
\hline\hline
\usebox{\savetype} & \usebox{\saveexcl}  \\
\hline
\usebox{\savetypea} & \usebox{\saveexcla}  \\
\hline
\end{tabular}\or
\begin{tabular}{|p{6cm}|p{6cm}|}
 \hline
$-K_X^3$=\usebox{\savedeg} & \usebox{\savekkk}\\
\hline\hline
\usebox{\savetype} & \usebox{\saveexcl} \\
\hline
\usebox{\savetypea} & \usebox{\saveexcla} \\
\hline
\usebox{\savetypeb} & \usebox{\saveexclb}  \\
\hline
\end{tabular}\or
\begin{tabular}{|p{6cm}|p{6cm}|}
 \hline
$-K_X^3$=\usebox{\savedeg} & \usebox{\savekkk}\\
\hline\hline
\usebox{\savetype} & \usebox{\saveexcl}  \\
\hline
\usebox{\savetypea} & \usebox{\saveexcla}  \\
\hline
\usebox{\savetypeb} & \usebox{\saveexclb}  \\
\hline
\usebox{\savetypec} & \usebox{\saveexclc}  \\
\hline
\end{tabular}\or
\begin{tabular}{|p{6cm}|p{6cm}|}
 \hline
$-K_X^3$=\usebox{\savedeg}&  \usebox{\savekkk}\\
\hline\hline
\usebox{\savetype} & \usebox{\saveexcl}  \\
\hline
\usebox{\savetypea} & \usebox{\saveexcla}  \\
\hline
\usebox{\savetypeb} & \usebox{\saveexclb} \\
\hline
\usebox{\savetypec} & \usebox{\saveexclc}  \\
\hline
\usebox{\savetyped} & \usebox{\saveexcld}  \\
\hline
\end{tabular}
\fi}

\author{Ivan Cheltsov and Jihun Park}

\title{Halphen Pencils on Weighted Fano Threefold Hypersurfaces}

\begin{document}

\begin{abstract}
We classify all  pencils on a general weighted hypersurface of
degree $\sum_{i=1}^{4}a_{i}$ in
$\mathbb{P}(1,a_{1},a_{2},a_{3},a_{4})$  whose general members are
surfaces of Kodaira dimension zero.
\vspace{5mm}

\noindent
Mathematics Subject Classifications (2000). 14E07, 14D06, 14J28, 14J45, 14J70.
\vspace{5mm}

\noindent Keywords. Birational automorphism, Center of canonical
singularities, Elliptic fibration, Halphen pencil, K3 surface,
Kawamata blow up, Noether-Fano inequality, Quotient singularity,
Weighted projective space, Weighted Fano threefold hypersurface.

 \end{abstract}
\address{\emph{Ivan Cheltsov}\newline \textnormal{School of
Mathematics, The University of Edinburgh,
Edinburgh, EH9 3JZ, UK;
 \texttt{cheltsov@yahoo.com}}}
\address{ \emph{Jihun Park}\newline \textnormal{Department of
Mathematics, POSTECH, Pohang, Kyungbuk 790-784, Korea;
\texttt{wlog@postech.ac.kr}}}

\maketitle\thispagestyle{empty}%
\tableofcontents

\newpage

\part{Introduction.} \label{section:introduction}

Throughout this article, all varieties are projective and  defined
over $\mathbb{C}$ and morphisms are proper, otherwise mentioned.

\section{Introduction.} Let $C$ be a smooth curve in
$\mathbb{P}^{2}$ defined by an cubic homogeneous equation
$f(x,y,z)=0$. Suppose that we have nine distinct points
$P_{1},\cdots ,P_{9}$ on $C$ such that the divisor
$$
\sum_{i=1}^{9}P_{i}-\mathcal{O}_{\mathbb{P}^{2}}\big(3\big)\Big\vert_{C}
$$
is a torsion divisor of order $m\geq 1$ on the curve $C$. Then,
there is a curve $Z\subset\mathbb{P}^{2}$ of degree $3m$  such
that $\mathrm{mult}_{P_{i}}(Z)=m$ for each point $P_{i}$. Let
$\mathcal{P}$ be the pencil given by the equation
$$
\lambda f^{m}\big(x,y,z\big)+\mu g\big(x,y,z\big)=0\subset\mathrm{Proj}\Big(\mathbb{C}[x,y,z]\Big)\cong\mathbb{P}^{2},%
$$
where $g(x,y,z)=0$ is a homogeneous equation of the curve $Z$ and
$(\lambda:\mu)\in\mathbb{P}^{1}$. Then, a general curve of the
pencil $\mathcal{P}$ is birational to an elliptic curve. The
pencil $\mathcal{P}$ is called a plane Halphen pencil
(\cite{Hal82}) and the construction of the pencil $\mathcal{P}$
can be generalized to the case when the curve $C$ has ordinary
double points and the points $P_{1},\cdots ,P_{9}$ are not
necessarily distinct (\cite{Do66}). In fact, every plane elliptic
pencil is birational to a Halphen pencil. Namely, the following
result is proved in \cite{Ber77} using the technique of plane
Cremona transformations but its rigorous proof is due to
\cite{Do66}.

\begin{theorem}
\label{theorem:Bertini-Dolgachev} Let $\mathcal{M}$ be a pencil on
$\mathbb{P}^{2}$ whose general curve is birational to an elliptic
curve. Then, there is a birational automorphism $\rho$ of
$\mathbb{P}^{2}$ such that $\rho(\mathcal{M})$ is a plane Halphen
pencil.
\end{theorem}

\begin{proof}
Every birational automorphism of $\mathbb{P}^{2}$ is a composition
of projective automorphisms and Cremona involutions (Section~2.5
in \cite{CoKoSm03}). Moreover, the arguments of Section~2.5
in \cite{CoKoSm03} together with the
two-dimensional analogue of Theorem~\ref{theorem:Noether-Fano}
imply that there is birational automorphism
$\rho\in\mathrm{Bir}(\mathbb{P}^{2})$ such that the singularities
of the log pair $(\mathbb{P}^{2}, \frac{3}{n}\mathcal{B})$ are
canonical, where $\mathcal{B}=\rho(\mathcal{M})$ and $n$ is the
natural number such that
$\mathcal{B}\sim_{\mathbb{Q}}\mathcal{O}_{\mathbb{P}^{2}}(n)$.

The singularities of the log pair $(\mathbb{P}^{2},
\frac{3}{n}\mathcal{B})$ are not terminal by
Theorem~\ref{theorem:Noether-Fano}.\footnote{Theorem~\ref{theorem:Noether-Fano}
can be generalized to each dimension $\geq 2$.} Therefore, there
is a birational morphism $\pi:S\to\mathbb{P}^{2}$ such that the
singularities of the log pair $(S, \frac{3}{n}\mathcal{B}_S)$ are
terminal and
$$
K_{S}+\frac{3}{n}\mathcal{B}_S\sim_{\mathbb{Q}}\pi^{*}\Big(K_{\mathbb{P}^{2}}+\frac{3}{n}\mathcal{B}\Big),%
$$
where $\mathcal{B}_S$ is the proper transform of the pencil
$\mathcal{B}$ by the birational morphism $\pi$. Then, $-K_{S}$ is
nef because the pencil $\mathcal{B}_S$ does not have any fixed
curve. The divisor $-K_{S}$ is not big by
Theorem~\ref{theorem:Noether-Fano}, which implies that
$K_{S}^{2}=0$. In particular, the complete linear system
$|-nK_{S}|$ does not have fixed points. Therefore,  the number $n$
is divisible by $3$ and
$$
\rho\big(\mathcal{M}\big)=\pi\Big(\big|-\frac{n}{3}K_{S}\big|\Big),
$$
which implies that $\rho(\mathcal{M})$ is a plane Halphen pencil.
\end{proof}

A problem similar to Theorem~\ref{theorem:Bertini-Dolgachev} can
be considered  for  Fano varieties  whose groups of birational
automorphisms are well understood\footnote{Elliptic pencils on
some del Pezzo surfaces defined over an algebraically non-closed
perfect field are birationally classified in \cite{Ch00a} using
the structure of their groups of birational automorphisms
(\cite{Ma66}).}. In particular, it would be very interesting to
classify pencils of K3 surfaces on three-dimensional weighted Fano
hypersurfaces.

\begin{definition}
\label{definition:Halphen-pencil}\index{Halphen pencil} A Halphen
pencil is a one-dimensional linear system whose general element is
birational to a smooth variety of Kodaira dimension zero.
\end{definition}

Let $X$ be a general quasismooth hypersurface of degree
$d=\sum_{i=1}^{4}a_{i}$ in $\mathbb{P}(1,a_{1},a_{2},a_{3},a_{4})$
that has terminal singularities, where $a_{1}\leq a_{2}\leq
a_{3}\leq a_{4}$. Then,
$$
-K_{X}\sim_{\mathbb{Q}}\mathcal{O}_{\mathbb{P}(1,\,a_{1},\,a_{2},\,a_{3},\,a_{4})}(1)\Big\vert_{X},
$$
which implies that $X$ is a Fano threefold. The divisor class
group $\mathrm{Cl}(X)$ is generated by the anticanonical divisor
$-K_{X}$ and there are exactly $95$ possibilities for the
quadruple $(a_{1},a_{2},a_{3},a_{4})$, which are found in
\cite{IF00}. We use the notation $\gimel$ for the  entry numbers
of these famous $95$ families. They are ordered in the same way of
\cite{IF00}, which is a standard way nowadays. We tabulate these
families together with their properties in
Part~\ref{section:weighted-Fanos}.

Birational geometry on such threefolds  is extensively studied in
\cite{Ch05}, \cite{ChPa05}, \cite{CPR}, \cite{Ry06}, and so forth.
The article \cite{CPR} describes the generators of the group
$\mathrm{Bir}(X)$ of birational automorphisms of $X$. Also, the
article \cite{ChPa05} shows the relations among these generators.
The former article proves the following result as well.

\begin{theorem}
\label{theorem:CPR} The threefold $X$ cannot be rationally fibred
by rational curves or surfaces.
\end{theorem}

 As for birational maps into elliptic fibrations, the hypersurface $X$ in each family except the families  of $\gimel =3$, $60$, $75$,
 $84$, $87$, and $93$ is birational to an elliptic fibration  (\cite{ChPa05}). Furthermore, all birational
transformations of the threefold $X$ into elliptic fibrations are
classified in \cite{Ch00a}, \cite{Ch03a}, \cite{Ch05},
\cite{ChPa05},
 and \cite{Ry06}.

It is known that Halphen pencils on the threefold $X$ always
exist, to be precise, the threefold $X$ can be always rationally
fibred by K3 surfaces (\cite{ChPa05}). In this article, we are to
classify all  Halphen pencils on the hypersurfaces in the $95$
families as what is done for their elliptic fibrations in
\cite{Ch05}.

Let us explain five examples of pencils on the threefold $X$. They
exhaust all the possible Halphen pencils on $X$. It follows from
\cite{CPR} that the pencils constructed below are
$\mathrm{Bir}(X)$-invariant
(Proposition~\ref{proposition:Invariance}). We will show,
throughout this article, that they  are indeed Halphen pencils.

\begin{example}
\label{example:Halphen-pencil-I} Suppose that $a_{2}=1$. Then,
every one-dimensional linear system in $|-K_{X}|$ is a Halphen
pencil. It follows from the Adjunction  that a general
surface in $|-K_X|$ is birational to a smooth K3 surface. In particular,
it belongs to Reid's $95$ codimension $1$ weighted K3 surfaces
(\cite{R79}).
\end{example}
Therefore, in the cases in Example~\ref{example:Halphen-pencil-I},
or equivalently $\gimel=1$, $2$, $3$, $4$, $5$, $6$, $8$, $10$,
$14$, there are infinitely many Halphen pencils on the
hypersurface $X$. Such cases will be studied in
Part~\ref{section:infinite-Halphens}, where we will prove that
every Halphen pencil is contained in $|-K_X|$.

\begin{example}
\label{example:Halphen-pencil-II} Suppose that $a_{1}\ne a_{2}$.
Then, the linear system $|-a_{1}K_{X}|$ is a pencil. If $a_1=1$,
then the linear system $|-K_X|$ is a Halphen pencil and its
general surface belongs to Reid's $95$ codimension $1$ weighted K3
surfaces as in Example~\ref{example:Halphen-pencil-I}. In fact, it
is a Halphen pencil if only if $a_1\ne a_2$ (\cite{ChPa05}). We
will see that it is a unique Halphen pencil except the cases with
$a_2=1$ and the cases in three Examples below.
\end{example}

  We will discuss
the cases with a unique Halphen pencil in
Parts~\ref{section:one-Halphen-1} and \ref{section:one-Halphen-2}.

Note that $a_{1}=a_{2}\ne 1$ exactly when $\gimel=18$, $22$, and
$28$.

\begin{example}
\label{example:Halphen-pencil-III} Suppose that $\gimel=18$, $22$,
or $28$. In such cases, $a_{1}=a_{2}\ne 1$ and $a_{3}=a_{1}+1$.
The threefold $X$ has singular points $O_{1}, \cdots , O_{r}$ of
type $\frac{1}{a_{1}}(1,1,a_{1}-1)$,  where
$r=\frac{3a_{1}+a_{4}+1}{a_{1}}$. There is a unique index $j\geq
3$ such that $a_{2}+a_{3}+a_{4}=ma_j$, where $m$ is a natural
number. In particular, the threefold $X$ is given by an equation
$$
\sum_{k=0}^{m}
x_j^{k}f_{k}\big(x_0,x_1,x_2,x_3,x_4\big)=0\subset\mathrm{Proj}\Big(\mathbb{C}[x_0,x_1,x_2,x_3,x_4]\Big),
$$
where $\mathrm{wt}(x_0)=1$, $\mathrm{wt}(x_l)=a_{l}$, and $f_{k}$
is a qua\-si\-ho\-mo\-ge\-ne\-ous polynomial of degree
$a_1+a_2+a_3+a_4-ka_j$ that is independent of the variable
$x_{j}$. Let $\mathcal{P}_{i}$ be the pencil of surfaces in
$|-a_{1}K_{X}|$ that pass through the point $O_{i}$ and
$\mathcal{P}$ be the pencil on the threefold $X$ that is cut out
by the pencil $\lambda x_0^{a_{1}}+\mu
f_{m}(x_0,x_1,x_2,x_3,x_4)=0$, where
$(\lambda:\mu)\in\mathbb{P}^{1}$. It will be proved that
$\mathcal{P}$ and $\mathcal{P}_{i}$ are Halphen pencils in
$|-a_{1}K_{X}|$.
\end{example}
The cases in Example~\ref{example:Halphen-pencil-III} are the only
cases that have  more than two finitely many but Halphen pencils.
These cases will be discussed in
Part~\ref{section:more-than-two-Halphens}.

\begin{example}
\label{example:Halphen-pencil-IV} Suppose that $\gimel=45$, $48$,
$55$, $57$, $58$, $66$, $69$, $74$, $76$, $79$, $80$, $81$, $84$,
$86$, $91$, $93$, or $95$. We then see $1\ne a_{1}\ne a_{2}$.
Moreover, there is a unique index $j\ne 2$ such that
$a_{1}+a_{3}+a_{4}=ma_j$, where $m$ is a natural number.
Therefore, the threefold $X$ is given by an equation
$$
\sum_{k=0}^m x^{k}_j
f_{k}\big(x_0,x_1,x_2,x_3,x_4\big)=0\subset\mathrm{Proj}\Big(\mathbb{C}[x_0,x_1,x_2,x_3,x_4]\Big),
$$
where $\mathrm{wt}(x_0)=1$, $\mathrm{wt}(x_i)=a_{i}$, and $f_{k}$
is a qua\-si\-ho\-mo\-ge\-ne\-ous polynomial of degree
$a_{1}+a_{2}+a_{3}+a_{4}-ka_{j}$ that is independent of the
variable $x_{j}$. Let $\mathcal{P}$ be the pencil on the threefold
$X$ that is cut out by the pencil $\lambda x_0^{a_{2}}+\mu
f_{m}(x_0,x_1,x_2,x_3,x_4)=0$, where
$(\lambda:\mu)\in\mathbb{P}^{1}$. It will be shown that
$\mathcal{P}$ is a Halphen pencil in $|-a_{2}K_{X}|$.
\end{example}

\begin{example}
\label{example:Halphen-pencil-V} Suppose that $\gimel=60$. Then,
$X$ is a general hypersurface of degree $24$ in
$\mathbb{P}(1,4,5,6,9)$. Hence, the threefold $X$ is given by an
equation
$$
w^{2}f_{6}\big(x,y,z,t\big)+wf_{15}\big(x,y,z,t\big)+f_{24}\big(x,y,z,t\big)=0\subset\mathrm{Proj}\Big(\mathbb{C}[x,y,z,t,w]\Big),
$$
where $\mathrm{wt}(x)=1$, $\mathrm{wt}(y)=4$, $\mathrm{wt}(z)=5$,
$\mathrm{wt}(t)=6$, $\mathrm{wt}(w)=9$, and $f_{k}(x,y,z,t)$ is a
general quasihomogeneous polynomial of degree $k$. Then, the
linear system on the threefold $X$ cut out by the pencil $\lambda
x^{6}+\mu f_{6}\big(x,y,z,t\big)=0$, where
$(\lambda:\mu)\in\mathbb{P}^{1}$, is a Halphen pencil in
$|-6K_{X}|$.
\end{example}

The cases in Examples~\ref{example:Halphen-pencil-IV} and
\ref{example:Halphen-pencil-V} have at least two Halphen pencils
because they also satisfy the condition for
Example~\ref{example:Halphen-pencil-II}. Furthermore, we will see
that  these are the only Halphen pencils on the hypersurface $X$
of each family. These cases are discussed in
Part~\ref{section:two-Halphens}.

 The main purpose of this
article is to prove the
following:\footnote{Theorem~\ref{theorem:main} is proved in
\cite{Ry05} and \cite{Ry06}  for the cases $\gimel=5$, $34$, $75$,
$88$, and $90$.}.

\begin{theorem}
\label{theorem:main} Let $X$ be a general hypersurface in the $95$
families. Then, the pencils constructed in
Examples~\ref{example:Halphen-pencil-I},
\ref{example:Halphen-pencil-II}, \ref{example:Halphen-pencil-III},
\ref{example:Halphen-pencil-IV}, and
\ref{example:Halphen-pencil-V} exhaust all possibilities for
Halphen pencils on the hypersurface $X$.
\end{theorem}

The following are immediate consequence of
Theorem~\ref{theorem:main}.
\begin{corollary} \label{corollary:main} Let $X$ be a general
hypersurface in the $95$ families with entry number $\gimel$.
\begin{enumerate}
\item There are finitely many Halphen pencils on the threefold $X$ if and only if $a_{2}\ne 1$.%
\item There are at most two Halphen pencils on $X$ in the case when  $a_{1}\ne a_{2}$.%
\item Every Halphen pencil on the threefold $X$ is contained in $|-K_{X}|$ if $a_{1}=1$.%
\item The linear system $|-K_{X}|$ is the only Halphen pencil on $X$ if $a_{1}=1$ and $a_{2}\ne 1$.%
\item The linear system $|-a_{1}K_{X}|$ is the only Halphen pencil
on the threefold $X$ if and only if $\gimel\ne 1$, $2$, $3$, $4$,
$5$, $6$, $8$, $10$, $14$, $18$, $22$, $28$, $45$, $48$, $55$,
$57$, $58$, $60$, $66$, $69$, $74$, $76$, $79$,
$80$, $81$, $84$, $86$, $91$, $93$, $95$.%
\end{enumerate}
\end{corollary}

Furthermore, Theorem~\ref{theorem:main} with
Proposition~\ref{proposition:Invariance} forces us to conclude
\begin{corollary} \label{corollary:main2} Let $X$ be a general
hypersurface in the $95$ families.
Then, every Halphen pencil on the threefold $X$ is invariant under the action of $\mathrm{Bir}(X)$.%
\end{corollary}

The proof of Theorem~\ref{theorem:main} is based on
Theorems~\ref{theorem:Noether-Fano}, \ref{theorem:main-tool} and
Lemmas~\ref{lemma:Kawamata}, \ref{lemma:Cheltsov-Kawamata}. We
prove the theorem case by case in order of the number of Halphen
pencils and the entry number $\gimel$.

In addition, we prove that general surfaces of the pencils
constructed in Examples~\ref{example:Halphen-pencil-I},
\ref{example:Halphen-pencil-II}, \ref{example:Halphen-pencil-III},
\ref{example:Halphen-pencil-IV}, \ref{example:Halphen-pencil-V}
are birational to  smooth K3 surfaces.

\begin{theorem}
\label{theorem:K3} Let $X$ be a general hypersurface in the $95$
families. Then, a general surface of every Halphen pencil on $X$ is birational to a smooth K3 surface.%
\end{theorem}

 It follows from
Proposition~\ref{proposition:Halphen-K3} that general surfaces of
the pencils constructed in Examples~\ref{example:Halphen-pencil-I}
and \ref{example:Halphen-pencil-II} are birational to smooth K3
surfaces. General surfaces of the other pencils will be discussed
in Part~\ref{section:two-Halphens}. The proof of
Theorem~\ref{theorem:K3} is based on
Corollaries~\ref{corollary:Halphen-K3} and
\ref{corollary:Halphen-K3-rational-curve}.

Theorems~\ref{theorem:main} and \ref{theorem:K3} tell us how a
general hypersurface in the $95$ families can be rationally fibred
by smooth surfaces of Kodaira dimension zero.

\begin{corollary} \label{corollary:fibration} Let $X$ be a general
hypersurface in the $95$ families and let $\pi:Y\to Z$ be a
morphism whose general fiber is birational to a smooth surface of
Kodaira dimension zero. If there is a birational map
$\alpha:X\dasharrow Y$, then there is an isomorphism $\phi:
\mathbb{P}^1\to Z$ such that  the following diagram commutes:

$$
\xymatrix{
&X\ar@{-->}[d]_{\psi}\ar@{-->}[rr]^{\alpha}&&Y\ar@{->}[d]^{\pi}&\\
&\mathbb{P}^1\ar@{->}[rr]_{\phi}&& Z,&}
$$
where the rational map $\psi:X\dasharrow \mathbb{P}^1 $ is induced
by one of the pencils in Examples~\ref{example:Halphen-pencil-I},
\ref{example:Halphen-pencil-II}, \ref{example:Halphen-pencil-III},
\ref{example:Halphen-pencil-IV}, and
\ref{example:Halphen-pencil-V}. In particular, a general fiber of
the morphism $\pi$ is birational to a smooth K3 surface.
\end{corollary}


\vspace{5mm} {\itshape Acknowledgments}. The authors would like to
thank I.\,Aliev, A.\,Co\-rti, M.\,Gri\-nen\-ko,
V.\,Is\-kov\-skikh, Yu.\,Pro\-kho\-rov, and V.\,Sho\-ku\-rov for
useful conversations. This work was initiated when the second
author visited University of Edinburgh in February, 2006 and they
could finish the article while the first author visited POSTECH in
Korea. The authors would like to thank POSTECH and University of
Edinburgh
 for their hospitality. The first author has been supported
by CRDF grant RUM1-2692MO-05 and the second author supported by
KRF Grant 2005-070-C00005. \vspace{5mm}

\section{Preliminaries.} \label{section:preliminaries}

Let $X$ be a threefold with $\mathbb{Q}$-factorial singularities
and $\mathcal{M}$ be a linear system on the threefold $X$ without
fixed components. We consider the log pair $(X, \mu\mathcal{M})$
for some nonnegative rational number $\mu$.

Let $\alpha:Y\to X$ be a proper birational morphism such that $Y$
is smooth and  the proper transform $\mathcal{M}_Y$ of the linear
system $\mathcal{M}$ by the birational morphism $\alpha$ is
base-point-free. Then, the rational equivalence
$$
K_{Y}+\mu\mathcal{M}_Y\sim_{\mathbb{Q}}\alpha^{*}\Big(K_{X}+\mu\mathcal{M}\Big)+\sum_{i=1}^{k}a_{i}E_{i}%
$$
holds, where $E_{i}$ is an exceptional divisor of the birational
morphism $\alpha$ and $a_{i}$ is a rational number.
\begin{definition}
\label{definition:canonical-terminal}\index{singularities} The
singularities of the log pair $(X, \mu\mathcal{M})$ are terminal
(canonical, log-terminal, respectively) if each rational number
$a_{i}$ is positive (nonnegative, greater than $-1$,
respectively). In case, we also say that the log pair $(X,
\mu\mathcal{M})$ is terminal (canonical, log-terminal,
respectively).
\end{definition}
It is convenient to specify where the log pair $(X,
\mu\mathcal{M})$ is not terminal.
\begin{definition}
\label{definition:center-of-canonical-singularities}\index{center
of canonical singularities} A proper irreducible subvariety
$Z\subset X$ is called a center of canonical singularities of the
log pair $(X, \mu\mathcal{M})$ if there is an exceptional divisor
$E_i$ such that $\alpha(E_{i})=Z$ and $a_{i}\leq 0$. The set of
all proper irreducible subvarieties of $X$ that are centers of
canonical singularities of the log pair $(X, \mu\mathcal{M})$ is
denoted by $\mathbb{CS}(X, \mu\mathcal{M})$.
\end{definition}

A curve not contained in the singular locus of the threefold $X$
is a center of canonical singularities of the log pair $(X,
\mu\mathcal{M})$  if and only if the multiplicity of a  general
surface of $\mathcal{M}$ along the curve is not smaller than
$\frac{1}{\mu}$. Furthermore, we obtain

\begin{lemma}
\label{lemma:curves} Let $C$ be a curve on the threefold $X$ that
is  not contained in the singular locus of $X$.  Suppose that the
curve $C$ is a center of canonical singularities of the log pair
$(X, \mu\mathcal{M})$ and the linear system $|-mK_{X}|$ is
base-point-free for some natural number $m>0$. If
$-K_{X}\sim_{\mathbb{Q}}\mu\mathcal{M}$, then $-K_{X}\cdot C\leq
-K_{X}^{3}$.
\end{lemma}

\begin{proof}
Let $M_{1}$ and $M_{2}$ be general surfaces in $\mathcal{M}$.
Then,
$$
\mathrm{mult}_{C}(M_{1}\cdot M_{2})\geq\mathrm{mult}_{C}(M_{1})\mathrm{mult}_{C}(M_{2})\geq \frac{1}{\mu^{2}}.%
$$
Let $H$ be a general surface in $|-mK_{X}|$. Then,
$$
\frac{m}{\mu^{2}}(-K_X^3)=H\cdot M_{1}\cdot M_{2}\geq
(-mK_{X}\cdot C)\mathrm{mult}_{C}(M_{1}\cdot
M_{2})\geq\frac{m}{\mu^{2}}(-K_{X}\cdot C),
$$
which implies $-K_{X}\cdot C\leq -K_{X}^{3}$.
\end{proof}

The following result is a generalization of so-called
Noether--Fano inequality (\cite{CoKoSm03}).

\begin{theorem}
\label{theorem:Noether-Fano}\index{Noether-Fano inequality}
Suppose that the linear system $\mathcal{M}$ is a pencil whose
general surface is birational to a smooth surface of Kodaira
dimension zero, the linear system $|-mK_{X}|$ is base-point-free
for some natural $m$, and $\mathcal{M}\sim_{\mathbb{Q}}-\mu K_X$.
If the linear system $|-mK_{X}|$ induces either a birational
morphism or an elliptic fibration, then the log pair $(X,
\mu\mathcal{M})$ is not terminal.
\end{theorem}

\begin{proof}
Let $M$ be a general surface in $\mathcal{M}$.  Suppose that the
log pair $(X, \mu M)$ is terminal. Then, for some positive
rational number $\epsilon>\mu$, the log pair $(X, \epsilon M)$ is
also ter\-mi\-nal and the divisor $K_X+\epsilon M$ is nef. We have
a resolution of indeterminacy of the rational map
$\rho:X\dasharrow \mathbb{P}^1$ induced by the pencil
$\mathcal{M}$ as follows:
$$
\xymatrix{
&&Y\ar@{->}[ld]_{\alpha}\ar@{->}[rd]^{\beta}&&\\%
&X\ar@{-->}[rr]_{\rho}&&\mathbb{P}^{1},&}
$$ %
where $Y$ is smooth and $\beta$ is a morphism. We consider the
linear equivalence
$$
K_{Y}+\epsilon M_Y\sim_{\mathbb{Q}} \alpha^{*}\Big(K_{X}+\epsilon M\Big)+\sum_{i=1}^{k}c_{i}E_{i},%
$$
where $M_Y$ is the proper transform of the surface $M$ and $c_{i}$
is a rational number. Then, each $c_i$ is positive. Also, we may
assume that the proper transform $\mathcal{M}_Y$ of the pencil
$\mathcal{M}$ by the birational morphism $\alpha$ is
base-point-free. In particular, the surface $M_Y$ is smooth.

Let $l$ be a sufficiently big and divisible natural number. Then,
the negativity property of the exceptional locus of a birational
morphism (Section~1.1 in \cite{Sho93}) implies that the linear
system $|l(K_{Y}+\epsilon M_Y)|$ gives a dominant rational map
$\xi:Y\dasharrow V$  with $\mathrm{dim}(V)\geq 2$. One the other
hand, since the proper transform  $\mathcal{M}_Y$ is a
base-point-free pencil, the Adjunction formula implies
$$
l\Big(K_{Y}+\epsilon M_Y\Big)\Big\vert_{M_Y}\sim lK_{M_Y}.
$$
However, the surface $M_Y$ has Kodaira dimension zero, which
implies that $\mathrm{dim}(V)\leq 1$. It is a contradiction.
\end{proof}

\begin{theorem-definition}\index{Kawamata blow up}\label{Kawamata blow up}
Let $(P\in U)$ be the germ of a threefold terminal quotient singularity $P$
of type $\frac{1}{r}(1, a, r-a)$, where $r\geq 2$, $r>a$, and $a$
is coprime to $r$. Suppose that $$f : (E\subset
W)\to (Z\subset U)$$ is a proper birational morphism such that
\begin{itemize}
\item the threefold $W$ has at worst terminal singularities;
\item the exceptional set $E$ of $f$ is an irreducible divisor of $W$;
\item the divisor $-K_W$ is $f$-ample;
\item the point $P$ is contained in the subvariety $Z$.
\end{itemize}
Then, it is the weighted blow up at the
point $P$ with weights $(1, a, r-a)$. In particular, $Z=P$. We
call such a birational morphism the Kawamata blow up at the point
$P$ with weights $(1, a, r-a)$ or simply the Kawamata blow up at
the point $P$.
\end{theorem-definition}
\begin{proof}
See \cite{Ka96}.
\end{proof}
Let $\pi: Y\to X$ be the Kawamata blow up at a quotient singular
point $P$ of type $\frac{1}{r}(1, a, r-a)$, where $r\geq 2$,
$r>a$, and $a$ is coprime to $r$. One can easily check that the
exceptional divisor $E$ of the birational morphism $\pi$ is
isomorphic to $\mathbb{P}(1, a, r-a)$. Furthermore, we see

\begin{equation*}
\left\{\aligned
&K_Y=\pi^*(K_X)+\frac{1}{r}E,\\
&K_Y^3=K_X^3+\frac{1}{ra(r-a)},\\
& E^3=\frac{r^2}{a(r-a)}.\\
\endaligned\right.
\end{equation*}

Otherwise mentioning, from this point throughout this section, we
always assume that the linear system $\mathcal{M}$ is a pencil
with $\mathcal{M}\sim_{\mathbb{Q}}-\mu K_X$. In addition, we
always assume that  a general surface of the pencil $\mathcal{M}$
is irreducible.

For our purpose Theorem-Definition~\ref{Kawamata blow up} can be modified as follows.

\begin{lemma}
\label{lemma:Kawamata} Let $P$ be a singular point  of a threefold
$X$ that is a quotient singularity of type $\frac{1}{r}(1,a,r-a)$,
$r\geq 2$, $r>a$, and $a$ is coprime to $r$. Suppose that the log
pair $(X, \mu\mathcal{M})$ is canonical but the set $\mathbb{CS}(X, \mu\mathcal{M})$ contains either  the
point $P$ or a curve passing through the point $P$. Let $\pi:Y\to X$ be the Kawamata blow up at the point
$P$ and let $\mathcal{M}_Y$ be the proper transform of
$\mathcal{M}$ by the birational morphism $\pi$. Then,
$$
\mu\mathcal{M}_Y\sim_{\mathbb{Q}}
\pi^{*}\Big(-K_{X}\Big)-\frac{1}{r}E\sim_{\mathbb{Q}} -K_{Y},
$$
where $E$ is the ex\-cep\-ti\-onal divisor of $\pi$.
\end{lemma}

\begin{proof}
We consider only the case $r=2$. Then, $E\cong\mathbb{P}^{2}$ and
$E\vert_{E}\sim_{\mathbb{Q}}\mathcal{O}_{\mathbb{P}^{2}}(-2)$. We
have
$$
\mathcal{M}_Y\sim_{\mathbb{Q}}\pi^{*}\Big(-\frac{1}{\mu}K_{X}\Big)-mE,
$$
where $m$ is a positive rational number. In particular, we have
$\mathcal{M}_Y\vert_{E}\equiv-mE\vert_{E}$.

Suppose that $m<\frac{1}{2\mu}$. Let $Q$ be a point of $E$.
Intersecting a general divisor of the pencil $\mathcal{M}_Y$ with
a general line on $E$ that passes through the point $Q$, we obtain
 the inequality $$\mathrm{mult}_{Q}(\mathcal{M}_Y)\leq
2m<\frac{1}{\mu}.$$

Suppose that the set $\mathbb{CS}(X, \mu\mathcal{M})$ contains  an
ir\-re\-du\-ci\-ble curve $Z$ that passes through the point $P$.
Then,
$$
\frac{1}{\mu}>\mathrm{mult}_{O}\big(\mathcal{M}_Y\big)\geq\mathrm{mult}_{Z_Y}\big(\mathcal{M}_Y\big)
\geq\mathrm{mult}_{Z}\big(\mathcal{M}\big)\geq\frac{1}{\mu},%
$$
where $Z_Y$ is the proper transform of the curve $Z$ and $O$ is a
point of the intersection of the curve $Z_Y$ and the exceptional
divisor $E$. Therefore, the singularities of the log pair $(X,
\mu\mathcal{M})$ are terminal in a punctured neighborhood of the
point $P$. The equivalence
$$
K_{Y}+\mu\mathcal{M}_Y\sim_{\mathbb{Q}}
\pi^{*}\Big(K_{X}+\mu\mathcal{M}\Big)+\Big(\frac{1}{2}-\mu m\Big)E,%
$$
shows that the set $\mathbb{CS}(Y, \mu\mathcal{M}_Y)$ contains a
proper subvariety $\Delta\subsetneq E$. It implies that the
inequality $\mathrm{mult}_{\Delta}(\mathcal{M}_Y)\geq
\frac{1}{\mu}$ holds, which is a contradiction.
\end{proof}

Lemma~\ref{lemma:Kawamata} can be generalized in the following way
(\cite{Ch05}).

\begin{lemma}
\label{lemma:Cheltsov-Kawamata} Under the assumptions and
notations of Lemma~\ref{lemma:Kawamata}, suppose that
 we have a proper subvariety
$Z\subset E\cong\mathbb{P}(1,a,r-a)$ that belongs to
$\mathbb{CS}(Y, \mu\mathcal{M}_Y)$. Then, the following hold:
\begin{enumerate}
\item The subvariety $Z$ is not a smooth point of the surface $E$;%
\item If the subvariety $Z$ is a curve, then it belongs to the
linear system $|\mathcal{O}_{\mathbb{P}(1,\,a,\,r-a)}(1)|$ defined
on the surface $E$ and all singular points of the surface $E$ are
contained in the set $\mathbb{CS}(Y, \mu\mathcal{M}_Y)$.
\end{enumerate}
\end{lemma}

\begin{proof}
We consider only the case when $r=5$ and $a=2$ because the proofs
for the other cases are very similar. Thus, we have
$E\cong\mathbb{P}(1,2,3)$. Let $Q_1$ and $Q_2$ be the singular
points of the surfaces $E$ and $L$ be the unique curve in
$|\mathcal{O}_{\mathbb{P}(1,\,2,\,3)}(1)|$ on the surface $E$.
Then, $L$ contains the singular points $Q_1$ and $Q_2$ but the
equivalence $\mu\mathcal{M}_Y\vert_{E}\equiv L$ holds by
Lemma~\ref{lemma:Kawamata}. Also, it follows from Lemma~\ref{lemma:Kawamata} that
the set  $\mathbb{CS}(Y, \mu\mathcal{M}_Y)$ contains both the points $Q_1$ and $Q_2$
if the curve $L$ is contained in the set  $\mathbb{CS}(Y, \mu\mathcal{M}_Y)$.

Suppose that the subvariety $Z$  is different from $L$, $Q_1$, and
$Q_2$. Let us show that this assumption gives us a contradiction.

Suppose that $Z$ is a point. Then, $Z$ is a smooth point of the
threefold $Y$, which implies the inequality
$\mathrm{mult}_{Z}(\mathcal{M}_Y)>\frac{1}{\mu}$. Let $C$ be a
general curve in $|\mathcal{O}_{\mathbb{P}(1,\,2,\,3)}(6)|$ on the
surface $E$ that passes through the point $Z$. Then, $C$ is not
contained in the base locus of $\mathcal{M}_Y$, which implies the
following contradictory inequality:
$$
\frac{1}{\mu}=C\cdot\mathcal{M}_Y\geq\mathrm{mult}_{Z}\big(\mathcal{M}_Y\big)>\frac{1}{\mu}.
$$

Therefore, the subvariety $Z$ must be a curve. Then,
$\mathrm{mult}_{Z}(\mathcal{M}_Y)\geq \frac{1}{\mu}$. Let $C$ be a
 general curve in the linear system
$|\mathcal{O}_{\mathbb{P}(1,\,2,\,3)}(6)|$ on the surface $E$.
Then, the curve $C$ is not contained in the base locus of the
pencil $\mathcal{M}_Y$. Therefore, we have
$$
\frac{1}{\mu}=C\cdot\mathcal{M}_Y\geq\mathrm{mult}_{Z}\big(\mathcal{M}_Y\big)C\cdot Z\geq \frac{1}{\mu}C\cdot Z,%
$$
which implies that $C\cdot Z=1$ on the surface  $E$. The equality
$C\cdot Z=1$ implies that the curve $Z$ is contained in the linear
system $|\mathcal{O}_{\mathbb{P}(1,\,2,\,3)}(1)|$ on the surface
$E$, which is impossible due to our assumption.
\end{proof}

\begin{lemma}
\label{lemma:negative-definite} Let $S$ be a normal surface and
$\Delta$ be an effective divisor on $S$ such that
$$
\Delta\equiv\sum_{i=1}^{r}a_{i}C_{i},
$$
where $C_{1},\cdots, C_{r}$ are irreducible curves on  $S$ and
$a_i$ is a rational number. If the intersection form of the curves
$C_1, \cdots, C_r$ on the surface $S$ is negative-definite, then
$\Delta=\sum_{i=1}^{r}a_{i}C_{i}$.
\end{lemma}

\begin{proof}
Let $\Delta=\sum_{i=1}^{k}c_{i}B_{i}$, where $B_{i}$ is an
irreducible curve on the surface $S$ and $c_{i}$ is a nonnegative
rational number. Suppose that
$$
\sum_{i=1}^{k}c_{i}B_{i}\ne \sum_{i=1}^{r}a_{i}C_{i}.
$$
Then, we may assume that each curve $B_{i}$ coincides with none of
the curves $C_{1},\cdots , C_{r}$. We have
\[\begin{split}
0&\geq\Big(\sum_{a_{i}>0}a_{i}C_{i}\Big)\cdot\Big(\sum_{a_{i}>0}a_{i}C_{i}\Big)\\
&=\Big(\sum_{i=1}^{k}c_{i}B_{i}\Big)\cdot\Big(\sum_{a_{i}>0}a_{i}C_{i}\Big)-
\Big(\sum_{a_{i}\leq 0}a_{i}C_{i}\Big)\cdot\Big(\sum_{a_{i}>0}a_{i}C_{i}\Big)\geq 0,%
\end{split}\]
which immediately implies $$\sum_{a_{i}>0}a_{i}C_{i}=0.$$
Therefore, we obtain the numerical equivalence
$$
\sum_{i=1}^{k}c_{i}B_{i}\equiv\sum_{a_{i}\leq 0}a_{i}C_{i}.%
$$
It then follows that $c_{i}=0$ and $a_{i}=0$ for every $i$.
\end{proof}

The following result is a generalization of Lemma~A.20 in
\cite{Ch05}.

\begin{theorem}
\label{theorem:main-tool} Let $\mathcal{B}$ be a linear system on
a threefold $X$ such that  general surface of the linear system
$\mathcal{B}$ is irreducible.  Then, the linear system
$\mathcal{B}$ coincides with the pencil $\mathcal{M}$ if one of
the following holds:
\begin{enumerate}
\item[(0)] There is a Zariski closed proper subset $\Sigma\subset
X$ such that
$$
\mathrm{Supp}\big(M\big)\cap\mathrm{Supp}\big(B\big)\subseteq\Sigma,%
$$
where $M$ and $B$ are general divisors of the pencil $\mathcal{M}$
and  the linear system $\mathcal{B}$, respectively. Note that the
general divisors $M$ and $B$ are chosen independently of the
proper subset $\Sigma$.
\end{enumerate}
For the below, let $B$ and $M$ be  general surfaces of the linear
system $\mathcal{B}$ and the pencil $\mathcal{M}$, respectively.
\begin{enumerate}

\item There is a nef and  big divisor $D$ on the threefold $X$ such that
$D\cdot M\cdot B=0$. \label{item:n-1}

\item The base locus of $\mathcal{B}$  consists of an irreducible
curve $C$ such that $M\cdot B\equiv\lambda C$ and $B\cdot C<0$ for
some positive rational number $\lambda$. \label{item:n-2}

\item The base locus of $\mathcal{M}$  consists of an irreducible
curve $C$ such that $M\cdot B\equiv\lambda C$ and $M\cdot C<0$ for
some positive rational number $\lambda$. \label{item:n-2-1}

\item\label{item:n-3} The equivalence $M\equiv \lambda B$ holds
for some positive rational number $\lambda$ and the base locus of
$\mathcal{B}$ consists of an irreducible curve $C$ such that
$B\cdot C<0$.

\item\label{item:n-3-1} The equivalence $M\equiv \lambda B$ holds
for some positive rational number $\lambda$ and the base locus of
$\mathcal{M}$ consists of an irreducible curve $C$ such that
$M\cdot C<0$.

\item\label{item:n-4} The surface $B$ is normal, the equivalence
$M\equiv\lambda B$ holds for some positive rational number
$\lambda$, and the base locus of $\mathcal{B}$ consists of
irreducible curves $C_{1},\cdots , C_{r}$ whose intersection form
on the surface $B$ is negative-definite.

\item\label{item:n-4-1} The surface $M$ is normal, the equivalence
$M\equiv \lambda B$ holds for some positive rational number
$\lambda$, and the base locus of $\mathcal{M}$ consists of
irreducible curves $C_{1},\cdots , C_{r}$ whose intersection form
on the surface $M$ is negative-definite.

\end{enumerate}
\end{theorem}
\begin{proof}
(0). Let $\rho:X\dasharrow\mathbb{P}^{1}$ be the rational map
induced by the pencil $\mathcal{M}$ and
$\xi:X\dasharrow\mathbb{P}^{r}$ be the rational map induced by the
linear system $\mathcal{B}$. We then consider a simultaneous
resolution of both the rational maps as follows:
$$
\xymatrix{
&&&W\ar@{->}[lld]_{\alpha}\ar@{->}[rrd]^{\beta}\ar@{->}[d]_{\pi}&&\\%
&\mathbb{P}^{r}&&X\ar@{-->}[rr]_{\rho}\ar@{-->}[ll]^{\xi}&&\mathbb{P}^{1},&}
$$ %
where $W$ is a smooth variety, $\pi$ is a birational morphism,
$\alpha$ and $\beta$ are morphisms.

Let $\Lambda$ be a Zariski closed subset of the variety $W$ such
that the morphism
$$
\pi\Big\vert_{W\setminus\Lambda}:W\setminus\Lambda\longrightarrow X\setminus\pi\big(\Lambda\big)%
$$
is an isomorphism and $\Delta$ be the union of the set $\Lambda$
and the closure of the proper transform of the set
$\Sigma\setminus\pi(\Lambda)$ on the variety $W$. Then, the set
$\Delta$ is a proper subvariety of $W$.

Suppose that the pencil $\mathcal{M}$ is different from the linear
system $\mathcal{B}$. Let $B_W$ be the pull-back of a  general
hyperplane of $\mathbb{P}^r$ by the morphism $\alpha$ and let
$M_W$ be a general fiber of the morphism $\beta$. Then, the
intersection $M_W\cap B_W$ is not empty  and the support
$\mathrm{Supp}(M_W\cap B_W)$ is not contained in $\Delta$. Hence,
we have
$$
\mathrm{Supp}\Big(\pi\big(M_W\big)\Big)\cap\mathrm{Supp}\Big(\pi\big(B_W\big)\Big)\not\subseteq\Sigma,%
$$
where $\pi(M_W)$ and $\pi(B_W)$ are general divisors in the linear
systems $\mathcal{M}$ and $\mathcal{B}$, respectively.

\ (\ref{item:n-1}). For some positive number $m$, there is an
ample divisor $A$ and an effective divisor $E$ such that $mD\sim
A+E$ since the divisor $D$ is nef and big. Then, the inequality $E\cdot
M\cdot B<0$ follows from $mD\cdot M\cdot B=0$ and $A\cdot M\cdot
B>0$. Therefore, the support of the cycle $M\cdot B$ must be
contained in the support of the effective divisor $E$, and hence
the statement~(0) completes the proof.

(\ref{item:n-2}). Let $H$ be an ample divisor on $X$. Then, there
is a positive rational number $\epsilon$ such that $(B+\epsilon
H)\cdot C=0$. Since $M\cdot B\equiv\lambda C$, we obtain
$(B+\epsilon H)\cdot M\cdot B=0$. Because the divisor $B +\epsilon
H$ is nef and big, (\ref{item:n-1}) implies $\mathcal{M}=\mathcal{B}$.

(\ref{item:n-2-1}). The proof is the same as (\ref{item:n-2}).

(\ref{item:n-3}). It is an immediate consequence of
(\ref{item:n-2}).

(\ref{item:n-3-1}). It is an immediate consequence of
(\ref{item:n-2-1}).

(\ref{item:n-4}). Since $M\equiv \lambda B$, the restricted
divisor  $M|_B$ of the effective divisor $M$ to the surface $B$ is
numerically equivalent to the divisor $\sum_{i=1}^{r}a_iC_i$ on
the normal surface $B$, where $a_i$ is a positive rational number.
It then follows from Lemma~\ref{lemma:negative-definite} that
\[M\Big|_B= \sum_{i=1}^{r}a_iC_i.\]
It implies that
$$\mathrm{Supp}\big(M\big)\cap\mathrm{Supp}\big(B\big)\subset\mathrm{Supp}\Big(\sum_{i=1}^{r}a_iC_i\Big).$$
Then, the statement~(0) completes the proof.

(\ref{item:n-4-1}). The proof is the same as (\ref{item:n-4}).
\end{proof}

\begin{theorem}
\label{theorem:Halphen} Suppose that a log pair $(X,
\mu\mathcal{M})$ is canonical with
$K_{X}+\mu\mathcal{M}\sim_{\mathbb{Q}} 0$. In addition, we suppose
that one of the following hold:
\begin{enumerate}
\item the base locus of the pencil $\mathcal{M}$ consists of
irreducible curves $C_{1},\cdots, C_{r}$ and  there is a nef and big
divisor $D$ on $X$ such that $D\cdot C_{i}=0$ for each $i$;

\item  the base locus of $\mathcal{M}$ consists of an irreducible
curve $C$ such that $\mathcal{M}\cdot C<0$;

\item  a general surface of the pencil $\mathcal{M}$ is normal,
the base locus of $\mathcal{M}$ consists of irreducible curves
$C_{1},\cdots, C_{r}$ whose intersection form is negative-definite
on a general surface in  the pencil $\mathcal{M}$.
\end{enumerate}
Then, the linear system $\mathcal{M}$ is a Halphen pencil and
there is a composition of antiflips $\xi:X\dasharrow X'$ along the
curves $C_{1},\cdots , C_{r}$ (or $C$) such that the proper
transform $\mathcal{M}_{X'}$ of the pencil $\mathcal{M}$ by $\xi$
is base-point-free.
\end{theorem}

\begin{proof}
The log pair $(X, \lambda\mathcal{M})$ is log-terminal for some
rational number $\lambda>\mu$. Hence, it follows from \cite{Sho93}
that there is a birational map $\xi:X\dasharrow X'$ such that
$\xi$ is an isomorphism in codimension one, the log pair $(X',
\lambda\mathcal{M}_{X'})$ is log-terminal, and the divisor
$K_{X'}+\lambda\mathcal{M}_{X'}$ is nef.

Let $H$ be a general surface in the pencil $\mathcal{M}_{X'}$.
Since
$$
H\equiv
\frac{1}{\lambda-\mu}\Big(K_{X'}+\lambda\mathcal{M}_{X'}-(K_X'+\mu\mathcal{M}_{X'})\Big),
$$
the divisor $H$ is nef. Hence, it follows from the log abundance
theorem (\cite{KMM}) that the linear system $|mH|$ is
base-point-free for some $m\gg 0$.

Let $\mathcal{B}$ be the proper transform of the linear system
$|mH|$ on $X$.  Also, let $B$ and $M$ be general surfaces of the
linear system $\mathcal{B}$ and the pencil $\mathcal{M}$,
respectively. Then, $B\equiv mM$ and one of the conditions in
Theorem~\ref{theorem:main-tool} is satisfied. Hence, we have
$\mathcal{M}=\mathcal{B}$, which implies that $m=1$ and
$\mathcal{M}_{X'}=|H|$ is base-point-free and induces a morphism
$\pi':X'\to \mathbb{P}^1$. Thus, every member of the pencil
$\mathcal{M}_{X'}$ is contracted to a point by the morphism
$\pi'$.

The log pair $(X', \mu\mathcal{M}_{X'})$ is canonical because the
map $\xi$ is a log flop with respect to the log pair $(X,
\mu\mathcal{M})$. In particular, the singularities of $X'$ are
canonical. Hence, the surface $H$ has at most Du Val singularities
because the pencil $\mathcal{M}_{X'}$ is base-point-free.
Moreover, the equivalence $K_{X'}+\mu H\sim_\mathbb{Q} 0$ and the Adjunction
formula imply that $K_{H}\sim 0$. Consequently, the linear system
$\mathcal{M}$ is a Halphen pencil.
\end{proof}

\begin{corollary}
\label{corollary:Halphen-K3} Under the assumption and notations of
Theorem~\ref{theorem:Halphen}, in addition, suppose that a general
surface of the pencil $\mathcal{M}$ is linearly equivalent to
$-nK_X$ for some natural number $n$. Then, a general element of
$\mathcal{M}$ is birational either to a smooth K3 surface or to an
abelian surface.
\end{corollary}
\begin{proof}
It immediately follows from the proof of Theorem~\ref{theorem:Halphen} and the classification of smooth surfaces
of Kodaira dimension zero.
\end{proof}

\begin{corollary}
\label{corollary:Halphen-K3-rational-curve} Under the assumption
and notations of Corollary~\ref{corollary:Halphen-K3},  suppose
that a general surface of the pencil $\mathcal{M}$ has a
rational curve not contained in the base locus of the pencil $\mathcal{M}$.
Then, a general
element of $\mathcal{M}$ is birational  to a smooth K3 surface.
\end{corollary}
\begin{proof}
In the proof of Theorem~\ref{theorem:Halphen}, suppose that the
surface $M$ has  a rational curve $L$ not contained in the base locus of the pencil $\mathcal{M}$.
Then, the surface $H$ contains a rational curve because the
birational map $\xi$ makes no change along the curve $L$. On the
other hand, the surface $H$ is birational either to a smooth K3
surface or to an abelian surface. However, an abelian surface
cannot contain a rational curve.
\end{proof}


\section{General results.}

Let $X$ be a general hypersurface of the $95$ families with entry
number $\gimel$ in $\mathbb{P}(1,a_1, a_2, a_3, a_4)$. In
addition, let $\mathcal{M}$ be a Halphen pencil on the threefold
$X$. Then, the log pair $(X, \frac{1}{n}\mathcal{M})$ is not
ter\-mi\-nal by Theorem~\ref{theorem:Noether-Fano}, where $n$ is
the natural number such that $\mathcal{M}\sim_{\mathbb{Q}}
-nK_{X}$. The following result is due to \cite{CPR}.

\begin{theorem}
\label{theorem:canonical-singularities} There is a birational
automorphism $\tau\in\mathrm{Bir}(X)$ such that the log pair $(X,
\frac{1}{m}\tau(\mathcal{M}))$ is canonical, where $m$ is the
natural number such that $\tau(\mathcal{M})\sim_{\mathbb{Q}}
-mK_{X}$.
\end{theorem}

To classify Halphen pencils on $X$ up to the action of
$\mathrm{Bir}(X)$, we may assume that  the log pair $(X,
\frac{1}{n}\mathcal{M})$ is canonical. However, it is not terminal
by Theorem~\ref{theorem:Noether-Fano}.

\begin{proposition}\label{proposition:Invariance}
The pencils constructed in
Examples~\ref{example:Halphen-pencil-I},
\ref{example:Halphen-pencil-II}, \ref{example:Halphen-pencil-III},
\ref{example:Halphen-pencil-IV}, and
\ref{example:Halphen-pencil-V} are invariant under the action of
the group $\mathrm{Bir}(X)$ of birational automorphisms of $X$.
\end{proposition}
\begin{proof}
Suppose that the hypersurface $X$ is defined by the equation
\[f_d(x,y,z,t,w)=0 \subset \mathrm{Proj} \big(\mathbb{C}[x,y,z,t,w]\big),\]
where $\mathrm{wt}(x)=1$, $\mathrm{wt}(y)=a_1$,
$\mathrm{wt}(z)=a_2$, $\mathrm{wt}(t)=a_3$, $\mathrm{wt}(w)=a_4$,
and $f_d$ is a general quasihomogeneous polynomial of degree
$d=\sum a_i$.

Since the hypersurface $X$ is general, it is not hard to see that the group $\mathrm{Aut}(X)$
of automorphisms of $X$ is either trivial or isomorphic to the
group of order $2$. The latter case happens when $2a_4=d$. In such
a case, the hypersurface $X$ can be defined by an equation of the
form
\[w^2=g_d(x,y,z,t),\]
where $g_d$ is a general quasihomogeneous polynomial of degree $d$
in variables $x$, $y$, $z$, and $t$. The group $\mathrm{Aut}(X)$
is generated by the involution $[x:y:z:t:w]\mapsto [x:y:z:t:-w]$.
Therefore, in both the cases,  we can see that the pencils
constructed in Examples~\ref{example:Halphen-pencil-I},
\ref{example:Halphen-pencil-II}, \ref{example:Halphen-pencil-III},
\ref{example:Halphen-pencil-IV}, and
\ref{example:Halphen-pencil-V} are invariant under the action of
the group $\mathrm{Aut}(X)$ of  automorphisms of $X$.

Suppose that the hypersurface $X$ is not superrigid, i.e., it has
a birational automorphism that is not biregular. Then, it is
either a quadratic involution or an elliptic involution that are
described  in \cite{CPR}. A quadratic involution has no effect on
things defined with the variables $x$, $y$, $z$, and $t$ (see
Theorem~4.9 in \cite{CPR}). On the other hand, an elliptic
involution has no effect on things defined with the variables $x$,
$y$, and $z$ (see Theorem~4.13 in \cite{CPR}). The pencils
constructed in Examples~\ref{example:Halphen-pencil-I},
\ref{example:Halphen-pencil-II}, \ref{example:Halphen-pencil-III},
and \ref{example:Halphen-pencil-IV} are defined by the variables
$x$, $y$, and $z$. Therefore, such pencils are invariant under the
action of the group $\mathrm{Bir}(X)$ of  birational automorphisms
of $X$. Meanwhile, the pencil constructed in
Example~\ref{example:Halphen-pencil-V} is contained in
$|-a_3K_X|$. However, in the case $\gimel=60$, the hypersurface
$X$ does not have an elliptic involution (see The Big Table in
\cite{CPR}), and hence the pencil is also
$\mathrm{Bir}(X)$-invariant.
\end{proof}

Every Halphen pencil on the threefold $X$ is, as shown throughout
the present article, birational to a pencil in
Examples~\ref{example:Halphen-pencil-I},
\ref{example:Halphen-pencil-II}, \ref{example:Halphen-pencil-III},
\ref{example:Halphen-pencil-IV}, and
\ref{example:Halphen-pencil-V} that is
$\mathrm{Bir}(X)$-invariant. It implies that every Halphen pencil
on $X$ is $\mathrm{Bir}(X)$-invariant.

\begin{lemma}
\label{lemma:smooth-points} Suppose that $\gimel\geq 3$. Then, the
set $\mathbb{CS}(X, \frac{1}{n}\mathcal{M})$  contains no smooth
point of $X$.
\end{lemma}
\begin{proof} It follows from the proof of Theorem~5.3.1 in
\cite{CPR}.
\end{proof}

\begin{corollary}
\label{corollary:curves} Suppose that the set $\mathbb{CS}(X,
\frac{1}{n}\mathcal{M})$ contains a curve $C$. Then, $-K_{X}\cdot
C\leq -K_{X}^{3}$.
\end{corollary}
\begin{proof}
It is an immediate consequence of Lemma~\ref{lemma:curves}.
\end{proof}

\begin{corollary}
\label{corollary:cheap-curves} The set $\mathbb{CS}(X,
\frac{1}{n}\mathcal{M})$ contains a singular point of $X$ whenever
$\gimel\geq 7$.
\end{corollary}
\begin{proof}
If $\gimel\geq 7$, then $-K_X^3<1$. Therefore, each curve $C$ with $-K_X\cdot C\leq -K_X^3$ passes through a singular point of $X$.
Then, the result follows from Lemmas~\ref{lemma:Kawamata}, \ref{lemma:smooth-points}, and  Corollary~\ref{corollary:curves}.
\end{proof}

In fact, we have a stronger result as follows:
\begin{theorem}
\label{theorem:Ryder} Suppose that $\gimel\geq 3$ and the set
$\mathbb{CS}(X, \frac{1}{n}\mathcal{M})$ contains a curve $C$.
Then,
$$
\mathrm{Supp}\big(C\big)\subset\mathrm{Supp}\Big(S_{1}\cdot S_{2}\Big),%
$$
where $S_{1}$ and $S_{2}$ are distinct surfaces of the linear
system $|-K_{X}|$.
\end{theorem}
\begin{proof}
See Section~3.1 in \cite{Ry06}.
\end{proof}

\begin{corollary}
\label{corollary:Ryder-a1} The set $\mathbb{CS}(X,
\frac{1}{n}\mathcal{M})$ contains no curves whenever $a_{1}\ne 1$.
\end{corollary}

\begin{corollary}
\label{corollary:Ryder-a2} If the set $\mathbb{CS}(X,
\frac{1}{n}\mathcal{M})$ contains a curve and $a_{2}\ne 1$, then
$\mathcal{M}=|-K_{X}|$.
\end{corollary}
\begin{proof}
It follows from Theorem~\ref{theorem:Ryder} and
Theorem~\ref{theorem:main-tool}.
\end{proof}

Suppose that $\gimel\geq 7$. Then,  the set $\mathbb{CS}(X,
\frac{1}{n}\mathcal{M})$ contains a singular point $P$ of type
$\frac{1}{r}(1,r-a,a)$, where $r\geq 2$, $r>a$, and $a$ is coprime
to $r$. Let $\pi:Y\to X$ be the Kawamata blow up at the singular
point $P$  and $E$ be its exceptional divisor. In addition, let
$\mathcal{M}_Y$ be the proper transform of the pencil
$\mathcal{M}$ by the birational morphism $\pi$. Then,
$\mathcal{M}_Y\sim_{\mathbb{Q}} -nK_{Y}$ by
Lemma~\ref{lemma:Kawamata}. The cone $\overline{\mathbb{NE}}(Y)$
of the threefold $Y$ contains two extremal rays $R_{1}$ and
$R_{2}$ such that $\pi$ is a contraction of the extremal ray
$R_{1}$. Moreover, the following result holds.

\begin{proposition}
\label{proposition:special-singular-points} Suppose that
$-K_{Y}^{3}\leq 0$ and $\gimel\ne 82$. Then, the threefold $Y$
contains irreducible surfaces $S\sim_{\mathbb{Q}}-K_{Y}$ and
$T\sim_{\mathbb{Q}} -bK_{Y}+cE$ whose scheme-theoretic
intersection is an irreducible reduced curve that generates
$R_{2}$, where $b>0$ and $c\geq 0$ are integer numbers.
\end{proposition}
\begin{proof}
See Lemma~5.4.3 in \cite{CPR}.
\end{proof}
The values of $b$ and $c$ for a given singular point in
Proposition~\ref{proposition:special-singular-points} appear in
The Table of Part~\ref{section:weighted-Fanos}.

\begin{lemma}
\label{lemma:special-singular-points-with-positive-c} Under the
assumptions and notations of
Proposition~\ref{proposition:special-singular-points}, suppose
that the inequality $-K_{Y}^{3}<0$ holds. Then, the number $c$ is
zero.
\end{lemma}

\begin{proof}
Let $M_{1}$ and $M_{2}$ be general surfaces of the pencil
$\mathcal{M}_Y$. Then, $M_{1}\cdot M_{2}\equiv n^{2}K^{2}_{Y}$,
which implies that $M_{1}\cdot
M_{2}\not\in\overline{\mathbb{NE}}(Y)$ in the case $c>0$ by
Proposition~\ref{proposition:special-singular-points}.
\end{proof}

\begin{lemma}
\label{lemma:special-singular-points-with-zero-c} Under the
assumptions and notations of
Proposition~\ref{proposition:special-singular-points}, suppose
that the inequality $-K_{Y}^{3}<0$ holds. Then, the pencil
$\mathcal{M}_Y$ is generated by the divisors $bS$ and $T$.
\end{lemma}

\begin{proof}
Let $M_{1}$ and $M_{2}$ be general surfaces in $\mathcal{M}_Y$.
Then, $M_{1}\cdot M_{2}\in \overline{\mathbb{NE}}(Y)$ and
$M_{1}\cdot M_{2}\equiv n^{2}K^{2}_{Y}$, which implies that
$M_{1}\cdot M_{2}\in \mathbb{R}^{+}R_2$ because $c=0$ by
Lemma~\ref{lemma:special-singular-points-with-positive-c}.
Moreover, we have
$$
\mathrm{Supp}\big(\Gamma\big)=\mathrm{Supp}\Big(M_{1}\cdot M_{2}\Big),%
$$
because $T\cdot R_2<0$ and $S\cdot R_2<0$. Therefore, the pencil
$\mathcal{M}_Y$ coincides with the pencil generated by the
divisors $bS$ and $T$ by Theorem~\ref{theorem:main-tool}, which
completes the proof.
\end{proof}

When $a_1=1$, a general surface of a pencil contained in the
linear system $|-K_X|$ is birational to a K3 surface. In
particular, it is one of  Reid's 95 codimension 1 weighted K3
surfaces. Furthermore, we have the following:
\begin{proposition}\label{proposition:Halphen-K3}
In every case, a general surface in a pencil contained in
$|-a_1K_X|$ is birational to a K3 surface.
\end{proposition}
\begin{proof}
See \cite{ChPa05}\footnote{The cases $\gimel=18$, $22$, $28$ is
not covered by  the article \cite{ChPa05}. Moreover, it has a
mistake in Lemma~3.1. So we reprove  Lemma~3.1 of \cite{ChPa05} in
this article. See K3-Propositions~\ref{K3proposition:n-18},
\ref{K3proposition:n-18-i}, \ref{K3proposiotn:n-22}, and
\ref{K3proposiotn:n-28}.}
\end{proof}
Therefore, general surfaces in the pencils of
Examples~\ref{example:Halphen-pencil-I} and
\ref{example:Halphen-pencil-II} are birational to K3 surfaces.

\section{Notations.}
Let us describe the notations we will use. Otherwise mentioned,
these notations are fixed from Part~\ref{section:one-Halphen-1} to
Part~\ref{section:weighted-Fanos}.

\begin{itemize}
\item In the weighted projective space $\mathbb{P}(1,a_1,
a_2,a_3,a_4)$, we assume that $a_1\leq a_2\leq a_3\leq a_4$. For
weighted homogeneous coordinates, we always use $x$, $y$, $z$,
$t$, and $w$ with $\mathrm{wt}(x)=1$, $\mathrm{wt}(y)=a_1$,
$\mathrm{wt}(z)=a_2$, $\mathrm{wt}(t)=a_3$, and
$\mathrm{wt}(w)=a_4$.

\item The number $\gimel$ always means the entry number of each
family of weighted Fano hypersurfaces in The  Table of
Part~\ref{section:weighted-Fanos}.

\item In each family, we always let $X$ be a general quasismooth
hypersurface of degree $d$ in the weighted projective space
$\mathbb{P}(1,a_1, a_2,a_3,a_4)$, where $d=\sum_{i=1}^4 a_i$.

\item On the threefold $X$, a given Halphen pencil is denoted by
$\mathcal{M}$.

\item For a given Halphen pencil $\mathcal{M}$, we always assume
that $\mathcal{M}\sim_{\mathbb{Q}} -nK_X$.

\item When a morphism $f:V\to W$ is given, the proper transforms
of  a curve $Z$, a surface $D$, and a linear system $\mathcal{D}$
on $W $ by the morphism $f$ will be always denoted by $Z_V$,
$D_V$, and $\mathcal{D}_V$, respectively, i.e., we use the ambient
space $V$ as their subscripts.

\item  $S$ : the surface on $X$ defined by the equation $x=0$.

 \item  $S^y$ : the surface on $X$ defined by the equation
$y=0$.

\item  $S^z$ : the surface on $X$ defined by the equation $z=0$.

\item $S^t$ : the surface on $X$ defined by the equation $t=0$.

\item $S^w$ : the surface on $X$ defined by the equation $w=0$.

\item $C$ : the curve on $X$ defined by the equations $x=y=0$.

\item $\bar{C}$ : the curve on $X$ defined by the equations
$x=z=0$.

\item $\tilde{C}$ : the curve on $X$ defined by the equations
$x=t=0$.

\item $\hat{C}$ : the curve on $X$ defined by the equations
$x=w=0$.
\end{itemize}

\section{Set-up.}
In each case, for a given general hypersurface $X$ and a given
Halphen pencil $\mathcal{M}$, we consider the log pair $(X,
\frac{1}{n}\mathcal{M})$. Note that the natural number $n$ is
given by $\mathcal{M}\sim_{\mathbb{Q}} -nK_X$. At the beginning of
each section, we state the degree of $X$, the weights of the
ambient weighted projective space, and  the singularities of $X$.
When the section contains a single case, we describe all the
singular points. Because we could not put the values of $b$ and
$c$ in Proposition~\ref{proposition:special-singular-points},
reader should refer to The  Table in
Part~\ref{section:weighted-Fanos} for each singular point.

After we describe these simple features, we present elliptic
fibrations into which the general hypersurface  $X$ can be
birationally transformed. For detail, reader should refer to
\cite{Ch05}. Note that for the cases $\gimel=3$, $60$, $75$, $84$,
$87$, $93$, the hypersurface $X$ cannot be birationally
transformed into an elliptic fibration (\cite{ChPa05}).

We try to show that the Halphen pencil $\mathcal{M}$ is one of the
pencils given in Examples~\ref{example:Halphen-pencil-I},
\ref{example:Halphen-pencil-II}, \ref{example:Halphen-pencil-III},
\ref{example:Halphen-pencil-IV}, and
\ref{example:Halphen-pencil-V}. To do so, we will consider the set
$$\mathbb{CS}\Big(X, \frac{1}{n}\mathcal{M}\Big).$$
Theorem~\ref{theorem:Noether-Fano} implies that  \emph{the set is always
non-empty} since the linear system $\mathcal{M}$ is a Halphen
pencil and the anticanonical divisor $-K_X$ is nef and big.
Furthermore, due to Theorem~\ref{theorem:canonical-singularities},
there is a birational automorphism $\rho\in\mathrm{Bir}(X)$ such
that the log pair $(X, \frac{1}{\bar{n}}\rho(\mathcal{M}))$ is
canonical for the natural number $\bar{n}$ with
$\rho(\mathcal{M})\sim_{\mathbb{Q}} -\bar{n}K_X$. It will turn out
that the pencil $\rho(\mathcal{M})$ is one of the pencils
constructed in Examples~\ref{example:Halphen-pencil-I},
\ref{example:Halphen-pencil-II}, \ref{example:Halphen-pencil-III},
\ref{example:Halphen-pencil-IV}, and
\ref{example:Halphen-pencil-V} that are
$\mathrm{Bir}(X)$-invariant
(Proposition~\ref{proposition:Invariance}). It implies
$\mathcal{M}=\rho(\mathcal{M})$.  For this reason, we may always assume
that \emph{the log pair $(X, \frac{1}{n}\mathcal{M})$ is canonical.}

At the first stage, we use Corollaries~\ref{corollary:Ryder-a1}
and \ref{corollary:Ryder-a2},  if applicable, in order to exclude
the case when the set $\mathbb{CS}(X, \frac{1}{n}\mathcal{M})$
contains a curve. Then,  \emph{usually the set contains only singular
points of $X$ due to Lemma~\ref{lemma:smooth-points}.}

Next, we apply
Lemmas~\ref{lemma:special-singular-points-with-positive-c} and
~\ref{lemma:special-singular-points-with-zero-c} with The
Table in order to minimize, as much as possible, the set
$\mathbb{CS}(X, \frac{1}{n}\mathcal{M})$ to be considered.

After such things to do, we start the game to identify the Halphen
pencil $\mathcal{M}$. We need to calculate base curves of proper
transforms of various linear systems by various Kawamata blow ups,
their multiplicities, and so on. Due to huge volume  of
calculations, we usually omit the calculations. However, from time
to time, we present them to show how to calculate.

In addition, \emph{unless otherwise mentioning, whenever we consider a log pair $(Y, \frac{1}{n}\mathcal{M}_Y)$
that is obtained from the log pair $(X, \frac{1}{n}\mathcal{M}_X)$ by a sequence of Kawamata blow ups, we always
assume that $\mathcal{M}_Y\sim_{\mathbb{Q}} -nK_Y$.} This condition is always satisfied when the Kawamata blow ups are obtained from
centers of canonical singularities due to Lemma~\ref{lemma:Kawamata}.

In Parts~\ref{section:two-Halphens} and
~\ref{section:more-than-two-Halphens}, we also show that  general
surfaces of   Halphen pencils of
types~\ref{example:Halphen-pencil-III},
\ref{example:Halphen-pencil-IV}, and
\ref{example:Halphen-pencil-V} are birational to smooth  K3
surfaces. It will be proved  directly or by using
Corollary~\ref{corollary:Halphen-K3-rational-curve}. Such
statements are titled by K3-Proposition to be simply distinguished
from the other works. The cases in
Examples~\ref{example:Halphen-pencil-I} and
\ref{example:Halphen-pencil-II} are covered by
Proposition~\ref{proposition:Halphen-K3}.

\newpage
\part{Fano threefold hypersurfaces with a single Halphen
pencil $|-K_X|$.}\label{section:one-Halphen-1}

\section{Case $\gimel=7$, hypersurface of degree $8$ in
$\mathbb{P}(1,1,2,2,3)$.} \index{$\gimel=07$}\label{section:n-7}

The threefold $X$ is a general hypersurface of degree $8$ in
$\mathbb{P}(1,1,2,2,3)$ with $-K_X^3=\frac{2}{3}$. The
singularities of the hypersurface $X$ consist of four points
$P_{1}$, $P_{2}$, $P_{3}$ and $P_{4}$ that are quotient
singularities of type $\frac{1}{2}(1,1,1)$ and a point $Q$ that is
a quotient singularity of type $\frac{1}{3}(1,1,2)$.

For each point $P_{i}$,  there is a commutative diagram
$$
\xymatrix{
&U_{i}\ar@{->}[d]_{\alpha_{i}}&&Y_{i}\ar@{->}[ll]_{\beta_{i}}\ar@{->}[d]^{\eta_{i}}&\\
&X\ar@{-->}[rr]_{\xi_{i}}&& \mathbb{P}(1,1,2),&}
$$
where \begin{itemize} \item $\xi_{i}$ is a rational map defined in
the outside of  the  points $P_{i}$ and $Q$,

\item $\alpha_{i}$ is the Kawamata blow up at the point $P_{i}$
with weights $(1,1,1)$,

\item $\beta_{i}$ is the Kawamata blow up with weights $(1,1,2)$
at the point $Q_i$ whose image to $X$ is the point $Q$,

\item $\eta_{i}$ is an elliptic fibration.

\end{itemize}
There is another rational map
$\xi_{0}:X\dasharrow\mathbb{P}(1,1,2)$ that gives us the following
commutative diagram:
$$
\xymatrix{
&U_{0}\ar@{->}[d]_{\alpha_{0}}&&Y_{0}\ar@{->}[ll]_{\beta_{0}}\ar@{->}[d]^{\eta_{0}}&\\
&X\ar@{-->}[rr]_{\xi_{0}}&& \mathbb{P}(1,1,2),&}
$$
where \begin{itemize} \item $\alpha_{0}$ is the Kawamata blow up
at the point $Q$ with weights $(1,1,2)$,

\item $\beta_{0}$ is the Kawamata blow up with weights $(1,1,1)$
at the singular point $O$ in the exceptional divisor of
$\alpha_0$,

\item  $\eta_{0}$ is an elliptic fibration.
\end{itemize}

The pencil $|-K_{X}|$ is invariant under the action of the group
of birational automorphisms $\mathrm{Bir}(X)$ and hence we may
assume that
$$
\varnothing\ne\mathbb{CS}\Big(X, \frac{1}{n}\mathcal{M}\Big)\subseteq\Big\{P_{1}, P_{2}, P_{3}, P_{4}, Q\Big\}%
$$
due to Lemmas~\ref{lemma:smooth-points},
\ref{lemma:special-singular-points-with-positive-c},
\ref{lemma:special-singular-points-with-zero-c} and
Corollary~\ref{corollary:Ryder-a2}. Note that the base locus of
the pencil $|-K_X|$ consists of the irreducible curve $C$ defined
by $x=y=0$.

\begin{lemma}
\label{lemma:n-7-two-singular-points-of-index-2} If the
$\mathbb{CS}(X, \frac{1}{n}\mathcal{M})$ contains distinct points
$P_i$ and $P_j$, then $\mathcal{M}=|-K_{X}|$.
\end{lemma}

\begin{proof}
 Let $\pi_j:W\to U_i$ be the Kawamata blow up with weights $(1,1,1)$ at the point whose image to $X$ is the point $P_j$
 and $D$ be a general
surface of the pencil $|-K_{X}|$. The base locus of the pencil
$|-K_{W}|$ consists of the irreducible curve $C_W$. The surface
$D_W$ is normal and  $C_W^{2}=-\frac{1}{3}$ on the surface $D_W$,
which implies that $\mathcal{M}_W=|-K_{W}|$ by
Theorem~\ref{theorem:main-tool}.
\end{proof}

The singularities of the log pair $(U_{i},
\frac{1}{n}\mathcal{M}_{U_i})$ are not terminal  because the
divisor $-K_{U_{i}}$ is nef and big.

\begin{lemma}
\label{lemma:n-7-singular-point-of-index-two} The set
$\mathbb{CS}(X, \frac{1}{n}\mathcal{M})$  cannot consist of a
single point $P_{i}$.
\end{lemma}

\begin{proof}
Suppose that the singularities of the log pair $(X,
\frac{1}{n}\mathcal{M})$ are terminal in the outside the singular
point $P_{i}$. Let $E_{i}$ be the exceptional divisor of
$\alpha_{i}$. Then, the set $\mathbb{CS}(U_{i},
\frac{1}{n}\mathcal{M}_{U_i})$ contains a line $Z\subset
E_{i}\cong\mathbb{P}^{2}$ by Lemmas~\ref{lemma:Kawamata} and
\ref{lemma:Cheltsov-Kawamata}. But this is a contradiction because
of Lemma~\ref{lemma:curves}.
\end{proof}

\begin{lemma}
\label{lemma:n-7-singular-points-of-index-two-and-three}  If the
set $\mathbb{CS}(X, \frac{1}{n}\mathcal{M})$ contains the points
$P_{i}$ and $Q$, then $\mathcal{M}=|-K_{X}|$.
\end{lemma}

\begin{proof}
Let $G_i$ be the exceptional divisor of  $\beta_{i}$. Then,
$\mathcal{M}_{Y_i}\sim_{\mathbb{Q}}-nK_{Y_{i}}$ by
Lemma~\ref{lemma:Kawamata}, which implies that every member of the
pencil $\mathcal{M}_{Y_i}$ is contracted to a curve by the
morphism $\eta_{i}$. In particular, the base locus of
$\mathcal{M}_{Y_i}$ does not contain curves that are not
contracted by $\eta_{i}$. Therefore, the set $\mathbb{CS}(Y_{i},
\frac{1}{n}\mathcal{M}_{Y_i})$ contains the singular point $O_i$
of the surface $G_i\cong\mathbb{P}(1,1,2)$ due to
Theorem~\ref{theorem:Noether-Fano} and
Lemma~\ref{lemma:Cheltsov-Kawamata}.

Let $\pi_i:W_i\to Y_{i}$ be the Kawamata blow up at the point
$O_i$ with weights $(1,1,1)$ and $D$ be a  general surface of the
pencil $|-K_{X}|$. The base locus of the pencil $|-K_{W_i}|$
consists of the irreducible curve $C_{W_i}$. The surface $D_{W_i}$
is normal and  $C_{W_i}^{2}=-\frac{1}{2}$ on the surface
$D_{W_i}$, which implies that $\mathcal{M}_{W_i}=|-K_{W_i}|$ by
Theorem~\ref{theorem:main-tool}.
\end{proof}

\begin{lemma}
\label{lemma:n-7-curves-on-second-floor} The log pair $(U_{0},
\frac{1}{n}\mathcal{M}_{U_0})$ is terminal in the outside of the
singular point $O$ of the exceptional divisor of the birational
morphism $\beta_{0}$.
\end{lemma}

\begin{proof}
Let $E_{0}$ be the exceptional divisor of the birational morphism
$\beta_{0}$. Suppose that
 the log pair $(U_{0},
\frac{1}{n}\mathcal{M}_{U_0})$ is not terminal in the outside of
the singular point $O$. Then, the set $\mathbb{CS}(U_{0},
\frac{1}{n}\mathcal{M}_{U_0})$ consists of the point $O$ and an
irreducible curve $L\subset E_{0}\cong\mathbb{P}(1,1,2)$ such that
$L\in|\mathcal{O}_{\mathbb{P}(1,\,1,\,2)}(1)|$.

The threefold $X$ can be given by the equation
$$
w^{2}z+wf_{5}\big(x,y,z,t\big)+f_{8}\big(x,y,z,t\big)=0,
$$
where  $f_{i}$ is a qua\-si\-ho\-mo\-ge\-ne\-ous polynomial of
degree $i$. Moreover, we may assume that the curve $L$ is cut out
on the surface $E_{0}$ by the surface $S_{U_{0}}$.

Let $\mathcal{P}$ be the pencil on $X$ cut out on the threefold
$X$ by the pencil $\lambda x^{2}+\mu z=0$, where
$(\lambda:\mu)\in\mathbb{P}^{1}$. Then, the base locus of
$\mathcal{P}$ consists of the irreducible curve $\bar{C}$ cut by
the equations $x=z=0$.

The base locus of the linear system $\mathcal{P}_{U_0}$ consists
of the irreducible curves $L$ and $\bar{C}_{U_0}$.  A general
surface $D$ in $\mathcal{P}_{U_0}$ is normal. Moreover, we have
$$
\mathcal{M}_{U_0}\Big\vert_{D}\equiv -nK_{U_{0}}\Big\vert_{D}\equiv nS_{U_0}\Big\vert_{D}\equiv n\big(L+\bar{C}_{U_0}\big)%
$$
by Lemma~\ref{lemma:Kawamata}. On the other hand, we have
$S_{U_0}\cdot D=L+\bar{C}_{U_0}$ and $E_0\cdot D=2L$, which
implies that $\bar{C}_{U_0}^{2}=-\frac{5}{4}$ on the surface $D$.
The latter together with the equality
$\mathrm{mult}_{L}(\mathcal{M}_{U_0})=n$ easily implies that
$\mathcal{M}_{U_0}=\mathcal{P}_{U_0}$ due to
Theorem~\ref{theorem:main-tool}. We obtain $n=2$, but $D$ is
normal and $\mathrm{mult}_{L}(D)\ne 2$, which is a contradiction.
\end{proof}

\begin{proposition}
\label{proposition:n-7} The linear system $|-K_{X}|$ is the only
Halphen pencil on $X$.
\end{proposition}
\begin{proof}
By the previous arguments, we have only to consider the case when
$$
\mathbb{CS}\Big(X, \frac{1}{n}\mathcal{M}\Big)=\Big\{ Q \Big\}%
$$
Then, $\mathcal{M}_{Y_0}\sim_{\mathbb{Q}}-nK_{Y_{0}}$, which
implies that every member of the pencil $\mathcal{M}_{Y_0}$ is
contracted to a curve by the morphism $\eta_{0}$. In particular,
the base locus of the pencil $\mathcal{M}_{Y_0}$ does not contain
any curves by Lemma~\ref{lemma:n-7-curves-on-second-floor}. Thus,
the singularities of the log pair $(Y_{0},
\frac{1}{n}\mathcal{M}_{Y_0})$ are terminal by
Lemma~\ref{lemma:Cheltsov-Kawamata}, which is impossible by
Theorem~\ref{theorem:Noether-Fano}.
\end{proof}

\section{Cases $\gimel=9$, $11$, and $30$.}\label{section:n-9}

We first consider the case $\gimel=9$.\index{$\gimel=09$} The
variety $X$ is a general hypersurface of degree $9$  in
$\mathbb{P}(1,1,2,3,3)$ with $-K_{X}^{3}=\frac{1}{2}$. The
singularities of the hypersurface $X$ consist of one point $O$
that is a quotient singularity of type $\frac{1}{2}(1,1,1)$ and
three points $P_{1}$, $P_{2}$ and $P_{3}$ that are quotient
singularities of type $\frac{1}{3}(1,1,2)$.

We have the following commutative diagram:
$$
\xymatrix{
&&W\ar@{->}[dl]_{\alpha}\ar@{->}[dr]^{\eta}&&\\
&X\ar@{-->}[rr]_{\psi}&&\mathbb{P}(1,1,2),&}
$$
where \begin{itemize} \item $\psi$ is the natural projection,

\item $\alpha$ is the composition of the Kawamata blow ups at the
points $P_{1}$, $P_{2}$, $P_{3}$ with weights $(1,1,2)$.

\item $\eta$ is an elliptic fibration.
\end{itemize}
There is another elliptic fibration as follows:
$$
\xymatrix{
&&V\ar@{->}[dl]_{\pi}\ar@{->}[dr]^{\eta_0}&&\\
&X\ar@{-->}[rr]_{\chi}&&\mathbb{P}(1,1,3),&}
$$
where \begin{itemize} \item $\pi$ is the Kawamata blow up at the
point $O$ with weights $(1,1,1)$,

\item $\eta_0$ is an elliptic fibration. \end{itemize}

It follows from Lemma~\ref{lemma:smooth-points} that the set
$\mathbb{CS}(X, \frac{1}{n}\mathcal{M})$ does not contain smooth
points of the hypersurface $X$. Moreover, if it contains a curve
then we obtain $\mathcal{M}=|-K_X|$ from
Corollary~\ref{corollary:Ryder-a2}. Therefore, we may assume that
$$
\varnothing\ne\mathbb{CS}\Big(X, \frac{1}{n}\mathcal{M}\Big)\subset \Big\{P_{1}, P_{2}, P_{3}, O\Big\}.%
$$

\begin{lemma}
\label{lemma:n-9-O} The set $\mathbb{CS}(X,
\frac{1}{n}\mathcal{M})$ does not consist of the point $O$.
\end{lemma}

\begin{proof}
Suppose that the set $\mathbb{CS}(X, \frac{1}{n}\mathcal{M})$
consists of the point $O$. Then,
$\mathcal{M}_{V}\sim_{\mathbb{Q}}-nK_{V}$ by
Lemma~\ref{lemma:Kawamata}, which implies that every member in the
pencil $\mathcal{M}_{V}$ is contracted to a curve by the morphism
$\eta_0$. In particular, the set $\mathbb{CS}(V,
\frac{1}{n}\mathcal{M}_{V})$ does not contain curves. On the other
hand, the set $\mathbb{CS}(V, \frac{1}{n}\mathcal{M}_{V})$ is not
empty by Theorem~\ref{theorem:Noether-Fano}. Hence, the set
$\mathbb{CS}(V, \frac{1}{n}\mathcal{M}_{V})$ contains a point of
the exceptional divisor of $\pi$, which is impossible by
Lemma~\ref{lemma:Cheltsov-Kawamata}.
\end{proof}

Note that the base locus of $|-K_X|$ consists of the irreducible
curve $C$ cut by $x=y=0$.
\begin{lemma}
\label{lemma:n-9-O-P} If the set $\mathbb{CS}(X,
\frac{1}{n}\mathcal{M})$ contains the points $O$ and $P_{i}$, then
$\mathcal{M}=|-K_{X}|$.
\end{lemma}

\begin{proof}
Let $\beta_i:W_i\to V$ be  the Kawamata blow up with weights
$(1,1,3)$ at the point whose image to $X$ is the point $P_{i}$.
Then, $|-K_{W_i}|$ is the proper transform of the pencil
$|-K_{X}|$ and the base locus of $|-K_{W_i}|$ consists of the
irreducible curve $C_{W_i}$. One can easily see that
$-K_{W_i}\cdot C_{W_i}<0$. On the other hand, we have
$\mathcal{M}_{W_i}\sim_{\mathbb{Q}}-nK_{W_i}$ by
Lemma~\ref{lemma:Kawamata}, which implies that
$\mathcal{M}_{W_i}=|-K_{W_i}|$ by Theorem~\ref{theorem:main-tool}.
\end{proof}

Therefore, we may further assume that
$$
\varnothing\ne\mathbb{CS}\Big(X, \frac{1}{n}\mathcal{M}\Big)\subset \Big\{P_{1}, P_{2}, P_{3}\Big\}.%
$$

Let $\alpha_{i}:U_{i}\to X$ be the Kawamata blow up at the point
$P_{i}$ and  $F_{i}$ be the exceptional divisor of $\alpha_{i}$.
The exceptional divisor $F_i$ contains a singular point $Q_i$ that
is a quotient singular point of the threefold $U_{i}$ of type
$\frac{1}{2}(1,1,1)$.

\begin{lemma}
\label{lemma:n-9-P-Q} If the set $\mathbb{CS}(U_{i},
\frac{1}{n}\mathcal{M}_{U_{i}})$ contains the point $Q_{i}$, then
$\mathcal{M}=|-K_{X}|$.
\end{lemma}

\begin{proof}
Let $\beta_i:Y_i\to U_{i}$ be the Kawamata blow up at the point
$Q_{i}$. Then, $|-K_{Y_i}|$ is the proper transform of the pencil
$|-K_{X}|$ and the base locus of $|-K_{V_i}|$ consists of the
irreducible curve $C_{Y_i}$. Also, we have  $-K_{Y_i}\cdot C<0$.
On the other hand,  $\mathcal{M}_{Y_i}\sim_{\mathbb{Q}}-nK_{Y_i}$
by Lemma~\ref{lemma:Kawamata}, which implies that
$\mathcal{M}_{Y_i}=|-K_{Y_i}|$ by Theorem~\ref{theorem:main-tool}.
Hence, we obtain $\mathcal{M}=|-K_{X}|$.
\end{proof}

Thus, it follows from Lemma~\ref{lemma:Cheltsov-Kawamata} that we
may assume that  the set $\mathbb{CS}(U_{i},
\frac{1}{n}\mathcal{M}_{U_{i}})$ does not contain subvarieties of
$F_{i}$ in the case when the set $\mathbb{CS}(X,
\frac{1}{n}\mathcal{M})$ contains the point $P_{i}$.

\begin{lemma}
\label{lemma:n-9-P} The set $\mathbb{CS}(X,
\frac{1}{n}\mathcal{M})$ does not consist of the point $P_{i}$.
\end{lemma}

\begin{proof}
Suppose that the set $\mathbb{CS}(X, \frac{1}{n}\mathcal{M})$
consists of the point $P_{i}$. Then,
$\mathcal{M}_{U_{i}}\sim_{\mathbb{Q}}-nK_{U_{i}}$ by
Lemma~\ref{lemma:Kawamata}, which implies that the set
$\mathbb{CS}(U_{i}, \frac{1}{n}\mathcal{M}_{U_{i}})$ is not empty
by Theorem~\ref{theorem:Noether-Fano}, because the divisor
$-K_{U_{i}}$ is nef and big. Therefore, the set
$\mathbb{CS}(U_{i}, \frac{1}{n}\mathcal{M}_{U_{i}})$ must contain
a subvariety of $G_{i}$, which is impossible by our assumption.
\end{proof}

\begin{proposition}
\label{proposition:n-9} If $\gimel=9$, then the linear system
$|-K_{X}|$ is the only Halphen pencil on $X$.
\end{proposition}
\begin{proof}
By the proof of Lemma~\ref{lemma:n-9-P}, we may assume that
$\mathbb{CS}(X, \frac{1}{n}\mathcal{M})=\{P_{1}, P_{2}, P_{3}\}$,
which implies that $\mathcal{M}_{W}\sim_{\mathbb{Q}}-nK_{W}$ by
Lemma~\ref{lemma:Kawamata}. Therefore, the set $\mathbb{CS}(W,
\frac{1}{n}\mathcal{M}_{W})$ must contain a subvariety of an
exceptional divisor of $\alpha$ by
Theorem~\ref{theorem:Noether-Fano}, which is impossible, because
we assumed that the set $\mathbb{CS}(U_{i},
\frac{1}{n}\mathcal{M}_{U_{i}})$ does not contain subvarieties of
$F_{i}$. The obtained contradiction concludes the proof.
\end{proof}

For the case $\gimel=30$\index{$\gimel=30$}, let $X$ be a general
hypersurface of degree $16$  in $\mathbb{P}(1,1,3,4,8)$ with
$-K_{X}^{3}=\frac{1}{6}$. The singularities of $X$ consist of one
point $O$ that is a quotient singularity of type
$\frac{1}{3}(1,1,2)$ and two points $P_{1}$, $P_{2}$ that are
quotient singularities of type $\frac{1}{4}(1,1,3)$.

We have the following commutative diagram:
$$
\xymatrix{
&&W\ar@{->}[dl]_{\alpha}\ar@{->}[dr]^{\eta}&&\\
&X\ar@{-->}[rr]_{\psi}&&\mathbb{P}(1,1,3),&}
$$
where \begin{itemize} \item $\psi$ is the natural projection,

\item $\alpha$ is the composition of the Kawamata blow ups at the
points $P_{1}$ and  $P_{2}$ with weights $(1,1,3)$.

\item $\eta$ is an elliptic fibration.
\end{itemize}
There is another elliptic fibration as follows:
$$
\xymatrix{
&&V\ar@{->}[dl]_{\pi}\ar@{->}[dr]^{\eta_0}&&\\
&X\ar@{-->}[rr]_{\chi}&&\mathbb{P}(1,1,4),&}
$$
where \begin{itemize} \item $\pi$ is the Kawamata blow up at the
point $O$ with weights $(1,1,2)$,

\item $\eta_0$ is an elliptic fibration. \end{itemize}

\begin{proposition}
\label{proposition:n-30} If $\gimel=30$, then the linear system
$|-K_{X}|$ is the only Halphen pencil on $X$.
\end{proposition}
\begin{proof}
The proof is the same as the case $\gimel=9$.
\end{proof}

In the case $\gimel=11$\index{$\gimel=11$}, the threefold $X$ is a
general hypersurface of degree $10$ in $\mathbb{P}(1,1,2,2,5)$
with $-K_X^3=\frac{1}{2}$.  Its singularities consist of five
points $P_1, \cdots, P_5$ that are quotient singularities of type
$\frac{1}{2}(1,1,1)$. For each  singular point $P_i$, we have an
elliptic fibration as follows:
$$
\xymatrix{
&&U_{i}\ar@{->}[dl]_{\pi_{i}}\ar@{->}[dr]^{\eta_{i}}&\\
&X\ar@{-->}[rr]_{\xi_{i}}&&\mathbb{P}(1,1,2),&&}
$$
where \begin{itemize}

\item $\pi_{i}$ is the Kawamata blow up at the point $P_{i}$ with
weights $(1,1,1)$,

\item $\eta_{i}$ is an elliptic fibration.

\end{itemize}

\begin{proposition}
\label{proposition:n-11} If $\gimel=11$, the linear system
$|-K_{X}|$ is a unique Halphen pencil on $X$.
\end{proposition}
\begin{proof}
If the set $\mathbb{CS}(X, \frac{1}{n}\mathcal{M})$ contains a
curve, then we obtain $\mathcal{M}=|-K_X|$ from
Corollary~\ref{corollary:Ryder-a2}. Thus, we may assume that
$$\mathbb{CS}\Big(X, \frac{1}{n}\mathcal{M}\Big)\subset\Big\{P_1, P_2, P_3,
P_4, P_5\Big\}$$ by Lemma~\ref{lemma:smooth-points}. Furthermore,
it cannot consist of a single point by Lemmas~\ref{lemma:curves}
and \ref{lemma:Cheltsov-Kawamata}. Therefore, it contains at least
two, say $P_i$ and $P_j$, of the five singular points. Let
$\pi:U\to U_i$ be the Kawamata blow up at the singular point whose
image to $X$ is the point $P_j$. Then, the pencil $|-K_U|$ is the
proper transform of the pencil $|-K_X|$ and its base locus
consists of the irreducible curve $C_U$. Because $-K_U\cdot C_U<0$
and $\mathcal{M}_U\sim_{\mathbb{Q}} -nK_U$,
Theorem~\ref{theorem:main-tool} completes the proof.
\end{proof}

\section{Case $\gimel=12$, hypersurface of degree $10$ in
$\mathbb{P}(1,1,2,3,4)$.}\index{$\gimel=12$} \label{section:n-12}

The threefold $X$ is a general hypersurface of degree $10$ in
$\mathbb{P}(1,1,2,3,4)$  with $-K_{X}^{3}=\frac{5}{12}$. The
singularities of the hypersurface $X$ consist of two singular
points that are quotient singularities of type
$\frac{1}{2}(1,1,1)$,  one point $P$ that is a quotient
singularity of type $\frac{1}{3}(1,1,2)$, and one point $Q$ that
is a quotient singularity of type $\frac{1}{4}(1,1,3)$.

We have the following commutative diagram:
$$
\xymatrix{
&&&Y\ar@{->}[dl]_{\gamma_{O}}\ar@{->}[dr]^{\gamma_{P}}\ar@{->}[drrrrr]^{\eta}&&&&&&\\
&&U_{PQ}\ar@{->}[dl]_{\beta_{Q}}\ar@{->}[dr]^{\beta_{P}}&&U_{QO}\ar@{->}[dl]^{\beta_{O}}&&&&\mathbb{P}(1,1,2),\\
&U_{P}\ar@{->}[dr]_{\alpha_{P}}&&U_{Q}\ar@{->}[dl]^{\alpha_{Q}}&&&&&&\\
&&X\ar@{-->}[rrrrrruu]_{\psi}&&&&&&&}
$$
where
\begin{itemize}
\item $\psi$ is the natural projection,

\item $\alpha_{P}$ is the Kawamata blow up at the point $P$ with
weights $(1,1,2)$,

\item $\alpha_{Q}$ is the Kawamata blow up at the point  $Q$ with
weights $(1,1,3)$,

\item $\beta_{Q}$ is the Kawamata blow up with weights $(1,1,3)$
at the point whose image by the birational morphism $\alpha_P$ is
the point $Q$,

\item $\beta_{P}$ is the Kawamata blow up with weights $(1,1,2)$
at the point whose image by the birational morphism $\alpha_Q$ is
the point $P$,

 \item $\beta_{O}$ is the Kawamata blow up with weights
$(1,1,2)$ at the singular point $O$ of the variety $U_{Q}$ that is
a quotient singularity of type $\frac{1}{3}(1,1,2)$ contained in
the exceptional divisor of the birational morphism $\alpha_{Q}$,

\item $\gamma_{P}$ is the Kawamata blow up with weights $(1,1,2)$
at the point whose image by the birational morphism
$\alpha_Q\circ\beta_O$ is the point $P$,

\item $\gamma_{O}$ is the Kawamata blow up with weights $(1,1,2)$
at the singular point of the variety $U_{PQ}$ that is a quotient
singularity of type $\frac{1}{3}(1,1,2)$ contained in the
exceptional divisor of the birational morphism $\beta_{Q}$,

\item  $\eta$ is an elliptic fibration.
\end{itemize}
The hypersurface $X$ can be given by the equation
$$
w^{2}z+f_{6}(x,y,z,t)w+f_{10}(x,y,z,t)=0,
$$
where  $f_{i}$ is a quasihomogeneous polynomial of degree $i$.
Moreover, there is a commutative diagram
$$
\xymatrix{
&&&U_{Q}\ar@{->}[dll]_{\pi}\ar@{->}[rr]^{\alpha_{Q}}&&X\ar@{-->}[ddrr]^{\psi}\ar@{-->}[ddll]_{\chi}&\\
&W\ar@{->}[rrd]_{\omega}&&&&&\\
&&&\mathbb{P}(1,1,2,3)\ar@{-->}[rrrr]_{\xi}&&&&\mathbb{P}(1,1,2),&}
$$
where \begin{itemize}

\item  $\xi$ and $\chi$ are the natural projections,

\item $\pi$ is a birational morphism,

\item $\omega$ is a double cover of $\mathbb{P}(1,1,2,3)$ ramified
along a surface $R$ of degree $12$.  \end{itemize} The surface $R$
is given by the equation
$$
f_{6}(x,y,z,t)^{2}-4zf_{10}(x,y,z,t)=0\subset\mathbb{P}(1,1,2,3)\cong\mathrm{Proj}\big(\mathbb{C}[x,y,z,t]\big),
$$
which implies that the surface $R$ has exactly $20$ ordinary
double points given by the equations $z=f_{6}=f_{10}=0$. Thus, the
morphism $\pi$  contracts $20$ smooth rational curves
$C_{1},\cdots , C_{20}$ to isolated ordinary double points of the
variety $W$ that dominate the singular points of $R$ in
$\mathbb{P}(1,1,2,3)$.

\begin{proposition}
\label{proposition:n-12} The linear system $|-K_X|$ is a unique
Halphen pencil on $X$.
\end{proposition}

Suppose that $\mathcal{M}\ne|-K_{X}|$. Let us show that this
assumption leads us to a contradiction.

It follows from Corollary~\ref{corollary:Ryder-a2} and
Lemmas~\ref{lemma:smooth-points},
\ref{lemma:special-singular-points-with-zero-c} that
$$\mathbb{CS}\Big(X, \frac{1}{n}\mathcal{M}\Big)\subset\Big\{P ,
Q\Big\}.$$

\begin{lemma}
\label{lemma:n-12-centers-on-U-4}  The set $\mathbb{CS}(U_{Q},
\frac{1}{n}\mathcal{M}_{U_{Q}})$ cannot contain a curve.
\end{lemma}

\begin{proof}
Suppose that the set $\mathbb{CS}(U_{Q},
\frac{1}{n}\mathcal{M}_{U_{Q}})$ contains an irreducible curve
$Z$. Let $G$ be the exceptional  divisor of the birational
morphism $\alpha_{Q}$. Then, $G\cong\mathbb{P}(1,1,3)$ and
$Z\subset G$. Moreover, it follows from
Lemma~\ref{lemma:Cheltsov-Kawamata} that $Z$ is a curve in the
linear system $|\mathcal{O}_{\mathbb{P}(1,1,3)}(1)|$ on the
surface $G$.

We consider the surface $S^{z}$ on $X$. Let $Z'$ be the curve
$S^z_{U_{Q}}\cap G$. Then, the surface $S^{z}_{U_{Q}}$ contains
every curve $C_i$, but not the curve $Z$, and its image $S^z_{W}$
by the morphism $\pi$ is isomorphic to $\mathbb{P}(1,1,3)$. The
curve $Z'$ is smooth and $\alpha_{Q}\vert_{Z'}$ is a double cover.

Because the hypersurface $X$ is general, the surface $S^z_{U_{Q}}$
is smooth along the curves $C_{i}$, the morphism
$\pi\vert_{S^z_{U_{Q}}}$ contracts the curve $C_{i}$ to a smooth
point of $S^z_{W}$, and either the intersection $Z\cap Z'$
consists of two points or the point $Z\cap Z'$ is not contained in
$\cup_{i=1}^{20}C_{i}$.

For  general surfaces $M_{U_{Q}}$, $M_{U_{Q}}^{\prime}$ in
$\mathcal{M}_{U_{Q}}$ and a general surface $D$ in
$|-6K_{U_{Q}}|$, we have
$$
2n^{2}=D\cdot M_{U_{Q}}\cdot M_{U_{Q}}^{\prime}\geq
2\mathrm{mult}_{Z}(M_{U_{Q}}\cdot M_{U_{Q}}^{\prime})\geq
2\mathrm{mult}_{Z}(M_{U_{Q}})\mathrm{mult}_{Z}(M_{U_{Q}}^{\prime})
\geq 2n^{2},%
$$
which immediately implies that the support of the cycle
$M_{U_{Q}}\cdot M_{U_{Q}}^{\prime}$ is contained in the union of
the curve $L$ and $\cup_{i=1}^{20}C_{i}$. Hence, we have
$$
\mathcal{M}_{U_{Q}}\Big\vert_{S^z_{U_{Q}}}=\mathcal{D}+\sum_{i=1}^{20}m_{i}C_{i},
$$
where $m_{i}$ is a natural number and $\mathcal{D}$ is a pencil
without fixed components.

Let $P'$ be a point of $Z\cap Z'$ and $D_{1}$, $D_{2}$ be general
curves in $\mathcal{D}$. Then,
$$
\mathrm{mult}_{P'}(D_{1})=\mathrm{mult}_{P'}(D_{2})\geq
\left\{\aligned
&n \phantom{-m_i}\ \text{in the case when}\ P'\not\in\cup_{i=1}^{20}C_{i},\\
&n-m_{i}\ \text{in the case when}\ P'\in C_{i},\\
\endaligned\right.
$$
which implies that
$$
\frac{n^{2}}{3}-\sum_{i=1}^{20}m_{i}^{2}=D_{1}\cdot D_{2}
\geq\sum_{P'\in Z\cap Z'}\mathrm{mult}_{P'}(D_{1})\mathrm{mult}_{P'}(D_{2}).%
$$
However, it is impossible because $m_{i}^{2}+(n-m_{i})^{2}\geq
\frac{n^{2}}{4}$.
\end{proof}

The exceptional divisor $E_P$ of $\alpha_P$ contains a singular
point $P_1$ of $U_P$ that is a quotient singularity of type
$\frac{1}{2}(1,1,1)$.
\begin{lemma}
\label{lemma:n-12-centers-on-U-3} The set $\mathbb{CS}(U_{P},
\frac{1}{n}\mathcal{M}_{U_{P}})$ cannot contain the singular point
$P_1$.
\end{lemma}

\begin{proof}
Let $\sigma_P:V_P\to U_{P}$ be the Kawamata blow up at the
singular point $P_1$ and $\mathcal{P}$ be the proper transform of
the pencil $|-K_X|$ via the birational morphism
$\alpha_{P}\circ\sigma_P$. Then,
$\mathcal{M}_{V_P}\sim_{\mathbb{Q}}-nK_{V_P}$ by
Lemma~\ref{lemma:Kawamata} and
$$
\mathcal{P}\sim_{\mathbb{Q}}-K_{V_P}\sim_{\mathbb{Q}}
(\alpha_P\circ\sigma_P)^*(
-K_{X})-\frac{1}{3}\sigma_P^*(E_P)-\frac{1}{2}F_P,
$$
where $F_P$ is the exceptional divisor of $\sigma_P$. Also, it has
a unique base curve $C_{V_P}$. On a general surface
$S_{V_P}\in\mathcal{P}$, the self-intersection number
$C^2_{V_P}=-K_{V_P}^3$ is negative. Also we have
$$
\mathcal{M}_{V_P}\Big\vert_{S_{V_P}}\equiv
-nK_{V_P}\Big\vert_{S_{V_P}}\equiv nC_{V_P},
$$
which implies that $\mathcal{M}$ is the pencil $|-K_X|$ by
Theorem~\ref{theorem:main-tool}. However, we assumed that
$\mathcal{M}\ne |-K_X|$.
\end{proof}
Therefore, the set $\mathbb{CS}(U_{P},
\frac{1}{n}\mathcal{M}_{U_{P}})$ must consist of the point
$\bar{Q}$ whose image to $X$ is the point $Q$ because it is not
empty by Theorem~\ref{theorem:Noether-Fano} and it cannot contain
a curve by Lemma~\ref{lemma:curves}.

Meanwhile, the exceptional divisor of $\beta_O$ contains a
singular point $O_1$ of $U_{QO}$ that is a quotient singularity of
type $\frac{1}{2}(1,1,1)$.

\begin{lemma}
\label{lemma:n-12-centers-on-U-O-1} The set $\mathbb{CS}(U_{QO},
\frac{1}{n}\mathcal{M}_{U_{QO}})$
 cannot contain the point
$O_1$.
\end{lemma}
\begin{proof}
Let $\sigma_O: V_O\to U_{QO}$ be the Kawamata blow up at the point
$O_1$. Then, the pencil $|-K_{V_O}|$ is the proper transform the
pencil $|-K_X|$. Its base locus consists of the irreducible curve
$C_{V_O}$. Because $\mathcal{M}_{V_O}\sim_{\mathbb{Q}} -nK_{V_O}$
and $-K_{V_O}\cdot C_{V_O}<0$, we obtain from
Theorem~\ref{theorem:main-tool} that $\mathcal{M}=|-K_X|$, which
contradicts our assumption.
\end{proof}

By the previous lemmas, we can see $\mathbb{CS}(X,
\frac{1}{n}\mathcal{M})=\{P, Q\}$. Furthermore, the set
$\mathbb{CS}(U_{PQ}, \frac{1}{n}\mathcal{M}_{U_{PQ}})$ consists of
the singular point of $U_{PQ}$ contained in the exceptional
divisor $\beta_Q$. Then, the set $\mathbb{CS}(Y,
\frac{1}{n}\mathcal{M}_{Y})$ must contain the singular point
contained in the exceptional divisor $\gamma_O$. In such a case,
Lemma~\ref{lemma:n-12-centers-on-U-O-1} shows a contradiction
$\mathcal{M}=|-K_X|$.

\section{Case $\gimel=13$, hypersurface of degree $11$ in $\mathbb{P}(1,1,2,3,5)$.}\index{$\gimel=13$}%
\label{section:n-13}

The threefold $X$ is a general hypersurface of degree $11$ in
$\mathbb{P}(1,1,2,3,5)$ with $-K_X^3=\frac{11}{30}$.  It has three
singular points. One is a quotient singularity of type
$\frac{1}{2}(1,1,1)$, another is a quotient singular point $P$ of
type $\frac{1}{3}(1,1,2)$, and the other is a quotient singular
point $Q$ of type $\frac{1}{5}(1,2,3)$.

We have the following commutative diagram:
$$
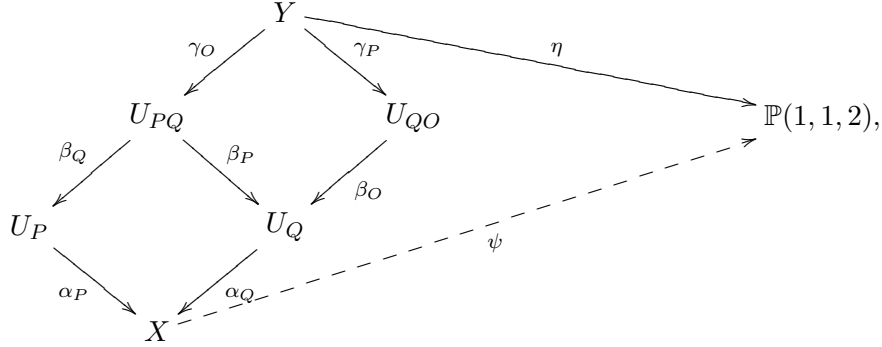
\xymatrix{
&&&Y\ar@{->}[dl]_{\gamma_{O}}\ar@{->}[dr]^{\gamma_{P}}\ar@{->}[drrrrr]^{\eta}&&&&&&\\
&&U_{PQ}\ar@{->}[dl]_{\beta_{Q}}\ar@{->}[dr]^{\beta_{P}}&&U_{QO}\ar@{->}[dl]^{\beta_{O}}&&&&\mathbb{P}(1,1,2),\\
&U_{P}\ar@{->}[dr]_{\alpha_{P}}&&U_{Q}\ar@{->}[dl]^{\alpha_{Q}}&&&&&&\\
&&X\ar@{-->}[rrrrrruu]_{\psi}&&&&&&&}
$$
where
\begin{itemize}
\item $\psi$ is the natural projection,

\item $\alpha_{P}$ is the Kawamata blow up at the point $P$ with
weights $(1,1,2)$,

\item $\alpha_{Q}$ is the Kawamata blow up at the point  $Q$ with
weights $(1,2,3)$,

\item $\beta_{Q}$ is the Kawamata blow up with weights $(1,2,3)$
at the point whose image by the birational morphism $\alpha_P$ is
the point $Q$,

\item $\beta_{P}$ is the Kawamata blow up with weights $(1,1,2)$
at the point whose image by the birational morphism $\alpha_Q$ is
the point $P$,

 \item $\beta_{O}$ is the Kawamata blow up with weights
$(1,2,1)$ at the singular point $O$ of the variety $U_{Q}$ that is
a quotient singularity of type $\frac{1}{3}(1,2,1)$ contained in
the exceptional divisor of the birational morphism $\alpha_{Q}$,

\item $\gamma_{P}$ is the Kawamata blow up with weights $(1,1,2)$
at the point whose image by the birational morphism
$\alpha_Q\circ\beta_O$ is the point $P$,

\item $\gamma_{O}$ is the Kawamata blow up with weights $(1,2,1)$
at the singular point of the variety $U_{PQ}$ that is a quotient
singularity of type $\frac{1}{3}(1,2,1)$ contained in the
exceptional divisor of the birational morphism $\beta_{Q}$,

\item  $\eta$ is an elliptic fibration.
\end{itemize}

Note that the base locus of the pencil $|-K_X|$ consists of the
irreducible curve $C$ defined by $x=y=0$.

It follows from Corollary~\ref{corollary:Ryder-a2} and
Lemmas~\ref{lemma:smooth-points},
\ref{lemma:special-singular-points-with-zero-c}  we may assume
that $\mathbb{CS}(X,\frac{1}{n}\mathcal{M})\subset \{P, Q\}$.

\begin{lemma}
\label{lemma:n-13-first} If
$\mathbb{CS}(X,\frac{1}{n}\mathcal{M})= \{P\}$, then
$\mathcal{M}=|-K_X|$.
\end{lemma}

\begin{proof}
It follows from Theorem~\ref{theorem:Noether-Fano} that the log
pair $(U_P, \mathcal{M}_{U_P})$ is not terminal at the singular
point $P_1$ contained in the exceptional divisor of $\alpha_P$.

Let $\beta_1 :W_P\to U_P$ be the Kawamata blow up at the singular
point $P_1$ with weights $(1,1,1)$. Then, the pencil $|-K_{W_P}|$
is the proper transform of $|-K_X|$. It has a unique irreducible
base curve $C_{W_P}$. We have $-K_{W_P}\cdot
C_{W_P}=-K^3_{W_P}=-\frac{3}{10}$ and
$\mathcal{M}_{W_P}\sim_{\mathbb{Q}} -nK_{W_P}$. Therefore, we
obtain $\mathcal{M}=|-K_X|$ from Theorem~\ref{theorem:main-tool}.
\end{proof}

The exceptional divisor $E\cong\mathbb{P}(1,2,3)$  of the
birational morphism $\alpha_Q$ contains two singular points $O$
and $O_2$ of types
 $\frac{1}{3}(1,2,1)$ and $\frac{1}{2}(1,1,1)$, respectively. We
 denote the unique irreducible curve in $|\mathcal{O}_{\mathbb{P}(1,2,3)}(1)|$ on
 $E$ by $L$.

\begin{lemma}
\label{lemma:n-13-second} If the set
$\mathbb{CS}(U_Q,\frac{1}{n}\mathcal{M}_{U_Q})$ contains the point
$O_2$, then $\mathcal{M}=|-K_X|$.
\end{lemma}

\begin{proof}
 Let $\beta_2:W_Q\to U_Q$ be the Kawamata blow up at the
point $O_2$ with weights $(1,1,1)$. The pencil $|-K_{W_Q}|$ is the
proper transform of $|-K_X|$. It  has exactly two irreducible base
curves $C_{W_Q}$ and $L_{W_Q}$.  On a general surface $D$ in the
pencil $|-K_{W_Q}|$, we have
$$
L_{W_Q}^2=-\frac{4}{3}, \ C_{W_Q}^2=-\frac{5}{6}, \ L_{W_Q}\cdot
C_{W_Q}=1.
$$
The surface $D$ is normal. The curves $C_{W_Q}$ and $L_{W_Q}$ form
a negative-definite intersection form on $D$. On the other hand,
$\mathcal{M}_{W_Q}|_D\equiv nC_{W_Q}+nL_{W_Q}$ by
Lemma~\ref{lemma:Kawamata}. Therefore, we obtain
$\mathcal{M}=|-K_X|$ by Theorem~\ref{theorem:main-tool}.
\end{proof}

\begin{lemma}
\label{lemma:n-13-3rd} If the set
$\mathbb{CS}(U_{QO},\frac{1}{n}\mathcal{M}_{U_{QO}})$ contains the
singular point $O_3$ contained in the exceptional of $\beta_O$,
then $\mathcal{M}=|-K_X|$.
\end{lemma}

\begin{proof} Let $\gamma_1: W_{QO}\to U_{QO}$ be the Kawamata blow up at
the point $O_3$ with weights $(1,1,1)$. Then, the pencil
$|-K_{W_{QO}}|$ is the proper transform of the pencil $|-K_X|$. It
 has exactly two irreducible base curves $C_{W_{QO}}$ and $L_{W_{QO}}$. On a
general surface $D$ in $|-K_{W_{QO}}|$, we have
$$
L_{W_{QO}}^2=-\frac{3}{2}, \ C_{W_{QO}}^2=-\frac{5}{6}, \
L_{W_{QO}}\cdot C_{W_{QO}}=1.
$$
The surface $D$ is normal. On the other hand,
$\mathcal{M}_{W_{QO}}|_D\equiv nC_{W_{QO}}+nL_{W_{QO}}$ by
Lemma~\ref{lemma:Kawamata}. Since the curves $C_{W_{QO}}$ and
$L_{W_{QO}}$ form a negative-definite intersection form on the
normal surface $D$, it follows from
Theorem~\ref{theorem:main-tool} that $\mathcal{M}=|-K_X|$.
\end{proof}

\begin{proposition}
\label{proposition:n-13} If $\gimel=13$, then
$\mathcal{M}=|-K_X|$.
\end{proposition}

\begin{proof}
Due to the previous lemmas, we may assume that
$$\mathbb{CS}\Big(X,\frac{1}{n}\mathcal{M}\Big)= \Big\{P, Q\Big\}.$$
Following the weighted blow ups $Y\to U_{QO}\to U_Q\to X$ and
using Lemmas~\ref{lemma:n-13-second} and \ref{lemma:n-13-3rd}, we
can furthermore assume that the set
$\mathbb{CS}(Y,\frac{1}{n}\mathcal{M}_{Y})$ contains the singular
point contained the exceptional divisor of the birational morphism
$\gamma_P$. The statement then follows from the proof
Lemma~\ref{lemma:n-13-first}.
\end{proof}

\section{Cases $\gimel=15$, $17$,  and $41$.}
\label{section:n-15-17-27-41-42-68}

We are to prove the following:
\begin{proposition}
If $\gimel=15$, $17$, or $41$, then the linear system $|-K_X|$ is
a unique Halphen pencil on $X$.
\end{proposition}

We consider the case of $\gimel=15$.\index{$\gimel=15$} Let $X$ be
the hypersurface given by a general quasihomogeneous equation of
degree $12$ in $\mathbb{P}(1,1,2,3,6)$ with $-K_X^3=\frac{1}{3}$.
Then, the singularities of $X$ consist of two singular points $P$
and $Q$ that are quotient singularities of type
$\frac{1}{3}(1,1,2)$ and two points of type $\frac{1}{2}(1,1,1)$.
We have a commutative diagram
$$
\xymatrix{
&&V\ar@{->}[dl]_{\sigma_Q}\ar@{->}[dr]^{\sigma_P}\ar@{->}[drrrrrrr]^{\eta}&&&&&&&\\
&Y_P\ar@{->}[dr]_{\pi_P}&&Y_Q\ar@{->}[dl]^{\pi_Q}&&&&&&\mathbb{P}(1,1,2),\\
&&X\ar@{-->}[rrrrrrru]_{\psi}&&&&&&&}
$$
where \begin{itemize}

\item $\psi$ is the natural projection,

\item $\pi_P$ is the Kawamata blow up at the point $P$ with
weights $(1,1,2)$,

\item $\pi_Q$ is the Kawamata blow up at $Q$ with weights
$(1,1,2)$,

\item $\sigma_Q$ is the Kawamata blow up with weights $(1,1,2)$ at
the point $Q_1$ whose image by the birational morphism $\pi_P$ is
the point $Q$,

\item $\sigma_P$ is the Kawamata blow up with weights $(1,1,2)$ at
the point $P_1$ whose image by the birational morphism $\pi_Q$ is
the point $P$,

\item $\eta$ is an elliptic fibration.

\end{itemize}

The set $\mathbb{CS}(X,\frac{1}{n}\mathcal{M})$ is nonempty by
Theorem~\ref{theorem:Noether-Fano}. If it contains either a
singular point of type $\frac{1}{2}(1,1,1)$ on $X$ or a curve,
then the identity $\mathcal{M}=|-K_X|$ follows from
Lemma~\ref{lemma:special-singular-points-with-zero-c} and
Corollary~\ref{corollary:Ryder-a2}, respectively. Furthermore, we
may assume that
$$\mathbb{CS}\Big(X,\frac{1}{n}\mathcal{M}\Big)\subset \Big\{P, Q\Big\}$$
due to Lemma~\ref{lemma:smooth-points}.
 Suppose
that the set $\mathbb{CS}(X,\frac{1}{n}\mathcal{M})$ contains the
point $P$.   The exceptional divisor $E_P$ of $\pi_P$ contains a
singular point $O$ of type $\frac{1}{2}(1,1,1)$.

\begin{lemma}\label{lemma:n-15-first}
If the set $\mathbb{CS}(Y_P,\frac{1}{n}\mathcal{M}_{Y_P})$
contains the singular point $O$, then $\mathcal{M}=|-K_X|$.
\end{lemma}

\begin{proof}
Let $\alpha:W\to Y_P$ be the Kawamata blow up at the point $O$
with weights $(1,1,1)$. The linear system $|-K_W|$ is the proper
transform of the linear system $|-K_X|$. It has a single base
curve $C_W$. On a general surface $D_W$ in $|-K_W|$, we have
$C_W^2=-\frac{1}{3}$. Note that the surface $D_W$ is normal.
Meanwhile, we have
\[\mathcal{M}_W\Big\vert_{D_W}\equiv -nK_W\Big\vert_{D_W}\equiv nC_W.\] Therefore, it follows
 from Theorem~\ref{theorem:main-tool} that the pencil
$\mathcal{M}$ coincides with $|-K_X|$.
\end{proof}

Therefore, we may assume that the set
$\mathbb{CS}(Y_P,\frac{1}{n}\mathcal{M}_{Y_P})$ consists of the
single point $Q_1$.  Then, we may assume that the set
$\mathbb{CS}(V,\frac{1}{n}\mathcal{M}_V)$ contains the singular
point of type $\frac{1}{2}(1,1,1)$ that is contained in the
exceptional divisor of $\sigma_Q$. We can then show
$\mathcal{M}=|-K_X|$ as in the proof of
Lemma~\ref{lemma:n-15-first}.

With the exactly same way as above, we can show that
$\mathcal{M}=|-K_X|$ if the set
$\mathbb{CS}(X,\frac{1}{n}\mathcal{M})$ contains the point $Q$.

Now, we consider the case of $\gimel=41$.\index{$\gimel=41$} Let
$X$ be the hypersurface given by a general quasihomogeneous
equation of degree $20$ in $\mathbb{P}(1,1,4,5,10)$
 with
$-K_X^3=\frac{1}{10}$. Then, it has two singular points $P$ and
$Q$ that are quotient singularities of type $\frac{1}{5}(1,1,4)$
and two singular points of type $\frac{1}{2}(1,1,1)$.

We have a commutative diagram
$$
\xymatrix{
&&V\ar@{->}[dl]_{\sigma_Q}\ar@{->}[dr]^{\sigma_P}\ar@{->}[drrrrrrr]^{\eta}&&&&&&&\\
&Y_P\ar@{->}[dr]_{\pi_P}&&Y_Q\ar@{->}[dl]^{\pi_Q}&&&&&&\mathbb{P}(1,1,4),\\
&&X\ar@{-->}[rrrrrrru]_{\psi}&&&&&&&}
$$
where \begin{itemize}

\item $\psi$ is  the natural projection,

\item $\pi_P$ is the Kawamata blow up at the point $P$ with
weights $(1,1,4)$,

\item $\pi_Q$ is the Kawamata blow up at $Q$ with weights
$(1,1,4)$,

\item $\sigma_Q$ is the Kawamata blow up with weights $(1,1,4)$ at
the point $Q_1$ whose image by the birational morphism $\pi_P$ is
the point $Q$,

\item $\sigma_P$ is the Kawamata blow up with weights $(1,1,4)$ at
the point $P_1$ whose image by the birational morphism $\pi_Q$ is
the point $P$,

\item $\eta$ is an elliptic fibration.
\end{itemize}

Using the exactly same method as in the case $\gimel=15$, one can
show that $\mathcal{M}=|-K_X|$.

For the case $\gimel=17$,\index{$\gimel=17$} let $X$ be the
hypersurface given by a general quasihomogeneous equation of
degree $12$ in $\mathbb{P}(1,1,3,4,4)$
 with
$-K_X^3=\frac{1}{4}$. Then, it has three singular points $P_1$,
$P_2$, and $P_3$ that are quotient singularities of type
$\frac{1}{4}(1,1,3)$.

In this case, we have the following commutative diagram:
$$
\xymatrix{
&&Y\ar@{->}[dl]_{\pi}\ar@{->}[dr]^{\eta}&&\\
&X\ar@{-->}[rr]_{\psi}&&\mathbb{P}(1,1,3),&}
$$
where \begin{itemize} \item $\psi$ is  the natural projection,

\item $\pi$ is the Kawamata blow up at the points $P_{1}$,
$P_{2}$, and $P_{3}$ with weights $(1,1,3)$,

\item $\eta$ is and elliptic fibration.
\end{itemize}

Even though three singular points are involved in this case, the
same method as in the previous cases can be applied to obtain
$\mathcal{M}=|-K_X|$.

\section{Case $\gimel=16$,  hypersurface of degree $12$ in $\mathbb{P}(1,1,2,4,5)$.}\index{$\gimel=16$}%
\label{section:n-16}

The threefold $X$ is a general hypersurface of degree $12$ in
$\mathbb{P}(1,1,2,4,5)$ with $-K_{X}^{3}=\frac{3}{10}$. Its
singularities
 consist of  three quotient singularities of type
$\frac{1}{2}(1,1,1)$ and one point $O$ that is a quotient
singularity of type $\frac{1}{5}(1,1,4)$.

There is a commutative diagram
$$
\xymatrix{
&U\ar@{->}[d]_{\alpha}&&W\ar@{->}[ll]_{\beta}&&Y\ar@{->}[ll]_{\gamma}\ar@{->}[d]^{\eta}&\\
&X\ar@{-->}[rrrr]_{\psi}&&&&\mathbb{P}(1,1,2)&}
$$
where\begin{itemize}

\item $\psi$ is the natural projection,

\item $\alpha$ is the Kawamata blow up  at the point $O$ with
weights $(1,1,4)$,

\item $\beta$ is the Kawamata blow up with weights $(1,1,3)$ at
the singular point of the variety $U$ that is contained in the
exceptional divisor of $\alpha$,

\item $\gamma$ is the Kawamata blow up with weights $(1,1,2)$ at
the singular point of $W$ that is contained in the exceptional
divisor of $\beta$,

\item $\eta$ is an elliptic fibration.
\end{itemize}
The hypersurface $X$ can be given by the equation
$$
w^{2}z+f_{7}(x,y,z,t)w+f_{12}(x,y,z,t)=0,
$$
where $f_{i}$ is a quasihomogeneous polynomial of degree $i$.
Moreover, there is commutative diagram
$$
\xymatrix{
&&&U\ar@{->}[dll]_{\pi}\ar@{->}[rr]^{\alpha}&&X\ar@{-->}[ddrr]^{\psi}\ar@{-->}[ddll]_{\chi}&\\
&V\ar@{->}[rrd]_{\omega}&&&&&\\
&&&\mathbb{P}(1,1,2,4)\ar@{-->}[rrrr]_{\xi}&&&&\mathbb{P}(1,1,2),&}
$$
where \begin{itemize} \item $\xi$ and $\chi$ are the natural
projections,

\item $\pi$ is a birational morphism,

\item $\omega$ is a double cover  of $\mathbb{P}(1,1,2,4)$
ramified along a surface $R$ of degree $12$.
\end{itemize}
The surface $R$ is given by the equation
$$
f_{7}(x,y,z,t)^{2}-4zf_{12}(x,y,z,t)=0\subset\mathbb{P}(1,1,2,4)\cong\mathrm{Proj}\big(\mathbb{C}[x,y,z,t]\big),
$$
which implies that $R$ has $21$ isolated ordinary double points,
given by the equations $z=f_{7}=f_{12}=0$. The morphism $\pi$
contracts $21$ smooth rational curves $C_{1}, C_{2}, \cdots ,
C_{21}$ to isolated ordinary double points of $V$ which dominate
the singular points of $R$.

\begin{proposition}
\label{proposition:n-16} The linear system $|-K_X|$ is a unique
Halphen pencil on $X$.
\end{proposition}

To prove Proposition~\ref{proposition:n-16}, due to
Corollary~\ref{corollary:Ryder-a2} and
Lemmas~\ref{lemma:smooth-points} and
\ref{lemma:special-singular-points-with-zero-c},  we may assume
that
$$\mathbb{CS}\Big(X,\frac{1}{n}\mathcal{M}\Big)=\Big\{O\Big\}.$$

Let $E$ be the exceptional divisor of the birational morphism
$\alpha$. It contains one singular point $P$ of $U$ that is a
quotient singularity of type $\frac{1}{4}(1,1,3)$. The surface $E$
is isomorphic to $\mathbb{P}(1,1,4)$. The set $\mathbb{CS}(U,
\frac{1}{n}\mathcal{M}_{U})$ contains the point $P$ by
Theorem~\ref{theorem:Noether-Fano} and
Lemma~\ref{lemma:Cheltsov-Kawamata}. Furthermore, the following
shows it consists of the point $P$.

\begin{lemma}
\label{lemma:n-16-point-P5} The set $\mathbb{CS}(U,
\frac{1}{n}\mathcal{M}_{U})$ cannot contain a curve.
\end{lemma}

\begin{proof}
Suppose that the set $\mathbb{CS}(U, \frac{1}{n}\mathcal{M}_{U})$
contains a curve $Z$. Then, $Z$ is contained in the surface $E$.
Furthermore, it follows from Lemma~\ref{lemma:Cheltsov-Kawamata}
that $Z$ is a curve in the linear system
$|\mathcal{O}_{\mathbb{P}(1,1,4)}(1)|$. Therefore, for a general
surface $M$ in $\mathcal{M}$,  we have
$$
\mathrm{Supp}\Big(M_{U}\cdot E\Big)=Z
$$
because
$M_{U}\vert_{E}\sim_{\mathbb{Q}}|\mathcal{O}_{\mathbb{P}(1,1,4)}(n)|$
and $\mathrm{mult}_{Z}(M_{U})\geq n$.

Let $M_{U}^{\prime}$ be a general surface in $\mathcal{M}_{U}$ and
$D$ be a general surface in $|-4K_{U}|$. Then,
$$
n^{2}=D\cdot M_{U}\cdot M_{U}^{\prime}\geq
\mathrm{mult}_{Z}(M_{U}\cdot M_{U}^{\prime})\geq
\mathrm{mult}_{Z}(M_{U})\mathrm{mult}_{Z}(M_{U}^{\prime})\geq n^{2},%
$$
which implies that $\mathrm{mult}_{L}(M_{U}\cdot
M_{U}^{\prime})=n^{2}$ and
$$
\mathrm{Supp}\Big(M_{U}\cdot M_{U}^{\prime}\Big)\subset Z\cup\bigcup_{i=1}^{21}C_{i}.%
$$

We consider the surface $S^z$.  The image $S^z_V$ of $S^z_U$ to
$V$  is isomorphic to $\mathbb{P}(1,1,4)$. The surface $S^z_{U}$
does not contain the curve $Z$ due to the generality in the choice
of $X$, but it contains every curve $C_{i}$. Moreover, the surface
$S^z_{U}$ is smooth along  the curves $C_{i}$ and the morphism
$\pi\vert_{S^z_{U}}$ contracts the curve $C_{i}$ to a smooth point
of $S^z_{V}$. Hence, we have
$$
\mathcal{M}_{U}\Big\vert_{S^z_{U}}=\mathcal{D}+\sum_{i=1}^{21}m_{i}C_{i},
$$
where $m_{i}$ is a natural number and $\mathcal{D}$ is a pencil
without fixed components. Therefore, the inequality $m_{i}>0$
implies that $C_{i}\cap Z\ne\varnothing$ and there is a point $P'$
of the intersection $Z\cap S^z_{U}$ that is different from the
singular point $P$. We may assume that $m_{1}>0$. Let $D_{1}$ and
$D_{2}$ be general curves in $\mathcal{D}$. Then,
$$
\mathrm{mult}_{P'}(D_{1})=\mathrm{mult}_{P'}(D_{2})\geq
\left\{\aligned
&n\phantom{-m_{i}}\ \text{in the case when}\ P'\not\in\cup_{i=1}^{21}C_{i},\\
&n-m_{i}\ \text{in the case when}\ P'\in C_{i},\\
\endaligned\right.
$$
and the curves $D_{1}$ and $D_{2}$ pass through the point $P$
because the point $P$ is a base point of the pencil
$\mathcal{M}_{U}$. Therefore, we have
$$
\frac{n^{2}}{4}-\sum_{i=1}^{21}m_{i}^{2}=D_{1}\cdot
D_{2}>\mathrm{mult}_{P}(D_{1})\mathrm{mult}_{P}(D_{2})\geq
(n-m_{1})^{2}\geq \frac{n^{2}}{4}-m_{1}^{2},
$$
which is a contradiction.
\end{proof}

Let $F$ be the exceptional divisor of the birational morphism
$\beta$. It contains the singular point $Q$ of $W$ that is a
quotient singularity of type $\frac{1}{3}(1,1,2)$.   The set
$\mathbb{CS}(W, \frac{1}{n}\mathcal{M}_{W})$ consists of the
singular point $Q$ by Theorem~\ref{theorem:Noether-Fano},
Lemmas~\ref{lemma:Cheltsov-Kawamata} and \ref{lemma:curves}.

Let $G$ be the exceptional divisor of $\gamma$ and $Q_1$ be the
unique singular point of $G$. The set $\mathbb{CS}(Y,
\frac{1}{n}\mathcal{M}_{Y})$ must consist of the point $Q_1$ by
Theorem~\ref{theorem:Noether-Fano} and
Lemma~\ref{lemma:Cheltsov-Kawamata} because every member in
$\mathcal{M}_{Y}$ is contracted to a curve by the morphism $\eta$.

Let $\sigma:V_1\to Y$ be the  Kawamata blow up at the point $Q_1$.
Then, $\mathcal{M}_{V_1}\sim_{\mathbb{Q}}-nK_{V_1}$ by
Lemma~\ref{lemma:Kawamata}, the linear system $|-K_{V_1}|$ is the
proper transform of the pencil $|-K_{X}|$, and the base locus of
the pencil $|-K_{V_1}|$ consist of the curve $C_{V_1}$. Therefore,
the inequality $-K_{V_1}\cdot C_{V_1}<0$ implies
$\mathcal{M}=|-K_X|$ by Theorem~\ref{theorem:main-tool}.

\section{Cases $\gimel=20$ and $31$.} \label{section:n-20}

First, we consider the case $\gimel=20$.\index{$\gimel=20$} The
threefold $X$ is a general hypersurface of degree $13$ in
$\mathbb{P}(1,1,3,4,5)$
 with $-K_{X}^{3}=\frac{13}{60}$. It has three
singular points. One is  a quotient singularity $P$ of type
$\frac{1}{4}(1,1,3)$, another is  a quotient singularity $Q$ of
type $\frac{1}{5}(1,1,4)$, and the other is   a quotient
singularity $Q'$ of type $\frac{1}{3}(1,1,2)$.

There is a commutative diagram
$$
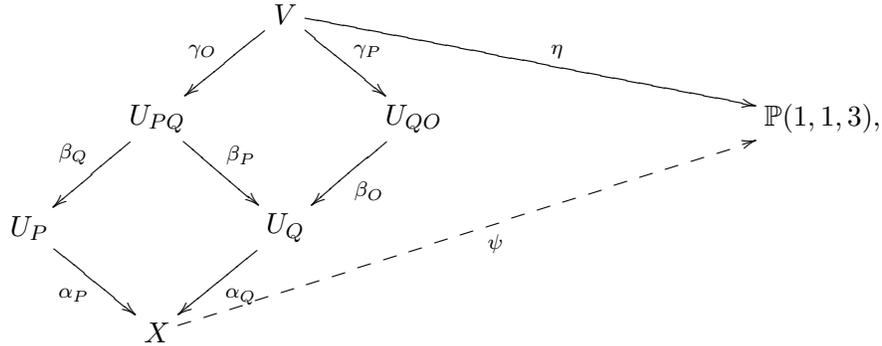
\xymatrix{
&&&V\ar@{->}[dl]_{\gamma_{O}}\ar@{->}[dr]^{\gamma_{P}}\ar@{->}[drrrrr]^{\eta}&&&&&&\\
&&U_{PQ}\ar@{->}[dl]_{\beta_{Q}}\ar@{->}[dr]^{\beta_{P}}&&U_{QO}\ar@{->}[dl]^{\beta_{O}}&&&&\mathbb{P}(1,1,3),\\
&U_{P}\ar@{->}[dr]_{\alpha_{P}}&&U_{Q}\ar@{->}[dl]^{\alpha_{Q}}&&&&&&\\
&&X\ar@{-->}[rrrrrruu]_{\psi}&&&&&&&}
$$
where \begin{itemize}

\item $\psi$ is  the natural projection,

\item $\alpha_{P}$ is the Kawamata blow up at the point $P$ with
weights $(1,1,3)$,

\item $\alpha_{Q}$ is the Kawamata blow up at the point $Q$ with
weights $(1,1,4)$,

\item $\beta_{Q}$ is the Kawamata blow up  with weights $(1,1,4)$
at the point whose image to $X$ is the point $Q$.

\item $\beta_{P}$ is the Kawamata blow up with weights $(1,1,3)$
at the point whose image to $X$ is the point $P$

\item $\beta_{O}$ is the Kawamata blow up with weights $(1,1,3)$
at the singular point $O$ of $U_{Q}$ contained in the exceptional
divisor of $\alpha_{Q}$,

\item $\gamma_{P}$ is the Kawamata blow up with weights $(1,1,3)$
at the point whose image to $X$ is the point $P$

\item $\gamma_{O}$ is the Kawamata blow up with weights $(1,1,3)$
at the singular point  of $U_{PQ}$ contained in the exceptional
divisor of $\beta_{Q}$,

\item $\eta$ is an elliptic fibration.

\end{itemize}
There is a rational map $\xi:X\dasharrow \mathbb{P}(1,1,4)$  that
gives us another commutative diagram
$$
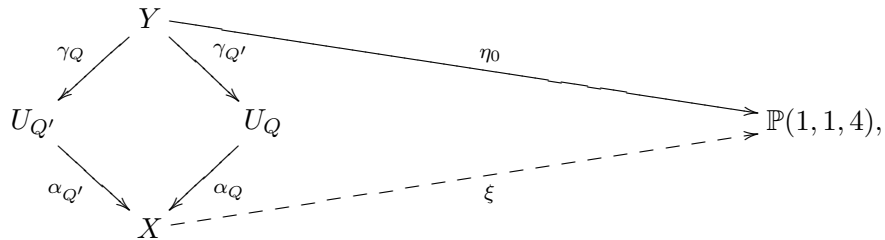
\xymatrix{
&&Y\ar@{->}[dl]_{\gamma_{Q}}\ar@{->}[dr]^{\gamma_{Q'}}\ar@{->}[drrrrrrr]^{\eta_0}&&&&&&&\\
&U_{Q'}\ar@{->}[dr]_{\alpha_{Q'}}&&U_{Q}\ar@{->}[dl]^{\alpha_{Q}}&&&&&&\mathbb{P}(1,1,4),\\
&&X\ar@{-->}[rrrrrrru]_{\xi}&&&&&&&}
$$
where \begin{itemize}

\item $\alpha_{Q'}$ is the Kawamata blow up at the point $Q'$ with
weights $(1,1,2)$,

\item $\alpha_{Q}$ is the Kawamata blow up at the point $Q$ with
weights $(1,1,4)$,

\item $\gamma_{Q}$ is the Kawamata blow up with weights $(1,1,4)$
at the point whose image to $X$ is the point $Q$,

\item $\gamma_{Q'}$ is the Kawamata blow up with weights $(1,1,2)$
at the point whose image to $X$ is the point $Q'$,

\item $\eta_0$ is an elliptic fibration.

\end{itemize}

It follows from Lemma~\ref{lemma:smooth-points} that the set
$\mathbb{CS}(X, \frac{1}{n}\mathcal{M})$ does not contain smooth
points of the hypersurface $X$. Moreover,
Corollary~\ref{corollary:Ryder-a2} implies that
$\mathcal{M}=|-K_{X}|$ in the case when the set  $\mathbb{CS}(X,
\frac{1}{n}\mathcal{M})$ contains a curve. We assume that the set
$\mathbb{CS}(X, \frac{1}{n}\mathcal{M})$ does not contain a curve,
which implies that $\mathcal{M}\ne |-K_{X}|$. Then,
$\mathbb{CS}(X, \frac{1}{n}\mathcal{M})\subseteq\{P,Q,Q'\}$.

\begin{lemma}
\label{lemma:n-20-points-P1-P3} The set $\mathbb{CS}(X,
\frac{1}{n}\mathcal{M})$ does not contain both the points $Q$ and
$Q'$.
\end{lemma}

\begin{proof}
Suppose that the set $\mathbb{CS}(X, \frac{1}{n}\mathcal{M})$
contains both the points $Q$ and $Q'$. Then,
$\mathcal{M}_{Y}\sim_{\mathbb{Q}} -nK_{Y}$, which implies that
every member in the pencil $\mathcal{M}_{Y}$ is contracted to a
curve by the elliptic fibration $\eta_0$.

Let $\bar{P}_{1}$, $\bar{P}_{2}$, $\bar{P}_{3}$ be the singular
points of $Y$ whose image to $X$ are the points $P_{1}$, $P_{2}$,
$P_{3}$, respectively. Then,
$$
\mathbb{CS}\Big(Y, \frac{1}{n}\mathcal{M}_{Y}\Big)\cap\Big\{\bar{P}_{1}, \bar{P}_{2}, \bar{P}_{3}\Big\}\ne\varnothing%
$$
by Theorem~\ref{theorem:Noether-Fano} and
Lemma~\ref{lemma:Cheltsov-Kawamata}.

Let $\pi_Q:W_Q\to Y$ be the Kawamata blow up of the point
$\bar{P}_{i}$ that is contained in the set $\mathbb{CS}(Y,
\frac{1}{n}\mathcal{M}_{Y})$ and $D$ be a general surface in
$|-K_{W_Q}|$. Then, $|-K_{W_Q}|$ is the proper transform of the
pencil $|-K_{X}|$ and the base locus of the pencil $|-K_{W_Q}|$
consists of the irreducible curve $C_{W_Q}$. We can easily check
$D\cdot C_{W_Q}<0$. Hence, we obtain $\mathcal{M}=|-K_{X}|$ by
Theorem~\ref{theorem:main-tool} because
$\mathcal{M}_{W_Q}\sim_{\mathbb{Q}} nD$ by
Lemma~\ref{lemma:Kawamata}. But it is impossible by our
assumption.
\end{proof}

\begin{lemma}
\label{lemma:n-20-points-P1-P2} The set $\mathbb{CS}(X,
\frac{1}{n}\mathcal{M})$ does not contain both the points $P$ and
$Q'$.
\end{lemma}

\begin{proof}
Suppose that the set $\mathbb{CS}(X, \frac{1}{n}\mathcal{M})$
contains both the points $P$ and $Q'$. Let $\pi_P:W_P\to U_P$ be
the Kawamata blow up at the point whose image to $X$ is the point
$Q'$. Then, $|-K_{W_P}|$ is the proper transform of the pencil
$|-K_{X}|$, the base locus of the pencil $|-K_{W_P}|$ consists of
the irreducible curve $C_{W_P}$,
 and
$D\cdot C_{W_P}=-\frac{1}{30}$. Hence,  $\mathcal{M}=|-K_{X}|$ by
Theorem~\ref{theorem:main-tool} because
$\mathcal{M}_{W_P}\sim_{\mathbb{Q}} nD$ by
Lemma~\ref{lemma:Kawamata}. But it is impossible by our
assumption.
\end{proof}

\begin{lemma}
\label{lemma:n-20-single-point} The set $\mathbb{CS}(X,
\frac{1}{n}\mathcal{M})$ cannot consist of a single point.
\end{lemma}

\begin{proof}
Suppose that the set $\mathbb{CS}(X, \frac{1}{n}\mathcal{M})$
consists of the point $P$.  Then,
$\mathcal{M}_{U_P}\sim_{\mathbb{Q}} -nK_{U_{P}}$ by
Lemma~\ref{lemma:Kawamata}. Hence, it follows from
Theorem~\ref{theorem:Noether-Fano},
Lemmas~\ref{lemma:Cheltsov-Kawamata} and \ref{lemma:curves} that
the set $\mathbb{CS}(U_{P}, \frac{1}{n}\mathcal{M}_{U_P})$
consists of the singular point $P_1$ of $U_{P}$ contained in the
exceptional divisor of $\alpha_P$.

Let $\sigma_P:V_P\to U_{P}$ be the Kawamata blow up at the
singular point $P_1$. Then, $\mathcal{M}_{V_P}\sim_{\mathbb{Q}}
-nK_{V_P}$ by Lemma~\ref{lemma:Kawamata}. Let $D$ be a general
surface of the pencil $|-K_{V_P}|$.  The proper transform
$C_{V_P}$ is   the unique base curve of the pencil $|-K_{V_P}|$
and $D\cdot C_{V_P}<0$. Hence, we have $\mathcal{M}=|-K_{X}|$ by
Theorem~\ref{theorem:main-tool} because
$\mathcal{M}_{V_P}\sim_{\mathbb{Q}} nD$ by
Lemma~\ref{lemma:Kawamata}. But it is impossible by our
assumption.

In the exactly same way, we can show that the set $\mathbb{CS}(X,
\frac{1}{n}\mathcal{M})$ cannot consist of the point $Q'$.

Suppose that the set $\mathbb{CS}(X, \frac{1}{n}\mathcal{M})$
consists of the point $Q$.  Then,
$\mathcal{M}_{U_Q}\sim_{\mathbb{Q}} -nK_{U_{Q}}$ by
Lemma~\ref{lemma:Kawamata}. It follows from
Theorem~\ref{theorem:Noether-Fano},
Lemmas~\ref{lemma:Cheltsov-Kawamata}, and \ref{lemma:curves} that
the set $\mathbb{CS}(U_{Q}, \frac{1}{n}\mathcal{M}_{U_Q})$
consists of the singular point $O$. The divisor $-K_{U_{QO}}$ is
nef and big, and hence it follows from
Theorem~\ref{theorem:Noether-Fano} and
Lemma~\ref{lemma:Cheltsov-Kawamata} that the set
$\mathbb{CS}(U_{QO}, \frac{1}{n}\mathcal{M}_{U_{QO}})$ contains
the singular point $O_1$ of the variety $U_{QO}$ contained in the
exceptional divisor of the $\beta_O$.

Let $\sigma_O:V_O\to U_{QO}$ be the Kawamata blow up at the
singular point $O_1$. Then, $\mathcal{M}_{V_O}\sim_{\mathbb{Q}}
-nK_{V_O}$ by Lemma~\ref{lemma:Kawamata}. Let $H$ be a general
surface of the pencil $|-K_{V_O}|$. The proper transform $C_{V_O}$
is the unique base curve of the pencil $|-K_{V_O}|$ and $H\cdot
C_{V_O}=-\frac{1}{24}$, which is impossible by
Theorem~\ref{theorem:main-tool} because $\mathcal{M}\ne|-K_{X}|$.
\end{proof}

Consequently, we see that $\mathbb{CS}(X,
\frac{1}{n}\mathcal{M})=\{P, Q\}$. Then,
$\mathcal{M}_{U_{PQ}}\sim_{\mathbb{Q}} -nK_{U_{PQ}}$ by
Lemma~\ref{lemma:Kawamata} and it follows from
Theorem~\ref{theorem:Noether-Fano} and
Lemma~\ref{lemma:Cheltsov-Kawamata} that the set
$\mathbb{CS}(U_{PQ}, \frac{1}{n}\mathcal{M}_{PQ})$ contains either
the singular point $Q_2$ of the variety $U_{PQ}$ contained in the
exceptional divisor of the birational morphism $\beta_{Q}$ or the
singular point $P_2$ of the variety $U_{PQ}$  contained in the
exceptional divisor of $\beta_{P}$.

\begin{lemma}
\label{lemma:n-20-exclusion-of-P6} The set $\mathbb{CS}(U_{PQ},
\frac{1}{n}\mathcal{M}_{U_{PQ}})$ does not contain the point
$P_2$.
\end{lemma}

\begin{proof}
Note that the point $P_2$ is a quotient singularity of type
$\frac{1}{3}(1,1,2)$ on $U_{PQ}$.

Suppose that the set $\mathbb{CS}(U_{PQ},
\frac{1}{n}\mathcal{M}_{U_{PQ}})$ contains the point $P_2$. Let
$\pi:V_{PQ}\to U_{PQ}$ be the Kawamata blow up at the point $P_2$
and $D$ be a general surface of the pencil $|-K_{V_{PQ}}|$.

Then, the proper transform $C_{V_{PQ}}$ is the unique  base curve
of the pencil $|-K_{V_{PQ}}|$. Because
$\mathcal{M}_{V_{PQ}}\sim_{\mathbb{Q}} -nK_{W}$ by
Lemma~\ref{lemma:Kawamata} and  $D\cdot C_{V_{PQ}}=-\frac{1}{24}$,
it follows from Theorem~\ref{theorem:main-tool} that
$\mathcal{M}=|-K_{X}|$.  But it is impossible by our assumption.
\end{proof}

Therefore, it follows from Lemma~\ref{lemma:curves} that the set
$\mathbb{CS}(U_{PQ}, \frac{1}{n}\mathcal{M}_{U_{PQ}})$ consists of
the singular point $Q_2$.  In particular, we have
$\mathcal{M}_{V}\sim_{\mathbb{Q}} -nK_{V}$, which implies that
each surface  in the pencil $\mathcal{M}_{V}$ is contracted to a
curve by the elliptic fibration $\eta_0$.

Let $Q_3$ be the singular point of $V$ that is contained in the
exceptional divisor of $\gamma_{O}$. It is a quotient singularity
of type $\frac{1}{3}(1,1,2)$. Theorem~\ref{theorem:Noether-Fano}
and Lemma~\ref{lemma:Cheltsov-Kawamata} imply that the set
$\mathbb{CS}(V, \frac{1}{n}\mathcal{M}_{V})$ contains the point
$Q_3$.

Let $\pi:W\to V$ be the Kawamata blow up at the point $Q_3$ and
$D$ be a general surface in $|-K_{W}|$. Then, $|-K_{W}|$ is the
proper transform of the pencil $|-K_{X}|$ and  the base locus of
the pencil $|-K_{W}|$ consists of the irreducible curve $C_W$.
Then, $\mathcal{M}=|-K_{X}|$ by Theorem~\ref{theorem:main-tool}
since $\mathcal{M}_{W}\sim_{\mathbb{Q}} nD$ by
Lemma~\ref{lemma:Kawamata} and $D\cdot C_W<0$. The obtained
contradiction concludes the following:
\begin{proposition}
\label{proposition:n-20} If $\gimel=20$, then the linear system
$|-K_X|$ is the only Halphen pencil on $X$.
\end{proposition}

From now, we consider the case $\gimel=31$.\index{$\gimel=31$} The
threefold $X$ is a hypersurface of degree $16$  in
$\mathbb{P}(1,1,4,5,6)$ with $-K_{X}^{3}=\frac{2}{15}$. The
singularities of $X$ consist of one quotient singularity of type
$\frac{1}{2}(1,1,1)$, one singular point $P$ that is a quotient
singularity of type $\frac{1}{5}(1,1,4)$, and one singular point
$Q$ that is a quotient singularity of type $\frac{1}{6}(1,1,5)$.

There is a commutative diagram
$$
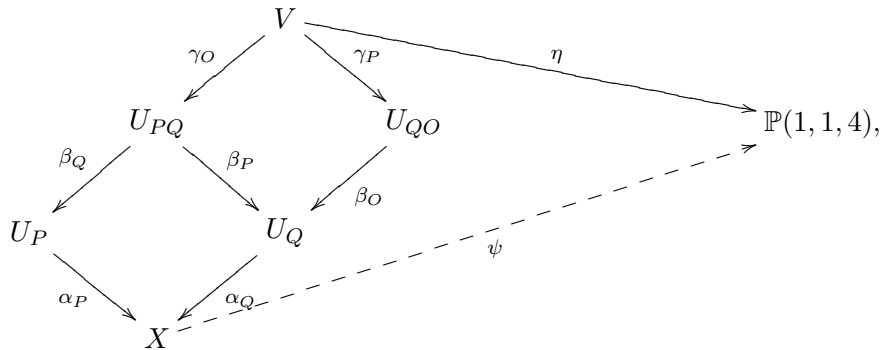
\xymatrix{
&&&V\ar@{->}[dl]_{\gamma_{O}}\ar@{->}[dr]^{\gamma_{P}}\ar@{->}[drrrrr]^{\eta}&&&&&&\\
&&U_{PQ}\ar@{->}[dl]_{\beta_{Q}}\ar@{->}[dr]^{\beta_{P}}&&U_{QO}\ar@{->}[dl]^{\beta_{O}}&&&&\mathbb{P}(1,1,4),\\
&U_{P}\ar@{->}[dr]_{\alpha_{P}}&&U_{Q}\ar@{->}[dl]^{\alpha_{Q}}&&&&&&\\
&&X\ar@{-->}[rrrrrruu]_{\psi}&&&&&&&}
$$
where \begin{itemize}

\item $\psi$ is  the natural projection,

\item $\alpha_{P}$ is the Kawamata blow up at the point $P$ with
weights $(1,1,4)$,

\item $\alpha_{Q}$ is the Kawamata blow up at the point $Q$ with
weights $(1,1,5)$,

\item $\beta_{Q}$ is the Kawamata blow up  with weights $(1,1,5)$
at the point whose image to $X$ is the point $Q$.

\item $\beta_{P}$ is the Kawamata blow up with weights $(1,1,4)$
at the point whose image to $X$ is the point $P$

\item $\beta_{O}$ is the Kawamata blow up with weights $(1,1,4)$
at the singular point $O$ of $U_{Q}$ contained in the exceptional
divisor of $\alpha_{Q}$,

\item $\gamma_{P}$ is the Kawamata blow up with weights $(1,1,4)$
at the point whose image to $X$ is the point $P$

\item $\gamma_{O}$ is the Kawamata blow up with weights $(1,1,4)$
at the singular point  of $U_{PQ}$ contained in the exceptional
divisor of $\beta_{Q}$,

\item $\eta$ is an elliptic fibration.

\end{itemize}
The hypersurface $X$ can be given by the equation
$$
t^{2}w+tf_{11}(x,y,z,w)+f_{16}(x,y,z,w)=0,
$$
where $f_{i}(x,y,z,w)$ is a general quasihomogeneous polynomial of
degree $i$. Consider the linear system on $X$ defined by the
equations
\[\mu w+\sum_{i=0}^{6}\lambda_ix^iy^{6-i}=0,\]
where $(\mu: \lambda_0: \lambda_1:\lambda_2:\lambda_3:
\lambda_4:\lambda_5:\lambda_6)\in\mathbb{P}^7$. It gives us a
dominant rational map $\xi :X\dasharrow\mathbb{P}(1,1,6)$ defined
in the outside of the point $P$. The normalization of a general
fiber is an elliptic curve. Therefore, we have another elliptic
fibration as follows:
$$
\xymatrix{
&U_P\ar@{->}[d]_{\alpha_P}&&Y\ar@{->}[ll]_{\beta}\ar@{->}[d]^{\eta_0}&\\%
&X\ar@{-->}[rr]_{\xi}&&\mathbb{P}(1,1,6),&}
$$
where \begin{itemize}

\item $\alpha_P$ is the Kawamata blow up at the point $P$ with
weights $(1,1,4)$,

\item $\beta$ is the Kawamata blow up with weights $(1,1,3)$ at
the singular point of the variety $U_P$ that is a quotient
singularity of type $\frac{1}{4}(1,1,3)$,

\item $\eta_0$ is an elliptic fibration.
\end{itemize}

\begin{proposition}
\label{proposition:n-31} If $\gimel=31$, then the linear system
$|-K_X|$ is the only Halphen pencil on $X$.
\end{proposition}

If the set $\mathbb{CS}(X, \frac{1}{n}\mathcal{M})$ contains
either the singular point of type $\frac{1}{2}(1,1,1)$ or a curve,
we obtain the identity $\mathcal{M}=|-K_X|$ from
Lemma~\ref{lemma:special-singular-points-with-zero-c} and
Corollary~\ref{corollary:Ryder-a2}. Therefore, due to
Lemma~\ref{lemma:smooth-points}, we may assume that
 $\mathbb{CS}(X,
\frac{1}{n}\mathcal{M})\subseteq\{P, Q\}$.

The proof of Proposition~\ref{proposition:n-31} is similar to that
of Proposition~\ref{proposition:n-20}. In the case $\gimel=31$, we
do not need Lemmas~\ref{lemma:n-20-points-P1-P3} and
~\ref{lemma:n-20-points-P1-P2}. The only different part is
Lemma~\ref{lemma:n-20-single-point}. However, it works for the
point $Q$ in the case $\gimel=31$ as well. Thus, the following
lemma will complete the proof.

\begin{lemma}
If the set $\mathbb{CS}(X, \frac{1}{n}\mathcal{M})$  consists of
the point $P$, then $\mathcal{M}=|-K_X|$.
\end{lemma}
\begin{proof}
Suppose that  the set $\mathbb{CS}(X, \frac{1}{n}\mathcal{M})$
consists of the point $P$. Then, the set $\mathbb{CS}(U_P,
\frac{1}{n}\mathcal{M}_{U_P})$ consists of the singular point
$P_1$ of $U_P$ contained in the exceptional divisor of $\alpha_P$
because of Lemmas~\ref{lemma:curves} and
\ref{lemma:Cheltsov-Kawamata}. Furthermore, the set
$\mathbb{CS}(Y, \frac{1}{n}\mathcal{M}_{Y})$ must contain the
singular point contained in the exceptional divisor of $\beta$ by
Theorem~\ref{theorem:Noether-Fano}. Consider the Kawamata blow up
at this point and apply the same method for
Lemma~\ref{lemma:n-20-single-point} to get $\mathcal{M}=|-K_X|$.
\end{proof}

\section{Cases $\gimel=21$, $35$, and $71$.}\index{$\gimel=21$}\index{$\gimel=35$}\index{$\gimel=71$} \label{section:n-35}

Suppose that $\gimel\in\{21, 35, 71 \}$. Then, the threefold $X
\subset \mathbb{P}(1,a_1,a_2,a_3,a_4)$ always contains the point
$O=(0:0:0:1:0)$. It  is a singular point of $X$ that is a quotient
singularity of type $\frac{1}{a_3}(1,1, a_4-a_3)$.

We also have a commutative diagram as follows:
$$
\xymatrix{
&U\ar@{->}[d]_{\alpha}&&W\ar@{->}[ll]_{\beta}\ar@{->}[d]^{\eta}&\\%
&X\ar@{-->}[rr]_{\psi}&&\mathbb{P}(1,1,a_2),&}
$$
where \begin{itemize} \item $\alpha$ is the Kawamata blow up at
the point $O$ with weights $(1,1,a_4-a_3)$,

\item $\beta$ is the Kawamata blow up with weights
$(1,1,a_4-a_3-1)$ at the point $P$ of $U$ that is a quotient
singularity of type $\frac{1}{a_4-a_3}(1,1,a_4-a_3-1)$,

\item  $\eta$ is an elliptic fibration.
\end{itemize}

By Lemma~\ref{lemma:special-singular-points-with-zero-c}, if the
set $\mathbb{CS}(X, \frac{1}{n}\mathcal{M})$ contains a singular
point of $X$ different from the singular point $O$, then the
identity $\mathcal{M}=|-K_X|$ holds. Moreover, if it contains a
curve, then Corollary~\ref{corollary:Ryder-a2} implies
$\mathcal{M}=|-K_X|$. Therefore, we may assume that
$$
\mathbb{CS}\Big(X, \frac{1}{n}\mathcal{M}\Big)=\Big\{O\Big\}.%
$$
due to Theorem~\ref{theorem:Noether-Fano} and
Lemma~\ref{lemma:smooth-points}

The exceptional divisor $E\cong\mathbb{P}(1,1,a_4-a_3)$ of
$\alpha$ contains one singular point $P$ that is a quotient
singularity of type $\frac{1}{a_4-a_3}(1,1, a_4-a_3-1)$.

\begin{lemma}
\label{lemma:n-35-no-ruling} The set $\mathbb{CS}(U,
\frac{1}{n}\mathcal{M}_U)$ cannot contain a curve.
\end{lemma}

\begin{proof}
Suppose that the set $\mathbb{CS}(U, \frac{1}{n}\mathcal{M}_U)$
contains a curve $Z$. Then, Lemma~\ref{lemma:Cheltsov-Kawamata}
implies that $Z\in|\mathcal{O}_{\mathbb{P}(1,1,a_4-a_3)}(1)|$,
which contradicts Lemma~\ref{lemma:curves}.
\end{proof}

\begin{proposition}
\label{proposition:n-35} The linear system $|-K_{X}|$ is a unique
Halphen pencil on $X$.
\end{proposition}
\begin{proof}
It follows from Corollary~\ref{corollary:Ryder-a1} and
Lemmas~\ref{lemma:smooth-points},
\ref{lemma:special-singular-points-with-zero-c} that  we may
assume that $\mathbb{CS}(X, \frac{1}{n}\mathcal{M})=\{O\}$. By
Lemma~\ref{lemma:n-35-no-ruling}, the set $\mathbb{CS}(U,
\frac{1}{n}\mathcal{M}_U)=\{P\}.$

The exceptional divisor $F$ of $\beta$ contains one singular point
$Q$ that is a quotient singularity of type
$\frac{1}{a_4-2a_3}(1,1, a_4-a_3-2)$. Because the set
$\mathbb{CS}(W, \frac{1}{n}\mathcal{M}_W)$ is not empty by
Theorem~\ref{theorem:Noether-Fano}, it must contain the point $Q$.

Let $\pi:Y\to W$ be the  Kawamata blow up at the point $Q$ with
weights $(1,1,a_4-a_3-2)$. Easy calculations show that the linear
system $|-K_{Y}|$ is the proper transform of the pencil $|-K_{X}|$
and the base locus of the pencil $|-K_{Y}|$ consists of the
irreducible curve $C_Y$ whose image to $X$ is the base curve of
the pencil $|-K_{X}|$. Also, we  can easily get
$$-K_{Y}\cdot C_Y=-K_{Y}^{3}=-\frac{1}{(a_4-a_3-1)(a_4-a_3-2)}<0.$$
The divisor
$(-K_X)^3(-K_{Y})+(-K_Y)^3(\alpha\circ\beta\circ\pi)^{*}(-K_{X})$
is nef and big. For a general surface $M$ in $\mathcal{M}_Y$ and a
general surface $D$ in $|-K_{Y}|$,
$$
\Big((-K_X)^3(-K_{Y})+(-K_Y)^3(\alpha\circ\beta\circ\pi)^{*}(-K_{X})\Big)
\cdot D\cdot M=0,%
$$
which implies that  $\mathcal{M}_Y=|-K_{Y}|$ by
Theorem~\ref{theorem:main-tool}. Therefore, we obtain the identity
$\mathcal{M}=|-K_{X}|$.
\end{proof}

\section{Cases $\gimel=24$ and $46$.} \label{section:n-46}

We first consider the case $\gimel=46$.\index{$\gimel=46$} The
threefold  $X$ is a general hypersurface of degree $21$ in
$\mathbb{P}(1,1,3,7,10)$  with $-K_{X}^{3}=\frac{1}{10}$. It has
only one singular point $P$ that is a quotient singularity of type
$\frac{1}{10}(1,3,7)$.

There is a commutative diagram
$$
\xymatrix{
&U\ar@{->}[d]_{\alpha}&&W\ar@{->}[ll]_{\beta}&&Y\ar@{->}[ll]_{\gamma}\ar@{->}[d]^{\eta}&\\%
&X\ar@{-->}[rrrr]_{\psi}&&&&\mathbb{P}(1,1,3),&}
$$
where \begin{itemize}

\item $\psi$ is the natural projection,

 \item $\alpha$ is the
Kawamata blow up at the point $P$ with weights $(1,3,7)$,

\item $\beta$ is the Kawamata blow up  with weights $(1,3,4)$ at
the singular point of $U$ that is a quotient singularity of type
$\frac{1}{7}(1,3,4)$,

\item $\gamma$ is the Kawamata blow up  with weights $(1,3,1)$ at
the singular point of the variety $W$ that is a quotient
singularity of type $\frac{1}{4}(1,3,1)$,

\item $\eta$ is an elliptic fibration.
\end{itemize}

The linear system $|-K_X|$ is a Halphen pencil. Furthermore, we
can obtain
\begin{proposition}
\label{proposition:n-46} The linear system $|-K_{X}|$ is the only
Halphen pencil on $X$.
\end{proposition}

\begin{proof}
Suppose that we have a Halphen pencil $\mathcal{M}$ different from
$|-K_{X}|$.  We are to show that this assumption leads us to a
contradiction.

Due to Lemma~\ref{lemma:smooth-points} and
Corollary~\ref{corollary:Ryder-a2}, we may assume that
$\mathbb{CS}(X, \frac{1}{n}\mathcal{M})=\{P\}$. Hence, we have
$\mathcal{M}_{U}\sim_{\mathbb{Q}} -nK_{U}$ by
Lemma~\ref{lemma:Kawamata}.

The exceptional divisor $E\cong\mathbb{P}(1,3,7)$ of $\alpha$
contains two singular points $P_1$ and $P_2$ of $U$ that are
quotient singularities of types $\frac{1}{3}(1,1,2)$ and
$\frac{1}{7}(1,3,4)$, respectively. The proof of
Lemma~\ref{lemma:n-29-second-floor} shows that the set
$\mathbb{CS}(U, \frac{1}{n}\mathcal{M}_{U})$ does not contain the
point  $P_{1}$ because $\mathcal{M}\ne |-K_{X}|$. Hence,
Lemma~\ref{lemma:Cheltsov-Kawamata} implies that the set
$\mathbb{CS}(U, \frac{1}{n}\mathcal{M}_{U})$ consists of the
singular point $P_2$.
 Therefore, we have $\mathcal{M}_{W}\sim_{\mathbb{Q}}
-nK_{W}$ by Lemma~\ref{lemma:Kawamata}, which implies that the set
$\mathbb{CS}(W, \frac{1}{n}\mathcal{M}_{W})$ is not empty.

The exceptional divisor $F\cong\mathbb{P}(1,3,4)$ of $\beta$
contains two singular points $Q_1$ and $Q_2$ of $W$ that are
quotient singularities of types $\frac{1}{3}(1,2,1)$ and
$\frac{1}{4}(1,3,1)$, respectively. It again follows from
Lemma~\ref{lemma:Cheltsov-Kawamata} that either the set
$\mathbb{CS}(W, \frac{1}{n}\mathcal{M}_{W})$ contains the point
$Q_1$ or the set $\mathbb{CS}(W, \frac{1}{n}\mathcal{M}_{W})$
consists of the point $Q_2$.

For the convenience, let $L$ be the unique curve  contained in
$|\mathcal{O}_{\mathbb{P}(1,3,7)}(1)|$ on the surface $E$ and
$\bar{L}$ be the unique curve  contained in
$|\mathcal{O}_{\mathbb{P}(1,3,4)}(1)|$ on the surface $F$.

Let $\sigma:V_1\to W$ be the Kawamata blow up at the point $Q_1$.
The pencil $|-K_{V_1}|$ is the proper transform of the pencil
$|-K_{X}|$ and the base locus of the pencil $|-K_{V_1}|$ consists
of three irreducible curves $C_{V_1}$, $L_{V_1}$ and
$\bar{L}_{V_1}$.

A general surface $D_{V_1}$ in $|-K_{V_1}|$ is normal. Moreover,
the intersection form of the curves $C_{V_1}$, $L_{V_1}$ and
$\bar{L}_{V_1}$ is negatively definite on the surface $D_{V_1}$.
Hence, the set $\mathbb{CS}(W, \frac{1}{n}\mathcal{M}_{W})$ does
not contain the point $Q_1$ by Lemma~\ref{lemma:Kawamata} and
Theorem~\ref{theorem:main-tool} because $\mathcal{M}_{V}\ne
|-K_{V}|$. Thus, we see that the set $\mathbb{CS}(W,
\frac{1}{n}\mathcal{M}_{W})$ must consist of the point $Q_2$,
which implies that $\mathcal{M}_{Y}\sim_{\mathbb{Q}} -nK_{Y}$ by
Lemma~\ref{lemma:Kawamata}. In particular, each member in the
pencil $\mathcal{M}_{Y}$ is contracted to a curve by the elliptic fibration
$\eta$. Theorem~\ref{theorem:Noether-Fano} and
Lemma~\ref{lemma:Cheltsov-Kawamata} imply that the set
$\mathbb{CS}(Y, \frac{1}{n}\mathcal{M}_{Y})$ contains the singular
point  $Q$ of $Y$ contained in the exceptional divisor of the
birational morphism $\gamma$.

Let $\pi:V_2\to W$ be the Kawamata blow up at the point $Q$ and
$D_{V_2}$ be a general surface in $|-K_{V_2}|$. Then, $D_{V_2}$ is
normal, the pencil $|-K_{V_2}|$ is the proper transform of the
pencil $|-K_{X}|$, and the base locus of the pencil $|-K_{V_2}|$
consists of three irreducible curves $C_{V_2}$, $L_{V_2}$, and
$\bar{L}_{V_2}$. The intersection form of the curves $C_{V_2}$,
$L_{V_2}$, and $\bar{L}_{V_2}$ is negative-definite on the surface
$D_{V_2}$ because the curves $C_{V_2}$, $L_{V_2}$, and
$\bar{L}_{V_2}$ are components of a single fiber of the elliptic
fibration $\eta\circ\pi\vert_{D_{V_2}}:D_{V_2}\to\mathbb{P}^{1}$
that consists of four components. On the other hand, we have
$$
\mathcal{M}_{V_2}\Big\vert_{D_{V_2}}\equiv
n(C_{V_2}+L_{V_2}+\bar{L}_{V_2})
$$
by Lemma~\ref{lemma:Kawamata}. Thus, we obtain the identity
$\mathcal{M}_{V_2}=|-K_{V_2}|$ from
Theorem~\ref{theorem:main-tool}, which is a contradiction to our
assumption.
\end{proof}

From now, we consider the case $\gimel=24$.\index{$\gimel=24$} The
variety $X$ is a general hypersurface of degree $15$ in
$\mathbb{P}(1,1,2,5,7)$
 with $-K_{X}^{3}=\frac{3}{14}$. Its singularities consist of
one point that is a quotient singularity of type
$\frac{1}{2}(1,1,1)$ and one point $P$ that is a quotient
singularity of type $\frac{1}{7}(1,2,5)$.

There is a commutative diagram
$$
\xymatrix{
&U\ar@{->}[d]_{\alpha}&&W\ar@{->}[ll]_{\beta}&&Y\ar@{->}[ll]_{\gamma}\ar@{->}[d]^{\eta}&\\%
&X\ar@{-->}[rrrr]_{\psi}&&&&\mathbb{P}(1,1,2),&}
$$
where \begin{itemize}

\item $\psi$ is the natural projection,

 \item $\alpha$ is the
Kawamata blow up at the point $P$ with weights $(1,2,5)$,

\item $\beta$ is the Kawamata blow up  with weights $(1,2,3)$ at
the singular point of $U$ that is a quotient singularity of type
$\frac{1}{5}(1,2,3)$,

\item $\gamma$ is the Kawamata blow up  with weights $(1,2,1)$ at
the singular point of the variety $W$ that is a quotient
singularity of type $\frac{1}{3}(1,2,1)$,

\item $\eta$ is an elliptic fibration.
\end{itemize}

If the set  $\mathbb{CS}(X, \frac{1}{n}\mathcal{M})$ contains a
curve, then we obtain the identity $\mathcal{M}=|-K_X|$ from
Corollary~\ref{corollary:Ryder-a2}. Therefore, due to
Lemmas~\ref{lemma:smooth-points} and
\ref{lemma:special-singular-points-with-zero-c}, we may assume
that $\mathbb{CS}(X, \frac{1}{n}\mathcal{M})=\{ P \}$.

\begin{proposition}
The linear system $|-K_{X}|$ is the only Halphen pencil on $X$.
\end{proposition}
\begin{proof}
The proof is  the same as that of
Proposition~\ref{proposition:n-46}.
\end{proof}

\section{Case $\gimel=25$, hypersurface of degree $15$ in
$\mathbb{P}(1,1,3,4,7)$.}\index{$\gimel=25$} \label{section:n-25}

The threefold $X$ is a general  hypersurface of degree $15$ in
$\mathbb{P}(1,1,3,4,7)$ with $-K_X^3=\frac{5}{28}$. It has two
singular points. One is
 a quotient singularity $P$ of type $\frac{1}{4}(1,1,3)$ and the other is a  quotient singularity $Q$ of type $\frac{1}{7}(1,3,4)$.

 There is a
commutative diagram
$$
\xymatrix{
&&&Y\ar@{->}[dl]_{\gamma_{O}}\ar@{->}[dr]^{\gamma_{P}}\ar@{->}[drrrrr]^{\eta}&&&&&&\\
&&U_{PQ}\ar@{->}[dl]_{\beta_{Q}}\ar@{->}[dr]^{\beta_{P}}&&U_{QO}\ar@{->}[dl]^{\beta_{O}}&&&&\mathbb{P}(1,1,3),\\
&U_{P}\ar@{->}[dr]_{\alpha_{P}}&&U_{Q}\ar@{->}[dl]^{\alpha_{Q}}&&&&&&\\
&&X\ar@{-->}[rrrrrruu]_{\psi}&&&&&&&}
$$
where \begin{itemize}

\item $\psi$ is  the natural projection,

\item $\alpha_{P}$ is the Kawamata blow up at the point $P$ with
weights $(1,1,3)$,

\item $\alpha_{Q}$ is the Kawamata blow up at the point $Q$ with
weights $(1,3,4)$,

\item $\beta_{Q}$ is the Kawamata blow up with weights $(1,3,4)$
at the point whose image to $X$ is the point $Q$,

\item $\beta_{P}$ is the Kawamata blow up with weights $(1,1,3)$
at the point whose image to $X$ is the point $P$,

\item $\beta_{O}$ is the Kawamata blow up with weights $(1,3,1)$
at the singular point $O$ of type $\frac{1}{4}(1,3,1)$ contained
in the ex\-cep\-ti\-onal divisor of the birational morphism
$\alpha_{Q}$,

\item $\gamma_{P}$ is the Kawamata blow up with weights $(1,1,3)$
at the point whose image to $X$ is the point $P$,

\item $\gamma_{O}$ is the Kawamata blow up with weights $(1,3,1)$
at the singular point of type $\frac{1}{4}(1,3,1)$ contained in
the exceptional divisor of the birational morphism $\beta_{Q}$,

\item  $\eta$ is an elliptic fibration.

\end{itemize}

\begin{proposition}
\label{proposition:n-25} The linear system $|-K_{X}|$ is a unique
Halphen pencil on $X$.
\end{proposition}
In what follows, we prove Proposition~\ref{proposition:n-25}.  For
the convenience,  let $D$ be a general surface in $|-K_X|$.

It follows from \cite{CPR} that $|-K_{X}|$ is invariant under the
action of the group $\mathrm{Bir}(X)$. Therefore, we may assume
that the log pair $(X, \frac{1}{n}\mathcal{M})$ is canonical. In
fact, we can assume that
$$
\varnothing\ne\mathbb{CS}\Big(X, \frac{1}{n}\mathcal{M}\Big)\subseteq\Big\{P, Q\Big\}%
$$
by Lemma~\ref{lemma:smooth-points} and
Corollary~\ref{corollary:Ryder-a2}.

\begin{lemma}
\label{lemma:n-25-P2-is-a-center} If the point $Q$ is not
contained in $\mathbb{CS}(X, \frac{1}{n}\mathcal{M})$, then
$\mathcal{M}=|-K_{X}|$.
\end{lemma}

\begin{proof}
The set $\mathbb{CS}(U_{P}, \frac{1}{n}\mathcal{M}_{U_P})$ is not
empty by Theorem~\ref{theorem:Noether-Fano} because
$\mathcal{M}_{U_P}\sim_{\mathbb{Q}}-nK_{U_{P}}$ by
Lemma~\ref{lemma:Kawamata}. Let $P_{1}$ be the singular point of
the variety $U_{P}$ contained in the exceptional divisor of the
birational morphism $\alpha_{P}$. It is a quotient singularity of
type $\frac{1}{3}(1,1,2)$. Lemma~\ref{lemma:Cheltsov-Kawamata}
implies that the set $\mathbb{CS}(U_{P},
\frac{1}{n}\mathcal{M}_{U_P})$ contains the point $P_{1}$.

Let $\pi_P:W_P\to U_{P}$ be the Kawamata blow up at the point
$P_{1}$ with weights $(1,1,2)$. We can easily check that
$|-K_{W_P}|$ is the proper transform of the pencil $|-K_X|$ and
the base locus of the pencil $|-K_{W}|$ consists of the
irreducible curve $C_{W_P}$. We can also see $D_{W_P}\cdot
C_{W_P}=-K_{W_P}^3=-\frac{1}{14}<0$. Hence,
Theorem~\ref{theorem:main-tool} implies the identity
$\mathcal{M}=|-K_{X}|$  because
$\mathcal{M}_{W_P}\sim_{\mathbb{Q}} nD_{W_P}$ by
Lemma~\ref{lemma:Kawamata}.
\end{proof}

 The exceptional divisor $E\cong\mathbb{P}(1,3,4)$ of the
 birational morphism
$\alpha_{Q}$ contains two singular points $O$ and $Q_1$ that are
quotient singularities of types $\frac{1}{4}(1,3,1)$ and
$\frac{1}{3}(1,2,1)$. Let $L$ be the unique curve of the linear
system $|\mathcal{O}_{\mathbb{P}(1,\,3,\,4)}(1)|$ on the surface
$E$.

Due to Lemma~\ref{lemma:n-25-P2-is-a-center}, we may assume that
the set  $\mathbb{CS}(X, \frac{1}{n}\mathcal{M})$ contains the
singular point $Q$. The proof of
Lemma~\ref{lemma:n-25-P2-is-a-center} also shows that the set
$\mathbb{CS}(U_{Q}, \frac{1}{n}\mathcal{M}_{U_Q})$ cannot consist
of the single point $\bar{P}$ whose image to $X$ is the point $P$.
It implies
$$
\mathbb{CS}\Big(U_{Q}, \frac{1}{n}\mathcal{M}_{U_Q}\Big)\cap\Big\{ O, Q_1 \Big\}\ne\varnothing%
$$
by Theorem~\ref{theorem:Noether-Fano} and
Lemma~\ref{lemma:Cheltsov-Kawamata}.

\begin{lemma}
\label{lemma:n-25-P3-and-P4} If the set $\mathbb{CS}(U_{Q},
\frac{1}{n}\mathcal{M}_{U_Q})$ contains both the point $O$ and the
point $Q_1$, then $\mathcal{M}=|-K_{X}|$.
\end{lemma}

\begin{proof}
Let $\gamma_Q:W_{Q}\to U_{QO}$ be the Kawamata blow up with
weights $(1,2,1)$ at the point whose image to $U_Q$ is the point
$Q_1$.

The proper transform $D_{W_{Q}}$ is irreducible and normal.  The
base locus of the pencil $|-K_{W_{Q}}|$ consists of the
irreducible curves $C_{W_{Q}}$ and $L_{W_{Q}}$. On the other hand,
we have
$$
\mathcal{M}_{W_{Q}}\Big\vert_{D_{W_{Q}}}\equiv-nK_{W_{Q}}\Big\vert_{D_{W_{Q}}}
\equiv nC_{W_{Q}}+nL_{W_{Q}},
$$
but the intersection form of the curves $L_{W_{Q}}$ and
$C_{W_{Q}}$ on the normal surface $D_{W_{Q}}$ is
negative-definite. Then, Theorem~\ref{theorem:main-tool} completes
the proof.
\end{proof}

It follows from Lemma~\ref{lemma:Cheltsov-Kawamata} that we may
assume  the following possibilities:
\begin{itemize}
\item $\mathbb{CS}(U_{Q}, \frac{1}{n}\mathcal{M}_{U_{Q}})=\{\bar{P}, O\}$;%

\item $\mathbb{CS}(U_{Q}, \frac{1}{n}\mathcal{M}_{U_{Q}})=\{ O\}$;%

\item $\mathbb{CS}(U_{Q}, \frac{1}{n}\mathcal{M}_{U_{Q}})=\{\bar{P}, Q_1\}$;%

\item $\mathbb{CS}(U_{Q}, \frac{1}{n}\mathcal{M}_{U_{Q}})=\{ Q_1\}$.%
\end{itemize}

 The exceptional divisor $F\cong\mathbb{P}(1,3,1)$ of
$\beta_{O}$ contains one singular point $Q_2$ that is a quotient
singularity of type $\frac{1}{3}(1,2,1)$.

\begin{lemma}
\label{lemma:n-25-P1-or-P4} If $\mathbb{CS}(U_{Q},
\frac{1}{n}\mathcal{M}_{U_Q})=\{ O \}$, then
$\mathcal{M}=|-K_{X}|$.
\end{lemma}

\begin{proof}
 The set $\mathbb{CS}(U_{QO},
\frac{1}{n}\mathcal{M}_{U_{QO}})$ contains the singular point
$Q_2$ by Theorem~\ref{theorem:Noether-Fano} and
Lemma~\ref{lemma:Cheltsov-Kawamata}.

Let $\gamma:W\to U_{QO}$ be the Kawamata blow up at the point
$Q_2$ with weights $(1,2,1)$. Then,
$\mathcal{M}_W\sim_{\mathbb{Q}}-nK_{W}$ by
Lemma~\ref{lemma:Kawamata} and the base locus of the pencil
$|-K_{W}|$ consists of the curves $C_W$ and $L_W$. The proper
transform $D_W$ is irreducible and normal, the equivalence
$\mathcal{M}_W\Big\vert_{D_W}\equiv nC_W+nL_W$ holds, but the
equalities
$$
C_W^{2}=-\frac{7}{12},\ \ L_W^{2}=-\frac{5}{6},\ \ C_W\cdot L_W=\frac{2}{3}%
$$
hold on the surface $D_W$. So, the intersection form of the curves
$C_W$ and $L_W$ on the normal surface $D_W$ is negative-definite,
which implies $\mathcal{M}=|-K_{X}|$ by
Theorem~\ref{theorem:main-tool}.
\end{proof}

\begin{lemma}
\label{lemma:n-25-P1-or-P3} If $\mathbb{CS}(U_{Q},
\frac{1}{n}\mathcal{M}_{U_Q})=\{\bar{P}, O\}$, then
$\mathcal{M}=|-K_{X}|$.
\end{lemma}

\begin{proof}
We have  $\mathcal{M}_Y\sim_{\mathbb{Q}}-nK_{Y}$, which implies
that every surface of the pencil $\mathcal{M}_Y$ is contracted to
a curve by the morphism  $\eta$. In particular, the set
$\mathbb{CS}(Y, \frac{1}{n}\mathcal{M}_Y)$ does not contain curves
because the exceptional divisors of $\beta_{O}\circ\gamma_{P}$ are
sections of $\eta$.

Because of Theorem~\ref{theorem:Noether-Fano} and
Lemmas~\ref{lemma:Cheltsov-Kawamata},  \ref{lemma:n-25-P1-or-P4},
we may assume that the set $\mathbb{CS}(Y,
\frac{1}{n}\mathcal{M}_Y)$ contains the singular point $P_2$ of
$Y$ contained in the exceptional divisor $\gamma_P$. Let
$\sigma_P:Y_P\to Y$ be the Kawamata blow up at the point $P_2$
with weights $(1,2,1)$. Then,
$\mathcal{M}_{Y_P}\sim_{\mathbb{Q}}-nK_{Y_P}$ but the base locus
of the pencil $|-K_{Y_P}|$ consists of the irreducible curves
$C_{Y_P}$ and $L_{Y_P}$.

The proper transform $D_{Y_P}$ is normal and
$\mathcal{M}_{Y_P}\vert_{D_{Y_P}}\equiv nC_{Y_P}+nL_{Y_P}$.  The
intersection form of the curves $C_{Y_P}$ and $L_{Y_P}$ on the
normal surface $D_{Y_P}$ is negative-definite because the curves
are contained in a fiber of $\eta\circ\sigma_P\vert_{D_{Y_P}}$
that consists of three irreducible components. Therefore, we
obtain the identity $\mathcal{M}=|-K_{X}|$ from
Theorem~\ref{theorem:main-tool}.
\end{proof}

Thus, to conclude the proof of Proposition~\ref{proposition:n-29},
we may assume  the following possibilities:
\begin{itemize}
\item $\mathbb{CS}(U_{Q}, \frac{1}{n}\mathcal{M}_{U_{Q}})=\{\bar{P}, Q_1\}$;%

\item $\mathbb{CS}(U_{Q}, \frac{1}{n}\mathcal{M}_{U_{Q}})=\{ Q_1\}$.%
\end{itemize}

The hypersurface $X$ can be given by the equation
$$
w^{2}y+wt^{2}+wtf_{4}(x,y,z)+wf_{8}(x,y,z)+tf_{11}(x,y,z)+f_{15}(x,y,z)=0,
$$
where  $f_{i}$ is a  general quasihomogeneous polynomial of degree
$i$.

\begin{lemma}
\label{lemma:n-25-P4} The case $\mathbb{CS}(U_{Q},
\frac{1}{n}\mathcal{M}_{U_Q})=\{ Q_1 \}$ never happens.
\end{lemma}
\begin{proof}
Suppose that $\mathbb{CS}(U_{Q}, \frac{1}{n}\mathcal{M}_{U_Q})=\{
Q_1 \}$.  Let $\pi:V\to U_{Q}$ be the Kawamata blow up at the
point $Q_1$ with weights $(1,2,1)$.

Let $G$ be the exceptional divisor of the birational morphism
$\pi$. The proof of Lemma~\ref{lemma:n-25-P1-or-P4} implies that
the set $\mathbb{CS}(V, \frac{1}{n}\mathcal{M}_V)$ does not
contain the singular point of $V$ contained in the exceptional
divisor $G$. So, the log pair $(V, \frac{1}{n}\mathcal{M}_V)$ is
terminal by Lemma~\ref{lemma:Cheltsov-Kawamata} and
Corollary~\ref{corollary:curves}.

We have
\begin{equation*}
\left\{\aligned
&(\alpha_{Q}\circ\pi)^{*}\big(-K_{X}\big)\sim_{\mathbb{Q}} S_V+\frac{3}{7}G+\frac{1}{7}E_V,\\
&(\alpha_{Q}\circ\pi)^{*}\big(-K_{X}\big)\sim_{\mathbb{Q}}S_V^y+\frac{10}{7}G+\frac{8}{7}E_V,\\
&(\alpha_{Q}\circ\pi)^{*}\big(-3K_{X}\big)\sim_{\mathbb{Q}}S_V^z+\frac{2}{7}G+\frac{3}{7}E_V,\\
&(\alpha_{Q}\circ\pi)^{*}\big(-4K_{X}\big)\sim_{\mathbb{Q}}S_V^t+\frac{5}{7}G+\frac{4}{7}E_V,\\
&(\alpha_{Q}\circ\pi)^{*}\big(-7K_{X}\big)\sim_{\mathbb{Q}} S_V^w.\\
\endaligned\right.
\end{equation*}
The equivalences 
imply
\[\left\{\aligned & (\alpha_{Q}\circ\pi)^{*}\Big(\frac{y}{x}\Big)\in |S_V|, \\
&(\alpha_{Q}\circ\pi)^{*}\Big(\frac{yz}{x^4}\Big)\in |4S_V|, \\
&(\alpha_{Q}\circ\pi)^{*}\Big(\frac{yt}{x^5}\Big)\in |5S_V|, \\
&(\alpha_{Q}\circ\pi)^{*}\Big(\frac{y^3w}{x^{10}}\Big)\in
|10S_V|,\\
\endaligned\right.\]
and hence the complete linear system $|-20K_{V}|$ induces a
birational map $\chi_1:V\dasharrow X^{\prime}$ such that
$X^{\prime}$ is a hypersurface in $\mathbb{P}(1,1,4,5,10)$, which
implies that the divisor $-K_{V}$ is big.

The base locus of the pencil $|-K_{V}|$ consists of the
irreducible curves $C_V$ and $L_V$. It follows from \cite{Sho93}
that there is a composition of antiflips $\zeta:V\dasharrow V'$
such that  $\zeta$ is regular in the outside of $C_V\cup L_V$ and
the anticanonical divisor  $-K_{V'}$ is nef and big. The
singularities of the log pair $(V', \frac{1}{n}\mathcal{M}_{V'})$
are terminal because the rational map $\zeta$ is a log flop with
respect to the log pair $(V, \frac{1}{n}\mathcal{M}_{V})$, which
contradicts Theorem~\ref{theorem:Noether-Fano}.
\end{proof}

Now, we suppose $\mathbb{CS}(U_{Q},
\frac{1}{n}\mathcal{M}_{U_{Q}})=\{\bar{P}, Q_1\}$.  Let
$\sigma:U\to U_{PQ}$ be the Kawamata blow up with weights
$(1,2,1)$ at the point $\bar{Q}_1$ whose image to $U_Q$ is the
point $Q_1$.  Let $\bar{E}$ and $\tilde{E}$ be the exceptional
divisors of $\alpha_P$ and $\sigma$, respectively. Then,
$$
\left\{\aligned
&S_U\sim_{\mathbb{Q}}(\alpha_{P}\circ\beta_Q\circ\sigma)^{*}\big(-K_{X}\big)-\frac{3}{7}\tilde{E}-\frac{1}{7}E_U-\frac{1}{4}\bar{E}_U\sim_{\mathbb{Q}}-K_{W},\\
&S_U^{y}\sim_{\mathbb{Q}}(\alpha_{P}\circ \beta_Q\circ\sigma)^{*}\big(-K_{X}\big)-\frac{10}{7}\tilde{E}-\frac{8}{7}E_U-\frac{1}{4}\bar{E}_U,\\
&S_U^{z}\sim_{\mathbb{Q}}(\alpha_{P}\circ\beta_Q\circ\sigma)^{*}\big(-3K_{X}\big)-\frac{2}{7}\tilde{E}-\frac{3}{7}E_U-\frac{3}{4}\bar{E}_U,\\
&S_U^{t}\sim_{\mathbb{Q}}(\alpha_{P}\circ\beta_Q\circ\sigma)^{*}\big(-4K_{X}\big)-\frac{5}{7}\tilde{E}-\frac{4}{7}E_U,\\
&S_U^{w}\sim_{\mathbb{Q}}(\alpha_{P}\circ\beta_Q\circ\sigma)^{*}\big(-7K_{X}\big)-\frac{11}{4}\bar{E}_U.\\
\endaligned\right.
$$
 The equivalences imply that the pull-backs of
rational functions
$$\frac{y}{x}, \ \frac{yz}{x}, \ \frac{y^3w}{x^{10}}, \
\frac{y^4tw}{x^{15}}$$ are contained in the linear system
$|aS_U|$, where $a=1$, $4$, $10$, and $15$, respectively.
Therefore, the linear system $|-60K_{U}|$ induces a birational map
$\chi_2:U\dasharrow X''$ such that the variety $X''$ is a
hypersurface of degree $30$ in $\mathbb{P}(1,1,4,10,15)$, which
implies that the anticanonical divisor $-K_{U}$ is big. However,
the proof of Lemma~\ref{lemma:n-25-P1-or-P4} shows that the
singularities of $(U, \frac{1}{n}\mathcal{M}_U)$ are terminal.
Then, we can obtain a contradiction in the same way as in the
proof of Lemma~\ref{lemma:n-25-P4}.

\section{Case $\gimel=26$, hypersurface of degree $15$ in
$\mathbb{P}(1,1,3,5,6)$.}\index{$\gimel=26$} \label{section:n-26}

The threefold $X$ is a general hypersurface of degree $15$ in
$\mathbb{P}(1,1,3,5,6)$ with $-K_{X}^{3}=\frac{1}{6}$. The
singularities of the hypersurface $X$ consist of two points
$P_{1}$ and $P_{2}$ that are quotient singularities of type
$\frac{1}{3}(1,1,2)$ and one point $Q$ that is a quotient
singularity of type $\frac{1}{6}(1,1,5)$.

For each of the singular points $P_{1}$ and $P_{2}$, we have a
commutative diagram
$$
\xymatrix{
&&U_{i}\ar@{->}[dl]_{\sigma_{i}}\ar@{->}[dr]^{\eta_{i}}&&\\
&X\ar@{-->}[rr]_{\xi_{i}}&&\mathbb{P}(1,1,6),&}
$$
where
\begin{itemize}

\item $\xi_{i}$ is a projection,

\item $\sigma_{i}$ is the Kawamata blow up at the point $P_{i}$
with weights $(1,1,2)$,

\item $\eta_{i}$ is an elliptic fibration.
\end{itemize}
Note that the threefold $X$ has a unique birational automorphism which is not biregular. It is a quadratic involution that is
defined in the outside of the point $Q$. Moreover, it interchanges the points $P_1$ and $P_2$. The linear system
that induces the rational map $\xi_1$ is transformed
into the linear system that induces the rational map $\xi_2$ by the quadratic involution and vice versa.

We have another elliptic fibration
$$
\xymatrix{
&U\ar@{->}[d]_{\alpha}&&W\ar@{->}[ll]_{\beta}&&Y\ar@{->}[ll]_{\gamma}\ar@{->}[d]^{\eta}&\\
&X\ar@{-->}[rrrr]_{\psi}&&&&\mathbb{P}(1,1,3)&}
$$
where \begin{itemize}

\item $\psi$ is  the natural projection,

\item $\alpha$ is the Kawamata blow up at the point $Q$ with
weights $(1,1,5)$,

\item $\beta$ is the Kawamata blow up with weights $(1,1,4)$ at
the singular point $Q_1$ of the variety $U$  contained in the
exceptional divisor of $\alpha$,

\item $\gamma$ is the Kawamata blow up with weights $(1,1,3)$ at
the singular point $Q_2$ contained in the exceptional divisor of
the birational morphism $\beta$,

\item $\eta$ is an elliptic fibration.

\end{itemize}

It follows from \cite{CPR} that the pencil $|-K_{X}|$ is invariant
under the action of the group $\mathrm{Bir}(X)$. Hence, we may
assume that the log pair $(X, \frac{1}{n}\mathcal{M})$ is
canonical. In fact, we can assume that
$$
\varnothing\ne\mathbb{CS}\Big(X, \frac{1}{n}\mathcal{M}\Big)\subseteq\Big\{P_{1}, P_{2}, Q\Big\},%
$$
due to Lemma~\ref{lemma:smooth-points} and
Corollary~\ref{corollary:Ryder-a2}.

\begin{lemma}
\label{lemma:n-26-points-P1-and-P2} If  the set $\mathbb{CS}(X,
\frac{1}{n}\mathcal{M})$ contains either the point $P_1$ or the
point $P_2$, then $\mathcal{M}=|-K_{X}|$.
\end{lemma}

\begin{proof}
Suppose that the set $\mathbb{CS}(X, \frac{1}{n}\mathcal{M})$
contains the point  $P_{1}$. Because
$\mathcal{M}_{U_1}\sim_{\mathbb{Q}}-nK_{U_{1}}$ by
Lemma~\ref{lemma:Kawamata}, each surface in the pencil
$\mathcal{M}_{U_1}$ is contracted to a curve by the morphism
$\eta_{1}$ and hence the set $\mathbb{CS}(U_{1},
\frac{1}{n}\mathcal{M}_{U_1})$ does not contain curves.

Let $\bar{P}_{2}$ and $\bar{Q}$ be the  points on $U_1$ whose
images to $X$ are the points  $P_{2}$ and $Q$, respectively, and
$O$ be the singular point of $U_{1}$ contained in the exceptional
divisor of $\sigma_{1}$. Then,
$$
\mathbb{CS}\Big(U_{1},\frac{1}{n}\mathcal{B}\Big)\cap\Big\{\bar{P}_{2}, \bar{Q}, O\Big\}\ne\varnothing%
$$
by Theorem~\ref{theorem:Noether-Fano} and
Lemma~\ref{lemma:Cheltsov-Kawamata}. We consider only the case
when the set $\mathbb{CS}(U_{1}, \frac{1}{n}\mathcal{M}_{U_1})$
contains the point $O$ because the other cases are similar.
Suppose that the set $\mathbb{CS}(U_{1},
\frac{1}{n}\mathcal{M}_{U_1})$ contains the point $O$.

Let $\pi:V\to U_{1}$ be the Kawamata blow up at the point $O$ with
weights $(1,1,1)$. Then, the base locus of the pencil $|-K_{V}|$
consists of the curve $C_V$. Let $D_V$ be a general surface in
$|-K_{V}|$. Then, $D_V$ is normal and $C_V^{2}=-\frac{1}{2}$ on
the surface $D_V$, which implies that $\mathcal{M}=|-K_{X}|$ by
Lemma~\ref{lemma:Kawamata} and Theorem~\ref{theorem:main-tool}.
\end{proof}

Therefore, we may assume that the set $\mathbb{CS}(X,
\frac{1}{n}\mathcal{M})$ consists of the point $Q$.

The excep\-ti\-onal divisor $E$ of the birational morphism
$\alpha$ contains a singular point  $Q_1$ that is  a quotient
singularity of type $\frac{1}{5}(1,1,4)$. The set $\mathbb{CS}(U,
\frac{1}{n}\mathcal{M}_U)$ is not empty by
Theorem~\ref{theorem:Noether-Fano}.

\begin{lemma}
\label{lemma:n-26-point-P4} The set $\mathbb{CS}(U,
\frac{1}{n}\mathcal{M}_U)$ consists of the point $Q_1$.
\end{lemma}

\begin{proof} Suppose that the set $\mathbb{CS}(U,
\frac{1}{n}\mathcal{M}_U)$ contains a subvariety $Z\subset U$ that
is different from the singular point $Q_1$. Then, the subvariety
$Z$ is a curve with  $-K_{U}\cdot Z=\frac{1}{5}$ by
Lemma~\ref{lemma:Cheltsov-Kawamata}, which is impossible by
Lemma~\ref{lemma:curves} because $-K_{U}^{3}=\frac{2}{15}$.
\end{proof}

Let $G$ be the excep\-ti\-onal divisor of the birational morphism
$\gamma$ and $Q_3$ be the singular point of $G$. Then,
Theorem~\ref{theorem:Noether-Fano}, Lemmas~\ref{lemma:Kawamata},
and \ref{lemma:Cheltsov-Kawamata} imply that
$\mathcal{M}_Y\sim_{\mathbb{Q}}-nK_{Y}$. Each member in the pencil
$\mathcal{M}_Y$ is contracted to a curve by the morphism $\eta$
and $\mathbb{CS}(Y, \frac{1}{n}\mathcal{M}_Y)\ne\varnothing$.

\begin{lemma}
\label{lemma:n-26-point-P6} The set $\mathbb{CS}(Y,
\frac{1}{n}\mathcal{H})$ consists of the point $Q_3$.
\end{lemma}

\begin{proof}
Suppose that the set $\mathbb{CS}(Y, \frac{1}{n}\mathcal{M}_Y)$
contains a subvariety $Z$ of the variety $Y$ that is different
from the point $Q_3$. Let $F$ be  the exceptional divisor of the
birational morphism $\beta$. The base locus of the pencil
$\mathcal{M}_Y$ does not contain curves  contained in $F$ since
$F$ is a section of $\eta$. Then, $Z$ is an irreducible curve such
that the curve $\gamma(Z)$ is a ruling of the cone $F$ by
Lemma~\ref{lemma:Cheltsov-Kawamata}. Thus, we have $-K_{W}\cdot
\gamma(C)=\frac{1}{4}$, which is impossible by
Lemma~\ref{lemma:curves}.
\end{proof}

Let $\gamma_1:Y_1\to Y$ be the Kawamata blow up at the point $Q_3$
with weights $(1,1,2)$. Then, the base locus of the pencil
$|-K_{Y_1}|$ consists of the curve $C_{Y_1}$. Let $D_{Y_1}$ be a
general surface in $|-K_{Y_1}|$. Then, $D_{Y_1}$ is normal and
$C_{Y_1}^{2}=-\frac{1}{6}$ on the surface $D_{Y_1}$, which implies
that $\mathcal{M}=|-K_{X}|$ by Lemma~\ref{lemma:Kawamata} and
Theorem~\ref{theorem:main-tool}.

Therefore, we have obtained
\begin{proposition}
\label{proposition:n-26} The linear system $|-K_{X}|$ is the only
Halphen pencil on $X$.
\end{proposition}

\section{Cases $\gimel=29$, $50$,  and $67$.}\index{$\gimel=29$}\index{$\gimel=50$}\index{$\gimel=67$} \label{section:n-29}

Suppose that $\gimel\in\{29, 50, 67 \}$. Then, the threefold $X
\subset \mathbb{P}(1,1,a_2,a_3,a_4)$ always contains the point
$O=(0:0:0:1:0)$. It  is a singular point of $X$ that is a quotient
singularity of type $\frac{1}{a_3}(1,a_2, a_3-a_2)$.

We also have a commutative diagram as follows:
$$
\xymatrix{
&U\ar@{->}[d]_{\alpha}&&W\ar@{->}[ll]_{\beta}\ar@{->}[d]^{\eta}&\\%
&X\ar@{-->}[rr]_{\psi}&&\mathbb{P}(1,1,a_2),&}
$$
where \begin{itemize}

\item $\psi$ is the natural projection,

\item $\alpha$ is the Kawamata blow up at the point $O$ with
weights $(1,a_2,a_3-a_2)$,

\item $\beta$ is the Kawamata blow up with weights $(1,a_2,1)$ at
the point $P$ of $U$ that is a quotient singularity of type
$\frac{1}{a_3-a_2}(1,a_2,1)$,

\item $\eta$ is an elliptic fibration.
\end{itemize}

The exceptional divisor $E$ of the birational morphism $\alpha$
contains two singular points $P$ and $Q$ that are quotient
singularity of types $\frac{1}{a_3-a_2}(1,a_2,1)$ and
$\frac{1}{a_2}(1,a_2-1, 1)$, respectively. The base locus of
$|-K_X|$ consists of the irreducible curve $C$ defined by $x=y=0$.
The base curve of $|-K_U|$ consists of the proper transform $C_U$
and the unique irreducible curve $L$ in
$|\mathcal{O}_{\mathbb{P}(1, a_2, a_3-a_2)}(1)|$ on the surface
$E$.

\begin{lemma}
\label{lemma:n-29-second-floor} If the set  $\mathbb{CS}(U,
\frac{1}{n}\mathcal{M}_U)$ contains the point $Q$, then
$\mathcal{M}=|-K_{X}|$.
\end{lemma}

\begin{proof}
Let $\pi:Y\to U$ be the Kawamata blow up at the point $Q$ with
weights $(1,a_2-1, 1)$ and $G_Q$ be its exceptional divisor. Then,
the base locus of the pencil $|-K_{Y}|$ consists of the
irreducible curves $C_Y$ and $L_Y$. Let $D$ be a general surface
in $|-K_{X}|$. We then have
$$
\left\{\aligned
&D_Y\sim_{\mathbb{Q}}\big(\alpha\circ\pi\big)^{*}\Big(-K_{X}\Big)-\frac{1}{a_3}\pi^{*}\big(E\big)-\frac{1}{a_2}G_Q,\\
&S_Y\sim_{\mathbb{Q}}\big(\alpha\circ\pi\big)^{*}\Big(-K_{X}\Big)-\frac{1}{a_3}\pi^{*}\big(E\big)-\frac{1}{a_2}G_Q,\\
&E_Y\sim_{\mathbb{Q}}\pi^{*}\big(E\big)-\frac{1}{a_2}G_Q\\
\endaligned
\right.
$$
and we also have $S_Y\cdot D_Y=C_Y+L_Y$ and $E_Y\cdot D_Y=L_Y$. It
then follows that
$$
\left\{\aligned
&D_Y\cdot L_Y=E_Y\cdot D_Y^2=E_Y\cdot K_Y^2=-\frac{a_3+1}{a_2(a_3-a_2)(a_2-1)}<0, \\
&D_Y\cdot C_Y=(-K_Y)^3-D_Y\cdot L_Y
\\ &\phantom{D\cdot C_Y}=\frac{a_4}{a_2a_3}-\frac{1}{a_2a_3(a_3-a_2)}-\frac{1}{a_2(a_2-1)}
+\frac{a_3+1}{a_2(a_3-a_2)(a_2-1)}\leq0.\\
\endaligned
\right.
$$
The divisors $-K_U$ and $-K_X$ are nef and big. Moreover,
$$
\left\{\aligned
&-K_U\cdot L = E\cdot K_U^2 =E\cdot K_U^2=\frac{1}{a_2(a_3-a_2)},\\
& -K_U\cdot C_U= S_U\cdot D_U^2+K_U\cdot L=-K_U^3+K_U\cdot L=\frac{1}{a_2a_3}\Big(a_4-\frac{1+a_3}{a_3-a_2}\Big), \\
& -K_X\cdot C= -K_X^3=\frac{a_4}{a_2a_3}. \\
\endaligned
\right.
$$
One can easily check that \[\lambda :=\frac{D_Y\cdot L_Y}{K_U\cdot
L}>0, \  \  \mu:= \frac{(D_Y\cdot C_Y)(K_U\cdot L)-(D_Y\cdot
L_Y)(K_U\cdot C_U)}{(K_X\cdot C)(K_U\cdot L)}\geq 0.\] Then,
$B=\pi^*(-\lambda K_U) + (\alpha\circ\pi)^*(-\mu K_X)+D_Y$ is a
nef and big divisor with $B\cdot L_Y=B\cdot C_Y=0$. Therefore,
$B\cdot D_Y\cdot M_Y=0$, where $M_Y$ is a general surface in
$\mathcal{M}_Y$, and hence $\mathcal{M}=|-K_X|$ by
Theorem~\ref{theorem:main-tool}.
\end{proof}

The exceptional divisor $F$ of the birational morphism $\beta$
contains one singular point $P_1$ of the threefold $W$ that is a
quotient singularity of type
 $\frac{1}{a_2}(1,a_2-1,
1)$.

\begin{proposition}
\label{proposition:n-29} The linear system $|-K_{X}|$ is the only
Halphen pencil on $X$.
\end{proposition}

\begin{proof}
By Lemma~\ref{lemma:special-singular-points-with-zero-c}, if the
set $\mathbb{CS}(X, \frac{1}{n}\mathcal{M})$ contains a singular
point of $X$ different from the singular point $O$, then the
identity $\mathcal{M}=|-K_X|$ holds. Moreover, if it contains a
curve, then Corollary~\ref{corollary:Ryder-a2} implies
$\mathcal{M}=|-K_X|$. Therefore, we may assume that
$$
\mathbb{CS}\Big(X, \frac{1}{n}\mathcal{M}\Big)=\Big\{O\Big\}.%
$$
due to Theorem~\ref{theorem:Noether-Fano} and
Lemma~\ref{lemma:smooth-points}.

 Furthermore,
Lemma~\ref{lemma:Cheltsov-Kawamata} implies that either
$\mathbb{CS}(U, \frac{1}{n}\mathcal{M}_U)=\{P\}$ or
$Q\in\mathbb{CS}(U, \frac{1}{n}\mathcal{M}_U)$. The latter case
implies $\mathcal{M}=|-K_X|$ by
Lemma~\ref{lemma:n-29-second-floor}. Suppose that the set
$\mathbb{CS}(U, \frac{1}{n}\mathcal{M}_U)$ consists of the point
$P$. Then, the set $\mathbb{CS}(W, \frac{1}{n}\mathcal{M}_W)$
contains the point $P_1$. Let $\sigma:V\to W$ be the Kawamata blow
up at the point $P_1$. Then, the base locus of the pencil
$|-K_{V}|$ consists of the irreducible curves $C_V$ and $L_V$.
Applying the same method as in
Lemma~\ref{lemma:n-29-second-floor}, we obtain the identity
$\mathcal{M}=|-K_{X}|$.
\end{proof}

\section{Cases $\gimel=34$,  $53$,   $70$,  and  $88$.}%
\label{section:n-34-53-70-88}

Suppose
$\gimel\in\{34,53,70,88\}$.\index{$\gimel=34$}\index{$\gimel=53$}\index{$\gimel=70$}\index{$\gimel=88$}
Then, the hypersurface $X$ has a singular point $P$ of    type
$\frac{1}{a_{2}+1}(1,1,a_{2})$.

We also have a commutative diagram as follows:

$$
\xymatrix{
&&&Y\ar@{->}[lld]_{\pi}\ar@{->}[rrd]^{\eta}&&&\\%
&X\ar@{-->}[rrrr]_{\psi}&&&&\mathbb{P}(1,1,a_{2})&}
$$
where \begin{itemize} \item $\psi$ is the natural projection,
\item $\pi$ is the Kawamata blow up at the point $P$ with weights
$(1, 1, a_2)$, \item $\eta$ is an elliptic fibration.\end{itemize}

We may assume that the set  $\mathbb{CS}(X,
\frac{1}{n}\mathcal{M})$ consists of the point $P$ by
Lemmas~\ref{lemma:smooth-points},
\ref{lemma:special-singular-points-with-zero-c} and
Corollary~\ref{corollary:Ryder-a2}.

\begin{proposition}
\label{proposition:n-34-53-70-88} If $\gimel\in\{34,53,70,88\}$,
then the linear system $|-K_{X}|$ a unique Halphen pencil.
\end{proposition}

\begin{proof}

The singularities of the log pair $(Y, \frac{1}{n}\mathcal{M}_Y)$
are not terminal by Theorem~\ref{theorem:Noether-Fano}. The base
locus of the pencil $\mathcal{M}_Y$ does not contain curves that
are not contained in a fiber of the elliptic fibration $\eta$.
Hence, it follows from Lemma~\ref{lemma:Cheltsov-Kawamata} that
the set $\mathbb{CS}(Y, \frac{1}{n}\mathcal{M}_Y)$ contains the
singular point $Q$ contained in the exceptional divisor of $\pi$,
which is a quotient singularity of type
$\frac{1}{a_{2}}(1,1,a_{2}-1)$ on $Y$.

Let $\alpha:U\to Y$ be the Kawamata blow up at the point $P$ with
weights $(1,1,a_{2}-1)$. Then, the pencil $|-K_U|$ is the proper
transform of the pencil $|-K_X|$. Its base locus consists of the
irreducible curve $C_U$. Moreover, $-K_U\cdot C_U=-K_U^3<0$ and
$\mathcal{M}_U\equiv -nK_U$. We then obtain $\mathcal{M}=|-K_{X}|$
from Theorem~\ref{theorem:main-tool}.
\end{proof}

\section{Case $\gimel=36$, hypersurface of degree $18$ in
$\mathbb{P}(1,1,4,6,7)$.}\index{$\gimel=36$} \label{section:n-36}

Let $X$ be the hypersurface given by a general quasihomogeneous
equation of degree $18$ in $\mathbb{P}(1,1,4,6,7)$
 with
$-K_X^3=\frac{3}{28}$. Then, the singularities of $X$ consist of
two singular points $P$ and $Q$ that are quotient singularities of
types  $\frac{1}{7}(1,1,6)$ and $\frac{1}{4}(1,1,3)$,
respectively, and one point of type $\frac{1}{2}(1,1,1)$.

 There
is a commutative diagram
$$
\xymatrix{
&U_{P}\ar@{->}[d]_{\alpha_{P}}&&V\ar@{->}[ll]_{\beta}&&Y\ar@{->}[ll]_{\gamma}\ar@{->}[d]^{\eta}&\\
&X\ar@{-->}[rrrr]_{\psi}&&&&\mathbb{P}(1,1,4)&}
$$
where \begin{itemize}

\item $\psi$ is  the natural projection,

\item $\alpha_P$ is the weighted blow up at the point $P$ with
weights $(1,1,6)$,

\item $\beta$ is the Kawamata blow up with weights $(1,1,5)$ at
the singular point $O_1$ of the variety $U_P$  contained in the
exceptional divisor of the birational morphism $\alpha_{P}$,

\item $\gamma$ is the Kawamata blow up with weights $(1,1,4)$ at
the singular point $O_2$ of the variety $V$ contained in the
exceptional divisor of the birational morphism $\beta$,

\item $\eta$ is an elliptic fibration.
\end{itemize}

The hypersurface $X$ can be given by the quasihomogeneous equation
of degree $18$
$$
z^{3}t+z^{2}f_{10}(x,y,t,w)+zf_{14}(x,y,t,w)+f_{18}(x,y,t,w)=0,
$$
where  $f_i$ is a quasihomogeneous polynomial of degree $i$. Let
$\xi:X\dasharrow\mathbb{P}^{7}$ be the rational map given by the
linear system of divisors cut on the hypersurface $X$ by the
equations $$\mu t+\sum_{i=0}^{6}\lambda_{i}x^{i}y^{6-i}=0$$, where
$(\mu:\lambda_0:\cdots: \lambda_6)\in\mathbb{P}^7$. Then, we
obtain another commutative diagram
$$
\xymatrix{
&&W\ar@{->}[dl]_{\beta_{Q}}\ar@{->}[dr]^{\beta_{P}}\ar@{->}[rrrrrrrd]^{\eta_1}&&&&&\\
&U_{P}\ar@{->}[dr]_{\alpha_{P}}&&U_{Q}\ar@{->}[dl]^{\alpha_{Q}}&&&&&&\mathbb{P}(1,1,6),\\
&&X\ar@{-->}[rrrrrrru]_{\xi}&&&&&&&}
$$
where \begin{itemize}

\item $\alpha_{P}$ is the Kawamata blow up at the singular point
$P$ with weights $(1,1,6)$,

\item $\alpha_{Q}$ is the Kawamata blow up at the singular point
$Q$ with weights $(1,1,3)$,

\item $\beta_{P}$ is the Kawamata blow up  with weights $(1,1,6)$
at the point $P_1$ whose image to $X$ is the point  $P$,

\item $\beta_{Q}$ is the Kawamata blow up with weights $(1,1,3)$
at the point $Q_1$ whose image to $X$ is the point  $Q$,

\item $\eta_1$ is an elliptic fibration.
\end{itemize}

If it contains either the singular point of type
$\frac{1}{2}(1,1,1)$ or a curve, then
Lemma~\ref{lemma:special-singular-points-with-zero-c} and
Corollary~\ref{corollary:Ryder-a2} imply $\mathcal{M}=|-K_X|$.
Therefore, we may also assume that
$$\varnothing\ne\mathbb{CS}\Big(X,\frac{1}{n}\mathcal{M}\Big)\subset \Big\{P, Q\Big\}.$$

The exceptional divisor $E_Q$ of the birational morphism
$\alpha_Q$ contains a unique singular point $O$ that is a quotient
singularity of type $\frac{1}{3}(1,1,2)$.

\begin{lemma}
If the set $\mathbb{CS}(U_Q,\frac{1}{n}\mathcal{M}_{U_Q})$
contains the point $O$, then  $\mathcal{M}=|-K_X|$.
\end{lemma}
\begin{proof}
Let $\beta_O:U_O\to U_Q$ be the Kawamata blow up at the point $O$
with weights $(1,1,2)$. Then, the proper transform
$\mathcal{D}_{U_O}$ of the pencil $|-K_X|$ by the birational
morphism $\alpha_Q\circ\beta_O$ has a unique base curve $C_{U_O}$.
Because a general surface in $\mathcal{D}_{U_O}$ is normal and the
self-intersection of $C_{U_O}$ on a general surface in
$\mathcal{D}_{U_O}$ is negative, we obtain $\mathcal{M}=|-K_X|$.
\end{proof}

\begin{lemma}
If $\mathbb{CS}(U_Q,\frac{1}{n}\mathcal{M}_{U_Q})=\{P_1\}$, then
$\mathcal{M}=|-K_X|$.
\end{lemma}
\begin{proof}
The proper transform $\mathcal{M}_W$ must contain the singular
point in the exceptional divisor $F_P$ of the birational morphism
$\beta_P$ that is a quotient singularity of type
$\frac{1}{6}(1,1,5)$. Let $\sigma_1:W_1\to W$ be the Kawamata blow
up at the singular point with weights $(1,1,5)$. Then, the proper
transform $\mathcal{D}_{W_1}$ of the pencil $|-K_X|$ has a unique
base curve $C_{W_1}$. Because a general surface in
$\mathcal{D}_{W_1}$ is normal and  the self-intersection of
$C_{W_1}$ on a general surface in $\mathcal{D}_{W_1}$ is negative,
we obtain $\mathcal{M}=|-K_X|$.
\end{proof}

Meanwhile, the exceptional divisor $E_P$ of the birational
morphism $\alpha_P$ contains one singular point $O_1$ of type
$\frac{1}{6}(1,1,5)$.
\begin{lemma}\label{lemma:n-36-curve}
The set $\mathbb{CS}(U_P,\frac{1}{n}\mathcal{M}_{U_P})$ cannot
contain a curve.
\end{lemma}
\begin{proof}
Because $-K_{U_P}^3=\frac{1}{12}$ and $E_P\cdot
K_{U_P}^2=\frac{1}{6}$, the statement  immediately follows from
Lemma~\ref{lemma:curves}.
\end{proof}

Therefore, we may assume that
$$\mathbb{CS}\Big(U_1,\frac{1}{n}\mathcal{M}_{U_1}\Big)= \Big\{O_1\Big\}.$$
 By the same method of
Lemma~\ref{lemma:n-36-curve}, we may assume that
$\mathbb{CS}(V,\frac{1}{n}\mathcal{M}_{V})=\{O_2\},$ and hence the
set $\mathbb{CS}(Y,\frac{1}{n}\mathcal{M}_{Y})$ must contain the
singular point contained in the exceptional divisor of $\gamma$
that is a quotient singularity of type $\frac{1}{4}(1,1,3)$. By
considering Kawamata blow up at this point with weights $(1,1,3)$
and the proper transform of the pencil $|-K_X|$, one can easily
conclude the following:

\begin{proposition}
\label{proposition:n-36} The linear system $|-K_X|$ is a unique
Halphen pencil on $X$.
\end{proposition}

\section{Cases $\gimel=47$, $54$, and $62$.} \label{section:n-47}

We first consider the case $\gimel=47$ which is more complicated
than the other.\index{$\gimel=47$} The threefold $X$ is a  general
hypersurface of degree $21$ in $\mathbb{P}(1,1,5,7,8)$ with
$-K_{X}^{3}=\frac{3}{40}$. The singularities of $X$ consist of one
point $P$ that is a quotient singularity of type
$\frac{1}{5}(1,2,3)$ and one point $Q$ that is a quotient
singularity of type $\frac{1}{8}(1,1,7)$.

We have the following  commutative diagram:
$$
\xymatrix{
&U\ar@{->}[d]_{\alpha}&&W\ar@{->}[ll]_{\beta}&&Y\ar@{->}[ll]_{\gamma}\ar@{->}[d]^{\eta}&\\%
&X\ar@{-->}[rrrr]_{\psi}&&&&\mathbb{P}(1,1,5),&}
$$
where \begin{itemize}

\item $\psi$ is the natural projection,

\item $\alpha$ is the Kawamata blow up at the point $Q$ with
weights $(1,1,7)$,

\item $\beta$ is the Kawamata blow up with weights $(1,1,6)$ at
the singular point of $U$ that is a quotient singularity of type
$\frac{1}{7}(1,1,6)$,

\item $\gamma$ is the Kawamata blow up  with weights $(1,1,5)$ at
the singular point of $W$ that is a quotient singularity of type
$\frac{1}{6}(1,1,5)$,

\item $\eta$ is an elliptic fibration.
\end{itemize}

The hypersurface $X$ can be given by the equation
$$
w^{2}z+\sum_{i=0}^{2}wz^{i}g_{13-5i}(x,y,t)+\sum_{i=0}^{3}z^{i}g_{21-5i}(x,y,t)=0,
$$
where $g_{i}(x,y,t)$ is a quasihomogeneous polynomial of degree
$i$. The base locus of the pencil $|-K_{X}|$ consists of the
irreducible curve $C$  cut out on $X$ by the equations $x=y=0$.
Let $D$ be a general surface in $|-K_{X}|$. Then, $D$ is smooth at
a generic  point of $C$. Note that $$\mathbb{CS}\Big(X,
|-K_X|\Big)=\Big\{P, Q, C\Big\}.$$

If the set $\mathbb{CS}(X, \frac{1}{n}\mathcal{M})$ contains a
curve, then the identity $\mathcal{M}=|-K_X|$ holds due to
Corollary~\ref{corollary:Ryder-a2}. Furthermore, by
Lemma~\ref{lemma:smooth-points}, we may assume that
$$\mathbb{CS}\Big(X, \frac{1}{n}\mathcal{M}\Big)\subseteq\Big\{P, Q\Big\}.$$

\begin{lemma}
\label{lemma:n-47-P2} The set $\mathbb{CS}(X,
\frac{1}{n}\mathcal{M})$ contains the point $Q$.
\end{lemma}

\begin{proof}
Suppose that the set $\mathbb{CS}(X, \frac{1}{n}\mathcal{M})$ does
not contain $Q$. Then, $\mathbb{CS}(X,
\frac{1}{n}\mathcal{M})=\{P\}$, and hence $\mathcal{M}\ne
|-K_{X}|$.

Let $\alpha_P:U_P\to X$ be the Kawamata blow up at the point $P$.
The exceptional divisor of $E_P\cong\mathbb{P}(1,2,3)$  contains
two singular points $P_1$ and $P_2$ of the threefold $U_P$ that
are quotient singularities of types $\frac{1}{2}(1,1,1)$ and
$\frac{1}{3}(1,2,1)$, respectively. For the convenience, let $L$
be the unique curve  in the linear system
$|\mathcal{O}_{\mathbb{P}(1,\,2,\,3)}(1)|$ on the surface $E_P$.

If the log pair $(U_P, \frac{1}{n}\mathcal{M}_{U_P})$ is not
terminal, then it is not terminal either at the point $P_1$ or at
the point $P_2$. In such cases, the proof of
Lemma~\ref{lemma:n-29-second-floor} leads us to the identity
$\mathcal{M}=|-K_X|$, which contradicts  our assumption.
Therefore, the log pair must be terminal.

The base locus of the pencil $|-K_{U_P}|$ consists of the curve
$C_{U_P}$ and the irreducible curve $L$. The inequalities
$-K_{U_P}\cdot C_{U_P}<0$ and $-K_{U_P}\cdot L>0$ implies that the
curve $C_{U_P}$ is the only curve on the variety $U_P$ that has
negative intersection with the divisor $-K_{U_P}$. It follows from
\cite{Sho93} that the antiflip $\zeta_1:U_P\dasharrow U_P'$ along
the curve $C_{U_P}$ exists. The divisor  $-K_{U_P'}$ is nef. The
 log pair $(U_P',
\frac{1}{n}\mathcal{M}_{U_P'})$ is terminal because the log pair
$(U_P, \frac{1}{n}\mathcal{M}_{U_P})$ is so.  On the other hand,
the pull-backs of the rational functions $1$, $\frac{y}{x}$,
$\frac{zy}{x^{6}}$, $\frac{ty}{x^{8}}$ and $\frac{yw}{x^{9}}$
induce a birational map $\chi_1$ of $U_P$ onto  a hypersurface
$X'$  of degree $24$ in $\mathbb{P}(1,1,6,8,9)$, which implies
that $-K_{U_P}$ is big. Hence, the divisor $-K_{U_P'}$ is big as
well, which contradicts Theorem~\ref{theorem:Noether-Fano}.
\end{proof}

The exceptional divisor $E$ of the birational morphism $\alpha$
contains the singular point $Q_1$ of the threefold $U$ that is a
quotient singularity of type $\frac{1}{7}(1,1,6)$.
Lemmas~\ref{lemma:curves} and \ref{lemma:Cheltsov-Kawamata}  show
that the set $\mathbb{CS}(U, \frac{1}{n}\mathcal{M}_{U})$ does not
contain  any other subvariety of the surface $E$ than the point
$Q_1$.

The pencil $|-K_{U}|$ is the proper transform of the pencil
$|-K_{X}|$ and the base locus of $|-K_{U}|$ consists of the curve
$C_{U}$.

\begin{lemma}
\label{lemma:n-47-P5} The set $\mathbb{CS}(U,
\frac{1}{n}\mathcal{M}_{U})$ contains the point $Q_1$.
\end{lemma}

\begin{proof}
Suppose that it does not contain the point $Q_1$. Then, it follows
from Theorem~\ref{theorem:Noether-Fano} and
Lemma~\ref{lemma:Cheltsov-Kawamata} that the set $\mathbb{CS}(U,
\frac{1}{n}\mathcal{M}_{U})$ consists of the point $\bar{P}$ whose
image to $X$ is the point $P$ because the divisor $-K_{U}$ is nef
and big.

Let $\beta_P:W_P\to U$ be the Kawamata blow up at the point
$\bar{P}$ with weights $(1,2,3)$ and
$F_P\cong\mathbb{P}(1,2,3)$ be its exceptional divisor. The proof
of Lemma~\ref{lemma:n-47-P2} implies that  the log pair $(W_P,
\frac{1}{n}\mathcal{M}_{W_P})$ is terminal.

The base locus of the pencil $|-K_{W_P}|$ consists of the
irreducible curve $C_{W_P}$ and the unique irreducible curve
$\bar{L}$ of the linear system
$|\mathcal{O}_{\mathbb{P}(1,\,2,\,3)}(1)|$ on the surface $F_P$.
Hence, there is an antiflip $\zeta_2:W_P\dasharrow W_P'$ along the
curve $C_{W_P}$, the divisor $-K_{W_P'}$ is nef, and  the  log
pair $(W_P', \frac{1}{n}\mathcal{M}_{W_P'})$ is terminal.

The pull-backs of the rational functions $1$, $\frac{y}{x}$,
$\frac{zy}{x^{6}}$, $\frac{ty}{x^{8}}$ and
$\frac{wzy^{2}}{x^{15}}$ induce a birational map
$\chi_2:W_P\dasharrow X''$ such that $X''$ is a hypersurface of
degree $30$ in $\mathbb{P}(1,1,6,8,15)$.  We have
$$
\left\{\aligned &S_{W_P}\sim_{\mathbb{Q}}
(\alpha\circ\beta_P)^*(-K_X)-\frac{1}{8}E_{W_P} -\frac{1}{5}F_P,\\
&S_{W_P}^y\sim_{\mathbb{Q}}
(\alpha\circ\beta_P)^*(-K_X)-\frac{1}{8}E_{W_P} -\frac{6}{5}F_P,\\
&S_{W_P}^z\sim_{\mathbb{Q}}
(\alpha\circ\beta_P)^*(-5K_X)-\frac{13}{8}E_{W_P},\\
&S_{W_P}^t\sim_{\mathbb{Q}}
(\alpha\circ\beta_P)^*(-7K_X)-\frac{7}{8}E_{W_P} -\frac{2}{5}F_P,\\
&S_{W_P}^w\sim_{\mathbb{Q}}
(\alpha\circ\beta_P)^*(-8K_X) -\frac{3}{5}F_P,\\
\endaligned
\right.
$$
which implies that $X''$ is the anticanonical model of $W_P'$. In
particular, the divisor $-K_{W_P'}$ is big, which is impossible by
Theorem~\ref{theorem:Noether-Fano}.
\end{proof}

The exceptional divisor $F$ of the birational  morphism $\beta$
contains a singular point $Q_2$ of $W$ that is a quotient
singularity of type $\frac{1}{6}(1,1,5)$. The divisor $-K_{W}$ is
nef and big, and hence  the set $\mathbb{CS}(W,
\frac{1}{n}\mathcal{M}_{W})$ is not empty by
Theorem~\ref{theorem:Noether-Fano}.

\begin{lemma}
\label{lemma:n-47-P6} The set $\mathbb{CS}(W,
\frac{1}{n}\mathcal{M}_{W})$ contains the point $Q_2$.
\end{lemma}

\begin{proof}
Suppose that it does not contain the point $Q_2$. Then, it follows
from Theorem~\ref{theorem:Noether-Fano},
Lemmas~\ref{lemma:curves}, and \ref{lemma:Cheltsov-Kawamata} that
the set  $\mathbb{CS}(W, \frac{1}{n}\mathcal{M}_{W})$ consists of
the point $\tilde{P}$ whose image to $X$ is the point $P$.

Let $\pi:Y_P\to W$ be the Kawamata blow up at the point
$\tilde{P}$ and $G_P$ be its exceptional divisor. Then, the proof
of Lemma~\ref{lemma:n-47-P2} shows  that  the log pair $(Y_P,
\frac{1}{n}\mathcal{M}_{Y_P})$ is terminal.

We have
$$
\left\{\aligned &S_{Y_P}\sim_{\mathbb{Q}}
(\alpha\circ\beta\circ\pi)^*(-K_X)-\frac{1}{8}E_{Y_P}-\frac{1}{4}F_{Y_P} -\frac{1}{5}G_P,\\
&S_{Y_P}^y\sim_{\mathbb{Q}}
(\alpha\circ\beta\circ\pi)^*(-K_X)-\frac{1}{8}E_{Y_P}-\frac{1}{4}F_{Y_P}  -\frac{6}{5}G_P,\\
&S_{Y_P}^z\sim_{\mathbb{Q}}
(\alpha\circ\beta\circ\pi)^*(-5K_X)-\frac{13}{8}E_{Y_P} -\frac{9}{4}F_{Y_P},\\
&S_{Y_P}^t\sim_{\mathbb{Q}}
(\alpha\circ\beta\circ\pi)^*(-7K_X)-\frac{7}{8}E_{Y_P} -\frac{3}{4}F_{Y_P} -\frac{2}{5}G_P,\\
&S_{Y_P}^w\sim_{\mathbb{Q}}
(\alpha\circ\beta\circ\pi)^*(-8K_X) -\frac{3}{5}G_P,\\
\endaligned
\right.
$$
which implies that the anticanonical model of $Y_P$ is a
hypersurface of degree $42$  in $\mathbb{P}(1,1,6,14,21)$ with
canonical singularities because the pull-backs of the rational
functions $\frac{y}{x}$, $\frac{zy}{x^{6}}$,
$\frac{tzy^{2}}{x^{14}}$ and $\frac{wz^{2}y^{3}}{x^{21}}$ are
contained in the linear systems $|S_{Y_P}|$, $|6S_{Y_P}|$,
$|14S_{Y_P}|$ and $|21S_{Y_P}|$, respectively. In particular, the
divisor $-K_{Y_P}$ is big.

The base locus of the linear system $|-42K_{Y_P}|$ consists of the
curve $C_{Y_P}$  and $-K_{Y_P}\cdot C_{Y_P}<0$. Therefore, the
existence of the antiflip $\zeta_3:Y_P\dasharrow Y_P'$ along the
curve $C_{Y_P}$ follows from \cite{Sho93}. The divisor $-K_{Y_P'}$
is nef and big. However, this contradicts
Theorem~\ref{theorem:Noether-Fano} because the log pair $(Y_P',
\frac{1}{n}\mathcal{M}_{Y_P'})$ is terminal.
\end{proof}

The set $\mathbb{CS}(Y, \frac{1}{n}\mathcal{M}_{Y})$ is not empty.
If it contains the singular point $Q_3$ contained in the
exceptional divisor of $\gamma$, then it easily follows from
Theorem~\ref{theorem:main-tool} that $\mathcal{M}=|-K_{X}|$. If
the set $\mathbb{CS}(Y, \frac{1}{n}\mathcal{M}_{Y})$ does not
contain the point $Q_3$, then Lemmas~\ref{lemma:curves} and
\ref{lemma:Cheltsov-Kawamata}  imply that the set $\mathbb{CS}(Y,
\frac{1}{n}\mathcal{M}_{Y})$ consists of the point $\hat{P}$ whose
image to $X$ is the point $P$.

Let $\sigma:V\to Y$ be the Kawamata blow up at the point
$\hat{P}$. Then,  the proof of Lemma~\ref{lemma:n-47-P2} shows
that the log pair $(V, \frac{1}{n}\mathcal{M}_{V})$ is terminal.
The proof of Lemma~\ref{lemma:n-47-P6} implies that the
anticanonical model of $V$ is the surface $\mathbb{P}(1,1,6)$. On
the other hand, there is an antiflip $\zeta:V\dasharrow V'$ along
the curve $C_{V}$, which implies that the linear system
$|-rK_{V'}|$ induces a surjective morphism
$V'\to\mathbb{P}(1,1,6)$. However, it contradicts
Theorem~\ref{theorem:Noether-Fano} because the singularities of
the log pair $(V', \frac{1}{n}\mathcal{M}_{V'})$ are terminal.

Consequently, we have shown
\begin{proposition}
\label{proposition:n-47} If $\gimel=47$, then the linear system
$|-K_{X}|$ is the only Halphen pencil on $X$.
\end{proposition}

Next, we consider the case $\gimel=54$.\index{$\gimel=54$} The
variety $X$ is a general hypersurface of degree $24$ in
$\mathbb{P}(1,1,6,8,9)$ with $-K_{X}^{3}=\frac{1}{18}$. The
singularities of $X$ consist of one point that is a quotient
singularity of type $\frac{1}{2}(1,1,1)$, one point that is a
quotient singularity of type $\frac{1}{3}(1,1,2)$, and a point $Q$
that is a quotient singularity of type $\frac{1}{9}(1,1,8)$.

There is a commutative diagram
$$
\xymatrix{
&U\ar@{->}[d]_{\alpha}&&W\ar@{->}[ll]_{\beta}&&Y\ar@{->}[ll]_{\gamma}\ar@{->}[d]^{\eta}&\\%
&X\ar@{-->}[rrrr]_{\psi}&&&&\mathbb{P}(1,1,6),&}
$$
where \begin{itemize}

\item $\psi$ is the natural projection,

\item $\alpha$ is the Kawamata blow up at the point $Q$ with
weights $(1,1,8)$,

\item $\beta$ is the Kawamata blow up with weights $(1,1,7)$ at
the singular point of $U$ that is a quotient singularity of type
$\frac{1}{8}(1,1,7)$,

\item $\gamma$ is the Kawamata blow up  with weights $(1,1,6)$ at
the singular point of $W$ that is a quotient singularity of type
$\frac{1}{7}(1,1,6)$,

\item $\eta$ is and elliptic fibration.
\end{itemize}

If the set  $\mathbb{CS}(X, \frac{1}{n}\mathcal{M})$ contains
either one of the singular points of types $\frac{1}{2}(1,1,1)$
and $\frac{1}{3}(1,1,2)$ or a curve, then we obtain
$\mathcal{M}=|-K_X|$ from
Lemma~\ref{lemma:special-singular-points-with-zero-c} and
Corollary~\ref{corollary:Ryder-a2}. Therefore, due to
Lemma~\ref{lemma:smooth-points}, we may assume that
$$\mathbb{CS}\Big(X, \frac{1}{n}\mathcal{M}\Big)=\Big\{Q\Big\}.$$

\begin{proposition}
\label{proposition:n-54} If $\gimel=54$, then the linear system
$|-K_{X}|$ is the only Halphen pencil on $X$.
\end{proposition}
\begin{proof}
Unlike the case $\gimel=47$, the set $\mathbb{CS}(X,
\frac{1}{n}\mathcal{M})$ does not contain any other point than the
point $Q$, which makes our situation far simpler than the case
$\gimel=47$. In this case, Lemmas~\ref{lemma:n-47-P2},
\ref{lemma:n-47-P5}, \ref{lemma:n-47-P6} are automatically
satisfied, so that the proof of Proposition~\ref{proposition:n-47}
proves this statement.
\end{proof}

From now, we consider the case $\gimel=62$. The variety $X$ is a general hypersurface of degree $26$  in
$\mathbb{P}(1,1,5,7,13)$ with $-K_{X}^{3}=\frac{2}{35}$\index{$\gimel=62$}.  It has
two singular points. One is  a quotient singularity $P$ of type
$\frac{1}{5}(1,2,3)$ and the other is a quotient singularity $Q$
of type $\frac{1}{7}(1,1,6)$.

There is a commutative diagram
$$
\xymatrix{
&U\ar@{->}[d]_{\alpha}&&Y\ar@{->}[ll]_{\beta}\ar@{->}[d]^{\eta}&\\%
&X\ar@{-->}[rr]_{\psi}&&\mathbb{P}(1,1,5),&}
$$
where \begin{itemize} \item $\psi$ is the natural projection,

\item $\alpha$ is the Kawamata blow up at the point $Q$ with
weights $(1,1,6)$,

\item $\beta$ is the Kawamata blow up with weights $(1,1,5)$ at
the singular point of the variety $U$ that is a quotient
singularity of type $\frac{1}{6}(1,1,5)$,

\item $\eta$ is an elliptic fibration.
\end{itemize}

If the set $\mathbb{CS}(X, \frac{1}{n}\mathcal{M})$ contains a
curve, then $\mathcal{M}=|-K_X|$ by
Corollary~\ref{corollary:Ryder-a2}. Therefore, we may assume that
$$\mathbb{CS}\Big(X,\frac{1}{n}\mathcal{M}\Big)\subset \Big\{P, Q\Big\}$$
by Lemma~\ref{lemma:smooth-points}.

\begin{lemma}\label{lemma:n-62-P}
The set $\mathbb{CS}(X, \frac{1}{n}\mathcal{M})$ contains the point $Q$.
\end{lemma}
\begin{proof} Suppose that  the set $\mathbb{CS}(X, \frac{1}{n}\mathcal{M})$ does not contain the point $Q$. Then, it must consist of
the point $P$.
Let $\alpha_P:U_P\to X$ be the Kawamata blow up at the point $P$
and let $G_P$ be its exceptional divisor. Then, the proof is the same as that of Lemma~\ref{lemma:n-47-P2}.
From the equivalences
\begin{equation*}
\left\{\aligned
&\alpha_P^{*}\big(-K_{X}\big)\sim_{\mathbb{Q}} S_{U_P}+\frac{1}{5}G_P,\\
&\alpha_P^{*}\big(-K_{X}\big)\sim_{\mathbb{Q}}S_{U_P}^y+\frac{6}{5}G_P,\\
&\alpha_P^{*}\big(-5K_{X}\big)\sim_{\mathbb{Q}}S_{U_P}^z,\\
&\alpha_P^{*}\big(-7K_{X}\big)\sim_{\mathbb{Q}}S_{U_P}^t+\frac{2}{5}G_P,\\
&\alpha_P^{*}\big(-13K_{X}\big)\sim_{\mathbb{Q}} S_{U_P}^w+\frac{3}{5}G_P,.\\
\endaligned\right.
\end{equation*}
we see that the pull-backs of the rational functions $1$, $\frac{y}{x}$,
$\frac{zy}{x^{6}}$, $\frac{ty}{x^{8}}$ and $\frac{y^2w}{x^{15}}$
induce a birational map of $U_P$ onto  a hypersurface
of degree $30$ in $\mathbb{P}(1,1,6,8,15)$.
\end{proof}

The exceptional
divisor $E$ of the birational morphism $\alpha$ contains the
singular point $O$ of $U$ that is a quotient singularity of type
$\frac{1}{6}(1,1,5)$.

\begin{lemma}\label{lemma:n-62-antiflip} The set $\mathbb{CS}(U, \frac{1}{n}\mathcal{M}_U)$ contains the point $O$.
\end{lemma}
\begin{proof}
Suppose that the set $\mathbb{CS}(U, \frac{1}{n}\mathcal{M}_U)$ does not contain the point $O$. Then, it must contain the point $\bar{P}$ whose image to $X$ is the point $P$. Let $\beta_P:W_P\to U$ be the Kawamata blow up at the point
$\bar{P}$ with weights $(1,2,3)$  and $F_P\cong\mathbb{P}(1,2,3)$
be its exceptional divisor.
The proof of Lemma~\ref{lemma:n-62-P} implies that the log pair $(W_P, \frac{1}{n}\mathcal{M}_{W_P})$ is terminal.
 We have
\begin{equation*}
\left\{\aligned
&(\alpha\circ\beta_P)^{*}\big(-K_{X}\big)\sim_{\mathbb{Q}} S_{W_P}+\frac{1}{7}E_{W_P}+\frac{1}{5}F_P,\\
&(\alpha\circ\beta_P)^{*}\big(-K_{X}\big)\sim_{\mathbb{Q}}S_{W_P}^y+\frac{1}{7}E_{W_P}+\frac{6}{5}F_P,\\
&(\alpha\circ\beta_P)^{*}\big(-5K_{X}\big)\sim_{\mathbb{Q}}S_{W_P}^z+\frac{12}{7}E_{W_P},\\
&(\alpha\circ\beta_P)^{*}\big(-7K_{X}\big)\sim_{\mathbb{Q}}S_{W_P}^t+\frac{2}{5}F_P,\\
&(\alpha\circ\beta_P)^{*}\big(-13K_{X}\big)\sim_{\mathbb{Q}} S_{W_P}^w+\frac{6}{7}E_{W_P}+\frac{3}{5}F_P,.\\
\endaligned\right.
\end{equation*}
The equivalences imply
\[\left\{\aligned & (\alpha\circ\beta_P)^{*}\Big(\frac{y}{x}\Big)\in |S_{W_P}|, \\
&(\alpha\circ\beta_P)^{*}\Big(\frac{yz}{x^6}\Big)\in |6S_{W_P}|, \\
&(\alpha\circ\beta_P)^{*}\Big(\frac{yzt}{x^{14}}\Big)\in |14S_{W_P}|, \\
&(\alpha\circ\beta_P)^{*}\Big(\frac{y^3zw}{x^{21}}\Big)\in
|21S_{W_P}|,\\
\endaligned\right.\]
and hence the pull-backs of rational functions $1$, $\frac{y}{x}$, $\frac{yz}{x^6}$, $\frac{yzt}{x^{14}}$, and $\frac{y^3zw}{x^{21}}$ induce  a
birational map $\chi:W_P\dasharrow X^{\prime}$ such that
$X^{\prime}$ is a hypersurface in $\mathbb{P}(1,1,6,14,21)$.
Therefore, the divisor $-K_{W_P}$ is big.

The base locus of the pencil $|-K_{W_P}|$ consists of the
irreducible curves $C_{W_P}$ and $\bar{L}_{W_P}$, where $\bar{L}_{W_P}$
is the unique curve in the linear system
$|\mathcal{O}_{\mathbb{P}(1,\,2,\,3)}(1)|$ on $F_P$. It follows
from \cite{Sho93} that there is an antiflip
$\zeta:W_P\dasharrow W_P'$ along the curve  $C_{W_P}$. The divisor $-K_{W_P'}$ is
nef and big. The singularities of the log pair $(W_P',
\frac{1}{n}\mathcal{M}_{W_P'})$ are terminal because the rational
map $\zeta$ is a log flop with respect to the log pair $(W_P,
\frac{1}{n}\mathcal{M}_{W_P})$, which contradicts
Theorem~\ref{theorem:Noether-Fano}.
\end{proof}

\begin{proposition}
The linear system $|-K_X|$ is the only Halphen pencil on $X$.
\end{proposition}
\begin{proof}
Suppose  that the set $\mathbb{CS}(Y, \frac{1}{n}\mathcal{M}_Y)$
contains the singular point contained in the exceptional divisor $F$
of the birational morphism $\beta$. Let $\gamma_1: W\to Y$ be the
Kawamata blow up at this point. Then, the pencil $|-K_W|$ is the
proper transform of the pencil $|-K_X|$. Its base locus consists
of the irreducible curve $C_W$. Because
$\mathcal{M}_W\sim_{\mathbb{Q}}-nK_W$, the inequality $-K_W\cdot
C_W<0$ implies that $\mathcal{M}=|-K_X|$ by Theorem~\ref{theorem:main-tool}.

Now, we suppose that the set $\mathbb{CS}(Y, \frac{1}{n}\mathcal{M}_Y)$
consists of the point $\hat{P}$ whose image to $X$ is the point $P$. Let $\gamma: V\to Y$ be the Kawamata blow up at the point $\hat{P}$ with weights $(1,2,3)$. Then, the curve $C_V$ is the only curve that intersects with $-K_V$ negatively.
The proof of Lemma~\ref{lemma:n-62-P} implies that the log pair $(V, \frac{1}{n}\mathcal{M}_{V})$ is terminal.
Therefore, there is an antiflip
$\hat{\zeta}:V\dasharrow V'$ along the curve  $C_{V}$. The log pair $(V', \frac{1}{n}\mathcal{M}_{V'})$ is also terminal. However, for some positive integer $m$, the linear system $|-mK_{V'}|$ induces an elliptic fibration onto $\mathbb{P}(1,1,6)$, which contradicts Theorem~\ref{theorem:Noether-Fano}.
\end{proof}

\section{Case $\gimel=51$, hypersurface of degree $22$ in
$\mathbb{P}(1,1,4,6,11)$.}\index{$\gimel=51$} \label{section:n-51}

The threefold  $X$ is a general hypersurface of degree $22$ in
$\mathbb{P}(1,1,4,6,11)$ with $-K_{X}^{3}=\frac{1}{12}$. The
singularities of $X$ consist of one point  that is a quotient
singularity of type $\frac{1}{2}(1,1,1)$, one point $P$ that is a
quotient singularity of type $\frac{1}{4}(1,1,3)$, and one point
$Q$ that is a quotient singularity of type $\frac{1}{6}(1,1,5)$.

There is a commutative diagram
$$
\xymatrix{
&U\ar@{->}[d]_{\alpha}&&Y\ar@{->}[ll]_{\beta}\ar@{->}[d]^{\eta}&\\%
&X\ar@{-->}[rr]_{\psi}&&\mathbb{P}(1,1,4),&}
$$
where \begin{itemize}

\item $\psi$ is the natural projection,

\item $\alpha$ is the Kawamata blow up at the point $Q$ with
weights $(1,1,5)$,

\item $\beta$ is the Kawamata blow up with weights $(1,1,4)$ at
the singular point of the variety $U$ contained in the exceptional
divisor of the birational morphism $\alpha$,

\item $\eta$ is an elliptic fibration.
\end{itemize}
By the generality of the hypersurface $X$, it can be given by an
equation of the form
$$
z^{4}t+z^{3}f_{10}(x,y,t)+z^{2}f_{14}(x,y,t,w)+zf_{18}(x,y,t,w)+f_{22}(x,y,t,w)=0,
$$
where  $f_i$ is a quasihomogeneous polynomial of degree $i$. Let
$\xi:X\dasharrow\mathbb{P}^{7}$ be the rational map induced by the
linear systems  on the hypersurface $X$ defined by the equations
$$
\mu t+\sum_{i=0}^{6}\lambda_{i}x^{i}y^{6-i}=0,
$$
where $(\mu:
\lambda_{0}:\lambda_{1}:\lambda_{2}:\lambda_{3}:\lambda_{4}:\lambda_{5}:\lambda_{6})\in\mathbb{P}^{7}$.
Then, the closure of the image of the rational map $\xi$ is the
surface $\mathbb{P}(1,1,6)$ and the normalization of a general
fiber of $\xi$ is an elliptic curve. We also have the following
commutative diagram:
$$
\xymatrix{
&&W\ar@{->}[dl]_{\gamma}\ar@{->}[dr]^{\eta_0}&&\\
&X\ar@{-->}[rr]_{\xi}&&\mathbb{P}(1,1,6),&}
$$
where \begin{itemize} \item $\gamma$ is the Kawamata blow up at
the point $P$ with weights $(1,1,3)$,

\item $\eta_0$ is an elliptic fibration.

\end{itemize}

If the set $\mathbb{CS}(X, \frac{1}{n}\mathcal{M})$ contains
either the singular point of type $\frac{1}{2}(1,1,1)$ or a curve,
then  $\mathcal{M}=|-K_{X}|$ by
Lemma~\ref{lemma:special-singular-points-with-zero-c} and
Corollary~\ref{corollary:Ryder-a2}.  Therefore,
Lemma~\ref{lemma:smooth-points} enables us to assume that
$\mathbb{CS}(X, \frac{1}{n}\mathcal{M})\subseteq\{P, Q\}$

\begin{lemma}
\label{lemma:n-51-P2-P3} If $\mathbb{CS}(X,
\frac{1}{n}\mathcal{M})=\{P, Q\}$, then $\mathcal{M}=|-K_{X}|$.
\end{lemma}

\begin{proof}
Let $\pi_1:V_1\to W$ be the Kawamata blow up with weights
$(1,1,5)$ at the point whose image to $X$ is the point $Q$. Then,
$\mathcal{M}_{V_1}\sim_{\mathbb{Q}} -nK_{V_1}$ by
Lemma~\ref{lemma:Kawamata}, the pencil $|-K_{V_1}|$ is the proper
transform of the pencil $|-K_{X}|$, and the base locus of the
pencil $|-K_{V_1}|$ consists of the irreducible curve $C_{V_1}$.
Because $-K_{V_1}\cdot C_{V_1}<0$, we obtain
$\mathcal{M}=|-K_{X}|$ from Theorem~\ref{theorem:main-tool}.
\end{proof}

\begin{lemma}
\label{lemma:n-51-P3} If $\mathbb{CS}(X,
\frac{1}{n}\mathcal{M})=\{P \}$, then  $\mathcal{M}=|-K_{X}|$.
\end{lemma}

\begin{proof}
Let $P_1$ be the singular point of $W$ contained in the
exceptional divisor of $\gamma$. Then, $P_{1}\in\mathbb{CS}(W,
\frac{1}{n}\mathcal{M}_{W})$ by
Theorem~\ref{theorem:Noether-Fano}, Lemmas~\ref{lemma:Kawamata},
and \ref{lemma:Cheltsov-Kawamata}.

Let $\pi_2:V_2\to W$ be the Kawamata blow up at the point $P_{1}$.
Then, $M_{V_2}\sim_{\mathbb{Q}} -nK_{V_2}$ by
Lemma~\ref{lemma:Kawamata}, the pencil $|-K_{V_2}|$ is the proper
transform of $|-K_{X}|$, and the base locus of $|-K_{V_2}|$
consists of the irreducible curve $C_{V_2}$.  The inequality
$-K_{V_2}\cdot C_{V_2}<0$ implies $\mathcal{M}=|-K_{X}|$ by
Theorem~\ref{theorem:main-tool}.
\end{proof}

\begin{proposition}
\label{proposition:n-51} The linear system $|-K_{X}|$ is the only
pencil on $X$.
\end{proposition}
\begin{proof}
Due to the previous arguments, we may assume that $\mathbb{CS}(X,
\frac{1}{n}\mathcal{M})=\{Q\}$, which implies that
$\mathcal{M}_{U}\sim_{\mathbb{Q}} -nK_{U}$ by
Lemma~\ref{lemma:Kawamata}. Hence, it follows from
Theorem~\ref{theorem:Noether-Fano},
Lemmas~\ref{lemma:Cheltsov-Kawamata}, and \ref{lemma:curves} that
the set $\mathbb{CS}(U, \frac{1}{n}\mathcal{M}_{U})$ consists of
the singular point of the threefold $U$ that is contained in the
exceptional divisor of $\alpha$ because  the divisor $-K_{U}$ is
nef and big.

Let $P_{2}$ be the singular point of $Y$ contained in the
exceptional divisor of the birational morphism $\beta$ and let
$\pi:V\to Y$ be the Kawamata blow up at the point $P_{2}$. Then,
$M_{V}\sim_{\mathbb{Q}} -nK_{V}$ by
Theorem~\ref{theorem:Noether-Fano}, Lemmas~\ref{lemma:Kawamata},
\ref{lemma:Cheltsov-Kawamata}, and \ref{lemma:curves}. On the
other hand, the pencil $|-K_{V}|$ is the proper transform of the
pencil $|-K_{X}|$ and the base locus of the pencil $|-K_{V}|$
consists of the irreducible curve $C_V$. Due to
Theorem~\ref{theorem:main-tool}, the inequality $-K_{V}\cdot
C_V<0$ completes our proof.
\end{proof}

\section{Case $\gimel=82$, hypersurface of degree $36$ in
$\mathbb{P}(1,1,5,12,18)$.}\index{$\gimel=82$}\label{section:n-82}

The threefold $X$ is a general hypersurface of degree $36$ in
$\mathbb{P}(1,1,5,12,18)$ with $-K_X^3=\frac{1}{30}$. Its
singularities consist of  two quotient singular points $P$ and $Q$
of types $\frac{1}{5}(1,2,3)$ and $\frac{1}{6}(1,1,5)$,
respectively. The hypersurface $X$ can be given by the equation
$$
z^{7}y+\sum_{i=0}^{6}z^{i}f_{36-5i}\big(x,y,z,t\big)=0,
$$
where  $f_{i}$ is a qua\-si\-ho\-mo\-ge\-ne\-ous polynomial of
degree $i$.

There is a commutative diagram
$$
\xymatrix{
&&U\ar@{->}[dl]_{\alpha}\ar@{->}[dr]^{\eta}&&\\
&X\ar@{-->}[rr]_{\psi}&&\mathbb{P}(1,1,5),&}
$$
where
\begin{itemize}
\item $\psi$ is the natural projection,

\item $\alpha$ is the Kawamata blow up of at the point $Q$ with
weights $(1,1,5)$,

\item $\eta$ is an elliptic fibration.
\end{itemize}

If the set $\mathbb{CS}(X,\frac{1}{n}\mathcal{M})$ contains a
curve, then we obtain $\mathcal{M}=|-K_X|$ from
Corollary~\ref{corollary:Ryder-a2}. Therefore, we may assume that
$$
\varnothing\ne\mathbb{CS}\Big(X,\frac{1}{n}\mathcal{M}\Big)\subseteq\Big\{P, Q\Big\}%
$$
due to Lemma~\ref{lemma:smooth-points}.

The exceptional divisor of the birational morphism $\alpha$
contains a singular point $O$ of the threefold $U$ that is a
quotient singularity of type $\frac{1}{5}(1,1,4)$.
\begin{lemma}
\label{lemma:n-82-good-point} If the set
$\mathbb{CS}(U,\frac{1}{n}\mathcal{M}_U)$ contains  the point $O$,
then $\mathcal{M}=|-K_{X}|$.
\end{lemma}
\begin{proof}
Let $\alpha_O:U_O\to U$ be the Kawamata blow up at the point $O$.
Then, the pencil $|-K_{U_O}|$ is the proper transform of the
pencil $|-K_X|$. Its base locus consists of $C_{U_O}$.
Furthermore, $-K_{U_O}\cdot C_{U_O}<0$ and
$\mathcal{M}_{U_O}\sim_{\mathbb{Q}} -nK_{U_O}$. Therefore,
Theorem~\ref{theorem:main-tool} implies the statement.
\end{proof}

Therefore, if the set the set
$\mathbb{CS}(X,\frac{1}{n}\mathcal{M})$ consists of the point $Q$,
then $\mathcal{M}=|-K_{X}|$, and hence we may assume that
$P\in\mathbb{CS}(X,\frac{1}{n}\mathcal{M})$.

 Let $\pi:W\to X$ be the Kawamata blow up at the point $P$ with
weights $(1,2,3)$ and $E\cong\mathbb{P}(1,2,3)$ be its exceptional
divisor. Then, the exceptional divisor $E$ contains two singular
points $P_1$ and $P_2$ that are quotient singularities of the
threefold $W$ of types $\frac{1}{2}(1,1,1)$ and
$\frac{1}{3}(1,2,1)$, respectively. Let $L$ be the unique curve of
the linear system $|\mathcal{O}_{\mathbb{P}(1,\,2,\,3)}(1)|$ on
the surface $E$.

The linear system $|-K_{W}|$ is the proper transform of $|-K_{X}|$
and its base locus  consists of the irreducible curves $L_W$ and
$C_W$. The divisor $-K_{W}$ is not nef because $-K_{W}\cdot
C_W<0$. The curve $C_W$ is the only curve that has negative
intersection with the divisor $-K_{W}$ because $-K_{W}\cdot
L_W>0$.

Let $\mathcal{D}$ be the proper transform of the linear system on
the threefold $X$ that is cut by
$$
h_{21}\big(x, y\big)+zyh_{15}\big(x, y\big)+zty^{3}h_{1}\big(x, y\big)+z^{2}y^{2}h_{9}\big(x, y\big)+
$$
$$+z^{3}y^{3}h_{3}\big(x, y\big)+ty^{2}h_{7}\big(x, y\big)+h_{0}wy^{3}=0%
$$
where $h_{i}$ is a homogeneous polynomial of degree $i$. Then,
$\mathcal{D}\sim_{\mathbb{Q}} -21K_{W}$ but  $\mathcal{D}$ induces
a bi\-ra\-ti\-o\-nal map $\gamma:W\dasharrow X^{\prime}$ such that
$X^{\prime}$ is a hypersurface of degree $42$ in
$\mathbb{P}(1,1,6,14,21)$ with canonical singularities (see the
proof of Theorem~5.5.1 in \cite{CPR}).

There is an antiflip $\xi:W\dasharrow W'$ along the curve $L_W$.
The divisor $-K_{W'}$ is nef and big and the linear system
$|-mK_{W'}|$ induced a birational morphism $\beta:W'\to
X^{\prime}$ for some natural number $m\gg 0$. Therefore, we have a
commutative diagram
$$
\xymatrix{
&W\ar@{-->}[rrrr]^{\xi}\ar@{->}[dl]_{\pi}\ar@{-->}[drrrrr]^{\gamma}&&&&W'\ar@{->}[dr]^{\beta}&\\
X\ar@{-->}[rrrrrr]^{\omega}\ar@{^{(}->}[dr]&&&&&&X^{\prime}\ar@{^{(}->}[dl]&\\
&\mathbb{P}\big(1,1,5,12,18\big)\ar@{-->}[rrrr]^{(x,\ y,\ zy,\ ty^{2},\ wy^{3})}&&&&\mathbb{P}\big(1,1,6,14,21\big),&}%
$$
where $\omega$ is the rational map induced by the linear system
$\pi(\mathcal{D})$.

\begin{lemma}
\label{lemma:n-82-non-terminal-singularities} The log pair
$(W,\frac{1}{n}\mathcal{M}_W)$ is not terminal.
\end{lemma}

\begin{proof}
Suppose that the log pair $(W,\frac{1}{n}\mathcal{M}_W)$ is
terminal. Then, the log pair $(W, \lambda\mathcal{M}_W)$  is
terminal for some rational number $\lambda>\frac{1}{n}$.  We have
$\mathcal{M}_{W'}\sim_{\mathbb{Q}} -nK_{W'}$, but $\xi$ is a log
flip for $(W, \lambda\mathcal{M}_W)$. Therefore, the
sin\-gu\-la\-ri\-ties of the log pair $(W',
\frac{1}{n}\mathcal{M}_{W'})$ are terminal, which contradicts
Theorem~\ref{theorem:Noether-Fano}.
\end{proof}

\begin{lemma}
\label{lemma:n-82-terminal-singularities} If the set
$\mathbb{CS}(W,\frac{1}{n}\mathcal{M}_W)$ contains the point
$P_{i}$, then $\mathcal{M}=|-K_{X}|$.
\end{lemma}

\begin{proof}
Let $\alpha_i:V_i\to W$ be the Kawamata blow up  at the point
$P_i$ with weights $(1,i,1)$ and $F$ be the exceptional divisor of
the birational morphism $\alpha_i$. The pencil $|-K_{V_i}|$ is the
proper transform of the pencil $|-K_X|$. Its base locus of the
pencil $|-K_{V_i}|$ consists of two irreducible curves $C_{V_i}$
and $L_{V_i}$.

On a general surface $D$ in $|-K_{V_i}|$, we have
$$
C_{V_i}^2=-\frac{7}{6},\ L_{V_1}^2=-\frac{4}{3}, \  L_{V_2}^2=-1,
\ C_{V_i}\cdot L_{V_i}=1.%
$$
 The surface $D$ is normal and the intersection form of the curves
$C_{V_i}$ and $L_{V_i}$ on the surface $D$ is negative-definite.
On the other hand, we have
$$
\mathcal{M}_{V_i}\Big\vert_{D}\equiv-nK_{V_i}\Big\vert_{D}\equiv nC_{V_i}+nL_{V_i}.%
$$
Therefore, $\mathcal{M}=|-K_X|$ by
Theorem~\ref{theorem:main-tool}.
\end{proof}

From now, we may assume that $\mathbb{CS}(X,
\frac{1}{n}\mathcal{M})=\{P,Q\}$.

Let $\bar{Q}$ be the point of the threefold $W$ such that
$\pi(\bar{Q})=Q$. Then,
$$
\mathbb{CS}\Big(W,\frac{1}{n}\mathcal{M}_W\Big)=\Big\{\bar{Q}\Big\}
$$
by Lemmas~\ref{lemma:n-82-non-terminal-singularities},
\ref{lemma:n-82-terminal-singularities}, and
\ref{lemma:Cheltsov-Kawamata}.

Let $\gamma:V\to W$ be the blow up of $\bar{Q}$ with weights
$(1,1,5)$ and $F$ be its exceptional divisor. Then, the base locus
of the pencil $|-K_{V}|$ consists of two irreducible curves $C_V$
and $L_V$.

The arguments used in the proof of
Lemma~\ref{lemma:n-82-terminal-singularities} together with
Lemma~\ref{lemma:Cheltsov-Kawamata} imply that either the
singularities of the log pair $(V, \frac{1}{n}\mathcal{M}_V)$ are
terminal or  $\mathcal{M}=|-K_X|$. We may assume that the log pair
$(V, \frac{1}{n}\mathcal{M}_V)$ is terminal, which implies that
the log pair $(V, \lambda\mathcal{M}_V)$ is still terminal for
some rational number $\lambda>\frac{1}{n}$.

Let $\mathcal{L}$ be the proper transform on the threefold $V$ of
the linear system
$$
\mu_{0} x^{6}+\mu_{1} y^{6}+\mu_{2} yz=0,
$$
where $(\mu_{0}:\mu_{1}:\mu_{2})\in\mathbb{P}^{2}$. Then,
$\mathcal{L}\sim_{\mathbb{Q}} -6K_{V}$ and  the base locus of the
linear system $\mathcal{L}$ consists of the curve $C_V$. So, the
curve $C_V$ is the only curve having negative intersection with
$-K_{V}$.

The linear system $\mathcal{L}$ induces a ra\-ti\-o\-nal map
$\upsilon:V\dasharrow\mathbb{P}(1,1,6)$ whose general fiber is an
elliptic curve. Moreover, it follows from \cite{Sho93} that there
is a log flip $\chi:V\dasharrow V^{\prime}$ along the curve $C_V$
with respect to the log pair $(V, \lambda\mathcal{M}_V)$. The log
pair $(V^{\prime}, \frac{1}{n}\chi(\mathcal{M}_V))$ is terminal
and
$$
\chi\big(\mathcal{M}_V\big)\equiv -nK_{V^{\prime}},
$$
which implies that $-K_{V^{\prime}}$ is nef. Therefore, it follows
from Theorem~3.1.1 in \cite{KMM} that the linear system
$|-mK_{V^{\prime}}|$ is base-point-free for some natural number
$m\gg 0$.  However, the equivalence $\chi(\mathcal{L})\equiv
-6K_{V^{\prime}}$ implies that the linear system
$|-mK_{V^{\prime}}|$ induces an elliptic fibration, which
contradicts  Theorem~\ref{theorem:Noether-Fano}.

Consequently, we have proved
\begin{proposition}
\label{proposition:n-82} The linear system $|-K_{X}|$ is a unique
Halphen pencil on $X$.
\end{proposition}

\newpage
\part{Fano threefold hypersurfaces with a single Halphen
pencil.}\label{section:one-Halphen-2}

\section{Case $\gimel=19$, hypersurface of degree $12$ in
$\mathbb{P}(1,2,3,3,4)$.}\index{$\gimel=19$} \label{section:n-19}

The threefold $X$ is a general hypersurface of degree $12$ in
$\mathbb{P}(1,2,3,3,4)$ with $-K_X^3=\frac{1}{6}$.  It has seven
singular points. Four points $P_1$, $P_2$, $P_3$, $P_4$ of them
are quotient singularities of type $\frac{1}{3}(1,2,1)$ and the
others are quotient singularities of type $\frac{1}{2}(1,1,1)$.

For each point $P_i$, we have the following diagram:
$$
\xymatrix{
&&Y_{i}\ar@{->}[dl]_{\pi_{i}}\ar@{->}[dr]^{\eta_{i}}&\\
&X\ar@{-->}[rr]_{\psi_{i}}&&\mathbb{P}(1,2,3),&&}
$$
where \begin{itemize}

\item $\psi_i$ is  a projection,

\item $\pi_{i}$ is the Kawamata blow up at the point $P_{i}$ with
weights $(1,2,1)$, \item $\eta_{i}$ is an elliptic fibration.
\end{itemize}
It follows from Corollary~\ref{corollary:Ryder-a1} and
Lemmas~\ref{lemma:smooth-points},
\ref{lemma:special-singular-points-with-positive-c} that
$$
\mathbb{CS}\Big(X, \frac{1}{n}\mathcal{M}\Big)\subseteq\Big\{P_{1}, P_{2}, P_3, P_4\Big\}.%
$$

\begin{lemma}
The set $\mathbb{CS}(X, \frac{1}{n}\mathcal{M})$ cannot consist of
a single point.
\end{lemma}
\begin{proof}
Suppose that it contains only the point $P_i$. Then, the set
$\mathbb{CS}(Y_i, \frac{1}{n}\mathcal{M}_{Y_i})$ must contain the
singular point $O_i$ contained in the exceptional divisor $E_i$ of
the birational  morphism $\pi_i$. Let $\alpha_i: V_i\to Y_i$ be
the Kawamata blow up at the point $O_i$. We consider the linear
system $\mathcal{D}$ on $X$ defined by the equations
\[\lambda_0x^4+\lambda_1x^2y+\lambda_2w=0,\]
where $(\lambda_0: \lambda_1: \lambda_2)\in\mathbb{P}^2$. The base
locus of the linear system $\mathcal{D}_{V_i}$ consists of the
irreducible curve $\hat{C}_{V_i}$ whose image to $X$ is the base
curve of $\mathcal{D}$. A general surface $D_i$ in
$\mathcal{D}_{V_i}$ is normal and $\hat{C}_{V_i}^2=-\frac{1}{6}$
on the surface $D_i$, which implies $\mathcal{M}=\mathcal{D}$. It
is a contradiction since $\mathcal{D}$ is not a pencil.
\end{proof}

\begin{proposition}
The linear system $|-2K_X|$ is the only Halphen pencil on $X$.
\end{proposition}

\begin{proof}
Suppose that the set $\mathbb{CS}(X, \frac{1}{n}\mathcal{M})$
contains  the points $P_i$ and $P_j$, where $i\ne j$. Then, the
set $\mathbb{CS}(Y_i, \frac{1}{n}\mathcal{M}_{Y_i})$ must contain
the singular point $Q_j$ whose image to $X$ is the point $P_j$.
Let $\beta_j: W_j\to Y_i$ be the Kawamata blow up at the point
$Q_j$. Then, we see that the linear system $|-2K_{W_j}|$ is the
proper transform of the pencil $|-2K_X|$.  The base locus of the
pencil $|-2K_{W_j}|$ consists of the irreducible curve $C_{W_j}$
whose image to $X$ is the base curve of $|-2K_X|$. A general
surface $D$ in $|-2K_{W_j}|$ is normal and
$C_{W_j}^2=-\frac{1}{3}$ on the surface $D$, which implies
$\mathcal{M}_{W_j}=|-2K_{W_j}|$. Therefore, $\mathcal{M}=|-2K_X|$.
\end{proof}

\section{Cases $\gimel=23$ and $44$.}\index{$\gimel=23$} \label{section:n-23}

For the case $\gimel=23$, let $X$ be a general hypersurface of
degree $14$ in $\mathbb{P}(1,2,3,4,5)$ with $-K_X^3=\frac{7}{60}$.
It has  three quotient singularities of type $\frac{1}{2}(1,1,1)$,
one quotient singularity of type $\frac{1}{3}(1,2,1)$, one
quotient singularity $P$ of type $\frac{1}{5}(1,2,3)$, and one
quotient singularity $Q$ of type $\frac{1}{4}(1,3,1)$.

We have  the following elliptic fibration:
$$
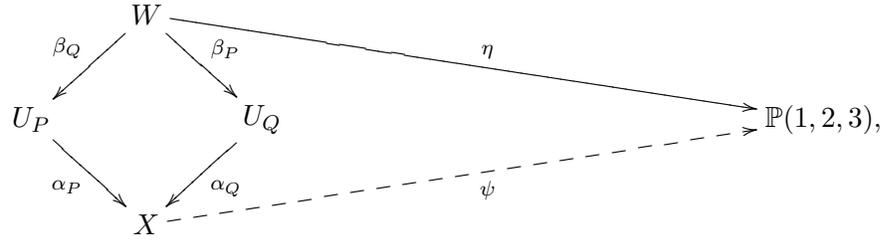
\xymatrix{
&&W\ar@{->}[dl]_{\beta_{Q}}\ar@{->}[dr]^{\beta_{P}}\ar@{->}[drrrrrrr]^{\eta}&&&&&&&\\
&U_{P}\ar@{->}[dr]_{\alpha_{P}}&&U_{Q}\ar@{->}[dl]^{\alpha_{Q}}&&&&&&\mathbb{P}(1,2,3),\\
&&X\ar@{-->}[rrrrrrru]_{\psi}&&&&&&&}
$$
where \begin{itemize}

\item $\psi$ is the natural projection,

 \item $\alpha_{P}$ is the Kawamata blow up at the point
$P$ with weights $(1,2,3)$,

\item $\alpha_{Q}$ is the Kawamata blow up at the point $Q$ with
weights $(1,3,1)$,

\item $\beta_{Q}$ is the Kawamata blow up with weights $(1,3,1)$
at the point whose image to $X$ is the point $Q$,

\item $\beta_{P}$ is the Kawamata blow up with weights $(1,2,3)$
at the point whose image to $X$ is the point $P$,

\item $\eta$ is an elliptic fibration.

\end{itemize}

Because we mainly consider the pencil $|-2K_X|$, we let $D$ be a
general surface in $|-2K_X|$.

 It follows from
Corollary~\ref{corollary:Ryder-a2} and
Lemmas~\ref{lemma:smooth-points},
\ref{lemma:special-singular-points-with-positive-c}, and
\ref{lemma:special-singular-points-with-zero-c} that  we may
assume that $\mathbb{CS}(X, \frac{1}{n}\mathcal{M})\subset \{P,
Q\}$.

\begin{lemma}
\label{lemma:n-23-first} If the log pair $(X,
\frac{1}{n}\mathcal{M})$ is terminal at the point $P$, then
$\mathcal{M}=|-2K_X|$.
\end{lemma}
\begin{proof}
Suppose that the log pair $(X, \frac{1}{n}\mathcal{M})$ is
terminal at the point $P$. Then, the set $\mathbb{CS}(U_Q,
\frac{1}{n}\mathcal{M}_{U_Q})$ contains the quotient singular
point $Q_1$ of type $\frac{1}{2}(1,1,1)$ on $U_Q$  contained in
the exceptional divisor $E$ of $\alpha_Q$. Let $\gamma_Q:W_Q\to
U_Q$ be the Kawamata blow up at the singular point $Q_1$ with
weights $(1,1,1)$.

Then, $|-2K_{W_Q}|$ is the proper transform of the pencil
$|-2K_{X}|$ and the base locus of the pencil $|-2K_{W_Q}|$
consists of the irreducible curve $C_{W_Q}$. We can easily check
$D_{W_Q}\cdot C_{W_Q}=-4K_{W_Q}^3<0$. Hence, we have
$\mathcal{M}=|-2K_{X}|$ by Theorem~\ref{theorem:main-tool} since
$\mathcal{M}_{W_Q}\sim_{\mathbb{Q}} -nK_{W_Q}$ by
Lemma~\ref{lemma:Kawamata}.
\end{proof}

Therefore, we may assume that the set $\mathbb{CS}(X,
\frac{1}{n}\mathcal{M})$ contains the point $P$.  Let $F$ be the
exceptional divisor of  the birational morphism $\alpha_P$. Then,
it contains two singular points $O_1$ and $P_1$ of types
$\frac{1}{2}(1,1,1)$ and $\frac{1}{3}(1,2,1)$, respectively.
\begin{lemma}
\label{lemma:n-23-second} The set $\mathbb{CS}(U_P,
\frac{1}{n}\mathcal{M}_{U_P})$ cannot contain the point $O_1$.
\end{lemma}
\begin{proof}
Suppose that  the set $\mathbb{CS}(U_P,
\frac{1}{n}\mathcal{M}_{U_P})$  contains the point $O_1$. Let
$\gamma_P:W_P\to U_P$ be the Kawamata blow up at the point $O_1$
with weights $(1,1,1)$. Let $G_1$ be the exceptional divisor of
the birational morphism  $\gamma_P$. Then, we have
\[\mathcal{M}_{W_P}\sim_{\mathbb{Q}} -nK_{W_P}\sim_{\mathbb{Q}} \gamma_P^*(-nK_{U_P})-\frac{n}{2}G_1.\]
Let $\mathcal{B}$ be the proper transform of the linear system
$|-3K_X|$ by the birational morphism $\alpha_P\circ\gamma_P$. We
then see that $\mathcal{B}\sim_{\mathbb{Q}}
\gamma_P^*(-3K_{U_P})-\frac{1}{2}G_1$. The base curve of
$\mathcal{B}$ consists of the proper transform $C_{W_P}$.
Furthermore, $B\cdot C_{W_P}=\frac{1}{12}$, where $B$ is a general
surface in $\mathcal{B}$. However, we obtain a contradictory
inequalities
\[\begin{split}
-\frac{n^2}{4}&=(\gamma_P^*(-3K_{U_P})-\frac{1}{2}G_1)\cdot(\gamma_P^*(-nK_{U_P})-\frac{n}{2}G_1)^2\\
&=(\gamma_P^*(-3K_{U_P})-\frac{1}{2}G_1)\cdot M_1\cdot M_2\geq
0,\end{split}\] where $M_1$ and $M_2$ are general surfaces in
$\mathcal{M}_{W_P}$.
\end{proof}

\begin{lemma}
\label{lemma:n-23-third} If the set $\mathbb{CS}(U_P,
\frac{1}{n}\mathcal{M}_{U_P})$  contains the point $P_1$, then
$\mathcal{M}=|-2K_X|$.
\end{lemma}

\begin{proof}
The proof of Lemma~\ref{lemma:n-23-first} implies the statement.
\end{proof}

\begin{proposition}
\label{proposition:n-23} If $\gimel=23$, then the linear system
$|-2K_X|$ is a unique Halphen pencil on $X$.
\end{proposition}
\begin{proof}
By the previous arguments, we may assume that the set
$\mathbb{CS}(U_P, \frac{1}{n}\mathcal{M}_{U_P})$ consists of the
single point $O$ such that $\alpha_P(O)=Q$. The set
$\mathbb{CS}(W, \frac{1}{n}\mathcal{M}_W)$ must contain the
singular point in the exceptional divisor of the birational
morphism $\beta_Q$. Then, the proof of
Lemma~\ref{lemma:n-23-first} completes the proof.
\end{proof}

From now on, we consider the case $\gimel=44$.\index{$\gimel=44$}
The threefold $X$ is a  general hypersurface of degree $20$ in
$\mathbb{P}(1,2,5,6,7)$ with $-K_{X}^{3}=\frac{1}{21}$. The
singularities of the hypersurface $X$ consist of three quotient
singular point of type $\frac{1}{2}(1,1,1)$,  one quotient
singular point $P$  of type $\frac{1}{7}(1,2,5)$, and one quotient
singular point $Q$ of type $\frac{1}{6}(1,5,1)$.

We have  the following elliptic fibration:
$$
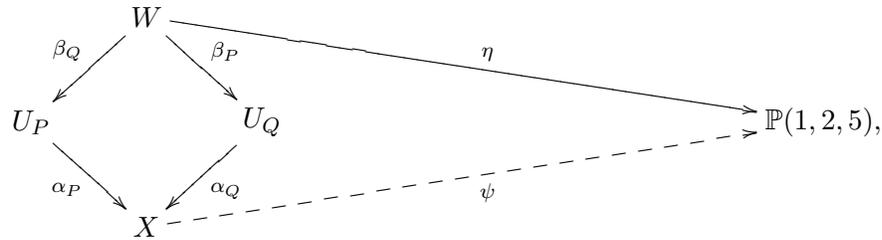
\xymatrix{
&&W\ar@{->}[dl]_{\beta_{Q}}\ar@{->}[dr]^{\beta_{P}}\ar@{->}[drrrrrrr]^{\eta}&&&&&&&\\
&U_{P}\ar@{->}[dr]_{\alpha_{P}}&&U_{Q}\ar@{->}[dl]^{\alpha_{Q}}&&&&&&\mathbb{P}(1,2,5),\\
&&X\ar@{-->}[rrrrrrru]_{\psi}&&&&&&&}
$$
where \begin{itemize}

\item $\psi$ is the natural projection,

 \item $\alpha_{P}$ is the Kawamata blow up at the point
$P$ with weights $(1,2,5)$,

\item $\alpha_{Q}$ is the Kawamata blow up at the point $Q$ with
weights $(1,5,1)$,

\item $\beta_{Q}$ is the Kawamata blow up with weights $(1,5,1)$
at the point whose image to $X$ is the point $Q$,

\item $\beta_{P}$ is the Kawamata blow up with weights $(1,2,5)$
at the point whose image to $X$ is the point $P$,

\item $\eta$ is an elliptic fibration.

\end{itemize}
The hypersurface $X$ can be given by the quasihomogeneous equation
of degree $20$ as follows:
\[w^2t+wf_{13}(x,y,z,t)+f_{20}(x,y,z,t)=0,\]
where $f_i$ is a quasihomogeneous polynomial of degree $i$.
Consider the linear subsystem of the linear system $|-6K_{X}|$
defined by equations
$$
\lambda_{0}t+\lambda_{1}y^{3}+\lambda_{2}y^{2}x^{2}+\lambda_{3}yx^{4}+\lambda_{4}x^{6}=0,
$$
where
$(\lambda_{0}:\lambda_{1}:\lambda_{2}:\lambda_{3}:\lambda_{4})\in\mathbb{P}^{4}$.
It gives us a rational map $\xi:X\dasharrow\mathbb{P}^{4}$ that is
 defined in the outside of  the point $P$. The closure of the image of
the rational map $\xi$ is the surface $\mathbb{P}(1,1,3)$, which
can be identified with a cone over a smooth rational curve  of
degree $3$ in $\mathbb{P}^{3}$. Moreover, the normalization of a
general fiber of the rational map $\xi$ is an elliptic curve.
Therefore, we have another elliptic fibration as follows:
$$
\xymatrix{
&U_{P}\ar@{->}[d]_{\alpha_{P}}&&Y\ar@{->}[ll]_{\gamma}\ar@{->}[d]^{\eta_0}&\\%
&X\ar@{-->}[rr]_{\xi}&&\mathbb{P}(1,1,3),&}
$$
where \begin{itemize}

\item $\gamma$ is the Kawamata blow up with weights $(1,2,3)$ at
the singular point of $U_{P}$ that is a quotient singularity of
type $\frac{1}{5}(1,2,3)$  contained in the exceptional divisor of
$\alpha_{P}$,

\item $\eta_0$ is an elliptic fibration.
\end{itemize}

It follows from Lemmas~\ref{lemma:smooth-points} and
\ref{lemma:special-singular-points-with-positive-c} that the set
$\mathbb{CS}(X, \frac{1}{n}\mathcal{M})\subseteq\{P, Q\}$.

The only different part from the proof for the case $\gimel=23$ is
Lemma~\ref{lemma:n-23-third}. Therefore, it is enough to show that
if the set $\mathbb{CS}(U_P, \frac{1}{n}\mathcal{M}_{U_P})$
contains the quotient singularity $P_1$ of type
$\frac{1}{5}(1,2,3)$ contained in the exceptional divisor of
$\alpha_P$, then $\mathcal{M}=|-2K_X|$. Suppose that the set
contains the point $P_1$. Then, the set $\mathbb{CS}(Y,
\frac{1}{n}\mathcal{M}_{Y})$ contains a singular point. The
exceptional divisor $F$ of the birational morphism $\gamma$
contains two singular points $P_1'$ and $P_2'$ of $Y$ that are
quotient singularities of types $\frac{1}{2}(1,1,1)$ and
$\frac{1}{3}(1,2,1)$.

By Lemma~\ref{lemma:n-23-first}, we may assume that the set
$\mathbb{CS}(Y, \frac{1}{n}\mathcal{M}_{Y})\subset \{ P_1',
P_2'\}$.

\begin{lemma}\label{lemma:n-44}
The set $\mathbb{CS}(Y, \frac{1}{n}\mathcal{M}_{Y})$ cannot
contain the point $P_1'$.
\end{lemma}
\begin{proof}
Suppose that the set $\mathbb{CS}(Y, \frac{1}{n}\mathcal{M}_{Y})$
 contains the point $P_1'$. Then, we consider the Kawamata blow up
$\sigma_1:Y_1\to Y$ at the point $P_1'$ with weights $(1,1,1)$.
Let $\mathcal{T}$ be the linear subsystem of the linear system
$|-6K_X|$ defined by the equations
\[\lambda_0t+\lambda_1x^6+\lambda_2 x^4y+\lambda_3x^2y^2=0,\]
where $(\lambda_0: \lambda_1: \lambda_2:
\lambda_3)\in\mathbb{P}^3$. Then, for a general surface $T$ in the
linear system $\mathcal{T}$, we have
\[T_{Y_1}\sim_{\mathbb{Q}}
(\alpha_P\circ\gamma\circ\sigma_1)^*(-6K_X)
-\frac{6}{7}(\gamma\circ\sigma_1)^*(E) -\frac{6}{5}\sigma_1^*(F)
-G,\] where $E$ and $G$ are the exceptional divisors of the
birational morphisms $\alpha_P$ and $\sigma_1$, respectively. In
addition, we see
\[F_{Y_1}\sim_{\mathbb{Q}} (\gamma\circ\sigma_1)^*(E) -\frac{3}{5}\sigma_1^*(F) -\frac{1}{2}G.\]

 Let $L$
be the unique curve in the linear system
$|\mathcal{O}_{\mathbb{P}(1,2,5)}(1)|$ on the surface $E$. The
base locus of the proper transform $\mathcal{T}_{Y_1}$ contains
the irreducible curve $\tilde{C}_{Y_1}$ and the irreducible curve
$L_{Y_1}$ such that
\[T_{Y_1}\cdot S_{Y_1}=\tilde{C}_{Y_1}+ L_{Y_1}, \  \
T_{Y_1}\cdot E_{Y_1}=2L_{Y_1}. \] The surface $T_{Y_1}$ is normal
and
\[\tilde{C}_{Y_1}^2=L_{Y_1}^2=-\frac{7}{4}, \ \
\tilde{C}_{Y_1}\cdot L_{Y_1}=\frac{5}{4}\] on the surface
$T_{Y_1}$. The intersection form of the curves $\tilde{C}_{Y_1}$
and $L_{Y_1}$ is negative-definite on the surface $T_{Y_1}$.
Because $\mathcal{M}_{Y_1}|_{T_{Y_1}}\equiv
nS_{Y_1}|_{T_{Y_1}}\equiv n\tilde{C}_{Y_1}+nL_{Y_1}$, we obtain an
contradictory identity $\mathcal{M}=\mathcal{T}$ from
Theorem~\ref{theorem:main-tool}.
\end{proof}

\begin{proposition}
\label{proposition:n-44} If $\gimel=44$, then the linear system
$|-2K_X|$ is a unique Halphen pencil on $X$.
\end{proposition}
\begin{proof}
By Lemma~\ref{lemma:n-44}, we may assume that the set
$\mathbb{CS}(Y, \frac{1}{n}\mathcal{M}_{Y})$  contains the point
$P_2'$. Let $\sigma_2: Y_2\to Y$ be the Kawamata blow up at the
point $P_2'$. Then, the pencil $|-2K_{Y_2}|$ is the proper
transform of the pencil $|-2K_X|$ and its base locus consists of
the irreducible curve $C_{Y_2}$.  We can easily check
$-K_{Y_2}\cdot C_{Y_2}=-2K_{Y_2}^3<0$ . Hence, we obtain the
identity  $\mathcal{M}=|-2K_{X}|$ from
Theorem~\ref{theorem:main-tool} because
$\mathcal{M}_{Y_2}\sim_{\mathbb{Q}} -nK_{Y_2}$ by
Lemma~\ref{lemma:Kawamata}.
\end{proof}

\section{Cases $\gimel=27$,  $42$, and $68$.}
\label{section:n-27-42-68}

The aim of this section is to prove the following:
\begin{proposition}
If $\gimel=27$, $42$, or $68$, then the linear system $|-a_1K_X|$
is a unique Halphen pencil on $X$.
\end{proposition}

We first consider the case $\gimel=68$.\index{$\gimel=68$} Let $X$
be the hypersurface given by a general quasihomogeneous equation
of degree $28$ in $\mathbb{P}(1,3,4,7,14)$  with
$-K_X^3=\frac{1}{42}$. The singularities consist of two quotient
singular points $P$ and $Q$ of type $\frac{1}{7}(1,3,4)$, one
point of type $\frac{1}{3}(1,1,2)$, and one point of type
$\frac{1}{2}(1,1,1)$.

We have a commutative diagram
$$
\xymatrix{
&&V\ar@{->}[dl]_{\sigma_Q}\ar@{->}[dr]^{\sigma_P}\ar@{->}[drrrrrrr]^{\eta}&&&&&&&\\
&Y_P\ar@{->}[dr]_{\pi_P}&&Y_Q\ar@{->}[dl]^{\pi_Q}&&&&&&\mathbb{P}(1,3,4),\\
&&X\ar@{-->}[rrrrrrru]_{\psi}&&&&&&&}
$$
where \begin{itemize}

\item $\psi$ is  the natural projection,

\item $\pi_P$ is the Kawamata blow up at the point $P$ with
weights $(1,3,4)$,

\item $\pi_Q$ is the Kawamata blow up at the point $Q$ with
weights $(1,3,4)$,

\item $\sigma_Q$ is the Kawamata blow up with weights $(1,3,4)$ at
the point $Q_1$ whose image by the birational morphism $\pi_P$ is
the point $Q$,

\item $\sigma_P$ is the Kawamata blow up with weights $(1,3,4)$ at
the point $P_1$ whose image by the birational morphism $\pi_Q$ is
the point $P$,

\item $\eta$ is an elliptic fibration. \end{itemize}

The set $\mathbb{CS}(X,\frac{1}{n}\mathcal{M})$ is nonempty by
Theorem~\ref{theorem:Noether-Fano}. We may also assume that
$$
\mathbb{CS}\Big(X,\frac{1}{n}\mathcal{M}\Big)\subset \Big\{P,
Q\Big\}
$$
due to Lemma~\ref{lemma:smooth-points} and
Corollary~\ref{corollary:Ryder-a1}. Suppose that the set
$\mathbb{CS}(X,\frac{1}{n}\mathcal{M})$ contains the point $P$.
The exceptional divisor $E\cong\mathbb{P}(1,3,4)$ of the
birational morphism $\pi_P$ contains two quotient singular points
$O_1$ and $O_2$ of types $\frac{1}{4}(1,3,1)$ and
$\frac{1}{3}(1,2,1)$, respectively.

\begin{lemma}
\label{lemma:n-68-first} The set
$\mathbb{CS}(Y_P,\frac{1}{n}\mathcal{M}_{Y_P})$ cannot contain the
point $O_2$.
\end{lemma}

\begin{proof}
Suppose so. We then consider the Kawamata blow up $\alpha: U\to
Y_P$ at the point $O_2$ with weights $(1,2,1)$. Let $\mathcal{D}$
be the proper transform of the linear system $|-4K_X|$ on $X$ by
the birational morphism $\pi_P\circ\alpha$. Its base locus
consists of the irreducible curve $\bar{C}_{U}$. For a general
surface $D_U$ in $\mathcal{D}$, we have
\[D_U\sim_{\mathbb{Q}} (\pi_P\circ
\alpha)^*(-4K_X)-\frac{4}{7}\alpha^*(E)-\frac{1}{3}F,\] where $F$
is the exceptional divisor of $\alpha$.  The surface $D_U$ is
normal. We also have
\[S_U\sim_{\mathbb{Q}} (\pi\circ
\sigma)^*(-K_X)-\frac{1}{7}\alpha^*(E)-\frac{1}{3}F\] and
$S_U\cdot D_U=\bar{C}_U$. On the normal surface $D_{U}$, the curve
$\bar{C}_U$ has negative self-intersection number
$\bar{C}_U^2=-\frac{5}{42}$. However, we have
\[\mathcal{M}_U\Big\vert_{D_U}\equiv -nK_U\Big\vert_{D_U}\equiv n\bar{C}_U.\]  Therefore,
Theorem~\ref{theorem:main-tool} implies $\mathcal{M}=|-4K_X|$. It
is a contradiction because the linear system $|-4K_X|$ is not a
pencil.
\end{proof}

\begin{lemma}
\label{lemma:n-68-second} If the set
$\mathbb{CS}(Y_P,\frac{1}{n}\mathcal{M}_{Y_P})$ contains the point
$O_1$, then $\mathcal{M}=|-3K_X|$.
\end{lemma}

\begin{proof}
 Let $\beta: W\to Y_P$ be the Kawamata blow up at the point $O_1$
with weights $(1,3,1)$. We consider the proper transform
$\mathcal{P}$ of the pencil $|-3K_X|$ by the birational morphism
$\pi_P\circ\beta$. Its base locus consists of the irreducible
curve $C_W$. Because the curve $C_W$ has negative
self-intersection on a general surface in $\mathcal{P}$ and a
general surface in $\mathcal{P}$ is normal, we can obtain
$\mathcal{M}=|-3K_X|$ from Theorem~\ref{theorem:main-tool}.
\end{proof}

Therefore, we may assume that the set
$\mathbb{CS}(Y_P,\frac{1}{n}\mathcal{M}_{Y_P})$ consists of the
single point $Q_1$.  The exceptional divisor of $\sigma_Q$
contains two singular points $Q_2$ and $Q_3$ of types
$\frac{1}{4}(1,3,1)$ and $\frac{1}{3}(1,2,1)$, respectively. Then,
the set $\mathbb{CS}(V,\frac{1}{n}\mathcal{M}_V)$ is nonempty.
Furthermore, it must contain either the point $Q_2$ or the point
$Q_3$. However, the proof of Lemma~\ref{lemma:n-68-first} shows it
cannot contain the point $Q_3$. Also,
Lemma~\ref{lemma:n-68-second} shows that $\mathcal{M}=|-3K_X|$ if
the set $\mathbb{CS}(V,\frac{1}{n}\mathcal{M}_V)$ contains the
point $Q_2$.

Following the same way, we  can also conclude that
$\mathcal{M}=|-3K_X|$ if the set
$\mathbb{CS}(X,\frac{1}{n}\mathcal{M})$ contains the point $Q$.

In the case $\gimel=42$,\index{$\gimel=42$} the hypersurface $X$
is given by a general quasihomogeneous equation of degree $20$ in
$\mathbb{P}(1,2,3,5,10)$
 with
$-K_X^3=\frac{1}{15}$. The singularities consist of two quotient
singular points $P$ and $Q$ of type $\frac{1}{5}(1,2,3)$, one
point of type $\frac{1}{3}(1,2,1)$, and two points of type
$\frac{1}{2}(1,1,1)$.

We have a commutative diagram
$$
\xymatrix{
&&V\ar@{->}[dl]_{\sigma_Q}\ar@{->}[dr]^{\sigma_P}\ar@{->}[drrrrrrr]^{\eta}&&&&&&&\\
&Y_P\ar@{->}[dr]_{\pi_P}&&Y_Q\ar@{->}[dl]^{\pi_Q}&&&&&&\mathbb{P}(1,2,3),\\
&&X\ar@{-->}[rrrrrrru]_{\psi}&&&&&&&}
$$
where \begin{itemize}

\item $\psi$ is  the natural projection,

\item $\pi_P$ is the Kawamata blow up at the point $P$ with
weights $(1,2,3)$,

\item $\pi_Q$ is the Kawamata blow up at $Q$ with weights
$(1,2,3)$,

\item $\sigma_Q$ is the Kawamata blow up with weights $(1,2,3)$ at
the point $Q_1$ whose image by the birational morphism $\pi_P$ is
the point $Q$,

\item $\sigma_P$ is the Kawamata blow up with weights $(1,2,3)$ at
the point $P_1$ whose image by the birational  morphism $\pi_Q$ is
the point $P$,

\item $\eta$ is an elliptic fibration.

\end{itemize}

Using the exactly same method as in the case $\gimel=68$, one can
show that $\mathcal{M}=|-2K_X|$.

From now, we consider the case $\gimel=27$.\index{$\gimel=27$} Let
$X$ be the hypersurface given by a general quasihomogeneous
equation of degree $15$ in $\mathbb{P}(1,2,3,5,5)$
 with
$-K_X^3=\frac{1}{10}$. The singularities consist of three quotient
singular points $P_1$, $P_2$, and $P_3$ of type
$\frac{1}{5}(1,2,3)$ and one point of type $\frac{1}{2}(1,1,1)$.

And we  have a commutative diagram
$$
\xymatrix{
&&Y\ar@{->}[dl]_{\pi}\ar@{->}[dr]^{\eta}&&\\
&X\ar@{-->}[rr]_{\psi}&&\mathbb{P}(1,2,3),&}
$$
where \begin{itemize}

\item $\psi$ is  the natural projection,

\item

$\pi$ is the Kawamata blow ups at the points $P_{1}$, $P_{2}$ and
$P_{3}$ with weights $(1,2,3)$,

\item $\eta$ is an elliptic fibration. \end{itemize}

Even though three singular points are involved in this case, the
same method as in the previous cases can be applied to obtain
$\mathcal{M}=|-2K_X|$.

\section{Case $\gimel=32$, hypersurface of degree $16$ in
$\mathbb{P}(1,2,3,4,7)$.}\index{$\gimel=32$} \label{section:n-32}

The hypersurface $X$ is given by a general quasihomogeneous
polynomial of degree $16$ in $\mathbb{P}(1,2,3, 4,7)$ with
$-K_X^3=\frac{2}{21}$. The singularities of the threefold $X$
consist of four quotient singular points of type
$\frac{1}{2}(1,1,1)$, one quotient  singular point of type
$\frac{1}{3}(1,2,1)$, and one quotient  singular point $P$ of type
$\frac{1}{7}(1,3,4)$. There is a commutative diagram
$$
\xymatrix{
&U\ar@{->}[d]_{\alpha}&&Y\ar@{->}[ll]_{\beta}\ar@{->}[d]^{\eta}&\\%
&X\ar@{-->}[rr]_{\psi}&&\mathbb{P}(1,2,3),&}
$$
where \begin{itemize} \item $\psi$ is the natural projection,
\item $\alpha$ is the Kawamata blow up at the point $P$ with
weights $(1,3,4)$,

\item $\beta$ is the Kawamata blow up with weights $(1,1,3)$ at
the singular point $Q$ of the variety $U$ that is a quotient
singularity of type $\frac{1}{4}(1,1,3)$ contained in the
exceptional divisor of $\alpha$,

\item  $\eta$ is an elliptic fibration.

\end{itemize}

The hypersurface $X$ can be given by the quasihomogeneous equation
$$
w^{2}y+wf_{9}\big(x,y,z,t\big)+f_{16}\big(x,y,z,t\big)=0
$$
where $f_{9}$ and $f_{16}$ are qua\-si\-ho\-mo\-ge\-ne\-ous
polynomials of degrees $9$ and $16$, respectively. Let $D$ be a
general surface in $|-2K_{X}|$. It is cut out on the threefold $X$
by the equation \[\lambda x^2+\mu y=0,\] where $(\lambda :
\mu)\in\mathbb{P}^1$.  The surface $D$ is irreducible and normal.
The base locus of the pencil $|-2K_{X}|$ consists of the curve
$C$, which implies that $C=D\cdot S$.

If the set $\mathbb{CS}(X, \frac{1}{n}\mathcal{M})$ contains the
singular point of type $\frac{1}{3}(1,2,1)$, we obtain
$\mathcal{M}=|-2K_{X}|$ from
Lemma~\ref{lemma:special-singular-points-with-zero-c}. It then
follows from Corollary~\ref{corollary:Ryder-a1} and
Lemma~\ref{lemma:smooth-points} that
$$
\mathbb{CS}\Big(X, \frac{1}{n}\mathcal{M}\Big)=\Big\{P\Big\}.%
$$
Furthermore, the set $\mathbb{CS}(U, \frac{1}{n}\mathcal{M}_U)$ is
not empty by Theorem~\ref{theorem:Noether-Fano} because $-K_{U}$
is nef and big.

The exceptional divisor $E\cong \mathbb{P}(1,3,4)$ of the
birational morphism $\alpha$ contains two singular points $O$ and
$Q$ that are quotient singularities of types $\frac{1}{3}(1,1,2)$
and $\frac{1}{4}(1,1,3)$, respectively.  Let $L$ be the unique
curve  contained in the linear system
$|\mathcal{O}_{\mathbb{P}(1,\,3,\,4)}(1)|$ on the surface $E$. Let
$F$ be the exceptional divisor of $\beta$. It contains a singular
point $Q_1$ that is quotient singularity of type
$\frac{1}{3}(1,1,2)$.

 Then, it follows from
Lemma~\ref{lemma:Cheltsov-Kawamata} that either
$Q\in\mathbb{CS}(U, \frac{1}{n}\mathcal{M}_U)$ or $\mathbb{CS}(U,
\frac{1}{n}\mathcal{M}_U)=\{O\}$.

\begin{lemma}
\label{lemma:n-32-points-P8} If the set $\mathbb{CS}(U,
\frac{1}{n}\mathcal{M}_U)$ consists of the point $Q$, then
$\mathcal{M}=|-2K_{X}|$.
\end{lemma}

\begin{proof}
It follows from Lemma~\ref{lemma:Kawamata} that
$\mathcal{M}_Y\sim_{\mathbb{Q}}-nK_{Y}$, which implies that every
surface in the pencil $\mathcal{M}_Y$ is contracted to a curve by
the morphism $\eta$ and the set $\mathbb{CS}(Y,
\frac{1}{n}\mathcal{M}_Y)$ contains the point $Q_1$.

Let $\pi:V\to Y$ be the Kawamata blow up at the point $Q_1$ with
weights $(1,1,2)$.  Then, the transform $D_V$ is normal but the
base locus of the pencil $|-2K_{V}|$ consists of the irreducible
curves $C_V$ and $L_V$.

The intersection form of the curves $C_V$ and $L_V$ on the surface
$D_V$ is negative-definite because the curves $C_V$ and $L_V$ are
components of a fiber of the elliptic fibration
$\eta\circ\pi\vert_{D_V}$ that contains three irreducible
components. On the other hand, we have
$$
\mathcal{M}_V\Big\vert_{D_V}\equiv-nK_{V}\Big\vert_{D_V}\equiv
nC_V+nL_V.
$$
Therefore, it follows from Theorem~\ref{theorem:main-tool} that
$\mathcal{M}=|-2K_{X}|$.
\end{proof}

From now on, we may assume that the set $\mathbb{CS}(U,
\frac{1}{n}\mathcal{M}_U)$ contains the point $O$ due to
Lemma~\ref{lemma:Cheltsov-Kawamata}. Let $\gamma:W\to U$ be the
Kawamata blow up at the point $O$  with weights $(1,1,2)$ and $G$
be the exceptional  divisor of the birational morphism $\gamma$.
Then, the surface $G\cong\mathbb{P}(1,1,2)$ and
$$
\mathcal{M}_W\sim_{\mathbb{Q}}-nK_{W}
\sim_{\mathbb{Q}}\gamma^{*}(-nK_{U})-\frac{n}{3}G
\sim_{\mathbb{Q}}(\alpha\circ\gamma)^{*}(-nK_{X})-\frac{n}{7}\gamma^*(E)-\frac{n}{3}G.
$$

In a neighborhood of the point $P$, the monomials $x$, $z$, and
$t$ can be considered as weighted local coordinates on $X$ such
that $\mathrm{wt}(x)=1$, $\mathrm{wt}(z)=3$, and
$\mathrm{wt}(z)=4$. Then, in a neighborhood of the singular point
$P$, the surface $D$ can be given by equation
$$
\lambda x^{2}+\mu\Big(\epsilon_{1}x^{9}+\epsilon_{2}zx^{6}+\epsilon_{3}z^{2}x^{3}+\epsilon_{4}z^{3}+\epsilon_{5}t^{2}x+\epsilon_{6}tx^{5}+\epsilon_{7}tzx^{2}+h_{16}(x,z,t)+\mbox{higher terms}\Big)=0,%
$$
where $\epsilon_{i}\in\mathbb{C}$ and $h_{16}$ is a
quasihomogeneous polynomial of degree $16$. In a neighborhood of
the singular point $O$, the birational morphism $\alpha$ can be
given by the equations
$$
x=\tilde{x}\tilde{z}^{\frac{1}{7}},\ z=\tilde{z}^{\frac{3}{7}},\ t=\tilde{t}\tilde{z}^{\frac{4}{7}},%
$$
where $\tilde{x}$, $\tilde{y}$, and $\tilde{z}$ are weighted local
coordinates on the variety $U$ in a neighborhood of the singular
point $O$ such that $\mathrm{wt}(\tilde{x})=1$,
$\mathrm{wt}(\tilde{z})=2$, and $\mathrm{wt}(\tilde{t})=1$.

In a neighborhood of the point $O$, the surface $E$ is given by
$\tilde{z}=0$, the surface $D_U$ is given by
$$
\lambda \tilde{x}^{2}+\mu\Big(\epsilon_{1}\tilde{x}^{9}\tilde{z}+\epsilon_{2}\tilde{z}\tilde{x}^{6}+\epsilon_{3}\tilde{z}\tilde{x}^{3}+\epsilon_{4}\tilde{z}+\epsilon_{5}\tilde{t}^{2}\tilde{x}\tilde{z}+\epsilon_{6}\tilde{t}\tilde{x}^{5}\tilde{z}+\epsilon_{7}\tilde{t}\tilde{z}\tilde{x}^{2}+\mbox{higher terms}\Big)=0,%
$$
and the surface $S_U$ is given by the equation $\tilde{x}=0$.

In a neighborhood of the singular point of $G$, the birational
morphism $\gamma$ can be given by
$$
\tilde{x}=\bar{x}\bar{z}^{\frac{1}{3}},\ \tilde{z}=\bar{z}^{\frac{2}{3}},\ \tilde{t}=\bar{t}\bar{z}^{\frac{1}{3}},%
$$
where $\bar{x}$, $\bar{z}$ and $\bar{t}$ are weighted local
coordinates on the variety $W$ in a neighborhood of the singular
point of $G$ such that
$\mathrm{wt}(\bar{x})=\mathrm{wt}(\bar{z})=\mathrm{wt}(\bar{t})=1$.
The surface $G$ is given by the equation $\bar{z}=0$, the proper
transform  $D_W$ is given by
$$
\lambda \bar{x}^{2}+\mu\Big(\epsilon_{1}\bar{x}^{9}\bar{z}^{3}+\epsilon_{2}\bar{z}^{2}\bar{x}^{6}+\epsilon_{3}\bar{z}\bar{x}^{3}+\epsilon_{4}+\epsilon_{5}\bar{t}^{2}\bar{x}\bar{z}+\epsilon_{6}\bar{t}\bar{x}^{5}\bar{z}^{2}+\epsilon_{7}\bar{t}\bar{z}\bar{x}^{2}+\mbox{higher terms}\Big)=0,%
$$
the proper transform $S_W$ is given by the equation $\bar{x}=0$,
and the proper transform  $E_W$  is given by the equation
$\bar{z}=0$.

Let $\mathcal{P}$ be the proper transforms on the variety $W$ of
the pencil $|-2K_{X}|$. The curves $C_W$ and $L_W$ are contained
in the base locus of the pencil $\mathcal{P}$. Moreover, easy
calculations show that the base locus of the pencil $\mathcal{P}$
does not contain  any other curve than $C_W$ and $L_W$. We also
have
\begin{equation*}
 \left\{\aligned
&E_W\sim_{\mathbb{Q}}\gamma^{*}\big(E\big)-\frac{2}{3}F,\\
&D_W\sim_{\mathbb{Q}}\big(\alpha\circ\gamma\big)^{*}\big(-2K_{X}\big)-\frac{2}{7}\gamma^{*}\big(E\big)-\frac{2}{3}G,\\
&S_W\sim_{\mathbb{Q}}\big(\alpha\circ\gamma\big)^{*}\big(-K_{X}\big)-\frac{1}{7}\gamma^{*}\big(E\big)-\frac{1}{3}G,\\
\endaligned
\right.
\end{equation*}
Also, we have  $C_W+L_W=S_W\cdot D_W$ and  $2L_W=D_W\cdot E_W$.

The curves $C_W$ and $L_W$ can be considered as irreducible
effective divisors on the normal surface $D_W$. Then, it follows
from the equivalences above that
$$L_W^2=-\frac{5}{8},\ C_W^2=-\frac{7}{24},\ C_W\cdot L_W=\frac{3}{8},%
$$
which implies that the intersection form of $C_W$ and $L_W$ on
$D_W$ is negative-definite. Let $M$ be a  general surface of the
linear system $\mathcal{M}_W$. Then,
$$
M\Big\vert_{D_W}\equiv-nK_{W}\Big\vert_{D_W}\equiv
nS_W\Big\vert_{D_W}\equiv nC_W+nL_W,
$$
which implies that $\mathcal{M}=|-2K_{X}|$ by
Theorem~\ref{theorem:main-tool}.

Consequently, we have proved
\begin{proposition}
\label{proposition:n-32} The linear system $|-2K_{X}|$ is the only
Halphen pencil on $X$.
\end{proposition}

\section{Cases   $\gimel=33$ and $38$} \label{section:n-38}

For the case $\gimel=38$,\index{$\gimel=38$} let $X$ be the
hypersurface given by a general quasihomogeneous equation of
degree $18$ in $\mathbb{P}(1,2,3,5,8)$ with $-K_X^3=\frac{1}{35}$.
Then, the singularities of $X$ consist of two singular points $P$
and $Q$ that are quotient singularities of types
$\frac{1}{5}(1,2,3)$ and $\frac{1}{8}(1,3,5)$, respectively, and
two points of type $\frac{1}{2}(1,1,1)$.

We have the following commutative diagram:
$$
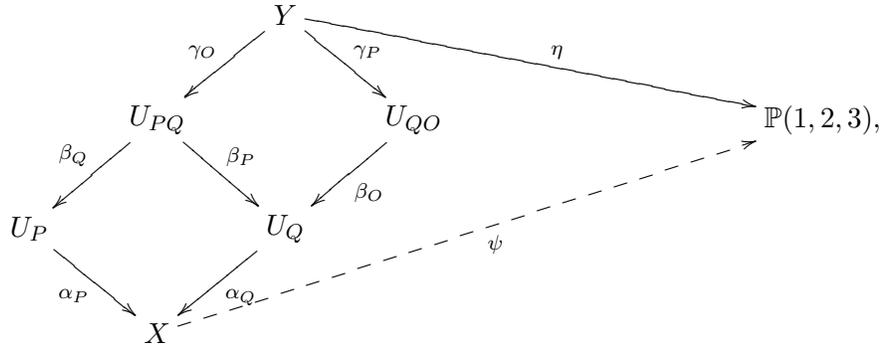
\xymatrix{
&&&Y\ar@{->}[dl]_{\gamma_{O}}\ar@{->}[dr]^{\gamma_{P}}\ar@{->}[drrrrr]^{\eta}&&&&&&\\
&&U_{PQ}\ar@{->}[dl]_{\beta_{Q}}\ar@{->}[dr]^{\beta_{P}}&&U_{QO}\ar@{->}[dl]^{\beta_{O}}&&&&\mathbb{P}(1,2,3),\\
&U_{P}\ar@{->}[dr]_{\alpha_{P}}&&U_{Q}\ar@{->}[dl]^{\alpha_{Q}}&&&&&&\\
&&X\ar@{-->}[rrrrrruu]_{\psi}&&&&&&&}
$$
where
\begin{itemize}
\item $\psi$ is  the natural projection,

\item $\alpha_{P}$ is the Kawamata blow up at the point $P$ with
weights $(1,2,3)$,

\item $\alpha_{Q}$ is the Kawamata blow up at the point  $Q$ with
weights $(1,3,5)$,

\item $\beta_{Q}$ is the Kawamata blow up with weights $(1,3,5)$
at the point whose image by the birational  morphism $\alpha_P$ is
the point $Q$,

\item $\beta_{P}$ is the Kawamata blow up with weights $(1,2,3)$
at the point whose image by the birational morphism $\alpha_Q$ is
the point $P$,

 \item $\beta_{O}$ is the Kawamata blow up with weights
$(1,3,2)$ at the singular point $O$ of the variety $U_{Q}$ that is
a quotient singularity of type $\frac{1}{5}(1,3,2)$ contained in
the exceptional divisor of the birational  morphism $\alpha_{Q}$,

\item $\gamma_{P}$ is the Kawamata blow up with weights $(1,2,3)$
at the point whose image by the birational morphism
$\alpha_Q\circ\beta_O$ is the point $P$,

\item $\gamma_{O}$ is the Kawamata blow up with weights $(1,3,2)$
at the singular point of the variety $U_{PQ}$ that is a quotient
singularity of type $\frac{1}{5}(1,3,2)$ contained in the
exceptional divisor of the birational morphism $\beta_{Q}$,

\item  $\eta$ is an elliptic fibration.
\end{itemize}
 By
Lemma~\ref{lemma:smooth-points} and
Corollary~\ref{corollary:Ryder-a1}, we may  assume that
$$\mathbb{CS}\Big(X,\frac{1}{n}\mathcal{M}\Big)\subset \Big\{P, Q\Big\}.$$

The exceptional divisor $E_P$ of the birational morphism
$\alpha_P$ contains two quotient singular points $P_1$ and $P_2$
of types $\frac{1}{3}(1,2,1)$ and $\frac{1}{2}(1,1,1)$,
respectively.

\begin{lemma}\label{lemma:n-38-P-1}
If the set $\mathbb{CS}(U_P,\frac{1}{n}\mathcal{M}_{U_P})$
contains the point $P_1$, then  $\mathcal{M}=|-2K_X|$.
\end{lemma}
\begin{proof}
Suppose it contains the point $P_1$. Let $\beta_1:W_1\to U_P$ be
the Kawamata blow up at the point $P_1$ with weights $(1,2,1)$.
The pencil $|-K_{W_1}|$ is the proper transform of the pencil
$|-2K_X|$. Its base locus consists of the irreducible curve
$C_{W_1}$.

For a general surface $D_{W_1}$ in $|-2K_{W_1}|$, we can easily
check $D_{W_1}\cdot C_{W_1}=-4K_{W_1}^3<0$. Hence, we obtain
$\mathcal{M}=|-2K_{X}|$ from Theorem~\ref{theorem:main-tool}
because $\mathcal{M}_{W_1}\sim_{\mathbb{Q}} nD_{W_1}$ by
Lemma~\ref{lemma:Kawamata}.
\end{proof}

\begin{lemma}\label{lemma:n-38-P-2}
The set $\mathbb{CS}(U_P,\frac{1}{n}\mathcal{M}_{U_P})$ cannot
contain the point $P_2$.
\end{lemma}
\begin{proof}
Suppose it contains the point $P_2$. Let $\beta_2:W_2\to U_P$ be
the Kawamata blow up at the point $P_2$ with weights $(1,1,1)$.
Also, let $\mathcal{D}_2$ be the proper transform of the linear
system $|-3K_X|$ by the birational morphism
$\alpha_P\circ\beta_2$. Its base locus consists of the irreducible
curve $\bar{C}_{W_2}$. A general surface $D_{W_2}$ in
$\mathcal{D}_2$ is normal and the self-intersection
$\bar{C}_{W_2}^2$ is negative on the surface $D_{W_2}$. Because
$\mathcal{M}_{W_2}|_{D_{W_2}}\equiv-n\bar{C}_{W_2}$, we obtain an
absurd equality $\mathcal{M}=|-3K_X|$ from
Theorem~\ref{theorem:main-tool}.
\end{proof}

Meanwhile, the exceptional divisor $E\cong\mathbb{P}(1,3,5)$ of
the birational morphism $\alpha_Q$ contains two singular points
$O$ and $Q_1$ of types $\frac{1}{5}(1,3,2)$ and
$\frac{1}{3}(1,1,2)$, respectively. For the convenience, let $L$
be the unique curve on the surface $E$ contained in the linear
system $|\mathcal{O}_{\mathbb{P}(1,3,5)}(1)|$

\begin{lemma}\label{lemma:n-38-Q-1}
If the set $\mathbb{CS}(U_Q,\frac{1}{n}\mathcal{M}_{U_Q})$
contains the point $Q_1$, then  $\mathcal{M}=|-2K_X|$.
\end{lemma}
\begin{proof}
Let $\pi_1:V_1\to U_Q$ be the Kawamata blow  up at the point $Q_1$
with weights $(1,1,2)$. The pencil $|-2K_{V_1}|$ is the proper
transform of the pencil $|-2K_X|$. Its base locus consists of the
irreducible curves $C_{V_1}$, $L_{V_1}$, and a curve $\bar{L}$ on
the exceptional divisor $F_Q\cong\mathbb{P}(1,1,2)$ of the
birational morphism $\pi_1$ contained in the linear system
$|\mathcal{O}_{\mathbb{P}(1,1,2)}(1)|$.

Let $D_{V_1}$ be a general surface in $|-2K_{V_1}|$. We see then
\[S_{V_1}\cdot D_{V_1}=C_{V_1}+L_{V_1}+\bar{L}, \ E_{V_1}\cdot D_{V_1}=2L_{V_1}, \ F_Q\cdot D_{V_1}=2\bar{L}.\]
 Using the following equivalences

\[\left\{\aligned
&E_{V_1}\sim_{\mathbb{Q}}\pi_1^{*}\big(E\big)-\frac{1}{3}F_Q,\\
&S_{V_1}\sim_{\mathbb{Q}}\big(\alpha_Q\circ\pi_1\big)^{*}\Big(-K_{X}\Big)-\frac{1}{8}\pi_1^{*}\big(E\big)-\frac{1}{3}F_Q,\\
&D_{V_1}\sim_{\mathbb{Q}}\big(\alpha_Q\circ\pi_1\big)^{*}\Big(-2K_{X}\Big)-\frac{2}{8}\pi_1^{*}\big(E\big)-\frac{2}{3}F_Q,\\
\endaligned
\right.\]  we obtain
\[C_{V_1}^2=-\frac{37}{20}, \ L_{V_1}^2=-\frac{7}{20}, \ \bar{L}^2=-\frac{3}{4}, \ C_{V_1}\cdot L_{V_1} =0,
\ C_{V_1}\cdot \bar{L}=1,\  L_{V_1}\cdot \bar{L}=\frac{1}{4}\] on
the normal surface $D_{V_1}$. One can see that the intersection
form of these three curves on $D_{V_1}$ is negative-definite. Let
$M$ be a  general surface in the linear system
$\mathcal{M}_{V_1}$. Then,
$$
M\Big\vert_{D_{V_1}}\equiv-nK_{V_1}\Big\vert_{D_{V_1}}\equiv
nS_{V_1}\Big\vert_{D_{V_1}}\equiv nC_{V_1}+nL_{V_1}+n\bar{L},
$$
which implies that $\mathcal{M}=|-2K_{X}|$ by
Theorem~\ref{theorem:main-tool}.
\end{proof}

The exceptional divisor $F_O$ of the birational morphism $\beta_O$
contains two quotient singular points $O_1$ and $O_2$ of types
$\frac{1}{3}(1,1,2)$ and $\frac{1}{2}(1,1,1)$, respectively.

\begin{lemma}\label{lemma:n-38-O-1}
If the set $\mathbb{CS}(U_{QO},\frac{1}{n}\mathcal{M}_{U_{QO}})$
contains the point $O_1$, then  $\mathcal{M}=|-2K_X|$.
\end{lemma}
\begin{proof}
Let $\sigma_1:U_1\to U_{QO}$ be the Kawamata blow up at the point
$O_1$ with weights $(1,1,2)$. The pencil $|-2K_{U_1}|$ is the
proper transform of the pencil $|-2K_X|$. Its base locus consists
of the irreducible curves $C_{U_1}$ and $L_{U_1}$.

Let $D_{U_1}$ be a general surface in  $|-2K_{U_1}|$.  We see then
\[S_{U_1}\cdot D_{U_1}=C_{U_1}+L_{U_1}, \ E_{U_1}\cdot D_{U_1}=2L_{U_1}. \]
 Using the same argument as in
Lemma~\ref{lemma:n-38-Q-1}, one can see that the intersection form
of these two curves on $D_{U_1}$ is negative-definite, and hence
$\mathcal{M}=|-2K_{X}|$.
\end{proof}

\begin{lemma}\label{lemma:n-38-O-2}
The set $\mathbb{CS}(U_{QO},\frac{1}{n}\mathcal{M}_{U_{QO}})$
cannot contain the point $O_2$.
\end{lemma}
\begin{proof}
Let $\sigma_2:U_2\to U_{QO}$ be the Kawamata blow up at the point
$O_2$ with weights $(1,1,1)$ and let $\mathcal{D}$ be the proper
transform of the linear system $|-3K_X|$ by the birational
morphism $\alpha_Q\circ\beta_O\circ\sigma_2$.

 The base locus of the linear system  $\mathcal{D}$
consists of the irreducible curve $C_{U_2}$. The same method as in
Lemma~\ref{lemma:n-38-P-1} shows that $\mathcal{M}=|-3K_X|$, which
is a contradiction.
\end{proof}

\begin{proposition}
\label{proposition:n-38} If $\gimel=38$, then the linear system
$|-2K_X|$ is the only Halphen pencil on $X$.
\end{proposition}

\begin{proof}
Due to the previous lemmas, we may assume that
$$\mathbb{CS}\Big(X,\frac{1}{n}\mathcal{M}\Big)= \Big\{P, Q\Big\}.$$
Following the Kawamata blow ups $Y\to U_{QO}\to U_Q\to X$ and
using Lemmas~\ref{lemma:n-38-Q-1}, \ref{lemma:n-38-O-1}, and
\ref{lemma:n-38-O-2}, we can furthermore assume that the set
$\mathbb{CS}(Y,\frac{1}{n}\mathcal{M}_{Y})$ contains one of
singular points contained in the exceptional divisor of the
birational  morphism $\gamma_P$. In this case,
Lemmas~\ref{lemma:n-38-P-1} and \ref{lemma:n-38-P-2} imply the
statement.
\end{proof}

From now, we consider the case $\gimel=33$.\index{$\gimel=33$} The
variety $X$ is a general hypersurface of degree $17$  in
$\mathbb{P}(1,2,3,5,7)$ with $-K_{X}^{3}=\frac{17}{210}$. The
singularities of $X$ consist of one quotient singularity of type
$\frac{1}{2}(1,1,1)$, one point that is a quotient singularity of
type $\frac{1}{3}(1,2,1)$, one point $P$ that is a quotient
singularity of type $\frac{1}{5}(1,2,3)$, and one point $Q$ that
is a quotient singularity of type $\frac{1}{7}(1,2,5)$.

We have a commutative diagram as follows:
$$
\xymatrix{
&&&Y\ar@{->}[dl]_{\gamma_{O}}\ar@{->}[dr]^{\gamma_{P}}\ar@{->}[drrrrr]^{\eta}&&&&&&\\
&&U_{PQ}\ar@{->}[dl]_{\beta_{Q}}\ar@{->}[dr]^{\beta_{P}}&&U_{QO}\ar@{->}[dl]^{\beta_{O}}&&&&\mathbb{P}(1,2,3),\\
&U_{P}\ar@{->}[dr]_{\alpha_{P}}&&U_{Q}\ar@{->}[dl]^{\alpha_{Q}}&&&&&&\\
&&X\ar@{-->}[rrrrrruu]_{\psi}&&&&&&&}
$$
where \begin{itemize}

\item $\psi$ is the natural projection,

\item $\alpha_{P}$ is the Kawamata blow up at the point $P$ with
weights $(1,2,3)$,

\item $\alpha_{Q}$ is the Kawamata blow up at the point $Q$ with
weights $(1,2,5)$,

\item $\beta_{Q}$ is the Kawamata blow up with weights $(1,2,5)$
at the point whose image to $X$ is the point $Q$,

\item $\beta_{P}$ is the Kawamata blow up with weights $(1,2,3)$
at the point whose image to $X$ is the point $P$,

\item $\beta_{O}$ is the Kawamata blow up with weights $(1,2,3)$
at the quotient singular point $O$ of type $\frac{1}{5}(1,2,3)$
contained in the ex\-cep\-ti\-onal divisor of the birational
morphism $\alpha_{Q}$,

\item $\gamma_{P}$ is the Kawamata blow up with weights $(1,2,3)$
at the point whose image to $X$ is the point $P$,

\item $\gamma_{O}$ is the Kawamata blow up with weights $(1,2,3)$
at the quotient singular point of type $\frac{1}{5}(1,2,3)$
contained in the exceptional divisor of the birational morphism
$\beta_{Q}$,

\item  $\eta$ is an elliptic fibration.

\end{itemize}
It follows from
Lemma~\ref{lemma:special-singular-points-with-positive-c} that the
set $\mathbb{CS}(X, \frac{1}{n}\mathcal{M})$ does not contain the
singular point of type $\frac{1}{2}(1,1,1)$.  Moreover,
Lemma~\ref{lemma:special-singular-points-with-zero-c} implies that
if the set $\mathbb{CS}(X, \frac{1}{n}\mathcal{M})$ contains the
singular point of type $\frac{1}{3}(1,2,1)$, then
$\mathcal{M}=|-2K_X|$.  Therefore, due to
Corollary~\ref{corollary:Ryder-a1} and
Lemma~\ref{lemma:smooth-points}, we may assume that
$$\mathbb{CS}\Big(X, \frac{1}{n}\mathcal{M}\Big)\subset\Big\{P, Q\Big\}.$$

\begin{proposition}
\label{proposition:n-33} If $\gimel=33$, then the linear system
$|-2K_X|$ is the only Halphen pencil on $X$.
\end{proposition}
The proof is almost same as the proof of the case $\gimel=38$.
Lemmas~\ref{lemma:n-38-P-1} and \ref{lemma:n-38-P-2} work for the
case $\gimel=33$.

The exceptional divisor $E\cong\mathbb{P}(1,2,5)$ of the
birational morphism $\alpha_Q$ contains two singular points $O$
and $Q_1$ of types $\frac{1}{5}(1,2,3)$ and $\frac{1}{2}(1,1,1)$,
respectively. For the case $\gimel=33$, Lemma~\ref{lemma:n-38-Q-1}
should be replaced by the following:

\begin{lemma}\label{lemma:n-33-Q-1}
The set $\mathbb{CS}(U_Q,\frac{1}{n}\mathcal{M}_{U_Q})$ cannot
contain the point $Q_1$.
\end{lemma}
\begin{proof}
Suppose that it contains the point $Q_1$. Let $\pi_1:V_1\to U_Q$
be the Kawamata blow up at the point $Q_1$. Consider the linear
system $\mathcal{T}$ on $X$ cut out by the equations
\[\lambda_0t+\lambda_1x^5+\lambda_2 x^3y=0,\]
where $(\lambda_0: \lambda_1: \lambda_2)\in \mathbb{P}^2$. For a
general surface $T$ in $\mathcal{T}$, we have
\[T_{V_1}\sim_{\mathbb{Q}}
(\alpha_Q\circ\pi_1)^*(-5K_X)-\frac{5}{7}(E)-\frac{1}{2}F_Q,\]
where $F_Q$ is the exceptional divisor of the birational morphism
$\pi_1$. The base locus of the proper transform
$\mathcal{T}_{V_1}$ consists of the irreducible curve
$\tilde{C}_{V_1}$. The surface $T_{V_1}$ is normal and
$\tilde{C}_{V_1}^2<0$ on $T_{V_1}$. Because
$\mathcal{M}_{V_1}|_{T_{V_1}}\equiv nS_{V_1}|_{T_{V_1}}\equiv
n\tilde{C}_{V_1}$, we obtain an contradictory identity
$\mathcal{M}=\mathcal{D}$ from Theorem~\ref{theorem:main-tool}.
\end{proof}
The exceptional divisor $F_O$ of the birational  morphism
$\beta_O$ contains two quotient  singular points $O_1$ and $O_2$
of types $\frac{1}{2}(1,1,1)$ and $\frac{1}{3}(1,2,1)$,
respectively. Lemma~\ref{lemma:n-38-O-1} should be also replaced
by the following:
\begin{lemma}\label{lemma:n-33-O-1}
The set $\mathbb{CS}(U_{QO},\frac{1}{n}\mathcal{M}_{U_{QO}})$
cannot contain the point $O_1$.
\end{lemma}
\begin{proof}
Suppose it contains the point $O_1$. Let $\sigma_1:U_1\to U_{QO}$
be the Kawamata blow up at the point $O_1$ with weights $(1,1,1)$.
Let $\mathcal{D}$ be the proper transform of the linear system
$|-3K_X|$ by the birational morphism
$\alpha_Q\circ\beta_O\circ\sigma_1$. Its base locus consists of
two irreducible curves. One is the curve $\bar{C}_{U_1}$ and the
other is the proper transform $L_{U_1}$ of the unique curve $L$ in
the linear system $|\mathcal{O}_{\mathbb{P}(1,2,5)}(1)|$ on the
surface $E$. For a general surface $D_{U_1}$ in $\mathcal{D}$,
\[\left\{\aligned
&E_{U_1}\sim_{\mathbb{Q}}(\beta_O\circ\sigma_1)^*(E)
-\frac{3}{5}\sigma_1^*(F_O)-\frac{1}{2}G,\\
&S_{U_1}\sim_{\mathbb{Q}}
(\alpha_Q\circ\beta_O\circ\sigma_1)^*(-K_X)-\frac{1}{7}(\beta_O\circ\sigma_1)^*(E)
-\frac{1}{5}\sigma_1^*(F_O)-\frac{1}{2}G,\\
&D_{U_1}\sim_{\mathbb{Q}}
(\alpha_Q\circ\beta_O\circ\sigma_1)^*(-3K_X)-\frac{3}{7}(\beta_O\circ\sigma_1)^*(E)
-\frac{3}{5}\sigma_1^*(F_O)-\frac{1}{2}G,\\
\endaligned
\right.\] where $G$ is the exceptional divisor of $\sigma_1$.
 We see also
\[S_{U_1}\cdot D_{U_1}=\bar{C}_{U_1}+L_{U_1}, \ E_{U_1}\cdot D_{U_1}=2L_{U_1}. \]
 Using the same argument as in
Lemma~\ref{lemma:n-38-Q-1}, one can see that the intersection form
of these two curves on $D_{U_1}$ is negative-definite, and hence
$\mathcal{M}=|-3K_{X}|$. It is a contradiction.
\end{proof}
Finally, we complete the proof of
Proposition~\ref{proposition:n-33} by replacing
Lemma~\ref{lemma:n-38-O-2} by the following:
\begin{lemma}\label{lemma:n-33-O-2}
If the set $\mathbb{CS}(U_{QO},\frac{1}{n}\mathcal{M}_{U_{QO}})$
contains the point $O_2$, then $\mathcal{M}=|-2K_X|$.
\end{lemma}
\begin{proof}
Let $\sigma_2:U_2\to U_{QO}$ be the Kawamata blow up at the point
$O_2$ with weights $(1,2,1)$. The pencil $|-2K_{U_2}|$ is the
proper transform of the pencil $|-2K_X|$. Its base locus consists
of the irreducible curve $C_{U_2}$. Because
$\mathcal{M}_{U_2}\sim_{\mathbb{Q}} -nK_{U_2}$ and $-K_{U_2}\cdot
C_{U_2}<0$, we obtain $\mathcal{M}=|-2K_X|$ from
Theorem~\ref{theorem:main-tool}.
\end{proof}

\section{Cases $\gimel= 37$, $39$, $52$,    $59$,   $73$ , and  $78$.}%
\label{section:n-37-39-52-59-73-78}

Suppose that $\gimel\in\{37, 39, 52,  59,  73
\}$.\index{$\gimel=37$}\index{$\gimel=39$}\index{$\gimel=52$}\index{$\gimel=59$}\index{$\gimel=73$}
Then, the threefold $X \subset \mathbb{P}(1,a_1,a_2,a_3,a_4)$
always contains the point $O=(0:0:0:1:0)$. It  is a singular point
of $X$ that is a quotient singularity of type
$\frac{1}{a_3}(1,a_2, a_3-a_2)$.

We also have a commutative diagram as follows:
$$
\xymatrix{
&&&Y\ar@{->}[lld]_{\pi}\ar@{->}[rrd]^{\eta}&&&\\%
&X\ar@{-->}[rrrr]_{\psi}&&&&\mathbb{P}(1,a_1,a_{2})&}
$$
where \begin{itemize} \item $\psi$ is the natural projection,
\item $\pi$ is the Kawamata blow up at the point $O$ with weights
$(1, a_2, a_3-a_2)$, \item $\eta$ is an elliptic
fibration.\end{itemize}

\begin{proposition}
\label{proposition:n-37-39-52-59-73} If $\gimel\in\{37, 39, 52,
59, 73\}$, then $\mathcal{M}=|-a_{1}K_{X}|$.
\end{proposition}

\begin{proof}
We may assume that the set $\mathbb{CS}(X,
\frac{1}{n}\mathcal{M})$ consists of the point $O$ by
Lemmas~\ref{lemma:smooth-points},
\ref{lemma:special-singular-points-with-positive-c},
\ref{lemma:special-singular-points-with-zero-c} and
Corollary~\ref{corollary:Ryder-a1}.

 Let $P$ be the singular point contained in the exceptional divisor of the birational morphism  $\pi$ that is
 a quotient
singular point of type $\frac{1}{a_2}(1,a_2-1,1)$. Then, the set
$\mathbb{CS}(Y, \frac{1}{n}\mathcal{M}_Y)$ contains the point $P$
by Theorem~\ref{theorem:Noether-Fano} and
Lemma~\ref{lemma:Cheltsov-Kawamata}.

Let $\alpha:W\to Y$ be the Kawamata blow up at the point $P$ with
weights $(1,a_2-1,1)$, and $\mathcal{P}$ be the proper transforms
of the pencil $|-a_{1}K_{X}|$ by the birational morphism
$\pi\circ\alpha$. Then,
$$
\mathcal{P}\sim_{\mathbb{Q}}-a_{1}K_{W},\ \mathcal{M}_W\sim_{\mathbb{Q}}-nK_{W}.%
$$
One can easily check that the base locus of the pencil
$\mathcal{P}$ consists of the irreducible curve $C_W$. A general
surface  $D$ in the pencil $\mathcal{P}$ is normal and the
inequality $C_W^{2}<0$ holds on the surface $D$, which implies
that $ \mathcal{M}=|-a_{1}K_{X}| $ by
Theorem~\ref{theorem:main-tool}.
\end{proof}

\begin{proposition}\index{$\gimel=78$}
\label{proposition:n-78} If $\gimel=78$, then
$\mathcal{M}=|-a_{1}K_{X}|$.
\end{proposition}
\begin{proof}
The only different thing from the proof of
Proposition~\ref{proposition:n-37-39-52-59-73} is that the
exceptional divisor $E$ contains another singular point $Q$. It is
a quotient singularity of type $\frac{1}{2}(1,1,1)$. For the case
$\gimel=78$, we have to consider  the case when $\mathbb{CS}(Y,
\frac{1}{n}\mathcal{M}_Y)=\{Q\}$. In this case, applying the
method in Proposition~\ref{proposition:n-37-39-52-59-73} to the
Kawamata blow up at the point $Q$ and the linear system
$|-a_2K_X|$, we can easily obtain $\mathcal{M}=|-a_2K_X|$.
However, it is a contradiction because the linear system
$|-a_2K_X|$ is not a pencil. Therefore, the case when
$\mathbb{CS}(Y, \frac{1}{n}\mathcal{M}_Y)=\{Q\}$ never happens.
\end{proof}

\section{Cases  $\gimel=40$ and  $61$.} \label{section:n-40}

For the case $\gimel=40$,\index{$\gimel=40$} let $X$ be the
hypersurface given by a general quasihomogeneous equation of
degree $19$ in $\mathbb{P}(1,3,4,5,7)$ with
$-K_X^3=\frac{19}{420}$. Then, the singularities of $X$ consist of
one quotient singular point $P$ of type $\frac{1}{7}(1,3,4)$, one
quotient  singular point $Q$ of type $\frac{1}{5}(1,3,2)$, one
point of type $\frac{1}{4}(1,3,1)$, and one point of type
$\frac{1}{3}(1,1,2)$.

We have a commutative diagram
$$
\xymatrix{
&&V\ar@{->}[dl]_{\sigma_Q}\ar@{->}[dr]^{\sigma_P}\ar@{->}[drrrrrrr]^{\eta}&&&&&&&\\
&Y_P\ar@{->}[dr]_{\pi_P}&&Y_Q\ar@{->}[dl]^{\pi_Q}&&&&&&\mathbb{P}(1,3,4),\\
&&X\ar@{-->}[rrrrrrru]_{\psi}&&&&&&&}
$$
where \begin{itemize} \item$\psi$ is  the natural projection,

\item $\pi_P$ is the Kawamata blow up at the point $P$ with
weights $(1,3,4)$,

\item $\pi_Q$ is the Kawamata blow up at the point $Q$ with
weights $(1,3,2)$,

\item $\sigma_Q$ is the Kawamata blow up with weights $(1,3,2)$ at
the point $Q_1$ whose image by the birational morphism $\pi_P$ is
the point $Q$,

\item $\sigma_P$ is the Kawamata blow up with weights $(1,3,4)$ at
the point $P_1$ whose image by the birational morphism $\pi_Q$ is
the point $P$,

\item $\eta$ is an elliptic fibration.
\end{itemize}

 We may assume that
$$\mathbb{CS}\Big(X,\frac{1}{n}\mathcal{M}\Big)\subset \Big\{P, Q\Big\}$$
due to Lemmas~\ref{lemma:smooth-points},
\ref{lemma:special-singular-points-with-positive-c},
\ref{lemma:special-singular-points-with-zero-c}, and
Corollary~\ref{corollary:Ryder-a1}.

 The exceptional divisor $E_Q$ of the birational  morphism $\pi_Q$
contains two singular point $Q_1$ and $Q_2$ that are quotient
singularities of types $\frac{1}{3}(1,1,2)$ and
$\frac{1}{2}(1,1,1)$, respectively.

\begin{lemma}
The set $\mathbb{CS}(Y_Q,\frac{1}{n}\mathcal{M}_{Y_Q})$ contain
neither the point $Q_1$ nor the point $Q_2$.
\end{lemma}
\begin{proof}
Suppose that the set
$\mathbb{CS}(Y_Q,\frac{1}{n}\mathcal{M}_{Y_Q})$ contains the point
$Q_1$. We then consider the Kawamata blow up $\alpha_1: U_1\to
Y_Q$ at the point $Q_1$ with weights $(1,1,2)$. Let
$\mathcal{D}_1$ be the proper transform of the linear system
$|-7K_X|$  by the birational  morphism $\alpha_1\circ\pi_Q$. We
have
$$
\mathcal{M}_{U_1}\sim_{\mathbb{Q}} -nK_{U_1}\sim_{\mathbb{Q}}
(\pi_Q\circ\alpha_1)^{*}(-nK_{X})-\frac{n}{5}\alpha^*_1(E_Q)-\frac{n}{3}F_1,
$$
where $F_1$ is the exceptional divisor of $\alpha_1$. Let $D_1$ be
a general surface in $\mathcal{D}_1$. Because the base locus of
the linear system $\mathcal{D}_1$ does not contain  any curve, the
divisor $D_1$ is nef. We then see
$$
D_1\sim_{\mathbb{Q}}
(\pi_Q\circ\alpha_1)^{*}(-7K_{X})-\frac{2}{5}\alpha_1^*(E_Q)-\frac{2}{3}F_1,
$$
which shows
\begin{equation*}
\begin{split}
D_1\cdot M_{1}\cdot M_{2}&=
\Big((\pi_Q\circ\alpha_1)^{*}(-7K_{X})-\frac{2}{5}\alpha_1^*(E_Q)-\frac{2}{3}F_1\Big)
\Big((\pi_Q\circ\alpha_1)^{*}(-nK_{X})-\frac{n}{5}\alpha^*_1(E_Q)-\frac{n}{3}F_1\Big)^{2}\\
&=-\frac{1}{12}n^{2},\end{split}%
\end{equation*}  where $M_{1}$ and $M_{2}$ are general surfaces in
$\mathcal{M}_{U_1}$. It is a contradiction.

To exclude the point $Q_2$, we use the exactly same method.
\end{proof}

Meanwhile, the exceptional divisor $E_P$ of the birational
morphism $\pi_P$ contains two singular point $P_1$ and $P_2$ that
are quotient singularities of types $\frac{1}{4}(1,3,1)$ and
$\frac{1}{3}(1,2,1)$, respectively.

\begin{lemma}
If the set $\mathbb{CS}(Y_P,\frac{1}{n}\mathcal{M}_{Y_P})$
contains the point $P_1$, then $\mathcal{M}=|-3K_X|$.
\end{lemma}
\begin{proof}
Suppose that the set
$\mathbb{CS}(Y_P,\frac{1}{n}\mathcal{M}_{Y_P})$ contains the point
$P_1$. Let $\beta_1:W_1\to Y_P$ be the Kawamata blow up at the
point $P_1$ with weights $(1,3,1)$. Also, let $\mathcal{L}_1$  be
the proper transform of the linear system $|-3K_X|$   by the
birational  morphism $\beta_1\circ\pi_P$. Then, the base locus of
the linear system $\mathcal{L}_1$ consists of the irreducible
curve $C_{W_1}$. Let $H_1$ be a general surface in
$\mathcal{L}_1$. Then, we have
$$
H_1\sim_{\mathbb{Q}}
(\pi_P\circ\beta_2)^{*}(-3K_{X})-\frac{3}{7}\beta_2^{*}(E_P)-\frac{3}{4}G_1,
$$
where $G_1$ is the exceptional divisor of $\beta_1$. Then, the
inequality $-K_{W_1}\cdot C_{W_1}=-3K_{W_1}^3<0$ and the
equivalence  $\mathcal{M}_{W_1}\sim_{\mathbb{Q}} -nK_{W_1}$ imply
$\mathcal{M}=|-3K_X|$ by Theorem~\ref{theorem:main-tool}.
\end{proof}

\begin{lemma}
The set $\mathbb{CS}(Y_P,\frac{1}{n}\mathcal{M}_{Y_P})$ cannot
contain the point $P_2$.
\end{lemma}
\begin{proof}
Suppose that the set
$\mathbb{CS}(Y_P,\frac{1}{n}\mathcal{M}_{Y_P})$ contains the point
$P_2$. Let $\beta_2:W_2\to Y_P$ be the Kawamata blow up at the
point $P_2$ with weights $(1,2,1)$ and  $\mathcal{L}_2$  be the
proper transform of the linear system $|-4K_X|$ by the birational
morphism $\beta_2\circ\pi_P$. Then, the base locus of the linear
system $\mathcal{L}_2$ consists of the irreducible curve
$\bar{C}_{W_2}$.

Let $H_2$ be a general surface in $\mathcal{L}_2$. Then, we have
$$
H_2\sim_{\mathbb{Q}}
(\pi_P\circ\beta_2)^{*}(-4K_{X})-\frac{4}{7}\beta_2^{*}(E_P)-\frac{1}{3}G_2,
$$
where $G_2$ is the exceptional divisor of $\beta_2$. The
equivalence  $\mathcal{M}_{W_2}\sim_{\mathbb{Q}} -nK_{W_2}$
implies $M\vert_{H_2}\equiv n \bar{C}_{W_2}$, where $M$ is a
general surface in $\mathcal{M}_{W_2}$.  Therefore, we obtain the
identity $\mathcal{M}=|-4K_X|$ from
Theorem~\ref{theorem:main-tool} because
$\bar{C}_{W_2}^{2}=-\frac{1}{30}$ on the normal surface $H_2$.
However, it is a contradiction because $|-4K_X|$ is not a pencil.
\end{proof}

Consequently, we have proved
\begin{proposition}
The linear system $|-3K_X|$ is the only Halphen pencil on $X$.
\end{proposition}

In the case $\gimel=61$,\index{$\gimel=61$} the hypersurface $X$
is given by a general quasihomogeneous equation of degree $25$ in
$\mathbb{P}(1,4,5,7,9)$ with $-K_X^3=\frac{5}{252}$.  It has three
singular points. One is a quotient singular point $P$ of type
$\frac{1}{9}(1,4,5)$, another is a quotient singular point $Q$ of
type $\frac{1}{7}(1,5,2)$, and the other is a quotient singular
point of type $\frac{1}{4}(1,3,1)$.

We have a commutative diagram
$$
\xymatrix{
&&V\ar@{->}[dl]_{\sigma_Q}\ar@{->}[dr]^{\sigma_P}\ar@{->}[drrrrrrr]^{\eta}&&&&&&&\\
&Y_P\ar@{->}[dr]_{\pi_P}&&Y_Q\ar@{->}[dl]^{\pi_Q}&&&&&&\mathbb{P}(1,4,5),\\
&&X\ar@{-->}[rrrrrrru]_{\psi}&&&&&&&}
$$
where \begin{itemize} \item $\psi$ is  the natural projection,

\item $\pi_P$ is the Kawamata blow up at the point $P$ with
weights $(1,4,5)$,

\item $\pi_Q$ is the Kawamata blow up at the point  $Q$ with
weights $(1,5,2)$,

\item $\sigma_Q$ is the Kawamata blow up with weights $(1,5,2)$ at
the point $Q_1$ whose image by the birational morphism $\pi_P$ is
the point $Q$,

\item $\sigma_P$ is the Kawamata blow up with weights $(1,4,5)$ at
the point $P_1$ whose image by the birational morphism $\pi_Q$ is
the point $P$,

\item $\eta$ is an elliptic fibration.\end{itemize}

\begin{proposition}
The linear system $|-4K_X|$ is the only Halphen pencil on $X$.
\end{proposition}
\begin{proof}
The proof is exactly same as that of the case $\gimel=40$.
\end{proof}

\section{Case $\gimel=43$, hypersurface of degree $20$ in
$\mathbb{P}(1,2,4,5,9)$.}\index{$\gimel=43$} \label{section:n-43}

The threefold $X$ is a general hypersurface of degree $20$  in
$\mathbb{P}(1,2,4,5,9)$ with $-K_{X}^{3}=\frac{1}{18}$. The
singularities of the hypersurface $X$ consist of five points that
are quotient singularities of type $\frac{1}{2}(1,1,1)$ and the
point $O=(0:0:0:0:1)$ that is a quotient singularity of type
$\frac{1}{9}(1,4,5)$.

There is a commutative diagram
$$
\xymatrix{
&U\ar@{->}[d]_{\alpha}&&Y\ar@{->}[ll]_{\beta}\ar@{->}[d]^{\eta}&\\%
&X\ar@{-->}[rr]_{\psi}&&\mathbb{P}(1,2,4),&}
$$
where \begin{itemize}

\item $\psi$ is the natural projection,

\item $\alpha$ is the Kawamata blow up at the point  $O$ with
weights $(1,4,5)$,

\item $\beta$ is the Kawamata blow up with weights $(1,4,1)$ at
the singular point of the variety $U$ that is a quotient
singularity of type $\frac{1}{5}(1,4,1)$ contained in the
exceptional divisor of $\alpha$,

\item $\eta$ is an elliptic fibration.

\end{itemize}

Due to Theorem~\ref{theorem:Noether-Fano},
Lemmas~\ref{lemma:smooth-points},
\ref{lemma:special-singular-points-with-positive-c}, and
Corollary~\ref{corollary:Ryder-a1}, we may assume $\mathbb{CS}(X,
\frac{1}{n}\mathcal{M})=\{O\}$.

 Let $E\cong\mathbb{P}(1,4,5)$ be the exceptional  divisor of the
 birational
morphism $\alpha$ and $D$ be a general surface of the pencil
$|-2K_{X}|$. The linear system $|-2K_{U}|$ is the proper transform
of the pencil $|-2K_{X}|$. Its base locus consists of the curve
$C_{U}$ and  the unique curve $L$ in the linear system
$|\mathcal{O}_{\mathbb{P}(1,4,5)}(1)|$ on the surface $E$.

The exceptional divisor $E$ contains two singular points $P$ and
$Q$ of $U$ that are quotient singularities of types
$\frac{1}{4}(1,3,1)$ and $\frac{1}{5}(1,4,1)$, respectively.

\begin{lemma}
\label{lemma:face-n-43-P8} The set $\mathbb{CS}(U,
\frac{1}{n}\mathcal{M}_{U})$ cannot contain a curve.
\end{lemma}

\begin{proof}
Suppose that it contains a curve $Z$. Then, $-K_U\cdot
Z=\frac{1}{9}$ by Lemma~\ref{lemma:Cheltsov-Kawamata}. Because
$-K_U^3=\frac{1}{20}$, it contradicts Lemma~\ref{lemma:curves}.
\end{proof}

Therefore, the set $\mathbb{CS}(U, \frac{1}{n}\mathcal{M}_{U})$
must contain either the point $P$ or the point $Q$.

\begin{lemma}
\label{lemma:n-43-P8} If the set $\mathbb{CS}(U,
\frac{1}{n}\mathcal{M}_{U})$ consists of  the point $P$, then
$\mathcal{M}=|-2K_X|$
\end{lemma}

\begin{proof}

Let $\gamma:W\to U$ be the Kawamata blow up at the point $P$ with
weights $(1,3,1)$. The exceptional divisor of $\gamma$ contains
the singular point $P_1$ of the variety $W$ that is a quotient
singularity of type $\frac{1}{3}(1,2,1)$. The pencil $|-2K_{W}|$
is the proper transform of the pencil $|-2K_{X}|$.

We first suppose that the set $\mathbb{CS}(W,
\frac{1}{n}\mathcal{M}_{W})$ contains the point $P_1$.  Let
$\pi:V\to W$ be the Kawamata blow up at the singular point $P_{1}$
with weights $(1,2,1)$ and $\mathcal{H}$ be the proper transform,
on the threefold $V$, of the linear system $|-5K_{X}|$. The pencil
$|-2K_{V}|$ is the proper transform of the pencil $|-2K_{X}|$. The
base locus of the linear system $\mathcal{H}$ consists of the
irreducible curve $\tilde{C}_V$ whose image to $X$ is the base
locus of the linear system $|-5K_{X}|$. For a general surface $H$
in $\mathcal{H}$, we have $H\cdot\tilde{C}_V=1$ and $H^{3}=6$,
which implies that the divisor $H$ is nef and big. On the other
hand, we have $H\cdot M_{V}\cdot D_{V}=0$, where $M_V$ is general
surface in $\mathcal{M}_V$. Therefore, it follows from
Theorem~\ref{theorem:main-tool} that $\mathcal{M}=|-2K_{X}|$.

From now, we suppose that  the set $\mathbb{CS}(W,
\frac{1}{n}\mathcal{M}_{W})$ does not contain the point $P_{1}$,
which implies that  the log pair $(W, \frac{1}{n}\mathcal{M}_{W})$
is terminal by Lemma~\ref{lemma:Cheltsov-Kawamata}.

In a neighborhood of the point $O$, the monomials $x$, $z$, and
$t$ can be considered as weighted local coordinates on $X$ such
that $\mathrm{wt}(x)=1$, $\mathrm{wt}(z)=4$, and
$\mathrm{wt}(t)=5$. Then, in a neighborhood of the singular point
$P$, the Kawamata blow up $\alpha$ is given by the equations
$$
x=\tilde{x}\tilde{z}^{\frac{1}{9}},\ z=\tilde{z}^{\frac{4}{9}},\ t=\tilde{t}\tilde{z}^{\frac{5}{9}},%
$$
where $\tilde{x}$, $\tilde{z}$, and $\tilde{t}$ are weighted local
coordinated on the variety $U$ in a neighborhood of the singular
point $P$ such that $\mathrm{wt}(\tilde{x})=1$,
$\mathrm{wt}(\tilde{z})=3$, and $\mathrm{wt}(\tilde{t})=1$. The
surface $E$ is given by the equation $\tilde{z}=0$, and the
surface $S_{U}$ is given by the equation $\tilde{x}=0$. Moreover,
it follows from the local equation of the surface $D_{U}$ that
$D_{U}\cdot S_{U}=C_{U}+2L$, where the curve $L$ is locally given
by the equations $\tilde{z}=\tilde{x}=0$.

In a neighborhood of the point $P_{1}$, the birational morphism
$\gamma$ is given by the equations
$$
\tilde{x}=\bar{x}\bar{z}^{\frac{1}{4}},\ \tilde{z}=\bar{z}^{\frac{3}{4}},\ \tilde{t}=\bar{t}\bar{z}^{\frac{1}{4}},%
$$
where $\bar{x}$, $\bar{z}$, and $\bar{t}$ are weighted local
coordinates on the variety $W$ in a neighborhood of the point
$P_{1}$ such that $\mathrm{wt}(\bar{x})=1$,
$\mathrm{wt}(\bar{z})=2$, and $\mathrm{wt}(\bar{t})=1$. In
particular, the exceptional divisor of the birational  morphism
$\gamma$ is given by the equation $\bar{z}=0$ and the surface
$S_{W}$ is locally given by the equation $\bar{x}=0$.

Let $F$ be the exceptional divisor of the birational  morphism
$\gamma$ and $\bar{L}$ be the curve on the variety $W$ that is
locally given by the equations $\bar{z}=\bar{x}=0$. Then,
$$
D_{W}\cdot S_{W}=C_{W}+2L_{W}+\bar{L},\ D_{W}\cdot E_{W}=2L_{W},\ D_{W}\cdot F=2\bar{L},%
$$
while the base locus of $|-2K_{W}|$ consists of the curves
$C_{W}$, $L_{W}$, and $\bar{L}$. We have
$$
D_{W}\cdot C_{W}=0,\ D_{W}\cdot L_{W}=-\frac{2}{5},\ D_{W}\cdot\bar{L}=\frac{2}{3},%
$$
because
$$
\left\{\aligned
&D_{W}\sim_{\mathbb{Q}}(\alpha\circ\gamma)^*(-2K_X)-\frac{2}{9}\gamma^*(E)
-\frac{2}{4}F,\\
&S_{W}\sim_{\mathbb{Q}}(\alpha\circ\gamma)^*(-K_X)-\frac{1}{9}\gamma^*(E)
-\frac{1}{4}F,\\
&E_{W}\sim_{\mathbb{Q}}\gamma^*(E)
-\frac{3}{4}F.\\
\endaligned
\right.
$$
We also have
$$
\left\{\aligned &S_W^y\sim_{\mathbb{Q}}
(\alpha\circ\gamma)^*(-2K_X)-\frac{11}{9}\gamma^*(E) -\frac{6}{4}F
\sim_{\mathbb{Q}} (\alpha\circ\gamma)^*(-2K_X)-\frac{11}{9}E_W
-\frac{5}{3}F,\\
&S_W^t
\sim_{\mathbb{Q}}(\alpha\circ\gamma)^*(-5K_X)-\frac{5}{9}\gamma^*(E)
-\frac{1}{4}F
\sim_{\mathbb{Q}}(\alpha\circ\gamma)^*(-5K_X)-\frac{5}{9}E_W
-\frac{2}{3}F,\\
\endaligned
\right.
$$
which implies that $$-14K_{W}\sim 14D_W \sim
2S_W^y+2S_W^t+2E_{W}.$$ The support of the cycle $S_W^y\cdot
S_W^t$ does not contain the curves $\bar{L}$ and $C_{W}$.
Therefore, the base locus of the linear system $|-14K_{W}|$ does
not contain curves except the curve $L_{W}$.

The log pair $(W, \lambda |-14K_{W}|)$ is log-terminal for some
rational number $\lambda>\frac{1}{14}$ but the divisor
$K_{W}+\lambda |-14K_{W}|$ has non-negative intersection with all
curves on the variety $W$ except the curve $L_{W}$. Hence, it
follows from \cite{Sho93} that there is a log-flip
$\zeta:W\dasharrow W'$ along the curve $L_{W}$ with respect to the
log pair $(W, \lambda |-14K_{W}|)$. In particular, the divisor
$-K_{W'}$ is nef. Thus, the singularities of the log pair $(W',
\frac{1}{n}\mathcal{M}_{W'})$ are terminal because the
singularities of the log pair $(W, \frac{1}{n}\mathcal{M}_{W})$
are terminal but the rational map $\zeta$ is a log flop with
respect to the log pair $(W, \frac{1}{n}\mathcal{M}_{W})$. Hence,
the divisor $-K_{W'}$ is not big by
Theorem~\ref{theorem:Noether-Fano}, and hence the divisor $-K_{W}$
is not big either.

The rational functions  $(\alpha\circ\gamma)^*(\frac{y}{x^{2}})$
and $(\alpha\circ\gamma)^*(\frac{ty}{x^{7}})$ are contained in the
linear systems $|2S_{W}|$ and $|7S_{W}|$, respectively.  The
equivalences
$$
S_W^z\sim_{\mathbb{Q}}
(\alpha\circ\gamma)^*(-4K_X)-\frac{4}{9}\gamma^*(E)
\sim_{\mathbb{Q}} (\alpha\circ\gamma)^*(-4K_X)-\frac{4}{9}E_W
-\frac{1}{3}F $$ imply that $-6K_{W}\sim S_W^y+S_W^z+E_{W}$. Thus,
the rational function $(\alpha\circ\gamma)^*(\frac{zy}{x^{6}})$ is
contained in the linear system $|6S_{W}|$, which implies that the
linear system $|-42K_{W}|$ maps the variety $W$ dominantly on a
threefold\footnote{In fact, the linear system $|-210K_{W}|$
induces a birational map $W\dasharrow X^{\prime}$, where
$X^{\prime}$ is a hypersurface of degree $30$  in
$\mathbb{P}(1,2,6,7,15)$ with canonical singularities.}. Hence,
the divisor $-K_{W}$ is big, which is a contradiction.
\end{proof}

Therefore, we may assume that the set $\mathbb{CS}(U,
\frac{1}{n}\mathcal{M}_{U})$ contains the point $Q$. The
exceptional divisor $G$ of $\beta$ contains the singular point
$Q_1$ of $Y$ that is a quotient singularity of type
$\frac{1}{3}(1,2,1)$. However, we have the following statement.
\begin{lemma}
\label{lemma:n-43-P10} If the set $\mathbb{CS}(Y,
\frac{1}{n}\mathcal{M}_{Y})$ contains the point $Q_1$, then
$\mathcal{M}=|-2K_X|$.
\end{lemma}

\begin{proof}
Suppose that the set $\mathbb{CS}(Y, \frac{1}{n}\mathcal{M}_{Y})$
contains the point $Q_1$. Let $\sigma_1:Y_1\to Y$ be the Kawamata
blow up at the singular point $Q_1$. Then, the pencil
$|-2K_{Y_1}|$ is the proper transform of the pencil $|-2K_{X}|$
and its base locus consists of the curves $C_{Y_1}$ and $L_{Y_1}$.
Thus, we see that $$-K_{Y_1}\cdot L_{Y_1}=-\frac{1}{4}, \ \
-K_{Y_1}\cdot C_{Y_1}=0$$ because $D_{Y_1}\cdot
S_{Y_1}=C_{Y_1}+2L_{Y_1}$. Hence, the divisor
$B:=D_{Y_1}+(\beta\circ\sigma_1)^*(-10K_{U})$ is nef and big. On
the other hand,  $B\cdot M_{V_1}\cdot D_{V_1}=0$, where $M_{Y_1}$
is a general surface in $\mathcal{M}_{Y_1}$. Therefore, we obtain
$\mathcal{M}_{Y_1}=|-2K_{Y_1}|$ from
Theorem~\ref{theorem:main-tool}.
\end{proof}

Hence, we may assume that the set $\mathbb{CS}(Y,
\frac{1}{n}\mathcal{M}_{Y})$ does not contain subvarieties of $G$.
However, it is not empty by Theorem~\ref{theorem:Noether-Fano}.
Therefore, it must consist of the point $\bar{P}$ whose image to
$U$ is the point $P$.

 Let $\sigma_2:Y_2\to Y$ be the Kawamata blow
up of the singular point $\bar{P}$. The proof of
Lemma~\ref{lemma:n-43-P8} shows that the pencil $\mathcal{M}$
coincides with $|-2K_X|$ if the set $\mathbb{CS}(Y_2,
\frac{1}{n}\mathcal{M}_{Y_2})$ contains the singular point of the
variety $Y_2$ that is a quotient singularity of type
$\frac{1}{3}(1,2,1)$ contained in the exceptional divisor of the
birational  morphism $\sigma_2$. Therefore, we may assume that the
log pair $(Y_2, \frac{1}{n}\mathcal{M}_{Y_2})$ is terminal by
Lemma~\ref{lemma:Cheltsov-Kawamata}.

The pencil $|-2K_{Y_2}|$ is the proper transform of the pencil
$|-2K_{X}|$ and its  base locus  consists of the curves $C_{Y_2}$,
$L_{Y_2}$, and a curve $\tilde{L}$ contained in the exceptional
divisor of $\sigma_2$. Moreover, the curve $L_{Y_2}$ is the only
curve that has negative intersection with the divisor $-K_{Y_2}$
because $D_{Y_2}\sim_{\mathbb{Q}} -2K_{Y_2}$ and $-K_{Y_2}\cdot
C_{Y_2}=0$. Hence, it follows from \cite{Sho93} that there is a
log-flip $\chi:Y_2\dasharrow Y_2'$ along the curve $L_{Y_2}$ with
respect to the log pair $(Y_2, \lambda |-2K_{Y_2}|)$ for some
rational number $\lambda>\frac{1}{2}$. In particular, the divisor
$-K_{Y_2}$ is nef.

The rational map $\chi$ is a log flop with respect to the log pair
$(Y_2, \frac{1}{n}\mathcal{M}_{Y_2})$. Thus, the singularities of
the log pair $(Y_2, \frac{1}{n}\mathcal{M}_{Y_2})$ are terminal.
Hence, the divisor $-K_{Y_2}$ is not big by
Theorem~\ref{theorem:Noether-Fano}. On the other hand, the
abundance theorem (\cite{KMM}) implies that the linear system
$|-rK_{Y_2}|$ is base-point-free for $r\gg 0$. Moreover, the
pull-backs of the rational functions  $\frac{y}{x^{2}}$ and
$\frac{zy}{x^{5}}$ are contained in the linear systems
$|2S_{Y_2}|$ and $|6S_{Y_2}|$, respectively. Thus, the linear
system $|-rK_{Y_2}|$ induces an elliptic fibration, which is
impossible by Theorem~\ref{theorem:Noether-Fano}.

Consequently, we have proved
\begin{proposition}
\label{proposition:n-43} The linear system $|-2K_{X}|$ is the only
Halphen pencil on $X$.
\end{proposition}

\section{Cases $\gimel=49$ and $64$.} \label{section:n-49-64}

We first consider the case $\gimel=49$.\index{$\gimel=49$} Let $X$
be the hypersurface given by a general quasihomogeneous equation
of degree $21$ in $\mathbb{P}(1,3,5,6,7)$ with
$-K_X^3=\frac{1}{30}$. Then, the singularities of $X$ consist of
one quotient singular point $P$ of type $\frac{1}{6}(1,5,1)$, one
quotient singular point $Q$ of type $\frac{1}{5}(1,3,2)$,  and
three quotient singular  points of type $\frac{1}{3}(1,2,1)$.

We have the following commutative diagram:
$$
\xymatrix{
&&&Y_P\ar@{->}[lld]_{\pi_P}\ar@{->}[rrd]^{\eta_P}&&&\\%
&X\ar@{-->}[rrrr]_{\psi}&&&&\mathbb{P}(1,3,5)&}
$$
where \begin{itemize} \item $\psi$ is the natural projection,
\item $\pi_P$ is the Kawamata blow up at the point $P$ with
weights $(1,5,1)$, \item $\eta_P$ is an elliptic fibration.
\end{itemize}

We may assume that $X$ is given by a quasihomogeneous equation
$$
z^{3}t+z^{2}f_{11}(x,y,t,w)+zf_{16}(x,y,t,w)+f_{21}(x,y,t,w)=0,
$$
where $f_{i}$ is a quasihomogeneous polynomial of degree $i$. We
then see that there is another elliptic fibration as follows:
$$
\xymatrix{
&&&Y_Q\ar@{->}[lld]_{\pi_Q}\ar@{->}[rrd]^{\eta_Q}&&&\\%
&X\ar@{-->}[rrrr]_{\xi}&&&& G\subset\mathbb{P}^3&}
$$
where \begin{itemize} \item $\xi$ is the rational map given by the
linear system spanned by $\{ x^6, x^3y, y^2, t\}$

\item $\pi_Q$ is the Kawamata blow up at the point $Q$ with
weights $(1,3,2)$.

 \item  $G$
is the image of $\xi:X\dasharrow \mathbb{P}^3$ that is a quadratic
cone in $\mathbb{P}^3$,

\item $\eta_Q$ is an elliptic fibration.

\end{itemize}

It follows from Corollary~\ref{corollary:Ryder-a1} and
Lemmas~\ref{lemma:smooth-points},
\ref{lemma:special-singular-points-with-positive-c} that
$$
\Big\{P, Q\Big\}\supseteq\mathbb{CS}\Big(X, \frac{1}{n}\mathcal{M}\Big)\ne\varnothing.%
$$

The exceptional divisor $E_P$ of the birational morphism $\pi_P$
contains one singular point $P_1$ that is a quotient singularity
of type $\frac{1}{5}(1,4,1)$.
\begin{lemma}
\label{lemma:n-19-49-64-first-point} If the set $\mathbb{CS}
(Y_P,\frac{1}{n}\mathcal{M}_{Y_P})$ contains the point $P_1$, then
$\mathcal{M}=|-3K_{X}|$.
\end{lemma}

\begin{proof}

We can assume that
$$
f_{21}(x,y,t,w)=t^{3}y+t^{2}g_{9}(x,y,w)+tg_{15}(x,y,w)+g_{21}(x,y,w),
$$
where $g_{i}$ is a general quasihomogeneous polynomial of degree
$i$. Then, locally at the singular point $P$, the monomials $x$,
$z$, and $w$ can be considered as weighted local coordinates on
the threefold $X$ with weights $\mathrm{wt}(x)=1$,
$\mathrm{wt}(z)=5$, and $\mathrm{wt}(w)=7$, which implies that
locally at the singular point $P_1$, the birational morphism $\pi$
is given by the equations
$$
x=\tilde{x}\tilde{z}^{\frac{1}{6}},\ z=\tilde{z}^{\frac{5}{6}},\ w=\tilde{w}\tilde{z}^{\frac{1}{6}},%
$$
where $\tilde{x}$, $\tilde{z}$, and $\tilde{w}$ can be considered
as weighted local coordinates in a neighborhood of the singular
point $P_1$ with weights $\mathrm{wt}(\tilde{x})=1$,
$\mathrm{wt}(\tilde{z})=4$, and $\mathrm{wt}(\tilde{w})=1$.

Let $D$ be a general surface in the linear system $|-3K_X|$. Then,
it is given by an equation
$$
\lambda x^3+\mu y=0,
$$
where $(\lambda : \mu)\in\mathbb{P}^1$. Locally at the point
$P_1$, the proper transform $D_{Y_P}$ is given by
$$
\lambda \tilde{x}^3+\mu\Big(\epsilon_1\tilde{w}^3+\epsilon_2\tilde{w}^2\tilde{x}+\epsilon_3\tilde{w}\tilde{x}+\mbox{higher terms}\Big)=0,%
$$
which implies that $|-3K_{Y_P}|$ is the proper transform of
$|-3K_X|$  and that it has no base curves on the surface $E_P$.
The exceptional divisor $E_P$ is defined by $\tilde{z}=0$.

Let $\alpha:U\to Y_P$ be the Kawamata blow up at the point $P_1$
with weights $(1,4,1)$ and let $F_P$ be its   exceptional divisor.
Around the singular point of $F_P$, the birational morphism
$\alpha$ can be given by the equations
$$
\tilde{x}=\bar{x}\bar{z}^{\frac{1}{5}},\ \tilde{z}=\bar z^{\frac{4}{5}},\ \tilde{w}=\bar w\bar z^{\frac{1}{5}},%
$$
where $\bar x$, $\bar z$, and $\bar w$ are weighted local
coordinates on the variety $U$ in a neighborhood of the singular
point of the surface $F_P$ with weights
$\mathrm{wt}(\tilde{x})=1$, $\mathrm{wt}(\tilde z)=3$, and
$\mathrm{wt}(\tilde w)=1$.

The proper transform $D_U$  is given by an equation of the form
$$
\lambda \bar{x}^3+\mu\Big(\epsilon_1\bar{w}^3+\epsilon_2\bar{w}^2\bar{x}+\epsilon_3\bar{w}\bar{x}+\mbox{higher terms}\Big)=0,%
$$
which shows that $|-3K_{U}|$ is the proper transform of $|-3K_X|$
 and that $|-3K_{U}|$ does not have base curves on $F_P$. Hence, the base
locus of $|-3K_{U}|$ consists of the curve $C_U$ whose image to
$X$ is defined by the equations $x=y=0$.

Then, the inequality $-K_U\cdot C_U=-3K_U^3<0$ and the equivalence
$\mathcal{M}_U\sim_{\mathbb{Q}}-nK_U$ imply that the pencil
$\mathcal{M}_W$  coincides with the pencil $|-3K_U|$ by
Theorem~\ref{theorem:main-tool}.
\end{proof}

The exceptional divisor $E_Q$ contains two quotient singular
points $Q_1$ and $Q_2$ of types $\frac{1}{2}(1,1,1)$ and
$\frac{1}{3}(1, 1, 2)$, respectively.

\begin{lemma}
\label{lemma:n-19-49-64-two-points} If the set $\mathbb{CS}
(Y_Q,\frac{1}{n}\mathcal{M}_{Y_Q})$ contains the point $Q_1$, then
$\mathcal{M}=|-3K_{X}|$.
\end{lemma}
\begin{proof}
The same method for the proof of
Lemma~\ref{lemma:n-19-49-64-first-point} implies that
$\mathcal{M}=|-3K_X|$.
\end{proof}

\begin{lemma}
\label{lemma:n-19-49-64-second-point} The set $\mathbb{CS}
(Y_Q,\frac{1}{n}\mathcal{M}_{Y_Q})$ cannot contain the point
$Q_2$.
\end{lemma}
\begin{proof}
Let $\beta:W\to Y_Q$ be the Kawamata blow up at $Q_2$ with weights
$(1,1,2)$. Also, let  $F_Q$ be the exceptional divisor of the
birational morphism $\beta$. Let $\mathcal{R}$ be the linear
system given by the equations
$$
\lambda x^6+\mu x^3y+\nu t=0,
$$
where $(\lambda : \mu :\nu)\in\mathbb{P}^2$.  The base locus of
the linear system $\mathcal{R}$ consists of the curve $\tilde{C}$
given by the equations $x=t=0$, in other words, we have
$\tilde{C}=R\cdot S$, where $R$ is a general surface of the linear
system $\mathcal{R}$.

Around the point $Q$, the monomials $x$, $y$, and $w$ can be
considered as weighted local coordinates on $X$ with weights
$\mathrm{wt}(x)=1$, $\mathrm{wt}(y)=3$, and $\mathrm{wt}(w)=2$.
Also, around the singular point $Q_2$, the birational morphism
$\pi_Q$ is given by the equations
$$
x=\tilde{x}\tilde{y}^{\frac{1}{5}},\ y=\tilde{y}^{\frac{3}{5}},\ w=\tilde{w}\tilde{y}^{\frac{2}{5}},%
$$
where $\tilde{x}$, $\tilde{y}$, and $\tilde{w}$ are weighted local
coordinates around the singular point $Q_2$ with weights
$\mathrm{wt}(\tilde{x})=1$, $\mathrm{wt}(\tilde{y})=1$, and
$\mathrm{wt}(\tilde{w})=2$. The proper transform $R_Y$ is given by
an equation of the form
$$
\lambda\tilde x^6+\mu \tilde x^3+\nu\Big(\delta_1\tilde{w}^3+\delta_2\tilde{x}\tilde{w}+\delta_3\tilde{x}^2\tilde{w}^2+\delta_4\tilde{y}\tilde{x}^{2}+\delta_5\tilde{y}^3+\mbox{higher terms}\Big)=0,%
$$
where $\delta_{i}\in\mathbb{C}$, which tells us that the proper
transform  $\mathcal{R}_Y$  has no base curve on the exceptional
divisor $E_Q$.

Locally at the unique singular point of $F_Q$, the birational
morphism $\beta$ can be expressed by
$$
\tilde{x}=\bar x\bar w^{\frac{1}{3}},\ \tilde{y}=\bar y\bar w^{\frac{1}{3}},\ \tilde{w}=\bar{w}^{\frac{2}{3}},%
$$
where $\bar x$, $\bar y$, and $\bar w$ are local coordinates with
weight $1$.  Then, the surface $R_W$ is given by an equation of
the form
$$
\lambda\bar x^6\bar{w}+\mu \bar x^3+\nu\Big(\delta_1\bar{w}+\delta_2\bar{x}+\delta_3\bar{x}^2\bar{w}+\delta_4\bar{y}\bar{x}+\delta_5\bar{y}^3+\mbox{higher terms}\Big)=0,%
$$
which implies that the proper transform  $\mathcal{R}_W$ does not
have base curves on the surface $F_Q$ either.  The surface $R_W$
is normal. We see
$$
\left\{\aligned
&R_W\sim_{\mathbb{Q}}\big(\pi_Q\circ\beta\big)^{*}\Big(-6K_{X}\Big)-\frac{6}{5}\beta^{*}\big(E_Q\big)-F_Q,\\
&S_W\sim_{\mathbb{Q}}\big(\pi_Q\circ\beta\big)^{*}\Big(-K_{X}\Big)-\frac{1}{5}\beta^{*}\big(E_Q\big)-\frac{1}{3}F_Q.\\
\endaligned
\right.
$$
These equivalences show that $\tilde{C}_{W}^2<0$ on the normal
surface $R_W$ because $S_W\cdot R_W=\tilde{C}_W$.

Let $M$ be a general surface of the pencil $\mathcal{M}_W$. Then,
$$
M\Big\vert_{R_W}\equiv-nK_{W}\Big\vert_{R_W}\equiv nS_W\Big\vert_{R_W}\equiv n\tilde{C}_W%
$$
by Lemma~\ref{lemma:Kawamata}, which implies that
$\mathcal{M}=\mathcal{R}$ by Theorem~\ref{theorem:main-tool}. It
is a contradiction because the linear system $\mathcal{R}$ is not
a pencil. \end{proof}

\begin{proposition}
\label{proposition:n-49} If $\gimel=49$, then
$\mathcal{M}=|-a_{1}K_{X}|$.
\end{proposition}
\begin{proof}
By the previous lemmas, we may assume that
$$
\mathbb{CS}\Big(X, \frac{1}{n}\mathcal{M}\Big)=\{P, Q\}.%
$$
Furthermore, we also assume that the set $\mathbb{CS}(Y_P,
\frac{1}{n}\mathcal{M}_{Y_P})$ consists of the point $O$ whose
image to $X$ is the point $Q$.

 Let $\gamma:V\to Y_P$ be the Kawamata blow up
at the point $O$. Then, the proof of
Lemma~\ref{lemma:n-19-49-64-first-point} implies that $|-3K_{V}|$
is the proper transform of the pencil $|-3K_{X}|$ and the base
locus of $|-3K_{V}|$ consists of the curve $C_V$ whose image to
$X$ is the base curve of the pencil $|-3K_{X}|$. Then, we can
easily check that $\mathcal{M}=|-3K_X|$.
\end{proof}

We now consider the case $\gimel=64$.\index{$\gimel=64$} Let $X$
be the hypersurface given by a general quasihomogeneous equation
of degree $26$ in $\mathbb{P}(1,2,5,6,13)$ with
$-K_X^3=\frac{1}{30}$. Then, the singularities of $X$ consist of
one quotient singular point $P$ of type $\frac{1}{6}(1,5,1)$, one
quotient  singular point $Q$ of type $\frac{1}{5}(1,2,3)$,  and
four quotient singular points of type $\frac{1}{2}(1,1,1)$.

We see the following commutative diagram:
$$
\xymatrix{
&&&Y_P\ar@{->}[lld]_{\pi_P}\ar@{->}[rrd]^{\eta_P}&&&\\%
&X\ar@{-->}[rrrr]_{\psi}&&&&\mathbb{P}(1,2,5)&}
$$
where \begin{itemize} \item $\psi$ is  the natural projection,
\item $\pi_P$ is the Kawamata blow up at the point $P$ with
weights $(1,5,1)$, \item $\eta_P$ is an elliptic fibration.
\end{itemize}

One the other hand, we may assume that $X$ is given by a
quasihomogeneous equation
$$
z^{4}t+z^{3}f_{11}(x,y,t,w)+z^2f_{16}(x,y,t,w)+zf_{21}(x,y,t,w)+f_{26}(x,y,t,w)=0,
$$
where $f_{i}$ is a quasihomogeneous polynomial of degree $i$.
Therefore, there is another elliptic fibration as follows:
$$
\xymatrix{
&&&Y_Q\ar@{->}[lld]_{\pi_Q}\ar@{->}[rrd]^{\eta_Q}&&&\\%
&X\ar@{-->}[rrrr]_{\xi}&&&& G\subset\mathbb{P}^4&}
$$
where \begin{itemize} \item $\xi$ is the map given by the linear
system spanned by $\{ x^6, x^4y, x^2y^2, y^3, t\}$

\item $\pi_Q$ is the Kawamata blow up at the point $Q$ with
weights $(1,2,3)$.

 \item  $G$
is the image of $\xi:X\dasharrow \mathbb{P}^3$ that is isomorphic to $\mathbb{P}(1,1,3)$,

\item $\eta_Q$ is an elliptic fibration.

\end{itemize}

\begin{proposition}
\label{proposition:n-64} If $\gimel=64$, then
$\mathcal{M}=|-a_{1}K_{X}|$.
\end{proposition}
\begin{proof}
The proof is the same as that of
Proposition~\ref{proposition:n-49}.
\end{proof}

\section{Case $\gimel=56$, hypersurface of degree $24$ in
$\mathbb{P}(1,2,3,8,11)$.}\index{$\gimel=56$} \label{section:n-56}

 The threefold  $X$ is a
 general hypersurface of degree $24$ in $\mathbb{P}(1,2,3,8,11)$
with $-K_{X}^{3}=\frac{1}{22}$. Its singularities consist of three
points that are quotient singularities of type
$\frac{1}{2}(1,1,1)$ and the point $O=(0:0:0:0:1)$ that is a
quotient singularity of type $\frac{1}{11}(1,3,8)$.

Before we proceed, let us first describe some birational
transformations of the hypersurface $X$ with elliptic fibrations,
which are useful to explain the geometrical nature of our proof.
There is a commutative diagram
$$
\xymatrix{
&&&Z\ar@{->}[dll]_{\nu}\ar@{-->}[r]^{\zeta}&Z^{\prime}\ar@{->}[r]^{\sigma}&\bar{Z}^{\prime}\ar@{^{(}->}[r]
&\mathbb{P}(1,2,5,14,21)\ar@{-->}[ddl]^{\phi}\\%
U\ar@{->}[d]_{\alpha}&W\ar@{->}[l]_{\beta}&&Y\ar@{->}[d]^{\eta}\ar@{->}[ll]_{\gamma}
&V\ar@{->}[l]_{\xi}\ar@{-->}[r]^{\chi}\ar@{->}[lu]_{\omega}&V^{\prime}\ar@{->}[d]^{\upsilon}&\\%
X\ar@{-->}[rrr]_{\psi}&&&\mathbb{P}(1,2,3)\ar@{-->}[rr]_{\rho}&&\mathbb{P}(1,2,5)&
}
$$
where \begin{itemize}
\item $\psi$ and $\phi$ are natural
projections,

\item $\alpha$ is the Kawamata blow up at the point $O$ with
weights $(1,3,8)$,

\item $\beta$ is the Kawamata blow up with weights $(1,3,5)$ at
the singular point $Q$ contained in the exceptional divisor $E$ of
$\alpha$ that is a quotient singularity of type
$\frac{1}{8}(1,3,5)$,

 \item $\gamma$ is the Kawamata blow up with weights $(1,3,2)$ at the
singular point $Q_1$ of contained in the exceptional divisor $F$
of $\beta$ that is a quotient singularity of type
$\frac{1}{5}(1,3,2)$,

\item $\nu$ is the Kawamata blow up with weights $(1,1,2)$ at the
singular point $Q_2$ of contained in the exceptional divisor of
$\beta$ that is a quotient singularity of type
$\frac{1}{3}(1,1,2)$,

\item $\xi$ is the Kawamata blow up with weights $(1,1,2)$ at the
point  $\bar{Q}_2$ whose image to $W$ is the point $Q_2$,

\item $\omega$ is the Kawamata blow up with weights $(1,3,2)$ at
the point $\bar{Q}_1$ whose image to $W$ is the point $Q_1$,

 \item $\eta$ and $\upsilon$ are elliptic fibrations,

\item the maps $\zeta$ and $\chi$ are compositions of antiflips,

\item the birational morphism $\sigma$ is given by the
plurianticanonical linear system of $Z^{\prime}$,

\item the rational map $\rho$ is a  toric map,
\end{itemize}
The exceptional divisor $E$ of the birational morphism $\alpha$
contains two singular points $P$ and $Q$ of $U$ that are quotient
singularities of types $\frac{1}{3}(1,1,2)$ and
$\frac{1}{8}(1,3,5)$, respectively. Meanwhile, the exceptional
divisor $F$ of the birational morphism $\beta$ also contains two
singular points $Q_1$ and $Q_2$ of $W$ that are quotient
singularities of types $\frac{1}{5}(1,3,2)$ and
$\frac{1}{3}(1,1,2)$, respectively.

\begin{remark}
\label{remark:n-56-anticanonical-models} The divisors
$-K_{Z^{\prime}}$, $-K_{U}$, and $-K_{W}$ are nef and big. Thus,
the anticanonical models of the threefolds $Z^{\prime}$, $U$ and
$W$ are Fano threefolds with canonical singularities. The
anticanonical model of $Z^{\prime}$ is a hypersurface
$\bar{Z}^{\prime}$ of degree $42$ in $\mathbb{P}(1,2,5,14,21)$.
The anticanonical model of $U$ is a hypersurface  of degree $26$
in $\mathbb{P}(1,2,3,8,13)$ and the anticanonical model of $W$ is
a hypersurface  of degree $30$ in $\mathbb{P}(1,2,3,10,15)$.
\end{remark}

For the convenience, we denote the pencil $|-2K_X|$ by
$\mathcal{B}$. In addition, a general surface in $\mathcal{B}$ is
denoted by $B$ and a general surface in $\mathcal{M}$ by $M$

It follows from Corollary~\ref{corollary:Ryder-a1} and
Lemmas~\ref{lemma:smooth-points},
\ref{lemma:special-singular-points-with-positive-c} that we may
assume that $\mathbb{CS}(X, \frac{1}{n}\mathcal{M})=\{O\}$.

\begin{lemma}
\label{lemma:n-56-points-P5} If the set $\mathbb{CS}(U,
\frac{1}{n}\mathcal{M}_{U})$ contains the point $P$, then
$\mathcal{M}=\mathcal{B}$.
\end{lemma}

\begin{proof}
Let $\beta_P:U_P\to U$ be the Kawamata blow up at the point $P$
and $E_P$ be its exceptional divisor. For a general surface $D$ in
$|-8K_{X}|$, we have
\[D_{U_P}\sim_{\mathbb{Q}} (\alpha\circ\beta_P)^*(-8K_X)-\frac{8}{11}\beta_P^*(E)
-\frac{2}{3}E_P.\] Because the base locus of the proper transform
of the linear system $|-8K_X|$ on $U_P$ does not contain  any
curve, the divisor $D_{U_P}$ is nef and big.

Since $M_{U_P}\sim_{\mathbb{Q}} nS_{U_P}$ by
Lemma~\ref{lemma:Kawamata} and $B_{U_P}\sim_{\mathbb{Q}}
2S_{U_P}$, we obtain \[D_{U_P}\cdot B_{U_P}\cdot S_{U_P}
=2n\Big(\beta_P^*(-8K_U)
-\frac{2}{3}E_P\Big)\cdot\Big(\beta_P^*(-K_U)
-\frac{1}{3}E_P\Big)^{2}=0.
\]
It implies  $\mathcal{M}=\mathcal{B}$ by
Theorem~\ref{theorem:main-tool}.
\end{proof}

Due to Theorem~\ref{theorem:Noether-Fano} and
Lemma~\ref{lemma:Cheltsov-Kawamata}, we may assume that the set
$\mathbb{CS}(U, \frac{1}{n}\mathcal{M}_{U})$ consists of the
singular point $Q$. Thus, it follows from
Theorem~\ref{theorem:Noether-Fano}, Lemmas~\ref{lemma:curves},
\ref{lemma:Cheltsov-Kawamata} that
$$
\varnothing\ne\mathbb{CS}\Big(W,
\frac{1}{n}\mathcal{M}_{W}\Big)\subseteq\Big\{Q_1,\ Q_2\Big\}.
$$
Now, we consider some local computation. We may assume that $X$ is
given by the equation
$$
w^{2}y+wf_{13}(x,y,z,t)+f_{24}(x,y,z,t)=0,
$$
where $f_{i}(x,y,z,t)$ is a  general quasihomogeneous polynomial
of degree $i$. The surface $B$ is given by the equation $\lambda
x^{2}+\mu y=0$, where $(\lambda : \mu)\in\mathbb{P}^{1}$. The base
locus of $\mathcal{B}$ consists of the irreducible curve $C$ that
is given by $x=y=0$. We have $B\cdot S=C$.

In a neighborhood of $O$, the monomials $x$, $z$, and $t$ can be
considered as weighted local coordinates on $X$ such that
$\mathrm{wt}(x)=1$, $\mathrm{wt}(z)=3$, and $\mathrm{wt}(z)=8$.
Then, in a neighborhood of the singular point $O$, the surface $B$
can be given by equation
$$
\lambda x^{2}+\mu\Big(\epsilon_{1}x^{13}+\epsilon_{2}zx^{10}+\epsilon_{3}z^{2}x^{7}+\epsilon_{4}z^{3}x^{4}+\epsilon_{5}z^{4}x
+\epsilon_{6}tx^{5}+\epsilon_{7}tzx^{2}+\epsilon_{8}t^{3}+\epsilon_{9}z^{8}+\mbox{other terms}\Big)=0,%
$$
where $\epsilon_{i}\in\mathbb{C}$. In a neighborhood of the
singular point $Q$, the birational morphism $\alpha$ can be given
by the equations
$$
x=\bar{x}\bar{t}^{\frac{1}{11}},\ z=\bar{z}\bar{t}^{\frac{3}{11}},\ t=\bar{t}^{\frac{8}{11}},%
$$
where $\bar{x}$, $\bar{z}$, and $\bar{t}$ are weighted local
coordinates on $U$ in a neighborhood of the singular point $Q$
such that $\mathrm{wt}(\bar{x})=1$, $\mathrm{wt}(\bar{z})=3$, and
$\mathrm{wt}(\bar{t})=8$. Thus, in a neighborhood of the singular
point $Q$, the divisor $E$ is given by the equation $\bar{t}=0$,
the divisor $S_{U}$ is given by $\bar{x}=0$, and the divisor
$B_{U}$ is given by the equation
$$
\lambda
\bar{x}^{2}+\mu\Big(\epsilon_{1}\bar{x}^{13}\bar{t}+\cdots+\epsilon_{5}\bar{z}^{4}\bar{x}\bar{t}+\epsilon_{6}\bar{t}\bar{x}^{5}+\epsilon_{7}\bar{t}\bar{z}\bar{x}^{2}+\epsilon_{8}\bar{t}^{2}+\epsilon_{9}\bar{z}^{8}\bar{t}^{2}+\mbox{other terms}\Big)=0,%
$$
which implies that $B_{U}\sim_{\mathbb{Q}} 2S_{U}$ and  the base
locus of $\mathcal{B}_{U}$ is the union of $C_{U}$ and the curve
$L\subset E$ that is given by $\bar{x}=\bar{t}=0$. We have
$E\cong\mathbb{P}(1,3,8)$ and the curve $L$ is the unique curve in
$|\mathcal{O}_{\mathbb{P}(1,3,8)}(1)|$ on the surface $E$. The
surface $B_{U}$ is not normal. Indeed, $B_{U}$ is singular at a
generic point of $L$. We have $S_{U}\cdot B_{U}=C_{U}+2L$ and
$E\cdot B_{U}=2L$, which implies that $S_{U}\cdot C_{U}=0$ and
$S_{U}\cdot L=\frac{1}{24}$.

\begin{lemma}
\label{lemma:n-56-P7} If the set $\mathbb{CS}(W,
\frac{1}{n}\mathcal{M}_{W})$ consists of the point $Q_2$, then
$\mathcal{M}=\mathcal{B}$.
\end{lemma}

\begin{proof}
In a neighborhood of $Q_2$, the birational morphism $\beta$ can be
given by the equations
$$
\bar{x}=\tilde{x}\tilde{z}^{\frac{1}{8}},\ \bar{z}=\tilde{z}^{\frac{3}{8}},\ \bar{t}=\tilde{t}\tilde{z}^{\frac{5}{8}},%
$$
where $\tilde{x}$, $\tilde{z}$, and $\tilde{t}$ are weighted local
coordinates on $W$ in a neighborhood of $Q_2$ such that
$\mathrm{wt}(\tilde{x})=1$, $\mathrm{wt}(\tilde{z})=1$, and
$\mathrm{wt}(\tilde{t})=2$. Thus, in a neighborhood of the
singular point $Q_2$, the divisor $F$ is given by the equation
$\tilde{z}=0$, the divisor $S_{W}$ is given by $\tilde{x}=0$, the
divisor $E_{W}$ is given by $\tilde{t}=0$, and the divisor $B_{W}$
is given by the equation
$$
\lambda
\tilde{x}^{2}+\mu\Big(\epsilon_{7}\tilde{t}\tilde{z}\tilde{x}^{2}+\epsilon_{8}\tilde{t}^{2}\tilde{z}+\epsilon_{9}\tilde{z}^{4}\tilde{t}^{2}+\mbox{other terms}\Big)=0,%
$$
which implies that $B_{W}\sim_{\mathbb{Q}} 2S_{W}$, the base locus
of $\mathcal{B}_{W}$ is the union of $C_{W}$, $L_{W}$, and the
curve $L^{\prime}\subset F$ that is given by
$\tilde{x}=\tilde{z}=0$.

The surface $F$ is isomorphic to $\mathbb{P}(1,3,5)$ and the curve
$L^{\prime}$ is the unique curve of the linear system
$|\mathcal{O}_{\mathcal{P}(1,3,5)}(1)|$ on the surface $F$. The
surface $B_{W}$ is smooth at a generic  point of $L^{\prime}$. We
have
$$S_{W}\cdot B_{W}=C_{W}+2L_{W}+L^{\prime}, \ \ E_{W}\cdot
B_{W}=2L_{W}, \ \ F\cdot B_{W}=2L^{\prime},$$ which implies that
$$S_{W}\cdot C_{W}=0, \ \ S_{W}\cdot L_{W}=0, \ \ S_{W}\cdot
L^{\prime}=\frac{1}{15}$$ because
$$
\left\{\aligned
&S_{W}\sim_{\mathbb{Q}}(\alpha\circ\beta)^*(-K_X)-\frac{1}{11}\beta^*(E)
-\frac{1}{8}F,\\
&B_{W}\sim_{\mathbb{Q}}(\alpha\circ\beta)^*(-2K_X)-\frac{2}{11}\beta^*(E)
-\frac{2}{8}F,\\
&E_{W}\sim_{\mathbb{Q}} \beta^*(E)
-\frac{5}{8}F.\\
\endaligned
\right.
$$

Let $R$ be the exceptional  divisor of $\nu$. Let $O_2$ be the
singular point of $Z$ that is contained in $R$. Then,
$R\cong\mathbb{P}(1,1,2)$ and $O_2$ is a quotient singularity of
type $\frac{1}{2}(1,1,1)$ on the threefold $Z$. In a neighborhood
of $O_2$, the birational morphism $\nu$ can be given by the
equations
$$
\tilde{x}=\hat{x}\hat{t}^{\frac{1}{3}},\ \tilde{z}=\hat{z}\hat{t}^{\frac{1}{3}},\ \tilde{t}=\hat{t}^{\frac{2}{3}},%
$$
where $\hat{x}$, $\hat{z}$ and $\hat{t}$ are weighted local
coordinates on $Z$ in a neighborhood of $O_2$ with weight $1$.
Thus, in a neighborhood of the singular point $O_2$, the divisor
$R$ is given by the equation $\hat{t}=0$, the divisor $S_{Z}$ is
given by $\hat{x}=0$, the divisor $E_{Z}$ does not pass through
the point $O_2$, the divisor $F_{Z}$ is given by $\bar{z}=0$, and
the divisor $B_{Z}$ is given by the equation
$$
\lambda
\hat{x}^{2}+\mu\Big(\epsilon_{8}\hat{t}\hat{z}+\epsilon_{9}\hat{z}^{4}\hat{t}^{2}+\mbox{other terms}\Big)=0,%
$$
which implies that $B_{Z}\sim_{\mathbb{Q}} 2S_{Z}$, the base locus
of $\mathcal{B}_{Z}$ consists of $C_{Z}$, $L_{Z}$,
$L^{\prime}_{Z}$, and the curve $L^{\prime\prime}$ that is given
by the equations $\hat{x}=\hat{t}=0$. The curve $L^{\prime\prime}$
is the unique curve in $|\mathcal{O}_{\mathbb{P}(1,1,2)}(1)|$ on
the surface $R$. The surface $B_{Z}$ is smooth at a generic point
of $L^{\prime\prime}$. Therefore, we obtain
$$
S_{Z}\cdot B_{Z}=C_{Z}+2L_{Z}+L^{\prime}_{Z}+L^{\prime\prime},\
E_{Z}\cdot B_{Z}=2L_{Z},\ F_{Z}\cdot B_{Z}=2L^{\prime}_{Z},\
R\cdot B_{Z}=2L^{\prime\prime},
$$
which gives $$S_{Z}\cdot C_{Z}=0, \ \  S_{Z}\cdot
L_{Z}=-\frac{1}{3}, \ \ S_{Z}\cdot L^{\prime}_{Z}=-\frac{1}{10}, \
\ S_{Z}\cdot L^{\prime\prime}=\frac{1}{2}$$ because
$$
\left\{\aligned
&S_{Z}\sim_{\mathbb{Q}}(\alpha\circ\beta\circ\nu)^*(-K_X)
-\frac{1}{11}(\beta\circ\nu)^*(E)
-\frac{1}{8}\nu^*(F)-\frac{1}{3}R,\\
&B_{Z}\sim_{\mathbb{Q}}(\alpha\circ\beta\circ\nu)^*(-2K_X)
-\frac{2}{11}(\beta\circ\nu)^*(E)
-\frac{2}{8}\nu^*(F)-\frac{2}{3}R,\\
&E_{Z}\sim_{\mathbb{Q}}(\beta\circ\nu)^*(E)
-\frac{5}{8}\nu^*(F)-\frac{2}{3}R,\\
&F_{Z}\sim_{\mathbb{Q}}\nu^*(F)-\frac{1}{3}R.\\
\endaligned
\right.
$$
In particular, the
curves $L_{Z}$ and $L^{\prime}_{Z}$ are the only curves on the
variety $Z$ that have negative intersection with the divisor
$-K_{Z}$.

Due to Lemma~\ref{lemma:Cheltsov-Kawamata},  either the set
$\mathbb{CS}(Z, \frac{1}{n}\mathcal{M}_{Z})$ contains the point
$O_2$ or the log pair $(Z, \frac{1}{n}\mathcal{M}_{Z})$ is
terminal.

We first suppose that the log pair $(Z,
\frac{1}{n}\mathcal{M}_{Z})$ is not terminal. Then, the set
$\mathbb{CS}(Z, \frac{1}{n}\mathcal{M}_{Z})$ must  contain the
point $O_2$. Let $\pi_2:Z_2\to Z$ be the Kawamata blow up at the
point $O_2$ and $H$ be the exceptional divisor of $\pi_2$. Then,
our local calculations imply that $B_{Z_2}\sim_{\mathbb{Q}}
2S_{Z_2}$ and the base locus of $\mathcal{B}_{Z_2}$ consists of
the curves $C_{Z_2}$, $L_{Z_2}$, $L^{\prime}_{Z_2}$, and
$L^{\prime\prime}_{Z_2}$. Furthermore, we have
$$
S_{Z_2}\cdot
B_{T}=C_{Z_2}+2L_{Z_2}+L^{\prime}_{Z_2}+L^{\prime\prime}_{Z_2},\ \
E_{Z_2}\cdot B_{Z_2}=2L_{Z_2}, $$
$$F_{Z_2}\cdot B_{Z_2}=2L^{\prime}_{Z_2},\ \  R_{Z_2}\cdot
B_{Z_2}=2L^{\prime\prime}_{Z_2},
$$
which implies that
$$S_{Z_2}\cdot C_{Z_2}=0, \ \ S_{Z_2}\cdot
L_{Z_2}=-\frac{1}{3}, \ \ S_{Z_2}\cdot
L^{\prime}_{Z_2}=-\frac{3}{5}, \ \ S_{Z_2}\cdot
L^{\prime\prime}_{Z_2}=0,$$ because
$$
\left\{\aligned
&S_{Z_2}\sim_{\mathbb{Q}}(\alpha\circ\beta\circ\nu\circ\pi_2)^*(-K_X)
-\frac{1}{11}(\beta\circ\nu\circ\pi_2)^*(E)
-\frac{1}{8}(\nu\circ\pi_2)^*(F)-\frac{1}{3}\pi_2^*R-\frac{1}{2}H,\\
&E_{Z_2}\sim_{\mathbb{Q}}(\beta\circ\nu\circ\pi_2)^*(E)
-\frac{5}{8}(\nu\circ\pi_2)^*(F)-\frac{2}{3}\pi_2^*R,\\
&F_{Z_2}\sim_{\mathbb{Q}}(\nu\circ\pi_2)^*(F)-\frac{1}{3}\pi_2^*R-\frac{1}{2}H,\\
&R_{Z_2}\sim_{\mathbb{Q}}\pi_2^*R-\frac{1}{2}H,\\
\endaligned
\right.
$$
The curves $L_{Z_2}$ and $L^{\prime}_{Z_2}$ are the only curves on
the variety $Z_2$ that have negative intersection with the divisor
$-K_{Z_2}$. Moreover, we see
$$
\Big(B_{Z_2}+(\beta\circ\nu\circ\pi_2)^*(-16K_{U})+(\nu\circ\pi_2)^*(-18K_{W})\Big)\cdot L_{Z_2}=0,$$
$$\Big(B_{Z_2}+(\beta\circ\nu\circ\pi_2)^*(-16K_{U})+(\nu\circ\pi_2)^*(-18K_{W})\Big)\cdot L^{\prime}_{Z_2}=0,%
$$
and hence the divisor
$D_{Z_2}:=B_{Z_2}+(\beta\circ\nu\circ\pi_2)^*(-16K_{U})+(\nu\circ\pi_2)^*(-18K_{W})$
is nef and big because $-K_{U}$ and $-K_{W}$ are nef and big.
Therefore, we obtain
$$
D_{Z_2}\cdot B_{Z_2}\cdot M_{Z_2}=0,
$$
and hence $\mathcal{M}=\mathcal{B}$ by
Theorem~\ref{theorem:main-tool}.

For now, we suppose that  the log pair $(Z,
\frac{1}{n}\mathcal{M}_{Z})$ is terminal. We will derive a
contradiction from this assumption, so that the set
$\mathbb{CS}(Z, \frac{1}{n}\mathcal{M}_{Z})$ must contain the
point $O_2$.

The log pair $(Z, \epsilon B_{Z})$ is terminal for some rational
number $\epsilon>\frac{1}{2}$ but the divisor $K_{Z}+\epsilon
B_{Z}$ has nonnegative intersection with all curves on the variety
$Z$ except the curves $L_{Z}$ and $L^{\prime}_{Z}$. It follows
from \cite{Sho93} that there is a composition of antiflips
$\zeta:Z\dasharrow Z^{\prime}$ and the divisor $-K_{Z^{\prime}}$
is nef. Then, the singularities of the log pair $(Z^{\prime},
\frac{1}{n}\mathcal{M}_{Z^{\prime}})$ are terminal because the
singularities of the log pair $(Z, \frac{1}{n}\mathcal{M}_{Z})$
are terminal and the rational map $\zeta$ is a log flop with
respect to the log pair $(Z, \frac{1}{n}\mathcal{M}_{Z})$.

We obtain
$$\left\{\aligned
&S_{Z} \sim_{\mathbb{Q}} (\alpha\circ\beta\circ\nu)^*(-K_X)-\frac{1}{11}(\beta\circ\nu)^*(E)-\frac{1}{8}\nu^*(F)-\frac{1}{3}R \\
&\phantom{S_Z} \sim_{\mathbb{Q}} (\alpha\circ\beta\circ\nu)^*(-K_X)-\frac{1}{11}E_{Z}-\frac{2}{11}F_{Z}-\frac{5}{11}R,\\
&S^{y}_{Z} \sim_{\mathbb{Q}} (\alpha\circ\beta\circ\nu)^*(-2K_X)-\frac{13}{11}(\beta\circ\nu)^*(E)-\frac{5}{8}\nu^*(F)-\frac{2}{3}R \\
&\phantom{D^y_Z}\sim_{\mathbb{Q}}(\alpha\circ\beta\circ\nu)^*(-2K_X)-\frac{13}{11}E_{Z}-\frac{15}{11}F_{Z}-\frac{21}{11}R,\\
&S^{z}_{Z} \sim_{\mathbb{Q}} (\alpha\circ\beta\circ\nu)^*(-3K_X)-\frac{3}{11}(\beta\circ\nu)^*(E)-\frac{3}{8}\nu^*(F)-\frac{1}{3}R \\
&\phantom{D^z_Z}\sim_{\mathbb{Q}} (\alpha\circ\beta\circ\nu)^*(-3K_X)-\frac{3}{11}E_{Z}-\frac{6}{11}F_{Z}-\frac{4}{11}R,\\
&S^{t}_{Z} \sim_{\mathbb{Q}} (\alpha\circ\beta\circ\nu)^*(-8K_X)-\frac{8}{11}(\beta\circ\nu)^*(E)\\
&\phantom{D^t_Z}\sim_{\mathbb{Q}}(\alpha\circ\beta\circ\nu)^*(-8K_X)-\frac{8}{11}E_{Z}-\frac{5}{11}F_{Z}-\frac{7}{11}R.\\
\endaligned\right.
$$
from
$$F_{Z}\sim_{\mathbb{Q}}\nu^*(F)-\frac{1}{3}R,\ \
E_{Z} \sim_{\mathbb{Q}}(\beta\circ\nu)^*
(E)-\frac{5}{8}\nu^*(F)-\frac{2}{3}R.$$ Thus, the pull-backs of
the rational functions $\frac{y}{x^{2}}$, $\frac{zy}{x^{5}}$ and
$\frac{ty^{3}}{x^{14}}$ are contained in the linear systems
$|2S_{Z}|$, $|5S_{Z}|$ and $|14S_{Z}|$, respectively. In
particular, the complete linear system $|-70K_{Z}|$ induces a
dominant rational map $Z\dasharrow\mathbb{P}(1,2,5,14)$. Thus, the
anticanonical divisor $-K_{Z^{\prime}}$ is nef and big. It
contradicts Theorem~\ref{theorem:Noether-Fano} because the log
pair $(Z^{\prime}, \frac{1}{n}\mathcal{M}_{Z^{\prime}})$ is
terminal.
\end{proof}

Due to the lemma above, we may assume that the set $\mathbb{CS}(W,
\frac{1}{n}\mathcal{M}_{W})$ contains the point $Q_1$. In
particular, the set $\mathbb{CS}(Y, \frac{1}{n}\mathcal{M}_{Y})$
is not empty and each member of the linear system
$\mathcal{M}_{Y}$ is contracted to a curve by the morphism
$\eta$.

Let $G$ be the exceptional  divisor of $\gamma$. Then, $G$
contains two singular points $Q_1'$ and $Q_2'$ of $Y$ that are
quotient singularities of types $\frac{1}{2}(1,1,1)$ and
$\frac{1}{3}(1,1,2)$, respectively. Then,
$$
\mathbb{CS}\Big(Y, \frac{1}{n}\mathcal{M}_{Y}\Big)\subseteq\Big\{Q_{1}',\ Q_{2}',\ \bar{Q}_{2}\Big\},%
$$
where $\bar{Q}_{2}$ is the point on $Y$ whose image to $W$ by
$\gamma$ is the point $Q_2$.

In a neighborhood of $Q_1$, the birational morphism $\beta$ can be
given by the equations
$$
\bar{x}=\tilde{x}\tilde{z}^{\frac{1}{8}},\ \bar{z}=\tilde{z}\tilde{t}^{\frac{3}{8}},\ \bar{t}=\tilde{t}^{\frac{5}{8}},%
$$
where $\tilde{x}$, $\tilde{z}$, and $\tilde{t}$ are weighted local
coordinates on $W$ in a neighborhood of $Q_1$ such that
$\mathrm{wt}(\tilde{x})=1$, $\mathrm{wt}(\tilde{z})=3$, and
$\mathrm{wt}(\tilde{t})=2$. Thus, in a neighborhood of the
singular point $Q_1$, the divisor $F$ is given by the equation
$\tilde{t}=0$, the divisor $S_{W}$ is given by $\tilde{x}=0$, the
divisor $E_{W}$ does not pass though the point $Q_1$, and the
divisor $B_{W}$ is given by the equation
$$
\lambda
\tilde{x}^{2}+\mu\Big(\epsilon_{8}\tilde{t}+\epsilon_{9}\tilde{z}^{8}\tilde{t}^{4}+\mbox{other terms}\Big)=0.%
$$
Therefore,  $B_{W}\sim_{\mathbb{Q}} 2S_{W}$ and  the base locus of
$\mathcal{B}_{W}$ is the union of $C_{W}$, $L_{W}$ and the curve
$L^{\prime}$ that is given by the equations
$\tilde{x}=\tilde{t}=0$. We have
$$
S_{W}\cdot B_{W}=C_{W}+2L_{W}+L^{\prime},\ E_{W}\cdot B_{W}=2L_{W}, \ F\cdot B_{W}=2L^{\prime},%
$$
which gives us

$$S_{W}\cdot C_{W}=S_{W}\cdot L_{U}=0, \ \ S_{W}\cdot
L^{\prime}=\frac{1}{15}.$$

In a neighborhood of $Q_2$, the birational morphism $\gamma$ can
be given by the equations
$$
\tilde{x}=\hat{x}\hat{t}^{\frac{1}{5}},\ \tilde{z}=\hat{z}^{\frac{3}{5}},\ \tilde{t}=\hat{t}\hat{z}^{\frac{2}{5}},%
$$
where $\hat{x}$, $\hat{z}$, and $\hat{t}$ are weighted local
coordinates on $Y$ in the neighborhood of $Q_2'$  such that
$\mathrm{wt}(\hat{x})=1$, $\mathrm{wt}(\hat{z})=1$, and
$\mathrm{wt}(\hat{t})=2$. Thus, in a neighborhood of the singular
point $Q_2'$, the divisor $G$ is given by the equation
$\hat{z}=0$, the divisor $S_{Y}$ is given by $\hat{x}=0$, the
divisor $F_{Y}$ is given by the equation $\bar{t}=0$, and the
divisor $B_{Y}$ is given by the equation
$$
\lambda
\hat{x}^{2}+\mu\Big(\epsilon_{8}\hat{t}+\epsilon_{9}\hat{z}^{6}\hat{t}^{2}+\mbox{other terms}\Big)=0.%
$$
Thus, $B_{Y}\sim_{\mathbb{Q}} 2S_{Y}$ and that the base locus of
$\mathcal{B}_{Y}$ is the union of the irreducible curves $C_{Y}$,
$L_{Y}$, and $L^{\prime}_{Y}$. We have

$$\left\{\aligned
&S_{Y}\sim_{\mathbb{Q}}
(\alpha\circ\beta\circ\gamma)^*(-K_X)-\frac{1}{11}(\beta\circ\gamma)^*(E)-\frac{1}{8}\gamma^*(F)-\frac{1}{5}G,\\
&E_{Y}\sim_{\mathbb{Q}} (\beta\circ\gamma)^*(E)-\frac{5}{8}\gamma^*(F),\\
&F_{Y}\sim_{\mathbb{Q}} \gamma^*(F)-\frac{2}{5}G,\\
\endaligned\right.
$$
and $$S_{Y}\cdot C_{Y}=S_{Y}\cdot L_{Y}=S_{Y}\cdot
L^{\prime}_{Y}=0,$$ which simply means that $C_{Y}$, $L_{Y}$ and
$L^{\prime}_{Y}$ are components of a fiber of $\eta$.

\begin{lemma}
\label{lemma:n-56-P11} If the set  $\mathbb{CS}(Y,
\frac{1}{n}\mathcal{M}_{Y})$ contains $Q_2'$, then
$\mathcal{M}=\mathcal{B}$.
\end{lemma}

\begin{proof}
Let $\sigma_2:Y_2\to Y$ be the Kawamata blow up at the point
$Q_2'$  and let $H_2$ be the exceptional divisor of $\sigma_2$.
Then, our local calculations imply that $B_{Y_2}\sim_{\mathbb{Q}}
2S_{Y_2}$ and the base locus of $\mathcal{B}_{Y_2}$ is the union
of curves $C_{Y_2}$, $L_{Y_2}$, and $L^{\prime}_{Y_2}$. Thus, we
have
$$
S_{Y_2}\cdot B_{Y_2}=C_{Y_2}+2L_{Y_2}+L^{\prime}_{Y_2},\ \
E_{Y_2}\cdot B_{Y_2}=2L_{Y_2},\ \  F_{Y_2}\cdot
B_{Y_2}=2L^{\prime}_{Y_2},
$$
which implies that $$S_{Y_2}\cdot C_{Y_2}=0, \ \ S_{Y_2}\cdot
L_{Y_2}=0, \ \ S_{Y_2}\cdot L^{\prime}_{Y_2}=-\frac{1}{3}$$
because

$$\left\{\aligned
&S_{Y_2}\sim_{\mathbb{Q}}
(\alpha\circ\beta\circ\gamma\circ\sigma_2)^*(-K_X)-\frac{1}{11}(\beta\circ\gamma\circ\sigma_2)^*(E)
-\frac{1}{8}(\gamma\circ\sigma_2)^*(F)-\frac{1}{5}\sigma_2^*(G)-\frac{1}{3}H_2,\\
&E_{Y_2}\sim_{\mathbb{Q}} (\beta\circ\gamma\circ\sigma_2)^*(E)-\frac{5}{8}(\gamma\circ\sigma_2)^*(F),\\
&F_{Y_2}\sim_{\mathbb{Q}} (\gamma\circ\sigma_2)^*(F)-\frac{2}{5}\sigma_2^*(G)-\frac{2}{3}H_2,\\
&G_{Y_2}\sim_{\mathbb{Q}} \sigma_2^*(G)-\frac{1}{3}H_2,\\
\endaligned\right.
$$
The curve $L^{\prime}_{Y_2}$ is the only curve on $Y_2$ that has
negative intersection with $-K_{Y_2}$. Moreover, we have
$(S_{Y_2}+(\gamma\circ\sigma_2)^*(-5K_{W}))\cdot
L_{Y_2}^{\prime}=0$, which implies that the divisor
$S_{Y_2}+(\gamma\circ\sigma_2)^*(-5K_{W})$  is nef and big because
$-K_{W}$ is nef and big. Therefore,
$$(S_{Y_2}+(\gamma\circ\sigma_2)^*(-5K_{W}))\cdot B_{Y_2}\cdot
M_{Y_2}=0$$ by Lemma~\ref{lemma:Kawamata}, and hence
$\mathcal{M}=\mathcal{B}$ by Theorem~\ref{theorem:main-tool}
\end{proof}

\begin{lemma}
\label{lemma:n-56-P10} The set $\mathbb{CS}(Y,
\frac{1}{n}\mathcal{M}_{Y})$ cannot contain the point $Q_1'$.
\end{lemma}

\begin{proof}
Suppose that the set $\mathbb{CS}(Y, \frac{1}{n}\mathcal{M}_{Y})$
contains the point $Q_1'$. Let $\sigma_1:Y_1\to Y$ be the Kawamata
blow up at the point $Q_1'$.

Let $\mathcal{D}$ be a general pencil in the linear system
$|-3K_{X}|$. Then, the base curve of $\mathcal{D}$ is the curve
$\bar{C}$ given by $x=z=0$. Moreover, the base locus of
$\mathcal{D}_{Y_1}$ consists of the curve $\bar{C}_{Y_1}$. Thus,
we see that $\bar{C}_{Y_1}=S_{Y_1}\cdot D_{Y_1}$ for a general
surface $D_{Y_1}$ in $\mathcal{D}_{Y_1}$. On the other hand, we
have $D_{Y_1}\cdot C_{Y_1}<0$, which implies that $n=3$ and
$\mathcal{M}_{Y_1}=\mathcal{D}_{Y_1}$ by
Theorem~\ref{theorem:main-tool}. However,
$D_{Y_1}\not\sim_{\mathbb{Q}} -3K_{Y_1}$.
\end{proof}

Consequently,  we may assume that the set $\mathbb{CS}(Y,
\frac{1}{n}\mathcal{M}_{Y})$ consists of the point $\bar{Q}_2$
whose image to $W$ is the point $Q_2$. It implies that
$\mathbb{CS}(W, \frac{1}{n}\mathcal{M}_{W})=\{Q_1, Q_2\}$. We have
$\mathcal{M}_{V}\sim_{\mathbb{Q}} -nK_{V}$ by
Lemma~\ref{lemma:Kawamata}.

Let $H$ be the exceptional divisor of $\xi$. Then, it follows from
the local computations made during the proof of
Lemma~\ref{lemma:n-56-P7} that $B_{V}\sim_{\mathbb{Q}} 2S_{V}$.
The base locus of $\mathcal{B}_{V}$ is the union of the
irreducible curves $C_{V}$, $L_{V}$, $L^{\prime}_{V}$, and the
curve $L^{\prime\prime}$ such that $L^{\prime\prime}=S_{V}\cdot
H$. We have
$$S_{V}\cdot C_{V}=0, \ \ S_{V}\cdot L_{V}=-\frac{1}{3}, \ \ S_{V}\cdot
L^{\prime}_{V}=-\frac{1}{6}, \ \ S_{V}\cdot
L^{\prime\prime}=\frac{1}{2}.$$

Let $\bar{O}$ be the singular point of $V$ that is contained in
$H$. Then, $\omega(\bar{O})$ is the singular point of $Z$
contained in the exceptional divisor of $\nu$. It follows from
Lemma~\ref{lemma:Cheltsov-Kawamata} that either the set
$\mathbb{CS}(V, \frac{1}{n}\mathcal{M}_{V})$ contains the point
$\bar{O}$ or  the log pair $(V, \frac{1}{n}\mathcal{M}_{V})$ is
terminal.

Suppose that the set $\mathbb{CS}(V, \frac{1}{n}\mathcal{M}_{V})$
contains the point $\bar{O}$. Then, the set $\mathbb{CS}(Z,
\frac{1}{n}\mathcal{M}_{Z})$ contains the point $\omega(\bar{O})$.
The proof of Lemma~\ref{lemma:n-56-P7} shows that
$\mathcal{M}=\mathcal{B}$ if the set $\mathbb{CS}(Z,
\frac{1}{n}\mathcal{M}_{Z})$ contains the point
$\omega(\bar{O})=O_2$.

From now,  we  suppose that the singularities of the log pair $(V,
\frac{1}{n}\mathcal{M}_{V})$ are terminal. The singularities of
the log pair $(V, \epsilon B_{V})$ are log-terminal for some
rational number $\epsilon>\frac{1}{2}$ but the divisor
$K_{V}+\epsilon B_{V}$ has nonnegative intersection with all
curves on the variety $V$ except the curves $L_{V}$ and
$L^{\prime}_{V}$. Then, there is a composition of antiflips
$\chi:V\dasharrow V^{\prime}$ and the divisor $-K_{V^{\prime}}$ is
nef. Hence, the linear system $|-rK_{V^{\prime}}|$ is
base-point-free for $r\gg 0$ by the log abundance theorem
(\cite{KMM}).

It follows from the proof of Lemma~\ref{lemma:n-56-P7} that the
pull-backs of the rational functions $\frac{y}{x^{2}}$ and
$\frac{zy}{x^{5}}$ are contained in the linear systems $|2S_{V}|$
and $|5S_{V}|$, respectively. In particular, the complete linear
system $|-10K_{V}|$ induces a dominant rational map
$V\dasharrow\mathbb{P}(1,2,5)$, which implies that the linear
system $|-rK_{V^{\prime}}|$ induces a dominant morphism to a
surface. In fact, the linear system $|-rK_{V^{\prime}}|$ induces
the morphism $\upsilon$. The singularities of the log pair
$(V^{\prime}, \frac{1}{n}\mathcal{M}_{V^{\prime}})$ are terminal
because the singularities of the log pair $(V,
\frac{1}{n}\mathcal{M}_{V})$ are terminal and the rational map
$\chi$ is a log flop with respect to the log pair $(V,
\frac{1}{n}\mathcal{M}_{V})$. However,  the singularities of the
log pair $(V^{\prime}, \frac{1}{n}\mathcal{M}_{V^{\prime}})$
cannot be terminal by Theorem~\ref{theorem:Noether-Fano}. We have
obtained a contradiction.

Summing up, we have proved

\begin{proposition}
\label{proposition:n-56} The linear system $|-2K_{X}|$ is a unique
Halphen pencil on $X$.
\end{proposition}

\section{Cases $\gimel=63$, $77$, $83$, and $85$.}
\label{section:n-63-83}

Suppose that $\gimel\in\{63,
83\}$.\index{$\gimel=63$}\index{$\gimel=83$} Then, the threefold
$X \subset \mathbb{P}(1,a_1,a_2,a_3,a_4)$ always contains the
point $O=(0:0:0:1:0)$. It  is a singular point of $X$ that is a
quotient singularity of type $\frac{1}{a_3}(1,a_2, a_3-a_2)$.

We also have a commutative diagram as follows:
$$
\xymatrix{
&U\ar@{->}[d]_{\alpha}&&W\ar@{->}[ll]_{\beta}\ar@{->}[d]^{\eta}&\\%
&X\ar@{-->}[rr]_{\psi}&&\mathbb{P}(1,a_1,a_2),&}
$$
where \begin{itemize} \item $\alpha$ is the Kawamata blow up at
the point $O$ with weights $(1,a_2,a_3-a_2)$,

\item $\beta$ is the Kawamata blow up with weights
$(1,a_2,a_3-2a_2)$ at the point $P$ of $U$ that is a quotient
singularity of type $\frac{1}{a_3-a_2}(1,a_2,a_3-2a_2)$,

\item $\eta$ is an elliptic fibration.
\end{itemize}
We may assume that
$$
\mathbb{CS}\Big(X, \frac{1}{n}\mathcal{M}\Big)=\Big\{O\Big\}%
$$
due to Theorem~\ref{theorem:Noether-Fano},
Lemmas~\ref{lemma:smooth-points},
\ref{lemma:special-singular-points-with-positive-c},
\ref{lemma:special-singular-points-with-zero-c}, and
Corollary~\ref{corollary:Ryder-a1}.

The exceptional divisor $E$ of the birational morphism $\alpha$
contains two singular points $P$ and $Q$ that are quotient
singularity of types $\frac{1}{a_3-a_2}(1,a_2,a_3-2a_2)$ and
$\frac{1}{a_2}(1,1, a_3-2a_2)$, respectively. The base locus of
$|-a_1K_X|$ consists of the irreducible curve $C$ defined by
$x=y=0$. The base locus of $|-a_1K_U|$ consists of the proper
transform $C_U$ and the unique irreducible curve $L$ in
$|\mathcal{O}_{\mathbb{P}(1, a_2, a_3-a_2)}(1)|$ on the surface
$E$.

\begin{lemma}
\label{lemma:n-63-83-second-floor} If the set  $\mathbb{CS}(U,
\frac{1}{n}\mathcal{M}_U)$ contains the point $Q$, then
$\mathcal{M}=|-a_1K_{X}|$.
\end{lemma}

\begin{proof}
Let $\pi:Y\to U$ be the Kawamata blow up at the point $Q$ with
weights $(1,1, a_3-2a_2)$ and $G_Q$ be its exceptional divisor.
Then, the base locus of the pencil $|-a_1K_{Y}|$ consists of the
irreducible curves $C_Y$ and $L_Y$. Then, our situation is exactly
same as Lemma~\ref{lemma:n-29-second-floor}. Using the same proof,
we get $\mathcal{M}=|-a_1K_{X}|$.
\end{proof}

The exceptional divisor $F$ of the birational morphism $\beta$
contains two singular points $P_1$ and $P_2$  that are quotient
singularities of types
 $\frac{1}{a_2}(1,1,a_2-1)$ and $\frac{1}{a_3-2a_2}(1,1,a_3-2a_2-1)$, respectively.

\begin{lemma}
\label{lemma:n-63-83-third-floor} The set  $\mathbb{CS}(W
\frac{1}{n}\mathcal{M}_W)$ cannot contain the point $P_2$.
\end{lemma}
\begin{proof}
Let $\sigma_2:V_2\to W$ be the Kawamata blow up at the point
$P_2$.  The proper transform $\mathcal{D}$ of the linear system
$|-a_2K_X|$ consists of two irreducible curves $\bar{C}_{V_2}$ and
$L_{V_2}$. Applying the same method as in
Lemma~\ref{lemma:n-63-83-second-floor} to the linear system
$\mathcal{D}$, we obtain an absurd identity
$\mathcal{M}=|-a_2K_X|$.
\end{proof}

\begin{proposition}
 The linear system $|-a_1K_{X}|$ is the only Halphen
pencil on $X$.
\end{proposition}

\begin{proof}
Lemma~\ref{lemma:Cheltsov-Kawamata} implies that either
$\mathbb{CS}(U, \frac{1}{n}\mathcal{M}_U)=\{P\}$ or
$Q\in\mathbb{CS}(U, \frac{1}{n}\mathcal{M}_U)$. The latter case
implies $\mathcal{M}=|-a_1K_X|$ by
Lemma~\ref{lemma:n-63-83-second-floor}. Suppose that the set
$\mathbb{CS}(U, \frac{1}{n}\mathcal{M}_U)$ consists of the point
$P$. Then, the set $\mathbb{CS}(W, \frac{1}{n}\mathcal{M}_W)$ must
contain the point $P_1$ by Lemma~\ref{lemma:n-63-83-third-floor}.
Let $\sigma_1:V_1\to W$ be the Kawamata blow up at the point
$P_1$. Then, the base locus of the pencil $|-a_1K_{V_1}|$ consists
of the irreducible curves $C_{V_1}$ and $L_{V_1}$. Applying the
same method as in Lemma~\ref{lemma:n-63-83-second-floor}, we
obtain $\mathcal{M}=|-a_1K_{X}|$.
\end{proof}

We suppose that $\gimel=77$ or
$85$.\index{$\gimel=77$}\index{$\gimel=85$} The hypersurface
$X\subset\mathbb{P}(1, a_1, a_2, a_3, a_4)$ always contains the
point $O=(0:0:0:1:0)$ as a quotient singularity  of type
$\frac{1}{a_3}(1,a_1, a_3-a_1)$.

There is a commutative diagram
$$
\xymatrix{
&U\ar@{->}[d]_{\alpha}&&W\ar@{->}[ll]_{\beta}\ar@{->}[d]^{\eta}&\\%
&X\ar@{-->}[rr]_{\psi}&&\mathbb{P}(1,a_1,a_2),&}
$$
where \begin{itemize} \item $\psi$ is the natural projection,

 \item $\alpha$ is the Kawamata blow up at
the point $O$ with weights $(1,a_1,a_3-a_1)$,

\item $\beta$ is the Kawamata blow up with weights
$(1,a_1,a_3-2a_1)$ at the singular point of the variety $U$ that
is a quotient singularity of type
$\frac{1}{a_3-a_1}(1,a_1,a_3-2a_1)$,

\item $\eta$ is an elliptic fibration.
\end{itemize}

As in the previous case, we may assume that
$$
\mathbb{CS}\Big(X, \frac{1}{n}\mathcal{M}\Big)=\Big\{O\Big\}.%
$$
The exceptional divisor $E$ of the birational morphism $\alpha$
contains two singular points $P$ and $Q$ that are quotient
singularities of types $\frac{1}{a_3-a_1}(1, a_1, a_3-2a_1)$ and
$\frac{1}{a_1}(1, 1, a_1-1)$, respectively.

Unlike the previous case, we have the opposite statement for the
point $Q$ as follows:

\begin{lemma}
\label{lemma:n-77-85-second-floor} The set  $\mathbb{CS}(U,
\frac{1}{n}\mathcal{M}_U)$ cannot contain the point $Q$.
\end{lemma}
\begin{proof}
Suppose the set  $\mathbb{CS}(U, \frac{1}{n}\mathcal{M}_U)$
contains the point $Q$. Let $\pi:Y\to U$ be the Kawamata blow up
at the point $Q$. The base locus of the proper transform
$\mathcal{D}$ of the linear system $|-a_4K_{X}|$ does not contain
any curve. Therefore, a general surface $D$ in the linear system
$\mathcal{D}$ is nef. However, we can easily check that $D\cdot
M_1\cdot M_2<0$ for general surfaces $M_1$ and $M_2$ in
$\mathcal{M}_Y$. It is a contradiction.
\end{proof}

\begin{proposition}
\label{proposition:n-77-85} If $\gimel\in\{77, 85\}$, then the
linear system  $|-a_{1}K_{X}|$ is a unique Halphen pencil on $X$.
\end{proposition}

\begin{proof}
The proof is the same as the cases $\gimel=63$ and $83$. The only
difference is Lemma~\ref{lemma:n-63-83-second-floor}. We replace
it by Lemma~\ref{lemma:n-77-85-second-floor}.
\end{proof}

\section{Case $\gimel=65$, hypersurface of degree $27$ in
$\mathbb{P}(1,2,5,9,11)$.}\index{$\gimel=65$} \label{section:n-65}

The threefold $X$ is a general hypersurface  of degree $27$ in
$\mathbb{P}(1,2,5,9,11)$ with $-K_{X}^{3}=\frac{3}{110}$. The
singularities of $X$ consist of one singular point $O$ that is a
quotient singularity of type $\frac{1}{11}(1,2,9)$, one point of
type $\frac{1}{5}(1,4,1)$, and four  points  of type
$\frac{1}{2}(1,1,1)$.

There is a commutative diagram
$$
\xymatrix{
&U\ar@{->}[d]_{\alpha}&&W\ar@{->}[ll]_{\beta}&&Y\ar@{->}[ll]_{\gamma}\ar@{->}[d]^{\eta}&\\%
&X\ar@{-->}[rrrr]_{\psi}&&&&\mathbb{P}(1,2,5),&}
$$
where \begin{itemize}

\item $\psi$ is the natural projection,

\item $\alpha$ is the Kawamata blow up at the point $O$ with
weights $(1,2,9)$,

\item $\beta$ is the Kawamata blow up with weights $(1,2,7)$ at
the singular point of the variety $U$ that is a quotient
singularity of type $\frac{1}{9}(1,2,7)$ contained in the
exceptional divisor of the birational morphism $\alpha$,

\item $\gamma$ is the Kawamata blow up with weights $(1,2,5)$ at
the singular point of the variety $W$ that is a quotient
singularity of type $\frac{1}{7}(1,2,5)$ contained in the
exceptional divisor of the birational morphism $\beta$,

\item $\eta$ is an elliptic fibration.

\end{itemize}

If the set $\mathbb{CS}(X, \frac{1}{n}\mathcal{M})$ contains the
singular point of type $\frac{1}{5}(1,4,1)$, then
$\mathcal{M}=|-2K_X|$ by
Lemma~\ref{lemma:special-singular-points-with-zero-c}. Therefore,
we may assume that
 $\mathbb{CS}(X,
\frac{1}{n}\mathcal{M})=\{O\}$ by Lemma~\ref{lemma:smooth-points}
and Corollary~\ref{corollary:Ryder-a2}.

The hypersurface $X$ can be given by the equation
$$
w^{2}z+wf_{16}(x,y,z,t)+f_{27}(x,y,z,t)=0,
$$
where $f_{i}(x,y,z,t)$ is a quasihomogeneous polynomial of degree
$i$. Let $\mathcal{P}$ be the pencil of surfaces  cut on the
hypersurface $X$ by
$$
\lambda x^{5}+\mu z=0,
$$
where $(\lambda : \mu)\in\mathbb{P}^{1}$. Even though the linear
system $\mathcal{P}$ is not a Halphen pencil, it is helpful for
our proof.  Note that the base locus of the pencil $\mathcal{P}$
consists of the irreducible curve $\bar{C}$.

The exceptional divisor $E\cong\mathbb{P}(1,2,9)$ contains two
singular points $P$ and $Q$ of $U$ that are quotient singularities
of types $\frac{1}{2}(1,1,1)$ and $\frac{1}{9}(1,2,7)$,
respectively. Let $L$ be the unique curve contained in the linear
system $|\mathcal{O}_{\mathbb{P}(1,\,2,\,9)}(1)|$ on the surface
$E$.

 The set $\mathbb{CS}(U, \frac{1}{n}\mathcal{M}_U)$
is not empty by Theorem~\ref{theorem:Noether-Fano}. Hence, either
the set $\mathbb{CS}(U, \frac{1}{n}\mathcal{D})$ contains the
point $P$ or it consists of the point $Q$ by
Lemma~\ref{lemma:Cheltsov-Kawamata}.

\begin{lemma}
\label{lemma:n-65-P4} The set $\mathbb{CS}(U,
\frac{1}{n}\mathcal{M}_U)$ does not contain the point $P$.
\end{lemma}

\begin{proof}
Suppose that $P\in\mathbb{CS}(U, \frac{1}{n}\mathcal{M}_U)$. Let
$\pi_P:U_P\to U$ be the Kawamata blow up at the point $P$ with
weights $(1,1,1)$ and $G_P$ be its exceptional divisor. Then,
$\mathcal{M}_{U_P}\sim_{\mathbb{Q}}-nK_{U_P}$ holds by
Lemma~\ref{lemma:Kawamata}.

The base locus of the pencil $\mathcal{P}_{U_P}$ consists of the
irreducible curves $\bar{C}_{U_P}$, $L_{U_P}$, and a line
$\bar{L}$ on $G_P\cong\mathbb{P}^2$. For a general surface
$D_{U_P}$ of the pencil $\mathcal{P}_{U_P}$,
$$
S_{U_P}\cdot D_{U_P}=\bar{C}_{U_P}+L_{U_P},\ E_{U_P}\cdot
D_{U_P}=5L_{U_P},\
G_P\cdot D_{U_P}=\bar{L}. %
$$
The surface $D_{U_P}$ is normal. On the other hand, we have
$$
\left\{\aligned
&E_{U_P}\sim_{\mathbb{Q}}\pi_P^{*}(E)-\frac{1}{2}G_P,\\
&D_{U_P}\sim_{\mathbb{Q}}(\alpha\circ\pi_P)^{*}(-5K_{X})-\frac{5}{11}\pi_P^{*}(E)
-\frac{1}{2}G_P,\\
&S_{U_P}\sim_{\mathbb{Q}}(\alpha\circ\pi_P)^{*}(-K_{X})-\frac{1}{11}\pi_P^{*}(E)
-\frac{1}{2}G_P,\\
\endaligned
\right.
$$
which implies
$$
L_{U_P}\cdot L_{U_P}=-\frac{73}{450},\ \bar{C}_{U_P}\cdot L_{U_P}=\frac{4}{45},\ \bar{C}_{U_P}\cdot \bar{C}_{U_P}=-\frac{497}{550}%
$$
 on the surface $D_{U_P}$. Therefore, the intersection form of
the curves $\bar{C}_{U_P}$ and $L_{U_P}$ on the surface $D_{U_P}$
is negative-definite. On the other hand, we have
$$
\mathcal{M}_{U_P}\Big\vert_{D_{U_P}}\equiv-nK_{_{U_P}}\Big\vert_{D_{U_P}}\equiv
nS_{U_P}\Big\vert_{D_{U_P}}\equiv n\bar{C}_{U_P}+nL_{U_P},
$$
which implies $\mathcal{M}_{U_P}=\mathcal{P}_{U_P}$ by
Theorem~\ref{theorem:main-tool}. Hence, $n=5$, but
$D_{U_P}\not\sim_{\mathbb{Q}}-5K_{U_P}$, which is a contradiction.
\end{proof}

Hence, the set $\mathbb{CS}(U, \frac{1}{n}\mathcal{M}_U)$ consists
of the point $Q$. The exceptional divisor
$F\cong\mathbb{P}(1,2,7)$ of $\beta$ contains two singular points
$Q_1$ and $Q_2$ of  $W$ that are quotient singularities of types
$\frac{1}{2}(1,1,1)$ and $\frac{1}{7}(1,2,5)$, respectively. Let
$\tilde{L}$ be the unique curve contained in the linear system
$|\mathcal{O}_{\mathbb{P}(1,\,2,\,7)}(1)|$ on the surface $F$.

It follows from Theorem~\ref{theorem:Noether-Fano} that the set
$\mathbb{CS}(W, \frac{1}{n}\mathcal{M}_W)$ is not empty because
the divisor $-K_{W}$ is nef and big.

\begin{lemma}
\label{lemma:n-65-P6} The set $\mathbb{CS}(W,
\frac{1}{n}\mathcal{M}_W)$ does not contain the point $Q_1$.
\end{lemma}

\begin{proof}
Suppose that $Q_1\in\mathbb{CS}(W, \frac{1}{n}\mathcal{M}_W)$. Let
$\pi_1:W_1\to W$ be the Kawamata blow up at the point $Q_1$ with
weights $(1,1,1)$ and  $G_1$ be its exceptional divisor.

The base locus of the pencil $\mathcal{P}_{W_1}$ consists of the
irreducible curves $\bar{C}_{W_1}$, $L_{W_1}$, $\tilde{L}_{W_1}$,
$\Delta_{1}$, $\Delta_{2}$, and $\Delta$, where  the curves
$\Delta_{1}$ and $\Delta_{2}$ are the lines on
$G_1\cong\mathbb{P}^2$  cut out by the divisors $E_{W_1}$ and
$F_{W_1}$, respectively, and the curve $\Delta$ is a line on $G_1$
different from the lines $\Delta_{1}$ and $\Delta_{2}$.

For a general surface $D_{W_1}$ in  the pencil
$\mathcal{P}_{W_1}$, we have
$$
S_{W_1}\cdot D_{W_1}=C_{W_1}+L_{W_1}+\tilde{L}_{W_1},\
 E_{W_1}\cdot D_{W_1}=5L_{W_1}+\Delta_{1},\
 F_{W_1}\cdot D_{W_1}=5\tilde{L}_{W_1}+\Delta_{2}. %
$$
The surface $D_{W_1}$ is normal and it is smooth in a neighborhood
of $G_1$. In particular, it follows from the local computations
and the Adjunction formula that the equalities
\begin{equation*}
\Delta_{1}\cdot\Delta_{2}=\Delta_{1}\cdot
\tilde{L}_{W_1}=\Delta_{2}\cdot L_{W_1}=1, \ \Delta_{1}\cdot
C_{W_1}=\Delta_{2}\cdot C_{W_1}=0,\
 \Delta_{1}^{2}=\Delta_{2}^{2}=-4 %
\end{equation*}
hold on the surface $D_{W_1}$. However, we have
\begin{equation*}
\left\{\aligned
&F_{W_1}\sim_{\mathbb{Q}}\pi_1^{*}(F)-\frac{1}{2}G,\\
&E_{W_1}\sim_{\mathbb{Q}}(\beta\circ\pi_1)^{*}(E)-\frac{7}{9}\pi_1^{*}(F)-\frac{1}{2}G_1,\\
&D_{W_1}\sim_{\mathbb{Q}}(\alpha\circ\beta\circ\pi_1)^{*}(-5K_{X})-\frac{5}{11}(\beta\circ\pi_1)^{*}(E)-\frac{5}{9}\pi_1^{*}(F)-\frac{3}{2}G_1,\\
&S_{W_1}\sim_{\mathbb{Q}}(\alpha\circ\beta\circ\pi_1)^{*}(-K_{X})-\frac{1}{11}(\beta\circ\pi_1)^{*}(E)-\frac{1}{9}\pi_1^{*}(F)-\frac{1}{2}G_1.\\
\endaligned
\right.
\end{equation*}
We can then obtain
$$
C_{W_1}\cdot C_{W_1}=L_{W_1}\cdot L_{W_1}=-\frac{1}{2},\
\tilde{L}_{W_1}\cdot \tilde{L}_{W_1}=-\frac{3}{7},\
C_{W_1}\cdot L_{W_1}=C_{W_1}\cdot \tilde{L}_{W_1}=L_{W_1}\cdot \tilde{L}_{W_1}=0%
$$
on the surface $D_{W_1}$. Therefore, the intersection form of the
curves $C_{W_1}$, $L_{W_1}$, and $\tilde{L}_{W_1}$ on the surface
$D_{W_1}$ is negative-definite. On the other hand, we have
$$
\mathcal{M}_{W_1}\Big\vert_{D_{W_1}}\equiv-nK_{W_1}\Big\vert_{D_{W_1}}\equiv
nS_{W_1}\Big\vert_{D_{W_1}}\equiv
nC_{W_1}+nL_{W_1}+n\tilde{L}_{W_1},
$$
which implies $\mathcal{M}_{W_1}=\mathcal{P}_{W_1}$ by
Theorem~\ref{theorem:main-tool}. Hence, we have $n=5$, but $
D_{W_1}\not\equiv-5K_{W_1}, $ which is a contradiction.
\end{proof}

Thus, the set $\mathbb{CS}(W, \frac{1}{n}\mathcal{M}_W)$ consists
of the point $Q_2$ by Lemma~\ref{lemma:Cheltsov-Kawamata}.

The exceptional divisor of the birational morphism $\gamma$
contains two singular points $O_1$ and $O_2$ of $Y$ that are
quotient singularities of types $\frac{1}{2}(1,1,1)$ and
$\frac{1}{5}(1,2,3)$, respectively. Then, the set $\mathbb{CS}(Y,
\frac{1}{n}\mathcal{M}_Y)$ must contain either the point $O_1$ or
the point $O_2$ by Theorem~\ref{theorem:Noether-Fano} and
Lemma~\ref{lemma:Cheltsov-Kawamata}.

\begin{lemma}
\label{lemma:n-65-P8} The set $\mathbb{CS}(Y,
\frac{1}{n}\mathcal{M}_Y)$ does not contain the point $O_1$.
\end{lemma}

\begin{proof}
Suppose that $O_1\in\mathbb{CS}(Y, \frac{1}{n}\mathcal{M}_Y)$. Let
$\sigma_1:V_1\to Y$ be the Kawamata blow up at the point $O_1$
 with weights $(1,1,1)$ and  $H_1$ be the exceptional divisor of
the birational morphism $\sigma_1$.

A general surface $D_{V_1}$ in $\mathcal{P}_{V_1}$ is normal and
$$
\mathcal{M}_{V_1}\Big\vert_{D_{V_1}}\equiv-nK_{V_1}\Big\vert_{D_{V_1}}\equiv
nS_{V_1}\Big\vert_{D_{V_1}}.
$$
The intersection $D_{V_1}\cap H_1$ is a line on
$G_1\cong\mathbb{P}^{2}$ that is different from the line
$S_{V_1}\cap G_1$. However, the curve $D_{V_1}\cap H_1$ is a fiber
of the elliptic fibration $\eta\circ\sigma_1\vert_{D_{V_1}}$ over
the point $\eta\circ\sigma_1(S_{V_1}\cap H_1)$. The support of the
cycle $S_{V_1}\cdot D_{V_1}$ contains all components of the fiber
of the elliptic fibration $\eta\circ\sigma_1\vert_{D_{V_1}}$ over
the point $\eta\circ\sigma_1(S_{V_1}\cap H_1)$ that are different
from the curve $D_{V_1}\cap H_1$. Hence, the intersection form of
the components of the cycle $S_{V_1}\cdot D_{V_1}$ are
negative-definite on the surface $D_{V_1}$, which implies that
$\mathcal{M}_{V_1}=\mathcal{P}_{V_1}$ by
Theorem~\ref{theorem:main-tool}. Hence, we have $n=5$, but it
follows from explicit calculations that
$$
D_{V_1}\sim_{\mathbb{Q}}\sigma_1^{*}(-5K_{Y})-\frac{1}{2}H_1\not\sim_{\mathbb{Q}} -5K_{V_1},%
$$
which is a contradiction.
\end{proof}

\begin{proposition}
\label{proposition:n-65} The linear system $|-2K_{X}|$ is a unique
Halphen pencil on $X$.
\end{proposition}

\begin{proof}
By what we have proved so far, we may assume that the set
$\mathbb{CS}(Y, \frac{1}{n}\mathcal{M}_Y)$ contains the point
$O_2$. Let $\sigma_2:V_2\to Y$ be the Kawamata blow up at the
point $O_1$ with weights $(1,2,3)$. Then, $|-2K_{V_2}|$ is the
proper transform of the pencil $|-2K_{X}|$. Its base locus
consists of the irreducible curve $C_{V_2}$. Because
$\mathcal{M}_{V_1}\sim_{\mathbb{Q}} -nK_{V_1}$ and $-K_{V_1}\cdot
C_{V_2}<0$, Theorem~\ref{theorem:main-tool} implies that
$\mathcal{M}=|-2K_{X}|$.
\end{proof}

\section{Cases $\gimel=72$,  $89$, $90$, $92$, and $94$.}%
\index{$\gimel=72$}\index{$\gimel=89$}\index{$\gimel=90$}\index{$\gimel=92$}\index{$\gimel=94$}\label{section:n-89-90-92-94}

Suppose that $\gimel\in\{72, 89, 90,92, 94\}$. Then, the threefold
$X \subset \mathbb{P}(1,a_1,a_2,a_3,a_4)$ always contains a
quotient singularity  $O$ of type
$\frac{1}{a_{1}+a_2}(1,a_1,a_{2})$.

We also have a commutative diagram as follows:

$$
\xymatrix{
&&&U\ar@{->}[lld]_{\pi}\ar@{->}[rrd]^{\eta}&&&\\%
&X\ar@{-->}[rrrr]_{\psi}&&&&\mathbb{P}(1,a_1,a_{2})&}
$$
where \begin{itemize} \item $\psi$ is a natural projection, \item
$\pi$ is the Kawamata blow up at the point $O$ with weights $(1,
a_1, a_2)$, \item $\eta$ is an elliptic fibration.\end{itemize}

 We may
assume that the set $\mathbb{CS}(X, \frac{1}{n}\mathcal{M})$
consists of the singular point $O$ due to
Lemmas~\ref{lemma:smooth-points},
\ref{lemma:special-singular-points-with-positive-c},
\ref{lemma:special-singular-points-with-zero-c} and
Corollary~\ref{corollary:Ryder-a1}.

The exceptional divisor $E$ of the birational morphism $\pi$
contains two singular points $P$ and $Q$ of types
$\frac{1}{a_1}(1, a ,a_1-a)$ and $\frac{1}{a_2}(1, a_1, a_2-a_1)$,
respectively, where $a=|2a_1-a_2|$. Then, the set $\mathbb{CS}(U,
\frac{1}{n} \mathcal{M}_U)$ contains either the singular point $P$
or the singular point $Q$ by Lemma~\ref{lemma:Cheltsov-Kawamata}.

\begin{lemma}
\label{lemma:n-89-90-92-94} The set $\mathbb{CS}(U, \frac{1}{n}
\mathcal{M}_U)$ does not contain the point $P$.
\end{lemma}

\begin{proof}
Suppose that the set $\mathbb{CS}(U, \frac{1}{n} \mathcal{M}_U)$
contains the point $P$. Let $\alpha:U_P\to U$ be the Kawamata blow
up at the point $P$ with weights $(1, a, a_1-a)$ and $G$ be the
ex\-cep\-ti\-o\-nal divisor of the birational morphism $\alpha$.
Then,
$$
\mathcal{M}_{U_P}\sim_{\mathbb{Q}}\big(\pi\circ\alpha\big)^{*}\Big(-nK_{X}\Big)
-\frac{n}{a_1+a_2}\alpha^{*}\big( E\big)-\frac{n}{a_1}G\sim_{\mathbb{Q}}-nK_{U_P}%
$$
by Lemma~\ref{lemma:Kawamata}.

Let $\mathcal{D}$ be the proper transform  of the linear system
$|-a_2K_{X}|$ by the birational morphism $\pi\circ\alpha$. Then,
$S_{U_P}\sim_{\mathbb{Q}}-K_{U_P}$ and
$$
\mathcal{D}\sim_{\mathbb{Q}}\big(\pi\circ\alpha\big)^{*}\Big(-a_2K_{X}\Big)
-\frac{a_2}{a_1+a_2}\alpha^{*}\big(E\big)-\frac{1}{a_1}G,%
$$
but the base locus of the linear system $\mathcal{D}$ consists of
the irreducible curve $\bar{C}_{U_P}$ whose image to $X$ is the
base curve of the linear system $|-a_2K_{X}|$. Let $D$ be a
general surface in $\mathcal{D}$ and $M$ be a general surface in
$\mathcal{M}_{U_P}$.  The surface $D$ is normal. We have
$\bar{C}_{U_P}^2<0$ on the normal surface $D$ and
$M\vert_{D}\equiv n\bar{C}_{U_P}$, which implies that
$\mathcal{M}_{U_P}=\mathcal{D}$ by
Theorem~\ref{theorem:main-tool}. However, the linear system
$\mathcal{D}$ is not a pencil.
\end{proof}

\begin{proposition}
\label{proposition:n-89-90-92-94} If $\gimel=72,  89, 90, 92, 94$,
then $\mathcal{M}=|-a_{1}K_{X}|$.
\end{proposition}
\begin{proof}
By Lemma~\ref{lemma:n-89-90-92-94}, we may assume that the set
$\mathbb{CS}(U, \frac{1}{n}\mathcal{M}_U)$ contains the point $Q$.

Let $\beta:W\to U$ be the  Kawamata blow up of the point $Q$ with
weights $(1, a_1, a_2-a_1)$. The linear system $|-a_1K_W|$ is the
proper transform of the pencil $|-a_1K_X|$ and the base locus of
the pencil $|-a_1K_W|$ consists of the irreducible curve $C_W$
whose image to $X$ is the base curve of $|-a_1K_X|$. Then, the
inequality $-K_W\cdot C_W<0$ and the equivalence
$\mathcal{M}_W\sim_{\mathbb{Q}}-nK_W$ imply that the pencil
$\mathcal{M}_W$  coincides with the pencil $|-a_1K_W|$ by
Theorem~\ref{theorem:main-tool}.
\end{proof}

\section{Cases $\gimel=75$ and $87$.}
\label{section:n-75-87}

In the case of $\gimel=75$,\index{$\gimel=75$} the threefold $X$
is a general hypersurface of degree $30$ in
$\mathbb{P}(1,4,5,6,15)$ with $-K_X^3=\frac{1}{60}$. Its
singularities consist of one quotient singular point of type
$\frac{1}{4}(1,1,3)$, one quotient singular point of type
$\frac{1}{3}(1,1,2)$, two quotient singular points of type
$\frac{1}{2}(1,1,1)$, and two quotient singular points of type
$\frac{1}{5}(1,4,1)$.

In the case of  $\gimel=87$,\index{$\gimel=87$} the threefold $X$
is a general hypersurface of degree $40$ in
$\mathbb{P}(1,5,7,8,20)$ with $-K_X^3=\frac{1}{140}$. It has one
quotient singular point of type $\frac{1}{4}(1,1,3)$, two quotient
singular points of type $\frac{1}{5}(1,2,3)$, and one quotient
singular point of type $\frac{1}{7}(1,1,6)$.

In both cases, the threefold $X$ cannot be birationally
transformed to an elliptic fibration (\cite{ChPa05}). However, it
can be rationally fibred by K3 surfaces.

\begin{proposition}
\label{proposition:n-75-87} If  $\gimel\in\{75, 87\}$, then
$|-a_{1}K_{X}|$ is a unique Halphen pencil on $X$.
\end{proposition}

\begin{proof}
Theorem~\ref{theorem:Noether-Fano},
Lemmas~\ref{lemma:smooth-points},
\ref{lemma:special-singular-points-with-positive-c},
\ref{lemma:special-singular-points-with-zero-c}, and
Corollary~\ref{corollary:Ryder-a1} immediately imply the result.
\end{proof}


\newpage
\part{Fano  threefold hypersurfaces with two Halphen
pencils.}\label{section:two-Halphens}

\section{Case $\gimel=45$, hypersurface of degree $20$ in
$\mathbb{P}(1,3,4,5,8)$.}\index{$\gimel=45$} \label{section:n-45}

Let $X$ be the hypersurface given by a general quasihomogeneous
equation of degree $18$ in $\mathbb{P}(1,3,4,5,8)$ with
$-K_X^3=\frac{1}{24}$. Then, the singularities of $X$ consist of
one singular point $P$ that is a quotient singularity of type
$\frac{1}{8}(1,3,5)$, one singular point of type
$\frac{1}{3}(1,1,2)$,  and two points of type
$\frac{1}{4}(1,3,1)$.

There is a commutative diagram
$$
\xymatrix{
&U\ar@{->}[d]_{\alpha}&&W\ar@{->}[ll]_{\beta}\ar@{->}[d]^{\eta}&\\%
&X\ar@{-->}[rr]_{\psi}&&\mathbb{P}(1,3,4),&}
$$
where \begin{itemize} \item $\psi$ is the natural projection,

\item $\alpha$ is a blow up at the singular point $P$ with weights
$(1,3,5)$,

\item $\beta$ is the blow up with weights $(1,3,2)$ of the
singular point of the variety $U$ that is a quotient singularity
of type $\frac{1}{5}(1,3,2)$,

\item $\eta$ is an elliptic fibration.
\end{itemize}

Also, by the generality of the hypersurface, we may assume that
the hypersurface $X$ is defined by the equation
$$
w^{2}z+wf_{12}(x,y,z,t)+f_{20}(x,y,z,t)=0,
$$
 where $f_{i}(x,y,z,t)$ is a quasihomogeneous polynomial of degree
$i$. Proposition~\ref{proposition:Halphen-K3} implies that the
linear system $|-3K_X|$ is a Halphen pencil. Let $\mathcal{P}$ be
the pencil on $X$ given by the equations
$$
\lambda x^4+\mu z=0,
$$
where $(\lambda, \mu)\in \mathbb{P}^1$. We  see that the linear
system  $\mathcal{P}$ is another Halphen pencil on $X$.
\begin{K3proposition}\label{K3proposition:n-45}
A general member of the pencil $\mathcal{P}$ is birational to a
smooth K3 surface.
\end{K3proposition}
\begin{proof}
It is a compactification of a double cover of $\mathbb{C}^2$
ramified along a sextic curve. It cannot be a rational surface by
Theorem~\ref{theorem:CPR}. Therefore, it is birational to a smooth
K3 surface.
\end{proof}

If the set $\mathbb{CS}(X,\frac{1}{n}\mathcal{M})$ contains one of
the singular points of type $\frac{1}{4}(1,3,1)$, then
$\mathcal{M}=|-3K_X|$ by
Lemma~\ref{lemma:special-singular-points-with-zero-c}. Therefore,
we may assume that
$$\mathbb{CS}\Big(X,\frac{1}{n}\mathcal{M}\Big)=\Big\{P\Big\}$$
by Lemma~\ref{lemma:smooth-points} and
Corollary~\ref{corollary:Ryder-a1}.

The exceptional divisor $E\cong\mathbb{P}(1,3,5)$ of the
birational morphism $\alpha$ contains two singular points $P_1$
and $P_2$ of types $\frac{1}{3}(1,1,2)$ and $\frac{1}{5}(1,3,2)$,
respectively. For the convenience, let $L$ be the unique curve
contained in the linear system
$|\mathcal{O}_{\mathbb{P}(1,3,5)}(1)|$ on the surface  $E$.

\begin{lemma}
The set $\mathbb{CS}(U,\frac{1}{n}\mathcal{M}_U)$ cannot contain
the point $P_1$.
\end{lemma}
\begin{proof}
Suppose so. Then, we consider the Kawamata blow up $\alpha_1:V\to
U$ at the point $P_1$ with weights $(1,1,2)$. Let
$\mathcal{D}_{V}$ be the proper transform, by the birational
morphism $\alpha\circ\alpha_1$, of the linear system $\mathcal{D}$
on $X$ defined by the equations
$\lambda_0x^5+\lambda_1x^2y+\lambda_2 t=0$, where $(\lambda_0:
\lambda_1: \lambda_2)\in \mathbb{P}^2$. Then, its base locus
consists of the irreducible curve $\tilde{C}_V$ whose image to $X$
is the base curve of the linear system $\mathcal{D}$. We easily
see that $S_V\cdot D_V= \tilde{C}_V$ and
\[D_V\sim_{\mathbb{Q}}
\big(\alpha\circ\alpha_1\big)^{*}\Big(-5K_{X}\Big)-\frac{5}{8}\alpha_1^{*}\big(E\big)-\frac{1}{3}F_1,\]
\[S_V\sim_{\mathbb{Q}}
\big(\alpha\circ\alpha_1\big)^{*}\Big(-K_{X}\Big)-\frac{1}{8}\alpha_1^{*}\big(E\big)-\frac{1}{3}F_1,\]
where $D_V$ is a general surface in $\mathcal{D}_V$ and $F_1$ is
the exceptional divisor of $\alpha_1$. Since the curve
$\tilde{C}_V$ has negative self-intersection on the normal surface
$D_V$, Theorem~\ref{theorem:main-tool} implies
$\mathcal{M}=\mathcal{D}$. But this is absurd because
$\mathcal{D}$ is not a pencil.
\end{proof}

Therefore, we may assume that
$\mathbb{CS}(U,\frac{1}{n}\mathcal{M}_W)=\{P_2\}$. The exceptional
divisor $F$ of the birational morphism $\beta$ contains two
singular points $Q_1$ and $Q_2$ that are quotient singularities of
types $\frac{1}{3}(1,1,2)$ and $\frac{1}{2}(1,1,1)$, respectively.

\begin{lemma}
If the set $\mathbb{CS}(W,\frac{1}{n}\mathcal{M}_W)$  contains the
point $Q_1$, then $\mathcal{M}=\mathcal{P}$.
\end{lemma}
\begin{proof} Suppose that the set $\mathbb{CS}(W,\frac{1}{n}\mathcal{M}_W)$  contains the
point $Q_1$.  Let $\beta_1:W_1\to W$ be the Kawamata blow up at
the point $Q_1$ with weights $(1,1,2)$.  The base locus of
$\mathcal{P}_{W_1}$ consists of the irreducible curve $C_{W_1}$
and  the irreducible curve $L_{W_1}$. A general surface $D_{W_1}$
in $\mathcal{P}_{W_1}$ is normal and the intersection form of the
curves $C_{W_1}$ and $L_{W_1}$ is negative-definite on the surface
$D_{W_1}$. Because $\mathcal{M}_{W_1}|_{D_{W_1}}\equiv
nC_{W_1}+nL_{W_1}$, we obtain $\mathcal{M}=\mathcal{P}$ from
Theorem~\ref{theorem:main-tool}.
\end{proof}

\begin{lemma}
The set $\mathbb{CS}(W,\frac{1}{n}\mathcal{M}_W)$  cannot contain
the point $Q_2$.
\end{lemma}
\begin{proof}
Suppose that the set $\mathbb{CS}(W,\frac{1}{n}\mathcal{M}_W)$
contains the point $Q_2$. Let $\beta_2:W_2\to W$ be the Kawamata
blow up at the point $Q_2$. We then consider the proper transform
$\mathcal{L}_{W_2}$ of the pencil $|-3K_X|$ by the birational
morphism $\alpha\circ\beta\circ\beta_2$. For a general surface
$D_{W_2}$ in $\mathcal{L}_{W_2}$, we have
\[D_{W_2}\sim_{\mathbb{Q}}
\big(\alpha\circ\beta\circ\beta_2\big)^{*}\Big(-3K_{X}\Big)-\frac{3}{8}\big(\beta\circ\beta_2\big)^{*}\big(E\big)
-\frac{3}{5}\beta_2^{*}\big(F\big)-\frac{1}{2}G,\] where $G$ is
the exceptional divisor of $\beta_2$. Also,  we have
\[S_{W_2}\sim_{\mathbb{Q}}
\big(\alpha\circ\beta\circ\beta_2\big)^{*}\Big(-K_{X}\Big)-\frac{1}{8}\big(\beta\circ\beta_2\big)^{*}\big(E\big)
-\frac{1}{5}\beta_2^{*}\big(F\big)-\frac{1}{2}G.\] Since
$\mathcal{M}_{W_2}|_{D_{W_2}}\equiv n C_{W_2}$, the surface
$D_{W_2}$ is normal, and the self-intersection number $C_{W_2}^2$
on the normal surface $D_{W_2}$ is $-\frac{1}{2}$, we obtain the
identity  $\mathcal{M}_{W_2}=\mathcal{L}_{W_2}$ from
Theorem~\ref{theorem:main-tool}. However,
$D_{W_2}\not\sim_{\mathbb{Q}} -3K_{W_2}$. It is a contradiction.
\end{proof}

\begin{proposition}
The linear systems  $|-3K_X|$ or $\mathcal{P}$ are the only
Halphen pencils on $X$.
\end{proposition}
\begin{proof}
We apply Lemma~\ref{lemma:special-singular-points-with-zero-c} to
the singular points on $X$ of type $\frac{1}{4}(1,3,1)$ and the
lemmas above to the singular point $P$.
\end{proof}

\section{Case $\gimel=48$, hypersurface of degree $21$ in
$\mathbb{P}(1,2,3,7,9)$.}\index{$\gimel=48$} \label{section:n-48}

In the case of $\gimel=48$, the hypersurface $X$ is defined by a
general quasihomogeneous equation of degree $21$ in
$\mathbb{P}(1,2,3,7,9)$ with $-K_X^3=\frac{1}{18}$. It has a
quotient singularity of type $\frac{1}{9}(1,2,7)$ at the point
$O=(0:0:0:0:1)$. It also has one quotient singular point of type
$\frac{1}{2}(1,1,1)$ and two quotient singular points of type
$\frac{1}{3}(1,2,1)$.

We  have an elliptic fibration as follows:
$$
\xymatrix{
&Y\ar@{->}[d]_{\alpha}&&U\ar@{->}[ll]_{\beta}&&V\ar@{->}[ll]_{\gamma}\ar@{->}[d]^{\eta}&\\%
&X\ar@{-->}[rrrr]_{\psi}&&&&\mathbb{P}(1,2,3),&}
$$
where \begin{itemize}

\item $\psi$ is the natural projection,

 \item $\alpha$ is the
Kawamata blow up at the point $O$ with weights $(1,2,7)$,

\item $\beta$ is the Kawamata blow up  with weights $(1,2,5)$ at
the singular point of $Y$ that is a quotient singularity of type
$\frac{1}{7}(1,2,5)$,

\item $\gamma$ is the Kawamata blow up  with weights $(1,2,3)$ at
the singular point of the variety $U$ that is a quotient
singularity of type $\frac{1}{5}(1,2,3)$,

\item $\eta$ is an elliptic fibration.
\end{itemize}

Proposition~\ref{proposition:Halphen-K3} implies that the linear
system $|-2K_X|$ is a Halphen pencil. However, we have another
Halphen pencil. The hypersurface $X$ can be given by the equation
\[w^2z+wf_{12}(x,y,z,t)+f_{21}(x,y,z,t)=0,\]
where $f_i$ is a quasihomogeneous polynomial of degree $i$. Let
$\mathcal{P}$ be the pencil on $X$ given by the equations
\[\lambda x^3+\mu z=0,\]
where $(\lambda: \mu)\in \mathbb{P}^1$. We will see that the
linear system $\mathcal{P}$ is a Halphen pencil on $X$
(K3-Proposition~\ref{K3proposition:n-48}).

If the set $\mathbb{CS}(X,\frac{1}{n}\mathcal{M})$ contains a
singular point of type $\frac{1}{3}(1,2,1)$ on $X$, then the
identity $\mathcal{M}=|-2K_X|$ follows from
Lemma~\ref{lemma:special-singular-points-with-zero-c}. Therefore,
we may assume that
$$\mathbb{CS}\Big(X,\frac{1}{n}\mathcal{M}\Big)=\Big\{O\Big\}$$
due to Theorem~\ref{theorem:Noether-Fano},
Lemmas~\ref{lemma:smooth-points},
\ref{lemma:special-singular-points-with-positive-c},
 and Corollary~\ref{corollary:Ryder-a1}.

The exceptional divisor $E\cong\mathbb{P}(1,2,7)$ of the
birational morphism $\alpha$ contains two quotient singular points
$P_1$ and $Q_1$ of types $\frac{1}{7}(1,2,5)$ and
$\frac{1}{2}(1,1,1)$, respectively. For the convenience, let $L$
be the unique curve in the linear system
$|\mathcal{O}_{\mathbb{P}(1,2,7)}(1)|$ on the surface $E$.

\begin{lemma}
The set $\mathbb{CS}(Y,\frac{1}{n}\mathcal{M}_Y)$ cannot contain
the point $Q_1$.
\end{lemma}
\begin{proof}

Suppose that the set $\mathbb{CS}(Y,\frac{1}{n}\mathcal{M}_Y)$
contains the point $Q_1$. Let $\alpha_1:W_1\to Y$ be the Kawamata
blow up at the point $Q_1$ with weights $(1,1,1)$.  By
Lemma~\ref{lemma:Kawamata}, we have
$\mathcal{M}_{W_1}\sim_{\mathbb{Q}}-nK_{W_1}.$

We consider the linear system $|-7K_X|$ on $X$ that has no base
curves. Also, the proper transform $\mathcal{D}_{W_1}$ of the
linear system $|-7K_X|$ by the birational $\alpha\circ\alpha_1$
has no base curve. Let $D_1$ be a general member in
$\mathcal{D}_{W_1}$. The divisor $D_1$ is nef. Meanwhile, we have
\[D_1\sim_{\mathbb{Q}}
(\alpha\circ\alpha_1)^*(-7K_X)-\frac{7}{9}\alpha_1^*(E)-\frac{1}{2}E_1,\]
where  $E_1$ is the exceptional divisor of $\alpha_1$. Therefore,
$D_1\cdot M_1\cdot M_2= -\frac{n^2}{6}<0$, where $M_1$ and $M_2$
are general members in $\mathcal{M}_{W_1}$. It is a contradiction.
\end{proof}

Therefore, we may assume that
$\mathbb{CS}(Y,\frac{1}{n}\mathcal{M}_Y)=\{P_1\}$.  The
exceptional divisor $F\cong\mathbb{P}(1,2,5)$ of the birational
morphism $\beta$ contains two singular points $P_2$ and $Q_2$ that
are quotient singularities of types $\frac{1}{5}(1,2,3)$ and
$\frac{1}{2}(1,1,1)$, respectively. The set
$\mathbb{CS}(U,\frac{1}{n}\mathcal{M}_U)$ contains either the
point $P_2$ or the point $Q_2$. For the convenience, we denote the
unique curve in the linear system
$|\mathcal{O}_{\mathbb{P}(1,2,5)}(1)|$ on the surface $F$ by
$\bar{L}$.

\begin{lemma}
\label{lemma:n-48-first} If the set
$\mathbb{CS}(U,\frac{1}{n}\mathcal{M}_U)$  contains the point
$Q_2$, then $\mathcal{M}=\mathcal{P}$.
\end{lemma}

\begin{proof}
Suppose that the set $\mathbb{CS}(U,\frac{1}{n}\mathcal{M}_U)$
contains the point $Q_2$.  Let $\beta_2 : W_2\to U$ be the
Kawamata blow up at the point $Q_2$ with weights $(1,1,1)$. Let
$F_2$ be the exceptional divisor of $\beta_2$. Then,
\[\mathcal{M}_{W_2}\sim_{\mathbb{Q}}-nK_{W_2}\sim_{\mathbb{Q}}(\alpha\circ\beta\circ\beta_2)^*(-nK_X)-\frac{n}{9}(\beta\circ\beta_2)^*(E)
-\frac{n}{7}\beta_2^*(F)-\frac{n}{2}F_2.\] The base locus of the
pencil $\mathcal{P}_{W_2}$ consists of three irreducible curves
$\bar{C}_{W_2}$, $L_{W_2}$, and $\bar{L}_{W_2}$.

Let $D_{W_2}$ be a general surface in the pencil
$\mathcal{P}_{W_2}$. Then, the surface $D_2$ is normal and we have
$$
\left\{\aligned
&D_{W_2}\sim_{\mathbb{Q}}(\alpha\circ\beta\circ\beta_2)^*(-3K_X)-\frac{3}{9}(\beta\circ\beta_2)^*(E)
-\frac{3}{7}\beta_2^*(F)-\frac{3}{2}F_2,\\
&S_{W_2}\sim_{\mathbb{Q}}(\alpha\circ\beta\circ\beta_2)^*(-K_X)-\frac{1}{9}(\beta\circ\beta_2)^*(E)
-\frac{1}{7}\beta_2^*(F)-\frac{1}{2}F_2,\\
&E_{W_2}\sim_{\mathbb{Q}}(\beta\circ\beta_2)^*(E)
-\frac{5}{7}\beta_2^*(F)-\frac{1}{2}F_2,\\
&F_{W_2}\sim_{\mathbb{Q}}\beta_2^*(F)-\frac{1}{2}F_2.\\
\endaligned
\right.
$$
Furthermore, \[S_{W_2}\cdot
D_{W_2}=\bar{C}_{W_2}+L_{W_2}+\bar{L}_{W_2}, \ \ E_{W_2}\cdot
D_{W_2}=3L_{W_2}, \ \ F_{W_2}\cdot D_{W_2}=3\bar{L}_{W_2}.\]

Consider the curves $\bar{C}_{W_2}$, $L_{W_2}$, and
$\bar{L}_{W_2}$ as divisors on $D_{W_2}$. Then, the equivalences
above imply that
\[\bar{C}_{W_2}\cdot \bar{C}_{W_2}=L_{W_2}\cdot L_{W_2}=-\frac{1}{2}, \ \ \bar{L}_{W_2}\cdot
\bar{L}_{W_2}=-\frac{2}{5}, \ \ C_{W_2}\cdot L_{W_2}=C_{W_2}\cdot
\bar{L}_{W_2}=L_{W_2}\cdot \bar{L}_{W_2}=0.\] Therefore, the
intersection form of these curves on $D_{W_2}$ is
negative-definite. On the other hand, we have
\[M_{W_2}\Big\vert_{D_{W_2}}\equiv-nK_{W_2}\Big\vert_{D_{W_2}}\equiv nS_{W_2}\Big\vert_{D_{W_2}}\equiv n\bar{C}_{W_2}+nL_{W_2}+n\bar{L}_{W_2},\]
where $M_{W_2}$ is a general surface in the pencil
$\mathcal{M}_{W_2}$. Therefore, we obtain
$\mathcal{M}=\mathcal{P}$ from Theorem~\ref{theorem:main-tool}.
\end{proof}

\begin{K3proposition}\label{K3proposition:n-48}
A general surface in the pencil $\mathcal{P}$ is birational to a
K3 surface.
\end{K3proposition}
\begin{proof}
We use the same notations in the proof of
Lemma~\ref{lemma:n-48-first}. Note that the proof of
Lemma~\ref{lemma:n-48-first} and
Proposition~\ref{proposition:Halphen-K3} shows the linear system
$\mathcal{P}$ is a Halphen pencil.

Suppose that the intersection curve $\Delta=F_2\cdot D_{W_2}$ is a
smooth curve on $F_2\cong \mathbb{P}^2$. Because the degree of the
curve $\Delta$ on $F_2$ is three, it must be an elliptic curve.
One can see that the singularities of the image surface
$D_U:=\beta_2(D_{W_2})$ are rational except the point $Q_2$.
Moreover, the restricted morphism $\beta_2|_{D_{W_2}}: D_{W_2}\to
D_U$ resolves the singular point $Q_1$.  Therefore, the
singularities of the surface $D_{W_2}$ is rational. Because the
divisor $-K_U$ is nef and big, the Leray spectral sequence for the
morphism $\bar{\beta}:=\beta_2|_{D_{W_2}}: D_{W_2}\to D_U$ shows
\[0\to H^1(D_U, \mathcal{O}_{D_U})=0\to  H^1(D_{W_2},
\mathcal{O}_{D_{W_2}})\to H^0(D_U,
R^1\bar{\beta}_*\mathcal{O}_{D_{W_2}}) \to H^2(D_U,
\mathcal{O}_{D_U})=0,\] and hence the irregularity $h^1(D_{W_2},
\mathcal{O}_{D_{W_2}})=1$. Because the surface $D_{W_2}$ has only
rational singularities, a smooth surface birational to the surface
$D_{W_2}$ has the same irregularity. However, the surface
$D_{W_2}$ is birational to a K3 surface or an abelian surface by
Corollary~\ref{corollary:Halphen-K3}, which is a contradiction.
Therefore, the curve $\Delta$ must be a singular curve, and hence
a rational curve. Therefore,
Corollary~\ref{corollary:Halphen-K3-rational-curve} completes the
proof.
\end{proof}

Due to the lemma above, we may assume that
$\mathbb{CS}(U,\frac{1}{n}\mathcal{M}_U)=\{P_2\}$. Let
$\gamma:V\to U$ be the Kawamata blow up at the point $P_2$ with
weights $(1,2,3)$. Then, the exceptional divisor
$G\cong\mathbb{P}(1,2,3)$ of the birational morphism $\gamma$
contains two singular points $P_3$ and $Q_3$ that are quotient
singularities of types $\frac{1}{3}(1,2,1)$ and
$\frac{1}{2}(1,1,1)$, respectively. Again, the set
$\mathbb{CS}(V,\frac{1}{n}\mathcal{M}_V)$ must contain either the
point $P_3$ or the point $Q_3$.

\begin{lemma}
\label{lemma:n-48-second} The set
$\mathbb{CS}(V,\frac{1}{n}\mathcal{M}_V)$  cannot contain the
point $Q_3$.
\end{lemma}

\begin{proof}
Suppose that the set $\mathbb{CS}(V,\frac{1}{n}\mathcal{M}_V)$
contains the point $Q_3$. Let $\gamma_3:W_3\to V$ be the Kawamata
blow up at the point $Q_3$ with weights $(1,1,1)$ and
 $G_3$ be the exceptional divisor of $\gamma_3$. Then,
$\mathcal{M}_{W_3}\sim_{\mathbb{Q}}-nK_{W_3}$ by
Lemma~\ref{lemma:Kawamata}. The base locus of the pencil
$\mathcal{P}_{W_3}$ consists of three irreducible curves
$\bar{C}_{W_3}$, $L_{W_3}$, and $\bar{L}_{W_3}$.

Let $D_{W_3}$ be a general surface in the pencil
$\mathcal{P}_{W_3}$. Then, the surface $D_{W_3}$ is normal and we
have
$$
\left\{\aligned
&D_{W_3}\sim_{\mathbb{Q}}(\alpha\circ\beta\circ\gamma\circ\gamma_3)^*(-3K_X)
-\frac{3}{9}(\beta\circ\gamma\circ\gamma_3)^*(E)-\frac{3}{7}(\gamma\circ\gamma_3)^*(F)
-\frac{3}{5}\gamma_3^*(G)-\frac{1}{2}G_3\\
&S_{W_3}\sim_{\mathbb{Q}}(\alpha\circ\beta\circ\gamma\circ\gamma_3)^*(-K_X)
-\frac{1}{9}(\beta\circ\gamma\circ\gamma_3)^*(E)-\frac{1}{7}(\gamma\circ\gamma_3)^*(F)
-\frac{1}{5}\gamma_3^*(G)-\frac{1}{2}G_3,\\
&E_{W_3}\sim_{\mathbb{Q}}(\beta\circ\gamma\circ\gamma_3)^*(E)
-\frac{5}{7}(\gamma\circ\gamma_3)^*(F)-\frac{3}{5}\gamma_3^*(G)-\frac{1}{2}G_3,\\
&F_{W_3}\sim_{\mathbb{Q}}(\gamma\circ\gamma_3)^*(F)-\frac{3}{5}\gamma_3^*(G)-\frac{1}{2}G_3,\\
&G_{W_3}\sim_{\mathbb{Q}}\gamma_3^*(G)-\frac{1}{2}G_3,\\
\endaligned
\right.
$$
Furthermore,
\[S_{W_3}\cdot
D_{W_3}=\bar{C}_{W_3}+L_{W_3}+\bar{L}_{W_3}, \ \ E_{W_3}\cdot
D_{W_3}=3L_{W_3}, \ \ F_{W_3}\cdot D_{W_3}=3\bar{L}_{W_3}.\]
Consider the curves $\bar{C}_{W_3}$, $L_{W_3}$, and
$\bar{L}_{W_3}$ as divisors on $D_{W_3}$. Then, the equivalences
above imply that
\[\bar{C}_{W_3}\cdot \bar{C}_{W_3}=L_{W_3}\cdot L_{W_3}=-\frac{44}{90}, \ \ \bar{L}_{W_3}\cdot
\bar{L}_{W_3}=-\frac{35}{90}, \]
\[ C_{W_3}\cdot
L_{W_3}=C_{W_3}\cdot \bar{L}_{W_3}=\frac{19}{90}, \ \ L_{W_3}\cdot
\bar{L}_{W_3}=\frac{1}{90}.\]
 It is easy to see that the intersection form
of these curves on $D_{W_3}$ is negative-definite. On the other
hand, we have
\[M_{W_3}\Big\vert_{D_{W_3}}\equiv-nK_{W_3}\Big\vert_{D_3}\equiv nS_{W_3}\Big\vert_{D_{W_3}}\equiv n\bar{C}_{W_3}+nL_{W_3}+n\bar{L}_{W_3},\]
where $M_{W_3}$ is a general surface in the pencil
$\mathcal{M}_{W_3}$. Therefore, we obtain
$\mathcal{M}_{W_3}=\mathcal{P}_{W_3}$ from
Theorem~\ref{theorem:main-tool}. However, it is a contradiction
because $D_3\not\sim_{\mathbb{Q}}-3K_{W_3}$.
\end{proof}

Therefore, the set  $\mathbb{CS}(V,\frac{1}{n}\mathcal{M}_V)$
consists of only one point $P_3$. Let $\delta:V_3\to V$ be the
Kawamata blow up at the point $P_3$ with weights $(1,2,1)$. The
pencil $|-2K_{V_3}|$ is the proper transform of the pencil
$|-2K_X|$. It has only one base curve $C_{V_3}$ whose image to $X$
is the base curve of the pencil $|-2K_X|$. Then, the inequality
$-K_{V_3}\cdot C_{V_3}=-2K_{V_3}^3<0$ and the equivalence
$\mathcal{M}_{V_3}\sim_{\mathbb{Q}} -nK_{V_3}$ imply
$\mathcal{M}=|-2K_X|$ by Theorem~\ref{theorem:main-tool}.

\begin{proposition}
If $\gimel= 48$, then the linear systems $|-2K_X|$ and
$\mathcal{P}$ are the only Halphen pencils on $X$.
\end{proposition}

\section{Cases $\gimel=55$, $80$, and $91$.}\index{$\gimel=55$}\index{$\gimel=80$}\index{$\gimel=91$}
\label{section:n-80-91}

The threefold $X \subset \mathbb{P}(1,a_1,a_2,a_3,a_4)$  of degree
$d=\sum a_i$ always contains the point $O=(0:0:0:1:0)$. It is a
singular point of $X$ that is a quotient singularity of type
$\frac{1}{a_3}(1,a_1, a_3-a_1)$. The threefold $X$ can be given by
$$
t^{3}z+\sum_{i=0}^{2}t^{i}f_{d-ia_{3}}\big(x,y,z,w\big)=0,
$$
where  $f_{i}$ is a general qua\-si\-ho\-mo\-ge\-ne\-ous
polynomial of degree $i$.

There is a commutative diagram
$$
\xymatrix{
&U\ar@{->}[d]_{\alpha}&&W\ar@{->}[ll]_{\beta}\ar@{->}[d]^{\eta}&\\%
&X\ar@{-->}[rr]_{\psi}&&\mathbb{P}(1,a_1,a_2),&}
$$
where \begin{itemize} \item $\alpha$ is the Kawamata blow up at
the point $O$ with weights $(1,a_1,a_3-a_1)$,

\item $\beta$ is the Kawamata blow up with weights
$(1,a_1,a_3-2a_1)$ at the singular point $Q$ of the variety $U$
that is a quotient singularity of type
$\frac{1}{a_3-a_1}(1,a_1,a_3-2a_1)$,

\item $\eta$ is an elliptic fibration.
\end{itemize}

Let $\mathcal{P}$ be the pencil defined by $$\lambda x^{a_{2}}+\mu
z=0,$$ where $(\lambda:\mu)\in\mathbb{P}^{1}$. Note that the base
curve $\bar{C}$ of $\mathcal{P}$ is given by the equations
$x=z=0$. We will see that the linear systems $|-a_1K_X|$ and
$\mathcal{P}$ are the only Halphen pencils on $X$.

We may assume that
$$
\mathbb{CS}\Big(X, \frac{1}{n}\mathcal{M}\Big)=\Big\{O\Big\}.%
$$
due to Theorem~\ref{theorem:Noether-Fano},
Lemmas~\ref{lemma:smooth-points},
\ref{lemma:special-singular-points-with-positive-c},
\ref{lemma:special-singular-points-with-zero-c}, and
Corollary~\ref{corollary:Ryder-a1}.

The exceptional divisor $E\cong\mathbb{P}(1,a_1,a_3-a_1)$ of the
birational morphism $\alpha$ contains two singular points $P$ and
$Q$ that are quotient singularities of types
$\frac{1}{a_1}(1,a_1-1,1)$ and
$\frac{1}{a_3-a_1}(1,a_1,a_3-2a_1)$, respectively. For the
convenience, let $L$ be the unique irreducible curve  contained in
$|\mathcal{O}_{\mathbb{P}(1,\,a_1,\,a_3-a_1)}(1)|$ on the surface
$E$.

\begin{lemma}
\label{lemma:n-80-81-91-P5} If the set $\mathbb{CS}(U,
\frac{1}{n}\mathcal{M}_U)$ contains the point $P$, then
$\mathcal{M}=\mathcal{P}$.
\end{lemma}

\begin{proof} We will consider only the case
$\gimel=91$. The other cases can be shown by the same method.

 Suppose that the set $\mathbb{CS}(U,
\frac{1}{n}\mathcal{M}_U)$ contains the point $P$. Let
$\gamma:V\to U$ be the Kawamata blow up at the point $P$ with
weights $(1,3,1)$ and let $G_P$ be the exceptional divisor of
$\gamma$.

Around the point $O$, the monomials $x$, $y$, and $w$ can be
considered as weighted local coordinates with weights
$\mathrm{wt}(x)=1$, $\mathrm{wt}(y)=4$, and $\mathrm{wt}(w)=9$.
Therefore, around the singular point $P$, the birational morphism
$\alpha$ is given by the equations
$$
x=\tilde{x}\tilde{y}^{\frac{1}{13}},\ y=\tilde{y}^{\frac{4}{13}},\ w=\tilde{w}\tilde{y}^{\frac{9}{13}},%
$$
where $\tilde{x}$, $\tilde{y}$, and $\tilde{w}$ are weighted local
coordinates with $\mathrm{wt}(\tilde{x})=\mathrm{wt}(\tilde{w})=1$
and $\mathrm{wt}(\tilde{y})=3$. Let $R$ be a general surface of
the pencil $\mathcal{P}$.  Then, $R$ is given by the equation of
the form
$$
\lambda x^5+\mu\Big(\delta_1 w^2+\sum_{i=0}^{4}\delta_{2+i} x^{18-4i}y^{i}+\delta_7 xy^{2}w+\delta_8 x^{5}yw+\delta_9 x^{9}w+\mbox{higher terms}\Big)=0%
$$
near the point  $O$, where $\delta_{i}\in\mathbb{C}$. The proper
transform $R_U$  is given by
$$
\lambda \tilde{x}^5+\mu\tilde{y}\Big(\delta_1 \tilde{w}^2+\delta_2\tilde{x}^{2}+\delta_7 \tilde{x}\tilde{w}+\mbox{higher terms}\Big)=0%
$$
near the point  $P$. The base locus of the pencil $\mathcal{P}_U$
consists of the irreducible curves $\bar{C}_U$ and  $L$. It shows
that  $|-5K_{V}|=\mathcal{P}_V$ and the base locus of the pencil
$|-5K_{V}|$ consists of the curves $\bar{C}_V$ and $L_V$.
Furthermore,
$$
R_V\cdot S_V=\bar{C}_V+L_V,\ R_V\cdot E_V=5L_V,%
$$
which implies that $R_V$ is normal. On the other hand, we have
$$
\left\{\aligned
&R_V\sim_{\mathbb{Q}}\big(\alpha\circ\gamma\big)^{*}\Big(-5K_{X}\Big)-\frac{5}{13}\gamma^{*}\big(E\big)-\frac{5}{4}G_P,\\
&S_V\sim_{\mathbb{Q}}\big(\alpha\circ\gamma\big)^{*}\Big(-K_{X}\Big)-\frac{1}{13}\gamma^{*}\big(E\big)-\frac{1}{4}G_P,\\
&E_V\sim_{\mathbb{Q}}\gamma^{*}\big(E\big)-\frac{3}{4}G_P,\\
\endaligned
\right.
$$
which implies that on the normal surface $R_V$ we have
$$
\left\{\aligned
&L_V\cdot L_V=\frac{E_V\cdot E_V\cdot R_V}{25}=-\frac{4}{45},\\
&\bar{C}_V\cdot \bar{C}_V=S_V\cdot S_V\cdot R_V-\frac{2}{5}S_V\cdot E_V\cdot R_V-\frac{E_V\cdot E_V\cdot R_V}{25}=-\frac{11}{30},\\
&\bar{C}_V\cdot L_V=\frac{S_V\cdot E_V\cdot R_V}{5}-\frac{E_V\cdot E_V\cdot R_V}{25}=\frac{1}{30}.\\
\endaligned
\right.
$$
It immediately implies that the intersection form of the curves
$\bar{C}_V$ and $L_V$ on the normal surface $R_V$ is
negative-definite. On the other hand, we have
$$
\mathcal{M}_V\Big\vert_{R_V}\equiv-nK_{V}\Big\vert_{R_V}\equiv nS_V\Big\vert_{R_V}\equiv n\bar{C}_V+nL_V,%
$$
which implies that $\mathcal{M}_V=|-5K_{V}|$ by
Theorem~\ref{theorem:main-tool}. In particular, we have
$\mathcal{M}=\mathcal{P}$.
\end{proof}
\begin{K3proposition}\label{K3proposition:n-55-80-91}
A general surface in $\mathcal{P}$ is birational to a smooth K3
surface.
\end{K3proposition}
\begin{proof}
We use the same notations in the proof of
Lemma~\ref{lemma:n-80-81-91-P5}. The exceptional
divisor $G_P$ of the Kawamata blow up $\gamma$ is isomorphic to
$\mathbb{P}(1, a_1-1, 1)$. Then, the intersection
$\Delta:=G_P\cdot R_V$ is a curve of degree $a_2=a_1+1$ on
$\mathbb{P}(1, a_1-1, 1)$.

In the case $\gimel=91$, the surface $G_P$ has a quotient singular
point. One can easily check the curve $G_P\cdot R_V$ passes
through the singular point, and hence it is a rational curve.
Therefore,
Corollary~\ref{corollary:Halphen-K3-rational-curve} implies that
the surface $R_V$ is birational to a smooth K3 surface.

In the cases $\gimel=55$ and $80$, the curve $\Delta$ does not
pass through a singular point of the surface $G_P$. We suppose
that the curve $\Delta$ is smooth.  Because  it does not pass
through any singular point of $G_P$ and its degree on $G_P$ is
$a_2=a_1+1$, it is an elliptic curve. One can see that the
singularities of the image surface $R_U:=\gamma(R_V)$ are rational
except the point $P$. Then, the same argument of
K3-Proposition~\ref{K3proposition:n-48} leads us to a
contradiction. Therefore, the surface $R_V$ is birational to a
smooth K3 surface.
\end{proof}

The exceptional divisor $F$ of the birational morphism $\beta$
contains two singular points $Q_1$ and $Q_2$ that are quotient
singularities of types $\frac{1}{a_1}(1,a_1-1,1)$ and
$\frac{1}{a_3-2a_1}(1,a_1,a_3-3a_1)$, respectively.

\begin{lemma}
\label{lemma:n-80-81-91-P7} The set $\mathbb{CS}(W,
\frac{1}{n}\mathcal{M}_{W})$ does not contain the point $Q_1$.
\end{lemma}

\begin{proof}
Suppose that the set $\mathbb{CS}(W, \frac{1}{n}\mathcal{M}_{W})$
contains the point $Q_1$. Let $\gamma_1:Y_1\to W$ be the Kawamata
blow up at the singular point  $Q_1$ with weights $(1,a_1-1,1)$
and let $G_1$ be the exceptional divisor of $\gamma_1$.

Simple calculations imply that the base locus of the proper
transform of the linear system $|-a_2K_{X}|$ on the threefold
$Y_1$ consists of the curves $\bar{C}_{Y_1}$ and $L_{Y_1}$.

Let $B$ be a general surface of the linear system $|-a_2K_{X}|$.
Then,
$$
B_{Y_1}\cdot S_{Y_1}=\bar{C}_{Y_1}+L_{Y_1},\ B_{Y_1}\cdot E_{Y_1}=L_{Y_1},%
$$
which implies that $B_{Y_1}$ is normal. On the other hand, we have
$$
\left\{\aligned
&B_{Y_1}\sim_{\mathbb{Q}}\big(\alpha\circ\beta\circ\gamma_1\big)^{*}\Big(-a_2K_{X}\Big)
-\frac{a_2}{a_3}\big(\beta\circ\gamma_1\big)^{*}\big(E\big)-\frac{a_2}{a_3-a_1}\gamma_1^{*}\big(F\big)-\frac{1}{a_1}G_1,\\
&S_{Y_1}\sim_{\mathbb{Q}}\big(\alpha\circ\beta\circ\gamma_1\big)^{*}\Big(-K_{X}\Big)
-\frac{1}{a_3}\big(\beta\circ\gamma_1\big)^{*}\big(E\big)-\frac{1}{a_3-a_1}\gamma_1^{*}\big(F\big)-\frac{1}{a_1}G_1,\\
&E_{Y_1}\sim_{\mathbb{Q}}\big(\beta\circ\gamma_1\big)^{*}\big(E\big)-\frac{a_3-2a_1}{a_3-a_1}\gamma_1^{*}\big(F\big)-\frac{1}{a_1}G_1,\\
&F_{Y_1}\sim_{\mathbb{Q}}\gamma_1^{*}\big(F\big)-\frac{a_1-1}{a_1}G_1.\\
\endaligned
\right.
$$
These equivalence shows that the intersection form of
$\bar{C}_{Y_1}$ and $L_{Y_1}$ on the surface $B_{Y_1}$ is
negative-definite.\footnote{The curves $\bar{C}_{Y_1}$, $L_{Y_1}$,
and $B_{Y_1}\cap F$ are components of a fiber of the elliptic
fibration $\eta\circ\gamma_1\vert_{B_{Y_1}}$, which implies that
the intersection form of $\bar{C}_{Y_1}$ and $L_{Y_1}$ on the
surface $B_{Y_1}$ is negative-definite.} Since
$\mathcal{M}_{Y_1}\vert_{B_{Y_1}}\equiv n\bar{C}_{Y_1}+nL_{Y_1}$,
the pencil $\mathcal{M}$  coincides with the linear system
$|-a_2K_{X}|$ by Theorem~\ref{theorem:main-tool}. However, the
linear system  $|-a_2K_{X}|$ is not a pencil.
\end{proof}

\begin{lemma}
\label{lemma:n-80-81-91-P6} If the set $\mathbb{CS}(W,
\frac{1}{n}\mathcal{M}_{W})$ contains the point $Q_2$, then
$\mathcal{M}=|-a_1K_{X}|$.
\end{lemma}

\begin{proof} Suppose that the set $\mathbb{CS}(W,
\frac{1}{n}\mathcal{M}_{W})$ contains the point $Q_2$. Let
$\gamma_2:Y_2\to W$ be the Kawamata blow up at the point $Q_2$
with weights $(1,a_1,a_3-3a_1)$. Then, $|-a_1K_{V}|$ is the proper
transform of the pencil $|-a_1K_{X}|$ and the base locus of
$|-a_1K_{V}|$ consists of the irreducible curve $C_V$ whose image
to $X$ is the base curve of the pencil $|-a_1K_{X}|$.

Let $D$ be a general surface of the pencil $|-a_1K_{V}|$. Then,
$D$ is normal and we can consider the curve $C_V$ as a divisor on
$D$. We have $C_V^{2}<0$ but $\mathcal{M}_{Y_2}\vert_{D}\equiv
nC_V$ by Lemma~\ref{lemma:Kawamata}, which implies that
$\mathcal{M}_{Y_2}=|-a_1K_{V}|$ by
Theorem~\ref{theorem:main-tool}.
\end{proof}

\begin{proposition}
\label{proposition:n-80-91} The linear systems $|-a_{1}K_{X}|$ and
$\mathcal{P}$ is the only Halphen pencils on $X$.
\end{proposition}
\begin{proof}
It immediately follows from the previous arguments.
\end{proof}

\section{Cases $\gimel=57$,  $66$, $81$,  and $86$.}%
\index{$\gimel=57$}\index{$\gimel=66$}\index{$\gimel=81$}\index{$\gimel=86$}\label{section:n-57-66-86}

Suppose that $\gimel\in\{ 57,  66, 81,   86\}$. Then, the
threefold $X \subset \mathbb{P}(1,a_1,a_2,a_3,a_4)$  of degree
$d=\sum_{i=1}^4 a_i$ always contains the point $O=(0:0:0:1:0)$. It
is a singular point of $X$ that is a quotient singularity of type
$\frac{1}{a_3}(1,a_1, a_3-a_1)$. The threefold $X$ can be given by
$$
t^{\frac{d-a_2}{a_3}}z+\sum_{i=0}^{\frac{d-a_2}{a_3}}t^{i}f_{d-ia_{3}}\big(x,y,z,w\big)=0,
$$
where $f_{i}$ is a general qua\-si\-ho\-mo\-ge\-ne\-ous polynomial
of degree $i$.

We also have a commutative diagram as follows:

$$
\xymatrix{
&&&Y\ar@{->}[lld]_{\pi}\ar@{->}[rrd]^{\eta}&&&\\%
&X\ar@{-->}[rrrr]_{\psi}&&&&\mathbb{P}(1,a_1,a_{2})&}
$$
where \begin{itemize} \item $\psi$ is the natural projection,
\item $\pi$ is the Kawamata blow up at the point $O$ with weights
$(1, a_1, a_3-a_1)$, \item $\eta$ is an elliptic
fibration.\end{itemize}

Let $\mathcal{P}$ be the pencil defined by $$\lambda x^{a_{2}}+\mu
z=0,$$ where $(\lambda:\mu)\in\mathbb{P}^{1}$.

We may assume that the set $\mathbb{CS}(X,
\frac{1}{n}\mathcal{M})$ consists of the point $O$ by
Lemmas~\ref{lemma:smooth-points},
\ref{lemma:special-singular-points-with-positive-c},
\ref{lemma:special-singular-points-with-zero-c}, and
Corollary~\ref{corollary:Ryder-a1}.

The exceptional divisor $E$ of the birational morphism $\pi$
contains two singular points $P$ and $Q$ that are quotient
singularities of types $\frac{1}{a_3-a_1}(1,a_3-a_2, a_2-a_1)$ and
$\frac{1}{a_1}(1,2a_1-a_3, a_3-a_1)$, respectively.

\begin{lemma}\label{lemma:n-57-66-81-86-P}
If the set $\mathbb{CS}(Y, \frac{1}{n}\mathcal{M}_Y)$ contains the
point $P$, then $\mathcal{M}=|-a_1K_X|$.
\end{lemma}
\begin{proof}
The proof of Proposition~\ref{proposition:n-34-53-70-88}
immediately implies that $\mathcal{M}=|-a_1K_{X}|$.
\end{proof}

\begin{lemma}
\label{lemma:n-57-66-81-86-Q} If the set $\mathbb{CS}(Y,
\frac{1}{n}\mathcal{M}_Y)$ contains the point $Q$, then
$\mathcal{M}=\mathcal{P}$.
\end{lemma}
\begin{proof}

Suppose that the set $\mathbb{CS}(Y, \frac{1}{n}\mathcal{M}_Y)$
contains the point $Q$. Let $\beta:W\to Y$ be the Kawamata blow up
at the point $Q$ with weights $(1,2a_1-a_3, a_3-a_1)$ and let $F$
be the exceptional divisor of the birational morphism $\beta$.
Then,
$$
\mathcal{P}_W\sim_{\mathbb{Q}}\big(\pi\circ\beta\big)^{*}\Big(-a_2K_{X}\Big)-\frac{a_2}{a_3}\beta^{*}\big(E\big)-\frac{a_2}{a_1}F,\\
$$
but the base locus of $\mathcal{P}_W$ consists of the irreducible
curve $\bar{C}_W$ whose image to $X$ is the unique base curve
$\bar{C}$ of the linear system $\mathcal{P}$.

Let $D$ and $M$ be general surfaces in $\mathcal{P}_W$ and
$\mathcal{M}_{W}$, respectively. Then, $D$ is normal and
$$
M\Big\vert_{D}\equiv-nK_{W}\Big\vert_{D}\equiv n\bar{C}_W,%
$$
but $\bar{C}_W^{2}<0$ on the surface $D$. Therefore, we obtain
$\mathcal{M}_{W}=\mathcal{P}_W$ from
Theorem~\ref{theorem:main-tool}.
\end{proof}
\begin{K3proposition}\label{K3proposition:n-57-66-81-86}
A general surface in the pencil $\mathcal{P}$ is birational to a
smooth K3 surface.
\end{K3proposition}
\begin{proof}
We use the same notation in the proof of
Lemma~\ref{lemma:n-57-66-81-86-Q}.  The exceptional divisor $E$ is
isomorphic to $\mathbb{P}(1, a_1, a_3-a_1)$. The curve $\Delta$
defined by the intersection of a general surface in
$\mathcal{P}_Y$ with $E$ is a curve of degree $a_2$. Because $a_2<
a_3+1$, the curve $\Delta$ is a rational curve. The result follows
from   Theorem~\ref{theorem:Halphen} since the curve $\Delta_W$ is not contained in the base locus of the pencil $\mathcal{P}_W$.
\end{proof}

\begin{proposition}
\label{proposition:n-57-66-81-86} The linear systems $|-a_1K_X|$
and $\mathcal{P}$ are the only Halphen pencils on $X$.
\end{proposition}

\begin{proof}
Because we may assume that the set $\mathbb{CS}(Y,
\frac{1}{n}\mathcal{M}_Y)$ contains either the point $P$ or the
point $Q$, it immediately follows from
Lemmas~\ref{lemma:n-57-66-81-86-P} and
\ref{lemma:n-57-66-81-86-Q}.
\end{proof}

\section{Case $\gimel=58$, hypersurface of degree $24$ in
$\mathbb{P}(1,3,4,7,10)$.}\index{$\gimel=58$} \label{section:n-58}

Let $X$ be the hypersurface given by a general quasihomogeneous
equation of degree $24$ in $\mathbb{P}(1,3,4,7,10)$
 with
$-K_X^3=\frac{1}{35}$. Then, the singularities of $X$ consist of
two singular points $P$ and $Q$ that are quotient singularities of
types $\frac{1}{7}(1,3,4)$ and $\frac{1}{10}(1,3,7)$,
respectively,  and one point of type $\frac{1}{2}(1,1,1)$. Also,
by the generality of the hypersurface, we may assume that the
hypersurface $X$ is defined by the equation
$$
w^{2}z+wf_{14}(x,y,z,t)+f_{24}(x,y,z,t)=0,
$$
 where $f_{i}(x,y,z,t)$ is a quasihomogeneous polynomial of degree
$i$. Proposition~\ref{proposition:Halphen-K3} implies that the
linear system $|-3K_X|$ is a Halphen pencil. Let $\mathcal{P}$ be
the pencil on $X$ given by the equations
$$
\lambda x^4+\mu z=0,
$$
where $(\lambda, \mu)\in \mathbb{P}^1$. We will see that the
linear system  $\mathcal{P}$ is another Halphen pencil on $X$
(K3-Proposition~\ref{K3proposition:n-58}).

We have the following commutative diagram:
$$
\xymatrix{
&&&Y\ar@{->}[dl]_{\gamma_{O}}\ar@{->}[dr]^{\gamma_{P}}\ar@{->}[drrrrr]^{\eta}&&&&&&\\
&&U_{PQ}\ar@{->}[dl]_{\beta_{Q}}\ar@{->}[dr]^{\beta_{P}}&&U_{QO}\ar@{->}[dl]^{\beta_{O}}&&&&\mathbb{P}(1,3,4),\\
&U_{P}\ar@{->}[dr]_{\alpha_{P}}&&U_{Q}\ar@{->}[dl]^{\alpha_{Q}}&&&&&&\\
&&X\ar@{-->}[rrrrrruu]_{\psi}&&&&&&&}
$$
where
\begin{itemize}
\item $\psi$ is  the natural projection,

\item $\alpha_{P}$ is the Kawamata blow up at the point $P$ with
weights $(1,3,4)$,

\item $\alpha_{Q}$ is the Kawamata blow up at the point  $Q$ with
weights $(1,3,7)$,

\item $\beta_{Q}$ is the Kawamata blow up with weights $(1,3,7)$
at the point whose image by the birational morphism $\alpha_P$ is
the point $Q$,

\item $\beta_{P}$ is the Kawamata blow up with weights $(1,3,4)$
at the point whose image by the birational morphism $\alpha_Q$ is
the point $P$,

 \item $\beta_{O}$ is the Kawamata blow up with weights
$(1,3,4)$ at the singular point $O$ of the variety $U_{Q}$ that is
a quotient singularity of type $\frac{1}{7}(1,3,4)$ contained in
the exceptional divisor of the birational morphism $\alpha_{Q}$,

\item $\gamma_{P}$ is the Kawamata blow up with weights $(1,3,4)$
at the point whose image by the birational morphism
$\alpha_Q\circ\beta_O$ is the point $P$,

\item $\gamma_{O}$ is the Kawamata blow up with weights $(1,3,4)$
at the singular point of the variety $U_{PQ}$ that is a quotient
singularity of type $\frac{1}{7}(1,3,4)$ contained in the
exceptional divisor of the birational morphism $\beta_{Q}$,

\item  $\eta$ is an elliptic fibration.
\end{itemize}
Because of Lemma~\ref{lemma:smooth-points} and
Corollary~\ref{corollary:Ryder-a1}, we may assume that
$$\mathbb{CS}\Big(X,\frac{1}{n}\mathcal{M}\Big)\subset \Big\{P, Q\Big\}.$$

The exceptional divisor $E_P$ of the birational morphism
$\alpha_P$ contains two quotient singular points $P_1$ and $P_2$
of types $\frac{1}{4}(1,3,1)$ and $\frac{1}{3}(1,2,1)$,
respectively.

\begin{lemma}\label{lemma:n-58-P-1}
If the set $\mathbb{CS}(U_P,\frac{1}{n}\mathcal{M}_{U_P})$
contains  the point $P_1$, then  $\mathcal{M}=|-3K_X|$.
\end{lemma}
\begin{proof}
Suppose it contains the point $P_1$. Let $\beta_1:W_1\to U_P$ be
the Kawamata blow up at the point $P_1$ with weights $(1,3,1)$.
Then, the pencil $|-3K_{W_1}|$ is the proper transform of the
pencil system $|-3K_X|$ . Its base locus consists of the
irreducible curve $C_{W_1}$ whose image to $X$ is the base curve
of the pencil $|-3K_X|$. Then, $-K_{W_1}\cdot C_{W_1}<0$ and
$\mathcal{M}_{W_1}\sim_{\mathbb{Q}} -nK_{W_1}$ imply
$\mathcal{M}=|-3K_X|$ by Theorem~\ref{theorem:main-tool}.
\end{proof}

\begin{lemma}\label{lemma:n-58-P-2}
The set $\mathbb{CS}(U_P,\frac{1}{n}\mathcal{M}_{U_P})$ cannot
contain the point $P_2$.
\end{lemma}
\begin{proof}
Suppose it contains the point $P_2$. Let $\beta_2:W_2\to U_P$ be
the Kawamata blow up at the point $P_2$ with weights $(1,2,1)$.
Also, let $\mathcal{D}_2$ be the proper transform of the linear
system $|-4K_X|$ by the birational morphism
$\alpha_P\circ\beta_2$. Its base locus consists of the irreducible
curve $\bar{C}_{W_2}$ whose image to $X$ is the base curve of the
linear system $|-4K_X|$. A general surface $D_2$ in
$\mathcal{D}_2$ is normal and the self-intersection
$\bar{C}_{W_2}^2$ is negative on the surface $D_2$. Because
$\mathcal{M}_{W_2}|_{D_2}\equiv-n\bar{C}_{W_2}$, we obtain an
absurd identity $\mathcal{M}=|-4K_X|$ from
Theorem~\ref{theorem:main-tool}.
\end{proof}

Meanwhile, the exceptional divisor $E\cong\mathbb{P}(1,3,7)$ of
the birational morphism $\alpha_Q$ contains two singular points
$O$ and $Q_1$ of types $\frac{1}{7}(1,3,4)$ and
$\frac{1}{3}(1,2,1)$. For the convenience, let $L$ be  the unique
curve contained in the linear system
$|\mathcal{O}_{\mathbb{P}(1,3,7)}(1)|$ on $E$ .

\begin{lemma}\label{lemma:n-58-Q-1}
If the set $\mathbb{CS}(U_Q,\frac{1}{n}\mathcal{M}_{U_Q})$
contains the point $Q_1$, then $\mathcal{M}=\mathcal{P}$.
\end{lemma}
\begin{proof}
Suppose that the set
$\mathbb{CS}(U_Q,\frac{1}{n}\mathcal{M}_{U_Q})$ contains  the
point $Q_1$. Let $\pi_1:V_1\to U_Q$ be the Kawamata blow  up at
the point $Q_1$ with weights $(1,2,1)$. The base locus of
$\mathcal{P}_{V_1}$ consists of two irreducible curves
$\bar{C}_{V_1}$ and $L_{V_1}$.

For a general surface $D_{V_1}$ in $\mathcal{P}_{V_1}$, we have
\[S_{V_1}\cdot D_{V_1}=\bar{C}_{V_1}+L_{V_1}, \ E_{V_1}\cdot D_{V_1}=4L_{V_1}.\]
 Using the following equivalences

\[\left\{\aligned
&E_{V_1}\sim_{\mathbb{Q}}\pi_1^{*}\big(E\big)-\frac{2}{3}F_Q,\\
&S_{V_1}\sim_{\mathbb{Q}}\big(\alpha_Q\circ\pi_1\big)^{*}\Big(-K_{X}\Big)-\frac{1}{10}\pi_1^{*}\big(E\big)-\frac{1}{3}F_Q,\\
&D_{V_1}\sim_{\mathbb{Q}}\big(\alpha_Q\circ\pi_1\big)^{*}\Big(-4K_{X}\Big)-\frac{4}{10}\pi_1^{*}\big(E\big)-\frac{4}{3}F_Q,\\
\endaligned
\right.\] where $F_Q$ is the exceptional divisor of $\pi_1$, we
can obtain
\[D_{V_1}\cdot \bar{C}_{V_1}=D_{V_1}\cdot L_{V_1}=-\frac{8}{7}.\]
Because $\pi_1^{*}(-K_{U_Q})\cdot L_1=\pi_1^{*}(-K_{U_Q})\cdot
\bar{C}_{V_1}=\frac{1}{21}$,  the divisor
$B:=24\pi_1^{*}(-K_{U_Q})+D_{V_1}$ is  nef and big  and $B\cdot
\bar{C}_{V_1}=B\cdot L_{V_1}=0$. Therefore,
Theorem~\ref{theorem:main-tool} implies $\mathcal{M}=\mathcal{P}$.
\end{proof}
\begin{K3proposition}\label{K3proposition:n-58}
A general surface in the pencil $\mathcal{P}$ is birational to a
smooth K3 surface.
\end{K3proposition}
\begin{proof}
We use the same notations in the proof of
Lemma~\ref{lemma:n-58-Q-1}. The exceptional divisor $F_Q$ is
isomorphic to $\mathbb{P}(1, 2, 1)$. Then, the intersection
$\Delta:=F_Q\cdot D_{V_1}$ is a curve of degree $4$ on
$\mathbb{P}(1, 2, 1)$.

Easy calculation shows that the curve $\Delta$ does not pass
through the singular point of the surface $F_Q$. We suppose that
the curve $\Delta$ is smooth. Because it does not pass through any
singular point of $F_Q$ and its degree on $F_Q$ is four, it is an
elliptic curve. One can see that the singularities of the image
surface $D_{U_Q}:=\pi_1(D_{V_1})$ are rational except the point
$Q_1$. Then, the same argument of
K3-Proposition~\ref{K3proposition:n-48} gives a contradiction.
Therefore, the surface $D_{V_1}$ is birational to a smooth  K3
surface.
\end{proof}

The exceptional divisor $F_O$ of the birational morphism $\beta_O$
contains two quotient singular points $O_1$ and $O_2$ of types
$\frac{1}{4}(1,3,1)$ and $\frac{1}{3}(1,2,1)$, respectively.
\begin{lemma}\label{lemma:n-58-O-1}
If the set $\mathbb{CS}(U_{QO},\frac{1}{n}\mathcal{M}_{U_{QO}})$
contains the point $O_1$, then $\mathcal{M}=|-3K_X|$.
\end{lemma}
\begin{proof}Suppose that the set $\mathbb{CS}(U_{QO},\frac{1}{n}\mathcal{M}_{U_{QO}})$
contains the point $O_1$.  Let $\sigma_1:U_1\to U_{QO}$ be the
Kawamata blow up at the point $O_1$ with weights $(1,3,1)$. Then,
the pencil $|-3K_{U_1}|$ is the proper transform  of the pencil
$|-3K_X|$. Its base locus consists of the irreducible curve
$C_{U_1}$. Because we have
$\mathcal{M}_{U_1}\sim_{\mathbb{Q}}-nK_{U_1}$ and $-K_{U_1}\cdot
C_{U_1}<0$, we obtain $\mathcal{M}=|-3K_X|$ from
Theorem~\ref{theorem:main-tool}.
\end{proof}

\begin{lemma}\label{lemma:n-58-O-2}
The set $\mathbb{CS}(U_{QO},\frac{1}{n}\mathcal{M}_{U_{QO}})$
cannot contain the point $O_2$.
\end{lemma}
\begin{proof}
Suppose that the set
$\mathbb{CS}(U_{QO},\frac{1}{n}\mathcal{M}_{U_{QO}})$ contains the
point $O_2$. Let $\sigma_2:U_2\to U_{QO}$ be the Kawamata blow up
at the point $O_2$ with weights $(1,2,1)$.

 The base locus of the pencil $\mathcal{P}_{U_2}$
consists of the irreducible curves $\bar{C}_{U_2}$ and $L_{U_2}$.
For a general surface $H$ in $\mathcal{P}_{U_2}$, we have
$$
H\sim_{\mathbb{Q}}(\alpha_{Q}\circ\beta_{O}\circ\sigma_2)^{*}(-4K_{X})-\frac{4}{10}(\beta_{O}\circ\sigma_2)^{*}(E)-
\frac{4}{7}\sigma_2^{*}(F_O)-\frac{1}{3}G,
$$
where $G$ is the exceptional divisor of $\sigma_2$. The general
surface $H$ is normal, $S_{U_2}\cdot H=\bar{C}_{U_2}+L_{U_2}$, and
$E_{U_2}\cdot H=4L_{U_2}$. Since
$$
E_{U_2}\sim_{\mathbb{Q}}(\beta_{O}\circ\sigma_2)^{*}(E)-\frac{4}{7}\sigma_2^{*}(F_O)-\frac{1}{3}G,\\
$$
we can see that the intersection form of the curves $L_{U_2}$ and
$\bar{C}_{U_2}$ on the surface $H$ is negative-definite. The
equivalence $\mathcal{M}_{U_2}\vert_{H}\equiv
n\bar{C}_{U_2}+nL_{U_2}$ holds. Therefore, we can obtain
$\mathcal{M}=\mathcal{P}$ from Theorem~\ref{theorem:main-tool}.
However, $H\not\sim_{\mathbb{Q}} -4K_{U_2}$.
\end{proof}

\begin{proposition}
The linear systems $|-3K_X|$ and $\mathcal{P}$ are the only
Halphen pencils on $X$.
\end{proposition}
\begin{proof}
Due to the previous lemmas, we may assume that
$$\mathbb{CS}\Big(X,\frac{1}{n}\mathcal{M}\Big)= \Big\{P, Q\Big\}.$$
Following the Kawamata blow ups $Y\to U_{QO}\to U_Q\to X$ and
using Lemmas~\ref{lemma:n-58-Q-1}, \ref{lemma:n-58-O-1}, and
\ref{lemma:n-58-O-2}, we can furthermore assume that the set
$\mathbb{CS}(Y,\frac{1}{n}\mathcal{M}_{Y})$ contains one of
singular points contained the exceptional divisor of the
birational morphism $\gamma_P$. In this case,
Lemmas~\ref{lemma:n-58-P-1} and \ref{lemma:n-58-P-2} imply the
statement.
\end{proof}

\section{Case $\gimel=60$, hypersurface of degree $24$ in
$\mathbb{P}(1,4,5,6,9)$.}\index{$\gimel=60$} \label{section:n-60}

 The threefold $X$ is a
 general hypersurface of degree $24$ in $\mathbb{P}(1,4,5,6,9)$
with $-K^3_X=\frac{1}{45}$. Its singularities consist of one
singular point $O$ that is a quotient singularity of type
$\frac{1}{9}(1,4,5)$, one quotient singular point of type
$\frac{1}{5}(1,4,1)$,  one  quotient singular point of type
$\frac{1}{3}(1,1,2)$, and two quotient singular points of type
$\frac{1}{2}(1,1,1)$.

It cannot be birationally transformed into an elliptic fibration
(\cite{ChPa05}). However, a general fiber of the natural
projection $\xi:X\dasharrow\mathbb{P}(1,4)$ is birational to a
smooth K3 surface by Proposition~\ref{proposition:Halphen-K3}, in
other words, the linear system $|-4K_X|$ is a Halphen pencil.

By coordinate change, we may assume that the threefold $X$ is
given by the equation
$$
w^{2}t+wf_{15}\big(x,y,z,t\big)+f_{24}\big(x,y,z,t\big)=0,
$$
where $f_{i}(x,y,z,t)$ is a general quasihomogeneous polynomial of
degree $i$. Let $\mathcal{P}$ be the pencil consisting of surfaces
 cut out on the threefold $X$ by the equations
$$
\lambda x^{6}+\mu t=0,
$$
where $(\lambda:\mu)\in\mathbb{P}^{1}$.

\begin{K3proposition}
A general surface in the pencil $\mathcal{P}$ is birational to a
smooth  K3 surface.
\end{K3proposition}
\begin{proof}
It is a compactification of a double cover of $\mathbb{C}^{2}$
branched over a curve of degree $6$. Therefore, the statement
follows from Theorem~\ref{theorem:CPR}.
\end{proof}

We are to show that the  pencils $|-4K_X|$ and $\mathcal{P}$ are
the only Halphen pencils on $X$.

If the set $\mathbb{CS}(X, \frac{1}{n}\mathcal{M})$ contains the
singular point of type $\frac{1}{5}(1,4,1)$, then
$\mathcal{M}=|-4K_X|$ by
Lemma~\ref{lemma:special-singular-points-with-zero-c}. Therefore,
due to Lemma~\ref{lemma:smooth-points} and
Corollary~\ref{corollary:Ryder-a1},
 we may assume that $\mathbb{CS}(X,
\frac{1}{n}\mathcal{M})=\{O\}$.

Let $\pi:Y\to X$ be the Kawamata blow up at the point $O$ with
weights $(1,4,5)$ and $E$ be its exceptional divisor. Then,
$\mathcal{M}_Y\sim_{\mathbb{Q}}-nK_{Y}$ by
Lemma~\ref{lemma:Kawamata}. Thus, the singularities of the log
pair $(Y, \frac{1}{n}\mathcal{M}_Y)$ are not terminal by
Theorem~\ref{theorem:Noether-Fano} because the divisor $-K_{Y}$ is
nef and big.

The exceptional divisor $E\cong\mathbb{P}(1,4,5)$ contains two
quotient singular points $P$ and $Q$  that are singularities of
types $\frac{1}{4}(1,3,1)$ and $\frac{1}{5}(1,4,1)$ on the
threefold $Y$, respectively. Then, the set $\mathbb{CS}(Y,
\frac{1}{n}\mathcal{M}_Y)$ contains either the singular point $P$
or the singular point $Q$ by Lemma~\ref{lemma:Cheltsov-Kawamata}.

\begin{lemma}
\label{lemma:n-60-center} If the log pair $(Y,
\frac{1}{n}\mathcal{M}_Y)$ is not terminal at the point $Q$, then
$\mathcal{M}=|-4K_{X}|$.
\end{lemma}

\begin{proof}
Suppose that  the set $\mathbb{CS}(Y, \frac{1}{n}\mathcal{M}_Y)$
contains  the point $Q$. Let $\alpha:U\to Y$ be the Kawamata blow
up at the point $Q$ with weights $(1,4,1)$.
 Then,
$\mathcal{M}_U\sim_{\mathbb{Q}}-nK_{U}$ by
Lemma~\ref{lemma:Kawamata}. The linear system $|-4K_{U}|$ is the
proper transform of the pencil $|-4K_{X}|$.  Its base locus
consists of the irreducible curve $C_U$. Since
$$-K_U\cdot C_U=-4K_U^3=-\frac{4}{30}$$ we obtain the identity
$\mathcal{M}=|-4K_{X}|$ from Theorem~\ref{theorem:main-tool}.
\end{proof}

We may assume that $\mathbb{CS}(Y,
\frac{1}{n}\mathcal{M}_Y)=\{P\}$. Let $\beta:W\to Y$ be the
Kawamata blow up at the singular point $P$ with weights $(1,3,1)$
and $F$ be its exceptional divisor. Then,
$$
\mathcal{M}_W\sim_{\mathbb{Q}}\big(\pi\circ\beta\big)^{*}\Big(-nK_{X}\Big)-\frac{n}{9}\beta^{*}\big(E\big)-\frac{n}{4}F\sim_{\mathbb{Q}}-nK_{W}
$$
by Lemma~\ref{lemma:Kawamata}. In addition, we see
$$
\mathcal{P}_W\sim_{\mathbb{Q}}\big(\pi\circ\beta\big)^{*}\Big(-6K_{X}\Big)-\frac{6}{9}\beta^{*}\big(E\big)-\frac{6}{4}F\sim_{\mathbb{Q}}-6K_{W}.
$$

Let $L$ be the unique curve  in the linear system
$|\mathcal{O}_{\mathbb{P}(1,\,4,\,9)}(1)|$ on the surface $E$. The
base locus of $\mathcal{P}_W$ consists of the irreducible curve
$\tilde{C}_W$ and  the irreducible curve $L_W$.

A general surface $D_W$ in  $\mathcal{P}_W$ is normal.  From the
equivalences
\begin{equation*}
\left\{\aligned
&E_W\sim_{\mathbb{Q}}\beta^{*}\big(E\big)-\frac{3}{4}F,\\
&S_W\sim_{\mathbb{Q}}\big(\pi\circ\beta\big)^{*}\Big(-K_{X}\Big)-\frac{1}{9}\beta^{*}\big(E\big)-\frac{1}{4}F,\\
&D_W\sim_{\mathbb{Q}}\big(\pi\circ\beta\big)^{*}\Big(-6K_{X}\Big)-\frac{6}{9}\beta^{*}\big(E\big)-\frac{6}{4}F,\\
\endaligned
\right.
\end{equation*}
we obtain
$$
L_W^{2}=\tilde{C}_W^{2}=-\frac{1}{5}, \ L_W\cdot \tilde{C}_W=0
$$
on the surface $D_W$ because $S_W\cdot D_W=\tilde{C}_W+L_W$ and
$D_W\cdot E_W=6L_W$.

For a general surface  $M_W$ of the linear system $\mathcal{M}_W$,
$$
M_W\Big\vert_{D_W}\equiv-nK_{W}\Big\vert_{D_W}\equiv
nS_W\Big\vert_{D_W}\equiv n\tilde{C}_W+nL_W,
$$
but the intersection form of the irreducible curves $L_W$ and
$\tilde{C}_W$ on the surface $D_W$ is negative-definite.
Therefore, Theorem~\ref{theorem:main-tool} implies that
$\mathcal{M}_W=\mathcal{P}_W$.

Consequently, we have obtained

\begin{proposition}
\label{proposition:n-60} The linear systems $|-4K_{X}|$ and
$\mathcal{P}$ are the only Halphen pencils on $X$.
\end{proposition}

\section{Case $\gimel=69$, hypersurface of degree $28$ in $\mathbb{P}(1,4,6,7,11)$.}%
\index{$\gimel=69$}\label{section:n-69}

The variety $X$ is a general hypersurface of degree $28$ in
$\mathbb{P}(1,4,6,7,11)$ with $-K_{X}^{3}=\frac{1}{66}$. The
singularities of the hypersurface $X$ consist of one singular
point $O$ that are quotient singularities of type
$\frac{1}{11}(1,4,7)$,  one point of type $\frac{1}{6}(1,1,5)$,
and two  points of type $\frac{1}{2}(1,1,1)$.

There is a commutative diagram
$$
\xymatrix{
&U\ar@{->}[d]_{\alpha}&&W\ar@{->}[ll]_{\beta}\ar@{->}[d]^{\eta}&\\%
&X\ar@{-->}[rr]_{\psi}&&\mathbb{P}(1,3,4),&}
$$
where \begin{itemize}

\item $\alpha$ is the Kawamata blow up at the point $O$ with
weights $(1,4,7)$,

\item $\beta$ is the Kawamata blow up with weights $(1,4,3)$ at
the singular point of the variety $U$ contained in the exceptional
divisor of the birational morphism $\alpha$ that is a quotient
singularity of type $\frac{1}{7}(1,4,3)$,

\item $\eta$ is an elliptic fibration.
\end{itemize}

The linear system $|-4K_{X}|$ is a Halphen pencil  on $X$ by
Proposition~\ref{proposition:Halphen-K3}. However, the pencil
$|-4K_{X}|$ is not a unique Halphen pencil on $X$. Indeed, the
hypersurface $X$ can be given by equation
$$
w^{2}z+wf_{17}(x,y,z,t)+f_{28}(x,y,z,t)=0,
$$
where $f_{i}(x,y,z,t)$ is a quasihomogeneous polynomial of degree
$i$. Let $\mathcal{P}$ be the pencil of surfaces cut out on the
hypersurface $X$ by the equations
$$
\lambda x^{6}+\mu z=0,
$$
where $(\lambda : \mu)\in\mathbb{P}^{1}$. We will see that the
 linear system $\mathcal{P}$ is another Halphen pencil.

First of all, if the set $\mathbb{CS}(X, \frac{1}{n}\mathcal{M})$
contains the singular point of type $\frac{1}{6}(1,1,5)$, the
identity $\mathcal{M}=|-4K_X|$ follows from
Lemma~\ref{lemma:special-singular-points-with-zero-c}. Moreover,
due to Lemmas~\ref{lemma:smooth-points},
\ref{lemma:special-singular-points-with-positive-c}, and
Corollary~\ref{corollary:Ryder-a2}, we may assume that
$\mathbb{CS}(X, \frac{1}{n}\mathcal{M})=\{O\}$.

It follows from Theorem~\ref{theorem:Noether-Fano} that the set
$\mathbb{CS}(U, \frac{1}{n}\mathcal{M}_U)$ is not empty.

The exceptional divisor $E\cong\mathbb{P}(1,4,7)$ of the
birational morphism $\alpha$ contains two singular point $P$ and
$Q$ of  $U$ that are quotient singularities of types
$\frac{1}{7}(1,3,4)$ and $\frac{1}{4}(1,1,3)$, respectively. For
the convenience, let $L$ be the unique curve in the linear system
$|\mathcal{O}_{\mathbb{P}(1,\,4,\,7)}(1)|$ on the surface $E$.

\begin{lemma}
\label{lemma:n-69-P6} The set $\mathbb{CS}(U,
\frac{1}{n}\mathcal{M}_U)$ does not contain the point $Q$.
\end{lemma}

\begin{proof}
Suppose that $Q\in\mathbb{CS}(U, \frac{1}{n}\mathcal{M}_U)$. Let
$\pi_Q:U_Q\to U$ be the Kawamata blow up at the point $Q$ with
weights $(1,1,3)$.

Let $\mathcal{D}_{U_Q}$ be the proper transform of the linear
system $|-7K_{X}|$ by the birational morphism $\alpha\circ\pi_Q$.
Its base locus consists of the irreducible curve
$\tilde{C}_{U_Q}$. A general surface $D_{U_Q}$ in
$\mathcal{D}_{U_Q}$ is normal. Moreover, the inequality
$\tilde{C}_{U_Q}^2<0$ holds on the surface $D_{U_Q}$. It implies
the identity  $\mathcal{M}_{U_Q}=\mathcal{D}_{U_Q}$ by
Theorem~\ref{theorem:main-tool} because
$\mathcal{M}_{U_Q}\vert_{D_{U_Q}}\equiv n\tilde{C}_{U_Q}$.
However, the linear system $\mathcal{D}_{U_Q}$ is not a pencil,
which is a contradiction.
\end{proof}

Therefore, by Theorem~\ref{theorem:Noether-Fano} and
Lemma~\ref{lemma:Cheltsov-Kawamata}, the set  $\mathbb{CS}(U,
\frac{1}{n}\mathcal{M}_U)$ consists of the point $P$.

The exceptional divisor $F$ of the birational morphism $\beta$
contains two singular points $P_1$ and $P_2$ of $W$ that are
quotient singularities of types $\frac{1}{3}(1,1,2)$ and
$\frac{1}{4}(1,3,1)$, respectively.

As usual, the set $\mathbb{CS}(W, \frac{1}{n}\mathcal{M}_W)$ must
contain either the point $P_1$ or the point $P_{2}$ by
Theorem~\ref{theorem:Noether-Fano} and
Lemma~\ref{lemma:Cheltsov-Kawamata}.

\begin{lemma}
\label{lemma:n-69-P7} If the set $\mathbb{CS}(W,
\frac{1}{n}\mathcal{M}_W)$ contains the point $P_1$, then
$\mathcal{M}=|-4K_{X}|$.
\end{lemma}

\begin{proof} Suppose that the set $\mathbb{CS}(W,
\frac{1}{n}\mathcal{M}_W)$ contains the point $P_1$. Let
$\pi_1:W_1\to W$ be the Kawamata blow up at the point $P_1$ with
weights $(1,1,2)$. Then, $|-4K_{W_1}|$ is the proper transform of
the pencil $|-4K_{X}|$. Its base locus consists of the irreducible
curve $C_{W_1}$. Because $\mathcal{M}_{W_1}\sim_{\mathbb{Q}}
-nK_{W_1}$, the inequality $-K_{W_1}\cdot C_{W_1}<0$ implies the
identity $\mathcal{M}=|-4K_X|$ by Theorem~\ref{theorem:main-tool}.
\end{proof}

\begin{lemma}
\label{lemma:n-69-P8} If the set $\mathbb{CS}(W,
\frac{1}{n}\mathcal{M}_W)$ contains the point $P_2$, then
$\mathcal{M}=\mathcal{P}$.
\end{lemma}

\begin{proof}
Suppose that the set $\mathbb{CS}(W, \frac{1}{n}\mathcal{M}_W)$
contains the point $P_2$. Let $\pi_2:W_2\to W$ be the Kawamata
blow up at the point $P_2$ with weights $(1,3,1)$ and let $G$ be
its exceptional divisor. Then, the base locus of the pencil
$\mathcal{P}_{W_2}$ consists of the irreducible curves
$\bar{C}_{W_2}$ and $L_{W_2}$.

Let $D$ be a general surface of the pencil $\mathcal{P}$. Then,
$D_{W_2}\sim_{\mathbb{Q}}-6K_{W_2}$, while $S_{W_2}\cdot
D_{W_2}=2L_{W_2}+\bar{C}_{W_2}$ and $E_{W_2}\cdot
D_{W_2}=6L_{W_2}$.

Let us find a divisor $B$ on the threefold $W_2$ such that $B$ is
nef and big but $R\cdot C_{W_2}=R\cdot L_{W_2}=0$. Namely,
consider a divisor $B$ such that
$$
B=(\alpha\circ\beta\circ\pi_2)^{*}(-\lambda K_{X})+
(\beta\circ\pi_2)^{*}(-\mu K_{U})+D_{W_2},%
$$
where $\lambda$ and $\mu$ are nonnegative rational numbers with
$(\lambda, \mu)\ne (0,0)$. Then, the equalities $B\cdot
\bar{C}_{W_2}=B\cdot L_{W_2}=0$ imply that $\lambda=0$ and
$\mu=42$ because
$$
-K_{X}\cdot\alpha\circ\beta\circ\pi_2(\bar{C}_{W_2})=\frac{1}{11},\
-K_{U}\cdot \beta\circ\pi_2(L_{W_2})=\frac{1}{28},\ -K_{U}\cdot
\beta\circ\pi_2(\bar{C}_{W_2})=0,$$
$$
D_{W_2}\cdot L_{W_2}=-\frac{3}{2},\ D_{W_2}\cdot \bar{C}_{W_2}=0.%
$$
The divisor $B$ is nef and big because the divisors $-K_{X}$ and
$-K_{U}$ are nef and big, while $\bar{C}_{W_2}$ and $L_{W_2}$ are
the only curves on the threefold $W_2$ that have negative
intersection with $D_{W_2}$. Let $M$ be a general surface in
$\mathcal{M}_{W_2}$. Then, $B\cdot D_{W_2}\cdot M=0$, which
implies that $\mathcal{M}_{W_2}=\mathcal{P}_{W_2}$ by
Theorem~\ref{theorem:main-tool}.
\end{proof}

\begin{remark}
The surface $D_{W_2}$ is not normal. Indeed, it follows from local
computations that the surface $D_{W_2}$ is singular along the
curve $C_{W_2}$, which is reflected by the fact that $S_{W_2}\cdot
D_{W_2}=2L_{W_2}+C_{W_2}$.
\end{remark}

\begin{K3proposition}
A general surface in the pencil $\mathcal{P}$ is birational to a
smooth  K3 surface.
\end{K3proposition}
\begin{proof}
Let $\Delta$ be the curve on the exceptional divisor $F$ of the
birational morphism $\beta$ defined by intersecting with a general
surface in $\mathcal{P}_W$. Then, the curve $\Delta$ is a curve of
degree $6$ in the surface $\mathbb{P}(1,4,3)$, and hence it is rational. The proper transform
$\Delta_{W_2}$ is a rational curve not contained in the base locus of the pencil $\mathcal{P}_{W_2}$. Therefore, a general
surface in the pencil $\mathcal{P}$ is birational to a smooth K3
surface by Corollary~\ref{corollary:Halphen-K3-rational-curve}.
\end{proof}

Consequently, we have shown
\begin{proposition}
\label{proposition:n-69}The linear systems $|-4K_{X}|$ and
$\mathcal{P}$ are the only pencils on $X$.
\end{proposition}

\section{Case $\gimel=74$, hypersurface of degree $30$ in $\mathbb{P}(1,3,4,10,13)$.}%
\index{$\gimel=74$}\label{section:n-74}

The hypersurface $X$ is defined by a general quasihomogeneous
equation of degree $30$ in  $\mathbb{P}(1,3,4,10,13)$ with
$-K_X^3=\frac{1}{52}$. Its singularities consist of  one quotient
singularity $O$ of type $\frac{1}{13}(1,3,10)$, one quotient
singular point of type $\frac{1}{4}(1,3,1)$, and one quotient
singular point of type $\frac{1}{2}(1,1,1)$.

There is a commutative diagram
$$
\xymatrix{
&U\ar@{->}[d]_{\alpha}&&W\ar@{->}[ll]_{\beta}&&Y\ar@{->}[ll]_{\gamma}\ar@{->}[d]^{\eta}&\\%
&X\ar@{-->}[rrrr]_{\psi}&&&&\mathbb{P}(1,3,4),&}
$$
where \begin{itemize}

\item $\psi$ is the natural projection,

\item $\alpha$ is the Kawamata blow up at the point $O$ with
weights $(1,3,10)$,

\item $\beta$ is the Kawamata blow up with weights $(1,3,7)$ at
the singular point of the variety $U$ that is a quotient
singularity of type $\frac{1}{10}(1,3,7)$ contained in the
exceptional divisor of the birational morphism $\alpha$,

\item $\gamma$ is the Kawamata blow up with weights $(1,3,4)$ at
the singular point of the variety $W$ that is a quotient
singularity of type $\frac{1}{7}(1,3,4)$ contained in the
exceptional divisor of the birational morphism $\beta$,

\item $\eta$ is an elliptic fibration.

\end{itemize}

Proposition~\ref{proposition:Halphen-K3} implies that the linear
system $|-3K_X|$ is a Halphen pencil. There is another Halphen
pencil as follows: The hypersurface $X$ can be given by the
equation
\[w^2z+wf_{17}(x,y,z,t)+f_{30}(x,y,z,t)=0,\]
where $f_i$ is a quasihomogeneous polynomial of degree $i$. Let
$\mathcal{P}$ be the pencil on $X$ given by the equations
$$
\lambda x^4+\mu z=0,
$$
where $(\lambda: \mu)\in \mathbb{P}^1$. Then, the linear system
$\mathcal{P}$ is another Halphen pencil on $X$
(K3-Proposition~\ref{K3proposition:n-74}).

If the set $\mathbb{CS}(X,\frac{1}{n}\mathcal{M})$ contains the
singular point of type $\frac{1}{4}(1, 3, 1)$, then the identity
$\mathcal{M}=|-3K_X|$ follows from
Lemma~\ref{lemma:special-singular-points-with-zero-c}. Therefore,
we may assume that
$$\mathbb{CS}\Big(X,\frac{1}{n}\mathcal{M}\Big)=\Big\{O\Big\}$$
due to Theorem~\ref{theorem:Noether-Fano},
Lemmas~\ref{lemma:smooth-points},
\ref{lemma:special-singular-points-with-positive-c},
 and Corollary~\ref{corollary:Ryder-a1}.

 The exceptional
divisor $E\cong\mathbb{P}(1,3,10)$ of the birational morphism
$\alpha$ contains two quotient singular points $P_1$ and $Q_1$ of
types $\frac{1}{10}(1,3,7)$ and $\frac{1}{3}(1,2,1)$,
respectively. For the convenience, let $L$ be
 the unique curve in the
linear system $|\mathcal{O}_{\mathbb{P}(1,3,10)}(1)|$ on the
surface $E$.

\begin{lemma}
\label{lemma:n-74-first} If the set
$\mathbb{CS}(U,\frac{1}{n}\mathcal{M}_U)$ contains the point
$Q_1$, then $\mathcal{M}=\mathcal{P}$.
\end{lemma}

\begin{proof}
Suppose that the set $\mathbb{CS}(U,\frac{1}{n}\mathcal{M}_U)$
contains the point $Q_1$. Let $\alpha_1:W_1\to U$ be the Kawamata
blow up at the point $Q_1$ with weights $(1,2,1)$.  By
Lemma~\ref{lemma:Kawamata}, we have
$\mathcal{M}_{W_1}\sim_{\mathbb{Q}}-K_{W_1}.$

The base locus of  the pencil $\mathcal{P}_{W_1}$ consists of the
curves $\bar{C}_{W_1}$ and $L_{W_1}$.

Let $D_{W_1}$ be a general surface in $\mathcal{P}_{W_1}$. We see
then that
\[S_{W_1}\cdot D_{W_1}=\bar{C}_{W_1}+2L_{W_1}, \
\  E_{W_1}\cdot D_{W_1}=4L_{W_1}.\] Note that the surface
$D_{W_1}$ is not normal and

\[\left\{\aligned & D_{W_1}\sim_{\mathbb{Q}}
(\alpha\circ\alpha_1)^*(-4K_X)-\frac{4}{13}\alpha_1^*(E)-\frac{4}{3}E_1,\\
& S_{W_1}\sim_{\mathbb{Q}}
(\alpha\circ\alpha_1)^*(-K_X)-\frac{1}{13}\alpha_1^*(E)-\frac{1}{3}E_1,\\
&E_{W_1}\sim_{\mathbb{Q}} \alpha_1^*(E)-\frac{2}{3}E_1,\\
\endaligned
\right.\] where $E_1$ is the exceptional divisor of the birational
morphism $\alpha_1$.

One can easily see that  $$-K_U\cdot L=\frac{1}{30}, \ \ -K_U\cdot
\bar{C}_U=0,\ \ D_{W_1}\cdot L_{W_1}=-\frac{6}{5}, \  \
D_{W_1}\cdot \bar{C}_{W_1}=0.$$ Because $-K_U$ is nef and big and
$L_{W_1}$ is the only curve intersecting  $D_{W_1}$ negatively,
$B:=36\alpha_1^*(-K_{U})+D_{W_1}$ is also a nef and big divisor
with $B\cdot L_{W_1}=B\cdot \bar{C}_{W_1}=0$. Let $M$ be a general
surface in $\mathcal{M}_{W_1}$. We then obtain $B\cdot M\cdot
D_{W_1}=0$, which implies that $\mathcal{M}=\mathcal{P}$.
\end{proof}
\begin{K3proposition}\label{K3proposition:n-74}
A general surface in the pencil $\mathcal{P}$ is birational to a
smooth  K3 surface.
\end{K3proposition}
\begin{proof}
We use the same notations in the proof of Lemma~\ref{lemma:n-74-first}.
Let $\Delta$ be the curve on the exceptional divisor $E_1$ defined by intersecting with the surface $D_{W_1}$.
Then, the curve $\Delta$ is a curve of degree $4$ on the surface $\mathbb{P}(1,2,1)$.

If the curve $\Delta$ is singular, then it is a rational curve on the surface $D_{W_1}$.
Suppose that the curve $\Delta$ is smooth. Then, it cannot pass through the singular point of the surface $E_1$, and hence it is
an elliptic curve. Then, the argument in the proof of K3-Proposition~\ref{K3proposition:n-48} gives a contradiction.
Therefore, the surface $D_{W_1}$ must have a rational curve not contained in the base locus of the pencil $\mathcal{P}_{W_1}$.
Then, Corollaries~\ref{corollary:Halphen-K3}
and \ref{corollary:Halphen-K3-rational-curve} complete the proof.
\end{proof}

Due to Lemma~\ref{lemma:n-74-first}, we may assume that
$\mathbb{CS}(U,\frac{1}{n}\mathcal{M}_U)=\{P_1\}.$

 The
exceptional divisor $F\cong\mathbb{P}(1,3,7)$ of the birational
morphism $\beta$ contains two singular points $P_2$ and $Q_2$ that
are quotient singularities of types $\frac{1}{7}(1,3,4)$ and
$\frac{1}{3}(1,2,1)$, respectively. We let $\bar{L}$ be the unique
curve in the linear system $|\mathcal{O}_{\mathbb{P}(1,3,7)}(1)|$
on surface $F$.

The set $\mathbb{CS}(W,\frac{1}{n}\mathcal{M}_W)$ contains either
the point $P_2$ or the point $Q_2$.

\begin{lemma}
\label{lemma:n-74-second}
 If the set $\mathbb{CS}(W,\frac{1}{n}\mathcal{M}_W)$ contains
the point $Q_2$, then $\mathcal{M}=\mathcal{P}$.
\end{lemma}

\begin{proof} Suppose that the set $\mathbb{CS}(W,\frac{1}{n}\mathcal{M}_W)$ contains
the point $Q_2$.  Let $\beta_2 : W_2\to W$ be the Kawamata blow up
at the point $Q_2$ with weights $(1,2,1)$ and $F_2$ be its
exceptional divisor. The base locus of the pencil
$\mathcal{P}_{W_2}$ consists of three irreducible curves
$\bar{C}_{W_2}$, $L_{W_2}$, and $\bar{L}_{W_2}$.

Let $D_{W_2}$ be a general surface in the pencil
$\mathcal{P}_{W_2}$. Then,  we have
$$
\left\{\aligned
&D_{W_2}\sim_{\mathbb{Q}}(\alpha\circ\beta\circ\beta_2)^*(-4K_X)-\frac{4}{13}(\beta\circ\beta_2)^*(E)
-\frac{4}{10}\beta_2^*(F)-\frac{4}{3}F_2,\\
&S_{W_2}\sim_{\mathbb{Q}}(\alpha\circ\beta\circ\beta_2)^*(-K_X)-\frac{1}{9}(\beta\circ\beta_2)^*(E)
-\frac{1}{7}\beta_2^*(F)-\frac{1}{3}F_2,\\
&E_{W_2}\sim_{\mathbb{Q}}(\beta\circ\beta_2)^*(E)
-\frac{7}{10}\beta_2^*(F)-\frac{1}{3}F_2,\\
&F_{W_2}\sim_{\mathbb{Q}}\beta_2^*(F)-\frac{2}{3}F_2.\\
\endaligned
\right.
$$
Furthermore,
\[S_{W_2}\cdot D_{W_2}=\bar{C}_{W_2}+2L_{W_2}+\bar{L}_{W_2},\ \
E_{W_2}\cdot D_{W_2}=4L_{W_2}, \ \ F_{W_2}\cdot
D_{W_2}=4\bar{L}_{W_2}.\] First of all, we can easily check that
$$D_{W_2}\cdot L_{W_2}=-\frac{2}{3}, \ \ D_{W_2}\cdot
\bar{L}_{W_2}=-\frac{8}{7}, \ \ D_{W_2}\cdot \bar{C}_{W_2}=0.$$ We
then consider the divisor
$$B_2=20(\beta\circ\beta_2)^*(-K_U)+24\beta_2^*(-K_W)+D_{W_2}$$
on $W_2$. Because
$$
\left\{\aligned
&-K_U\cdot (\beta\circ\beta_2)(L_{W_2})=\frac{1}{30}, \ -K_U\cdot (\beta\circ\beta_2)(\bar{C}_{W_2})=0\\
&-K_W\cdot\beta_2(L_{W_2})=0, \ -K_W\cdot\beta_2(\bar{L}_{W_2})=\frac{1}{21}, \ -K_W\cdot\beta_2(\bar{C}_{W_2})=0\\
\endaligned
\right.
$$
we see that $B_2\cdot L_{W_2}=B_2\cdot \bar{L}_{W_2}=B_2\cdot
\bar{C}_{W_2}=0$.

Because $-K_U$ and $-K_W$ are nef and big and $L_{W_2}$ and
$\bar{L}_{W_2}$ are the only curves intersecting $D_{W_2}$
negatively, the divisor $B_2$ is also nef and big. Let $M_2$ be a
general surface in $\mathcal{M}_{W_2}$. We then obtain $B_2\cdot
M_2\cdot D_{W_2}=0$, which implies that $\mathcal{M}=\mathcal{P}$.
\end{proof}

Due to the lemma above, we may assume that
$\mathbb{CS}(W,\frac{1}{n}\mathcal{M}_W)=\{P_2\}$. The exceptional
divisor $G\cong\mathbb{P}(1,3,4)$ of the birational morphism
$\gamma$ contains two singular points $P_3$ and $Q_3$ that are
quotient singularities of types $\frac{1}{4}(1,3,1)$ and
$\frac{1}{3}(1,2,1)$, respectively. Again, the set
$\mathbb{CS}(Y,\frac{1}{n}\mathcal{M}_Y)$ contains either the
point $P_3$ or the point $Q_3$.

\begin{lemma}
\label{lemma:n-74-third} The set
$\mathbb{CS}(Y,\frac{1}{n}\mathcal{M}_Y)$ cannot contain the point
$Q_3$.
\end{lemma}

\begin{proof}
Suppose that the set $\mathbb{CS}(Y,\frac{1}{n}\mathcal{M}_Y)$
contains the point $Q_3$. Let $\gamma_3:W_3\to Y$ be the Kawamata
blow up at the point $Q_3$ with weights $(1,2,1)$ and let $G_3$ be
the exceptional divisor of $\gamma_3$.
The base locus of the pencil $\mathcal{P}_{W_3}$ consists of three
irreducible curves $\bar{C}_{W_3}$, $L_{W_3}$, and
$\bar{L}_{W_3}$. Let $D_{W_3}$ be a general surface in the pencil
$\mathcal{P}_{W_3}$. We see
$$
\left\{\aligned
&D_{W_3}\sim_{\mathbb{Q}}(\alpha\circ\beta\circ\gamma\circ\gamma_3)^*(-4K_X)
-\frac{4}{13}(\beta\circ\gamma\circ\gamma_3)^*(E)-\frac{4}{10}(\gamma\circ\gamma_3)^*(F)
-\frac{4}{7}\gamma_3^*(G)-\frac{1}{3}G_3\\
&S_{W_3}\sim_{\mathbb{Q}}(\alpha\circ\beta\circ\gamma\circ\gamma_3)^*(-K_X)
-\frac{1}{13}(\beta\circ\gamma\circ\gamma_3)^*(E)-\frac{1}{10}(\gamma\circ\gamma_3)^*(F)
-\frac{1}{7}\gamma_3^*(G)-\frac{1}{3}G_3,\\
&E_{W_3}\sim_{\mathbb{Q}}(\beta\circ\gamma\circ\gamma_3)^*(E)
-\frac{7}{10}(\gamma\circ\gamma_3)^*(F),\\
&F_{W_3}\sim_{\mathbb{Q}}(\gamma\circ\gamma_3)^*(F)-\frac{4}{7}\gamma_3^*(G)-\frac{1}{3}G_3,\\
\endaligned
\right.
$$
Furthermore, \[S_{W_3}\cdot
D_{W_3}=\bar{C}_{W_3}+2L_{W_3}+\bar{L}_{W_3}, \ \ E_{W_3}\cdot
D_{W_3}=4L_{W_3}, \  \ F_{W_3}\cdot D_{W_3}=4\bar{L}_{W_3}.\] From
this, one can obtain
\[D_{W_3}\cdot L_{W_3}= 0, \ \ D_{W_3}\cdot \bar{L}_{W_3}=-\frac{1}{24}, \ \ D_{W_3}\cdot \bar{C}_{W_3}= -\frac{1}{8}.\]
Using the same method as in the previous lemma with nef and big
divisors $-K_X$, $-K_U$, and $-K_W$, we can find a nef and big
divisor $B_3$ on $W_3$ such that $B_3\cdot D_{W_3}\cdot M_3=0$,
where $M_3$ is a general surface in $\mathcal{M}_{W_3}$, which
implies that $\mathcal{M}=\mathcal{P}$. However it is a
contradiction because $D_3\not \sim_{\mathbb{Q}}-4K_{W_3}$.
\end{proof}

Therefore, the set  $\mathbb{CS}(Y,\frac{1}{n}\mathcal{M}_Y)$
contains the point $P_3$. Let $\delta:V\to Y$ be the Kawamata blow
up at the point $P_3$ with weights $(1,3,1)$. The pencil $|-3K_V|$
is the proper transform of the pencil $|-3K_X|$. It has only one
irreducible base curve $C_V$. For a general surface $M$ in
$\mathcal{M}_V$, we have
\[M\Big\vert_{D_V}\equiv-nK_V\Big\vert_{D_V}\equiv nC_V.\]
Therefore, the inequality $-K_C\cdot C_V<0$  implies that
$\mathcal{M}=|-3K_X|$ by Theorem~\ref{theorem:main-tool}.

Consequently, we have proved
\begin{proposition}
\label{proposition:n-74} The linear systems $|-3K_X|$ and
$\mathcal{P}$ are the only Halphen pencils on $X$.
\end{proposition}

\section{Case $\gimel=76$, hypersurface of degree $30$ in $\mathbb{P}(1,5,6,8,11)$.}%
\index{$\gimel=76$}\label{section:n-76}

The variety $X$ is a general hypersurface of degree $25$ in
$\mathbb{P}(1,5,6,8,11)$ with $-K_{X}^{3}=\frac{1}{88}$. The
singularities of $X$ consist of one point that is a quotient
singularity of type $\frac{1}{2}(1,1,1)$, one point $P$ that is a
quotient singularity of type $\frac{1}{8}(1,5,3)$, and one point
$Q$ that is a quotient singularity of type $\frac{1}{11}(1,5,6)$.

There is a commutative diagram
$$
\xymatrix{
&&Y\ar@{->}[dl]_{\beta_{Q}}\ar@{->}[dr]^{\beta_{P}}\ar@{->}[drrrrrrr]^{\eta}&&&&&&&\\
&U_{P}\ar@{->}[dr]_{\alpha_{P}}&&U_{Q}\ar@{->}[dl]^{\alpha_{Q}}&&&&&&\mathbb{P}(1,5,6),\\
&&X\ar@{-->}[rrrrrrru]_{\psi}&&&&&&&}
$$
where \begin{itemize}

\item $\psi$ is the natural projection,

 \item $\alpha_{P}$ is the Kawamata blow up at the point
$P$ with weights $(1,5,3)$,

\item $\alpha_{Q}$ is the Kawamata blow up at the point $Q$ with
weights $(1,5,6)$,

\item $\beta_{Q}$ is the Kawamata blow up with weights $(1,5,6)$
at the point whose image to $X$ is the point $Q$,

\item $\beta_{P}$ is the Kawamata blow up with weights $(1,5,3)$
at the point whose image to $X$ is the point $P$,

\item $\eta$ is an elliptic fibration.

\end{itemize}

The threefold $X$ can be given by the equation
$$
t^{3}z+t^{2}f_{14}\big(x,y,z,w\big)+t\Big(w^{2}+f_{22}\big(x,y,z\big)\Big)+
f_{30}\big(x,y,z,w\big)=0,
$$
where $f_{i}(x,y,z,t)$ is a general quasihomogeneous polynomial of
degree $i$. Let $\mathcal{P}$ be the pencil consisting of surfaces
 cut out on the threefold $X$ by the equations
$$
\lambda x^{6}+\mu z=0,
$$
where $(\lambda:\mu)\in\mathbb{P}^{1}$.
\begin{proposition}
\label{proposition:n-76} The linear systems $|-5K_{X}|$ and
$\mathcal{P}$ are the only Halphen pencils on $X$.
\end{proposition}

 By
Lemmas~\ref{lemma:special-singular-points-with-positive-c},
\ref{lemma:smooth-points}, and Corollary~\ref{corollary:Ryder-a1},
we have $\mathbb{CS}(X, \frac{1}{n}\mathcal{M})\subset\{P, Q\}$.

 The exceptional divisor $E_P$ of the birational morphism
 $\alpha_p$ contains two quotient singular points $P_1$ and $P_2$ of types
$\frac{1}{5}(1,2,3)$ and $\frac{1}{3}(1,2,1)$, respectively. The
divisor $-K_{U_{P}}$ is nef and big. Thus, the set
$\mathbb{CS}(U_P, \frac{1}{n}\mathcal{M}_{U_P})$ is not empty by
Theorem~\ref{theorem:Noether-Fano}.
\begin{lemma}
\label{lemma:n-76-P1} If the set $\mathbb{CS}(U_P,
\frac{1}{n}\mathcal{M}_{U_P})$ contains $P_1$, then
$\mathcal{M}=\mathcal{P}$.
\end{lemma}

\begin{proof}

Suppose that the set $\mathbb{CS}(U_P,
\frac{1}{n}\mathcal{M}_{U_P})$ contains the point $P_1$. Let
$\beta_{1}:W_{1}\to U_{P}$ be the Kawamata blow up at the point
$P_1$ with weights $(1,2,3)$. Then,
$\mathcal{M}_{W_1}\sim_{\mathbb{Q}} -nK_{W_{1}}$ by
Lemma~\ref{lemma:Kawamata}. Also, one can see that
$\mathcal{P}_{W_1}\sim_{\mathbb{Q}} -6K_{W_{1}}$ and the base
locus of $\mathcal{P}_{W_1}$ consists of the irreducible curve
$\bar{C}_{W_1}$. The inequality $-K_{W_1}\cdot \bar{C}_{W_1}<0$
implies $\mathcal{M}=\mathcal{P}$ by
Theorem~\ref{theorem:main-tool}.
\end{proof}

\begin{K3proposition}\label{K3proposition:n-76}
A general surface in the pencil $\mathcal{P}$ is birational to a
K3 surface.
\end{K3proposition}
\begin{proof}
We use the same notation in the proof of
Lemma~\ref{lemma:n-76-P1}. The exceptional divisor $E_P$ is
isomorphic to the weighted projective space $\mathbb{P}(1,3,5)$
and the curve $\Delta$ defined by the intersection of $E_P$ with a
general member in $\mathcal{P}_{U_P}$ has degree $6$ on $E_P$.
Therefore, the curve $\Delta_{W_1}$ is a rational curve not contained in the base locus of the pencil $\mathcal{P}_{W_1}$, and hence a general surface in the
pencil $\mathcal{P}$ is birational to a smooth  K3 surface by
Corollary~\ref{corollary:Halphen-K3-rational-curve}.
\end{proof}

\begin{lemma}
\label{lemma:n-76-P2} If the set $\mathbb{CS}(U_P,
\frac{1}{n}\mathcal{M}_{U_P})$ contains $P_2$, then
$\mathcal{M}=|-5K_X|$.
\end{lemma}
\begin{proof}
Suppose that the set $\mathbb{CS}(U_P,
\frac{1}{n}\mathcal{M}_{P})$ contains the point $P_2$. Let
$\beta_{2}:W_{2}\to U_{P}$ be the Kawamata blow up of the point
$P_2$ with weights $(1,2,1)$. Then,
$\mathcal{M}_{W_2}\sim_{\mathbb{Q}} -nK_{W_{2}}$ by
Lemma~\ref{lemma:Kawamata}. The pencil $|-5K_{W_2}|$ is the proper
transform of the pencil $|-5K_X|$ and its base locus consists of
the irreducible curve $C_{W_2}$. Then, the inequality
$-K_{W_2}\cdot C_{W_2}<0$ implies $\mathcal{M}=|-5K_X|$ by
Theorem~\ref{theorem:main-tool}.
\end{proof}

Meanwhile, the exceptional divisor $E_Q$ of the birational
morphism
 $\alpha_Q$ contains two quotient singular points $Q_1$ and $Q_2$ of types
$\frac{1}{5}(1,4,1)$ and $\frac{1}{6}(1,5,1)$, respectively. The
divisor $-K_{U_{Q}}$ is nef and big. Thus, the set
$\mathbb{CS}(U_Q, \frac{1}{n}\mathcal{M}_{U_Q})$ is not empty by
Theorem~\ref{theorem:Noether-Fano}.

\begin{lemma}
\label{lemma:n-76-Q1} If the set $\mathbb{CS}(U_Q,
\frac{1}{n}\mathcal{M}_{U_Q})$ contains $Q_1$, then
$\mathcal{M}=\mathcal{P}$.
\end{lemma}
\begin{proof}
The proof is similar to that of Lemma~\ref{lemma:n-76-P1}.
\end{proof}

\begin{lemma}
\label{lemma:n-76-Q2} If the set $\mathbb{CS}(U_P,
\frac{1}{n}\mathcal{M}_{U_P})$ contains $Q_2$, then
$\mathcal{M}=|-5K_X|$.
\end{lemma}
\begin{proof}
The proof is similar to that of Lemma~\ref{lemma:n-76-P2}.
\end{proof}

Therefore, for the proof of Proposition~\ref{proposition:n-76}, we
may assume that $\mathbb{CS}(X, \frac{1}{n}\mathcal{M})=\{P, Q\}$
by the previous lemmas.

Each member of the pencil $\mathcal{M}_Y$ is contracted to a curve
by  the elliptic fibration $\eta$ by Lemma~\ref{lemma:Kawamata}
but the set $\mathbb{CS}(Y, \frac{1}{n}\mathcal{M}_Y)$ is not
empty by Theorem~\ref{theorem:Noether-Fano}. Hence, it follows
from Lemma~\ref{lemma:Cheltsov-Kawamata} that the set
$\mathbb{CS}(Y, \frac{1}{n}\mathcal{M}_Y)$ contains a singular
point of the threefold $Y$ that is contained either in the
exceptional divisor of the birational morphism  $\alpha_{P}$ or in
the exceptional divisor of the birational morphism  $\beta_{Q}$.
Then, Lemmas~\ref{lemma:n-76-P1}, \ref{lemma:n-76-P2},
\ref{lemma:n-76-Q1}, and \ref{lemma:n-76-Q2} conclude the proof of
Proposition~\ref{proposition:n-76}.

\section{Case $\gimel=79$, hypersurface of degree $33$ in
$\mathbb{P}(1,3,5,11,14)$.}\index{$\gimel=79$}
\label{section:n-79}

The threefold $X$ is a general hypersurface of degree $33$ in
$\mathbb{P}(1,3,5,11,14)$ with $-K_X^3=\frac{1}{70}$. It has two
singular points. One is  a quotient singularity of type
$\frac{1}{5}(1,1,4)$ and the other is a quotient singularity $O$
of type $\frac{1}{14}(1,3,11)$. The hypersurface $X$ can be given
by the equation
$$
w^{2}z+wf_{19}\big(x,y,z,t\big)+f_{33}\big(x,y,z,t\big)=0,
$$
where  $f_{i}$ is a qua\-si\-ho\-mo\-ge\-ne\-ous polynomial of
degree $i$. Let $\mathcal{P}$ be the pencil cut out on $X$ by
$$\lambda x^{5}+\mu z=0,$$ where $(\lambda:\mu)\in\mathbb{P}^{1}$.

There is a commutative diagram
$$
\xymatrix{
&U\ar@{->}[d]_{\alpha}&&W\ar@{->}[ll]_{\beta}&&Y\ar@{->}[ll]_{\gamma}\ar@{->}[d]^{\eta}&\\%
&X\ar@{-->}[rrrr]_{\psi}&&&&\mathbb{P}(1,3,5),&}
$$
where \begin{itemize} \item $\psi$ is the natural projection,
\item $\alpha$ is the Kawamata blow up at the point $O$ with
weights $(1,3,11)$,

\item $\beta$ is the Kawamata blow up with weights $(1,3,8)$ at
the singular point of type $\frac{1}{11}(1,3,8)$ contained in the
exceptional divisor of the birational morphism $\alpha$,

\item $\gamma$ is the Kawamata blow up  with weights $(1,3,5)$ at
the singular point of type $\frac{1}{8}(1,3,5)$ contained in the
exceptional divisor of the birational morphism $\beta$, \item
$\eta$ is an elliptic fibration.
\end{itemize}

If the set $\mathbb{CS}(X, \frac{1}{n}\mathcal{M}_X)$ contains the
singular point of type $\frac{1}{5}(1,1,4)$, then
$\mathcal{M}=|-3K_X|$ by
Lemma~\ref{lemma:special-singular-points-with-zero-c}. Therefore,
we may~assume~that
$$
\mathbb{CS}\Big(X, \frac{1}{n}\mathcal{M}\Big)=\Big\{O\Big\}
$$
due to Lemma~\ref{lemma:smooth-points} and
Corollary~\ref{corollary:Ryder-a1}.

The exceptional divisor $E$ of $\alpha$ contains two quotient
singular points $P$ and $Q$  of types $\frac{1}{3}(1,1,2)$ and
$\frac{1}{11}(1,3,8)$, res\-pec\-tive\-ly.

\begin{lemma}
\label{lemma:n-79-points-P4} The set $\mathbb{CS}(U,
\frac{1}{n}\mathcal{M}_U)$ consists of the point $Q$.
\end{lemma}

\begin{proof}
 Suppose that $\mathbb{CS}(U,
\frac{1}{n}\mathcal{M}_U)\ne\{Q\}$. Then, the set $\mathbb{CS}(U,
\frac{1}{n}\mathcal{M}_U)$ contains the point $P$. Let
$\pi_P:U_P\to U$ be the Kawamata blow up at $P$ with weights
$(1,1,2)$. Then, $\mathcal{M}_{U_P}\sim_{\mathbb{Q}}-nK_{U_P}$ by
Lemma~\ref{lemma:Kawamata}.

Let $\mathcal{D}$ be the proper transforms of $|-11K_{X}|$ on the
threefold $U_P$ and $D$ be a general surface of the linear system
$\mathcal{D}$. Then, the base locus of the linear system
$\mathcal{D}$ does not contain curves, which implies that the
divisor $D$ is nef. Thus, we obtain an absurd inequality
$$
0\leq D\cdot M_{1}\cdot M_{2}=-\frac{n^{2}}{5},
$$
where $M_{1}$ and $M_{2}$ are general surfaces of the pencil
$\mathcal{M}_{U_P}$.
\end{proof}

The exceptional divisor $F$ of the birational morphism $\beta$
contains two singular points $Q_1$ and $Q_2$ that are quotient
singularities of types $\frac{1}{3}(1,1,2)$ and
$\frac{1}{8}(1,3,5)$ respectively.

\begin{lemma}
\label{lemma:n-79-points-P5} If the set $\mathbb{CS}(W,
\frac{1}{n}\mathcal{M}_W)$ contains the point $Q_1$, then
$\mathcal{M}=\mathcal{P}$.
\end{lemma}

\begin{proof} Suppose that the set $\mathbb{CS}(W,
\frac{1}{n}\mathcal{M}_W)$ contains the point $Q_1$. Let
$\pi:W_1\to W$ be the Kawamata blow up of $Q_1$ with weights
$(1,1,2)$ and $G$ be its exceptional divisor. Then,
$\mathcal{M}_{W_1}\sim_{\mathbb{Q}}-nK_{W_1}$ by
Lemma~\ref{lemma:Kawamata}.

Let $\mathcal{L}$ be the linear system on the hypersurface $X$ cut
out by
$$
\lambda_{0}x^{30}+\lambda_{1}y^{10}+\lambda_{2}z^{6}+\lambda_{3}t^{2}x^{8}+\lambda_{4}t^{2}y^{2}x^{2}+\lambda_{5}ty^{6}x+\lambda_{6}wtz=0
$$
where $(\lambda_{0}:\cdots:\lambda_{6})\in\mathbb{P}^{6}$. Then,
the base locus of $\mathcal{L}$ does not contain curves.
 Then, it follows from simple calculations that the base locus of the linear system $\mathcal{L}_{W_1}$ does
not contain any curve and for a general surface $B$ in
$\mathcal{L}$, we obtain
$$
B_{W_1}\sim_{\mathbb{Q}}\big(\alpha\circ\beta\circ\pi\big)^{*}\Big(-30K_{X}\Big)
-\frac{30}{14}\big(\beta\circ\pi\big)^{*}\big(E\big)-\frac{8}{11}\pi^{*}\big(F\big)-\frac{2}{3}G.
$$
 In
particular, the divisor $B_{W_1}$ is nef and big.

Let $M$ be a general surface of the pencil $\mathcal{M}_{W_1}$ and
$D$ be a general surface of the linear system $|-5K_{W_1}|$. Then,
$B_{W_1}\cdot M\cdot D=0$, which implies that
$\mathcal{M}=\mathcal{P}$ by Theorem~\ref{theorem:main-tool}
because the linear system $|-5K_{W_1}|$ is the proper transform of
the pencil $\mathcal{P}$.
\end{proof}
\begin{K3proposition}\label{K3proposition:n-79}
A general surface in the pencil $\mathcal{P}$ is birational to a
smooth  K3 surface.
\end{K3proposition}
\begin{proof}
We use the same notations in the proof of
Lemma~\ref{lemma:n-79-points-P5}. The surface $G$ is isomorphic to
the projective space $\mathbb{P}(1,1,2)$. Let $T$ be a general
surface in the pencil $\mathcal{P}$ and let $\Delta=G\cdot
T_{W_1}$. Then, by simple calculation, we see that the curve
$\Delta$ on $G$ is defined by the equation
\[\epsilon_1 \bar{x}^5+\epsilon_2 \bar{x}^2\bar{y}\bar{t}+\epsilon_3\bar{x}\bar{y}^2\bar{t}+\epsilon_4
\bar{y}\bar{t}^2=0\subset
\mathrm{Proj}(\mathbb{C}[\bar{x},\bar{y},\bar{t}])=\mathbb{P}(1,1,2),\]
where each $\epsilon_i$ is a general complex number. It has two
nodes at the points $(1:0:0)$ and $(0:1:0)$. But it is smooth at
the point $(0:0:1)$ which is a $\mathbb{A}_1$ singular point of
the surface $G$. Let $\tilde{G}$ be the blow up of the surface $G$
at these three points. The genus of the normalization
$\tilde{\Delta}$ of the curve $\Delta$ is
\[ p_g(\tilde{\Delta})=\frac{(K_{\tilde{G}}+
\tilde{\Delta})\cdot \tilde{\Delta}}{2} +1=\frac{(K_G+
\Delta)\cdot \Delta-\frac{1}{2}}{2} +1-2=0,\] and hence the curve
$\Delta$ is a rational curve not contained in the base locus of the pencil $\mathcal{P}_{W_1}$.
Therefore, Corollary~\ref{corollary:Halphen-K3-rational-curve}
completes the proof.
\end{proof}

We may assume that $\mathbb{CS}(W,
\frac{1}{n}\mathcal{M}_W)=\{Q_2\}$ due to
Theorem~\ref{theorem:Noether-Fano} and
Lemma~\ref{lemma:Cheltsov-Kawamata}. Let $O_1$ and $O_2$ be the
quotient  singular points of the threefold $Y$ contained in the
exceptional divisor of $\gamma$ that are of types
$\frac{1}{3}(1,1,2)$ and $\frac{1}{5}(1,3,2)$, respectively. Then,
$$
\varnothing\ne \mathbb{CS}\Big(Y, \frac{1}{n}\mathcal{M}_Y\Big)\subset\Big\{O_1, O_2\Big\}%
$$
by Theorem~\ref{theorem:Noether-Fano},
Lemmas~\ref{lemma:Kawamata}, and \ref{lemma:Cheltsov-Kawamata}.
The proof of Lemma~\ref{lemma:n-80-81-91-P7} implies that the
$\mathbb{CS}(Y, \frac{1}{n}\mathcal{M}_Y)$ does not contain the
point $O_1$. Now, the proofs of Lemma~\ref{lemma:n-80-81-91-P6}
shows $\mathcal{M}=|-3K_{X}|$.

Therefore, we have proved
\begin{proposition}
\label{proposition:n-79} The linear systems $|-3K_{X}|$ and
$\mathcal{P}$ are the only Halphen pencils on $X$.
\end{proposition}

\section{Cases $\gimel=84$ and  $93$.} \label{section:n-84-93}

In the case of $\gimel=84$,\index{$\gimel=84$} the threefold $X$
is a general hypersurface of degree $36$ in
$\mathbb{P}(1,7,8,9,12)$ with $-K_X^3=\frac{1}{168}$. Its
singularities  consist of one quotient singular point of type
$\frac{1}{4}(1,3,1)$, one quotient singular point of type
$\frac{1}{3}(1,1,2)$, one quotient singular point of type
$\frac{1}{7}(1,2,5)$, and one quotient singular point of type
$\frac{1}{8}(1,7,1)$.

In the case of  $\gimel=93$,\index{$\gimel=93$} the threefold $X$
is a general hypersurface of degree $50$ in $\mathbb{P}(1,7,8,10,
25)$ with $-K_X^3=\frac{1}{280}$. It has  one quotient singular
point of type $\frac{1}{2}(1,1,1)$, one quotient singular point of
type $\frac{1}{5}(1,2,3)$, one quotient singular point of type
$\frac{1}{7}(1,3,4)$, and one quotient singular point of type
$\frac{1}{8}(1,7,1)$.

In both cases, the threefold $X$ cannot be birationally
transformed to an elliptic fibration (\cite{ChPa05}). However, it
can be rationally fibred by K3 surfaces.

The threefold $X$ can be given by the equation
$$
y^{\frac{d-8}{7}}z+\sum_{i=0}^{\frac{d-15}{7}}y^{i}f_{d-7i}\big(x,z,t,w\big)=0,
$$
where $d$ is the degree of $X$ and $f_{i}$ is quasihomogeneous
polynomial of degree $i$. Let $\mathcal{P}$ be the pencil that is
cut out on $X$ by
$$\lambda x^{8}+\mu z=0,$$ where $(\lambda:\mu)\in\mathbb{P}^{1}$.

\begin{K3proposition}\label{K3proposition:n-84-93}
A general surface in the pencil $\mathcal{P}$ is birational to a
smooth  K3 surface.
\end{K3proposition}
\begin{proof}
If $\gimel=84$, then a general surface in the pencil $\mathcal{P}$
is a compactification of a quartic surface in $\mathbb{C}^{3}$ and
must be birational to a smooth K3 surface by
Theorem~\ref{theorem:CPR}. If $\gimel=93$, then a general surface
in the pencil $\mathcal{P}$ is a compactification of a double
cover of $\mathbb{C}^{2}$ ramified along a sextic curve that must
be birational to a smooth K3 surface.
\end{proof}

\begin{proposition}
\label{proposition:n-84-93} If  $\gimel\in\{84, 93\}$, then the
pencils $|-a_{1}K_{X}|$ and $\mathcal{P}$ are unique Halphen
pencils on $X$.
\end{proposition}
\begin{proof}
Theorem~\ref{theorem:Noether-Fano},
Lemmas~\ref{lemma:smooth-points},
\ref{lemma:special-singular-points-with-positive-c},
\ref{lemma:special-singular-points-with-zero-c}, and
Corollary~\ref{corollary:Ryder-a1} immediately imply the result.
\end{proof}

\section{Case $\gimel=95$, hypersurface of degree $66$ in $\mathbb{P}(1,5,6,22,33)$.}%
\index{$\gimel=95$}\label{section:n-95}

Let $X$ be a general hypersurface of degree $66$ in
$\mathbb{P}(1,5,6,22,33)$ with $-K_X^3=\frac{1}{330}$. Its
singularities consist of a quotient singular point of type
$\frac{1}{2}(1,1,1)$, a quotient singular point of type
$\frac{1}{3}(1,2,1)$, a quotient singular point $P$ of type
$\frac{1}{5}(1,2,3)$, and a quotient singular point $O$ of type
$\frac{1}{11}(1,5,6)$.

 We have an elliptic fibration as follows:
$$
\xymatrix{
&&&U\ar@{->}[lld]_{\pi}\ar@{->}[rrd]^{\eta}&&&\\%
&X\ar@{-->}[rrrr]_{\psi}&&&&\mathbb{P}(1,5,6)&}
$$
where \begin{itemize} \item $\psi$ is the natural projection,
\item $\pi$ is the Kawamata blow up at the point $O$ with weights
$(1, 5, 6)$, \item $\eta$ is an elliptic fibration.\end{itemize}

The threefold $X$ can be given by the equation
$$
y^{12}z+\sum_{i=0}^{11}y^{i}f_{66-5i}\big(x,z,t,w\big)=0
$$
in $ \mathbb{P}(1,5,6,22,33)$, where $f_{i}$ is a
qua\-si\-ho\-mo\-ge\-ne\-ous polynomial of degree $i$. Let
$\mathcal{P}$ be the pencil on the threefold $X$ that is cut out
by the pencil $\lambda x^{6}=\mu z$, where
$(\lambda:\mu)\in\mathbb{P}^{1}$.  Then,
Lemma~\ref{lemma:special-singular-points-with-zero-c} implies
$\mathcal{P}$ is a Halphen pencil as well.
\begin{proposition}
\label{proposition:n-95} If $\gimel=95$, then the linear systems
$|-5K_{X}|$ and $\mathcal{P}$ are the only Halphen pencils.
\end{proposition}
\begin{proof}
The proof is almost same as the cases $\gimel=89$, $90$, $92$,
$94$.  We have one thing different from these cases. It has a
singular point at $(0:1:0:0:0)$ such that it gives us $b=6$,
$c=0$, and $-K_Y^3<0$ for
Proposition~\ref{proposition:special-singular-points}. The result
therefore follows from Lemmas~\ref{lemma:n-89-90-92-94}
and~\ref{lemma:special-singular-points-with-zero-c}.
\end{proof}

\begin{K3proposition}
A general surface in the pencil $\mathcal{P}$ is birational to a
K3 surface.
\end{K3proposition}
\begin{proof}
Let $\alpha: Y\to X$ be the Kawamata blow up at the point $P$ with
weights $(1,2,3)$ and let $E$ be its exceptional divisor. Then,
the surface $E$ is isomorphic to $\mathbb{P}(1,2,3)$. Let $D$ be a
general surface in the pencil $\mathcal{P}$. Then, the
intersection $\Delta:= E\cdot D_Y$ is a curve of degree six on
$E$. It does not pass through any singular point of the surface
$E$.
 We suppose that the curve
$\Delta$ is smooth. Because it does not pass through any singular
point of $E$ and its degree on $E$ is $6$, it is an elliptic
curve. The singularities of the  surface $D$ are rational except
the point $P$. Then, the same argument of
K3-Proposition~\ref{K3proposition:n-48} gives a contradiction.
Therefore, the surface $D_{V_1}$ is birational to a smooth K3
surface.
\end{proof}

\newpage

\part{Fano threefold hypersurfaces with more than two Halphen
pencils.}\label{section:more-than-two-Halphens}

\section{Case $\gimel=18$, hypersurface of degree $12$ in $\mathbb{P}(1,2,2,3,5)$.}%
\index{$\gimel=18$}\label{section:n-18}

 The threefold $X$ is a general
hypersurface of degree $12$ in $\mathbb{P}(1,2,2,3,5)$ with
$-K_{X}^{3}=\frac{1}{5}$. The singularities of  $X$ consist of six
points $O_{1}$, $O_{2}$, $O_{3}$, $O_{4}$, $O_{5}$ and $O_{6}$
that are quotient singularities of type $\frac{1}{2}(1,1,1)$ and
one point $P$ that is a quotient singularity of type
$\frac{1}{5}(1,2,3)$.

There is a commutative diagram
$$
\xymatrix{
&U\ar@{->}[d]_{\alpha}&&W\ar@{->}[ll]_{\beta}\ar@{->}[d]^{\eta}&\\%
&X\ar@{-->}[rr]_{\psi}&&\mathbb{P}(1,2,2),&}
$$
where \begin{itemize}

\item $\psi$ is  the natural projection,

\item $\alpha$ is the Kawamata blow up at  the point $P$ with
weights $(1,2,3)$,

\item $\beta$ is the Kawamata blow up with weights $(1,2,1)$ of
the singular point of the variety $U$ that is a quotient
singularity of type $\frac{1}{3}(1,2,1)$,

\item $\eta$ is an elliptic fibration.
\end{itemize}

The hypersurface $X$ can be given by the equation
$$
w^{2}z+wf_{7}(x,y,z,t)+f_{12}(x,y,z,t)=0,
$$
where $f_{i}(x,y,z,t)$ is a general quasihomogeneous polynomial of
degree $i$. Let $\mathcal{P}$ be the pencil of surfaces that are
cut out on the hypersurface $X$ by the equations $\lambda
x^{2}+\mu z=0$, where $(\lambda :\mu)\in\mathbb{P}^{1}$.

\begin{K3proposition}\label{K3proposition:n-18}
A general surface of the pencil $\mathcal{P}$ is birational to a
K3 surface. In particular, the linear system $\mathcal{P}$ is a
Halphen pencil.
\end{K3proposition}
\begin{proof}
A general surface of the pencil $\mathcal{P}$ is not ruled because
$X$ is birationally rigid (\cite{CPR}). Hence, a general surface
of the pencil $\mathcal{P}$ is birational to a K3 surface because
it is a compactification of a double cover of $\mathbb{C}^{2}$
branched over a sextic curve.
\end{proof}

The hypersurface $X$ can also be given by the equation
$$
xg_{11}(x,y,z,t,w)+tg_{9}(x,y,z,t,w)+wg_{7}(x,y,z,t,w)+yg_{5}(y,z)=0
$$
such that the point $O_{1}$ is given by the equations $x=y=t=w=0$,
where $g_{i}$ is a general quasihomogeneous polynomial of degree
$i$. Let $\mathcal{P}_{1}$ be the pencil of surfaces that are cut
out on the hypersurface $X$ by the pencil $\lambda x^{2}+\mu y=0$,
where $(\lambda :\mu)\in\mathbb{P}^{1}$. We will see that the
linear system  $\mathcal{P}_{1}$ is a Halphen pencil. The base
locus of $\mathcal{P}_{1}$ does not contain the points $O_{2}$,
$O_{3}$, $O_{4}$, $O_{5}$ and $O_{6}$. Similarly, we can construct
a Halphen pencil $\mathcal{P}_{i}$ such that
$\mathcal{P}_{i}\subset|-2K_{X}|$ and the base locus of the pencil
$\mathcal{P}_{i}$ contains the point $O_{i}$.

\begin{proposition}
\label{proposition:n-18}  The linear systems $\mathcal{P}$,
$\mathcal{P}_{1}$, $\mathcal{P}_{2}$, $\mathcal{P}_{3}$,
$\mathcal{P}_{4}$, $\mathcal{P}_{5}$, and $\mathcal{P}_{6}$  are
the only Halphen pencils on $X$.
\end{proposition}

We may assume that the singularities of the log pair $(X,
\frac{1}{n}\mathcal{M})$ are canonical. Moreover, it follows from
Lemmas~\ref{lemma:smooth-points} and
Corollary~\ref{corollary:Ryder-a1} that
$$
\mathbb{CS}\Big(X, \frac{1}{n}\mathcal{M}\Big)\subset\Big\{O_{1}, O_{2}, O_{3}, O_{4}, O_{5}, O_{6}, P\Big\}.%
$$

\begin{lemma}
\label{lemma:n-18-O1-O6} If $O_{i}\in \mathbb{CS}(X,
\frac{1}{n}\mathcal{M})$, then $\mathcal{M}=\mathcal{P}_{i}$.
\end{lemma}

\begin{proof}
Let $\pi_i:V_i\to X$ be the Kawamata blow up at the point $O_{i}$
with weights $(1,1,1)$. Then,
$\mathcal{M}_{V_i}\sim_{\mathbb{Q}}-nK_{V_i}$  by
Lemma~\ref{lemma:Kawamata}.

The linear system $|-2K_{V_i}|$ is the proper transform of the
pencil $\mathcal{P}_{i}$ and the base locus of $|-2K_{V_i}|$
consists of the irreducible curve $C_{V_i}$ such that
$\pi(C_{V_i})$ is the base curve of the pencil $\mathcal{P}_{i}$.

Let $D$ be a general surface in $|-2K_{V_i}|$. Then, the surface
$D$ is normal and $C_{V_i}^{2}<0$  on the surface $D$. On the
other hand, we have $C_{V_i}\equiv -K_{V_i}\vert_{D}$, which
implies that $\mathcal{M}_{V_i}=|-2K_{V_i}|$ by
Theorem~\ref{theorem:main-tool}.
\end{proof}
\begin{K3proposition}\label{K3proposition:n-18-i}
A general surface of each pencil $\mathcal{P}_i$ is birational to
a K3 surface. In particular, $\mathcal{P}_i$ is a Halphen pencil.
\end{K3proposition}
\begin{proof}
We use the same notations as in the proof of
Lemma~\ref{lemma:n-18-O1-O6}. The pencil $|-2K_{V_i}|$ satisfies
the condition of Theorem~\ref{theorem:Halphen}. Therefore, it is a
Halphen pencil. The intersection of the surface $D$ and the
exceptional divisor $E_i\cong\mathbb{P}^2$ of the birational
morphism $\pi$ is a conic on $E_i$. An irreducible component of
the intersection $D\cdot E_i$ is a rational curve not contained in the base locus of the pencil $|-2K_{V_i}|$.
Therefore, the surface $D$ is
birational to a K3 surface by
Corollary~\ref{corollary:Halphen-K3-rational-curve}.
\end{proof}

Let $E$ be the exceptional divisor of the birational morphism
$\alpha$. It has two  singular points   $Q$ and $O$ that are
quotient singularities of types $\frac{1}{3}(1,2,1)$ and
$\frac{1}{2}(1,1,1)$, respectively.

Let $C$ be the base curve of the pencil $\mathcal{P}$ and $L$ be
the unique curve in of the linear system
$|\mathcal{O}_{\mathbb{P}(1,\,2,\,3)}(1)|$ on $E$.

\begin{lemma}
\label{lemma:n-18-O} If the set $\mathbb{CS}(U,
\frac{1}{n}\mathcal{M}_U)$ contains the point $O$, then
$\mathcal{M}=\mathcal{P}$.
\end{lemma}

\begin{proof}
Let $\pi:V\to U$ be the Kawamata blow up at the point  $O$ with
weights $(1,1,1)$ and $F$ be the exceptional divisor of the
birational morphism  $\pi$. Let $\mathcal{L}$ be the proper
transform of the linear system $|-3K_{U}|$ by the birational
morphism $\pi$. We have $\mathcal{M}_V\sim_{\mathbb{Q}}-nK_{V}$ by
Lemma~\ref{lemma:Kawamata},
$\mathcal{P}_V\sim_{\mathbb{Q}}-2K_{V}$, and
$$
\mathcal{L}\sim_{\mathbb{Q}}\pi^{*}(-3K_{U})-\frac{1}{2}F.
$$
The base locus of the linear system $\mathcal{L}$ consists of the
irreducible curve $\tilde{C}_V$. Moreover, for a general surface
$T$ of the linear system $\mathcal{L}$, the inequality $T\cdot
\tilde{C}_V>0$ holds, which implies that the divisor
$\pi^{*}(-6K_{U})-F$ is nef and big.

Let $M$ and $D$ be general surfaces of the pencils $\mathcal{M}_V$
and $\mathcal{P}_V$, respectively. Then,
$$
\Big(\pi^{*}(-6K_{U})-F\Big)\cdot M\cdot D=\Big(\pi^{*}(-6K_{U})-F\Big)\cdot\Big(\pi^{*}(-nK_{U})-\frac{n}{2}F\Big)\cdot\Big(\pi^{*}(-2K_{U})-F\Big)=0,%
$$
which implies that $\mathcal{M}_V=\mathcal{P}_V$ by
Theorem~\ref{theorem:main-tool}.
\end{proof}

For now, to prove Proposition~\ref{proposition:n-18}, we may
assume that $\mathbb{CS}(X, \frac{1}{n}\mathcal{M})=\{P\}$.
Because $\mathcal{M}_U\sim_{\mathbb{Q}}-nK_{U}$ by
Lemma~\ref{lemma:Kawamata}, the set $\mathbb{CS}(U,
\frac{1}{n}\mathcal{M})_U$ is not empty by
Theorem~\ref{theorem:Noether-Fano}.
 Therefore, Lemma~\ref{lemma:n-18-O} enables us to assume that the set
$\mathbb{CS}(U, \frac{1}{n}\mathcal{M}_U)$ consists of the point
$Q$. The equivalence $\mathcal{M}_W\sim_{\mathbb{Q}}-nK_{W}$  by
Lemma~\ref{lemma:Kawamata} implies that every surface in the
pencil $\mathcal{M}_W$ is contracted to a curve by the elliptic
fibration $\eta$. Moreover, the set $\mathbb{CS}(W,
\frac{1}{n}\mathcal{M}_W)$ is not empty by
Theorem~\ref{theorem:Noether-Fano}.

Let $G$ be the exceptional divisor of the birational morphism
$\beta$ and $Q_1$ be the singular point of the surface $G$. Then,
the point $Q_1$ is the quotient singularity of type
$\frac{1}{2}(1,1,1)$ on the variety $W$. Moreover, it follows from
Lemma~\ref{lemma:Cheltsov-Kawamata} that the set $\mathbb{CS}(W,
\frac{1}{n}\mathcal{M}_W)$ contains the point $Q_1$.

Let $\gamma:Y\to W$ be the Kawamata blow up at the point $Q_1$
with weights $(1,1,1)$. The base locus of the pencil
$\mathcal{P}_Y$ consists of the irreducible curves $C_Y$ and
$L_Y$. Let $D$ be a general surface of the pencil $\mathcal{P}_Y$.
Then, explicit local calculations  show that
$D\sim_{\mathbb{Q}}-2K_{Y}$. On the other hand, the surface $D$ is
normal and the intersection form of the curves $C_Y$ and $L_Y$ on
the surface $D$ is negative-definite. Hence, we obtain the
identity $\mathcal{M}_Y=\mathcal{P}_Y$ from
Theorem~\ref{theorem:main-tool} because
$\mathcal{M}_Y\vert_{D}\equiv n(C_Y+L_Y)$. Therefore, we see that
$\mathcal{M}=\mathcal{P}$, which completes our proof of
Proposition~\ref{proposition:n-18}.

\section{Case $\gimel=22$, hypersurface of degree $14$ in
$\mathbb{P}(1,2,2,3,7)$.}\index{$\gimel=22$} \label{section:n-22}

 The threefold $X$ is a general
hypersurface of degree $14$ in $\mathbb{P}(1,2,2,3,7)$ with
$-K_{X}^{3}=\frac{1}{6}$. The singularities of  $X$ consist of
seven points $O_{1}$, $O_{2}$, $O_{3}$, $O_{4}$, $O_{5}$, $O_6$,
and $O_{7}$ that are quotient singularities of type
$\frac{1}{2}(1,1,1)$ and one point $P$ that is a quotient
singularity of type $\frac{1}{3}(1,2,1)$.

There is a commutative diagram
$$
\xymatrix{
&&Y\ar@{->}[dl]_{\pi}\ar@{->}[dr]^{\eta}&\\
&X\ar@{-->}[rr]_{\psi}&&\mathbb{P}(1,2,2),&&}
$$
where \begin{itemize}

\item $\psi$ is the natural  projection,

\item $\pi$ is the Kawamata blow up at the point $P$ with weights
$(1,2,1)$,

\item $\eta$ is an elliptic fibration.
\end{itemize}
The hypersurface $X$ can be given by the equation
$$
t^{4}z+t^3f_{5}(x,y,z,w)+t^2f_{8}(x,y,z,w)+tf_{11}(x,y,z,w)+f_{14}(x,y,z,w)=0,
$$
where $f_{i}(x,y,z,t)$ is a general quasihomogeneous polynomial of
degree $i$. Let $\mathcal{P}$ be the pencil of surfaces that are
cut out on the hypersurface $X$ by the equations $\lambda
x^{2}+\mu z=0$, where $(\lambda :\mu)\in\mathbb{P}^{1}$. We will
see that the linear system $\mathcal{P}$ is a Halphen pencil.

The hypersurface $X$ can also be given by the equation
$$
xg_{13}(x,y,z,t,w)+tg_{11}(x,y,z,t,w)+wg_{7}(x,y,z,t,w)+yg_{12}(y,z)=0
$$
such that the point $O_{1}$ is given by the equations $x=y=t=w=0$,
where $g_{i}$ is a general quasihomogeneous polynomial of degree
$i$. Let $\mathcal{P}_{1}$ be the pencil of surfaces  cut out on
the hypersurface $X$ by the pencil $\lambda x^{2}+\mu y=0$, where
$(\lambda :\mu)\in\mathbb{P}^{1}$. Then, $\mathcal{P}_{1}$ is a
Halphen pencil. Indeed, a general surface in this pencil  is a
compactification of a double cover of $\mathbb{C}^{2}$ branched
over a sextic curve. The base locus of $\mathcal{P}_{1}$ does not
contain the points $O_{2}$, $O_{3}$, $O_{4}$, $O_{5}$, $O_6$, and
$O_{7}$. Similarly, we can construct a Halphen pencil
$\mathcal{P}_{i}$ such that $\mathcal{P}_{i}\subset|-2K_{X}|$ and
the base locus of the pencil $\mathcal{P}_{i}$ contains the point
$O_{i}$.

\begin{proposition}
\label{proposition:n-22} The linear systems $\mathcal{P}$,
$\mathcal{P}_{1}$, $\mathcal{P}_{2}$, $\mathcal{P}_{3}$,
$\mathcal{P}_{4}$, $\mathcal{P}_{5}$, $\mathcal{P}_{6}$,  and
$\mathcal{P}_{7}$ are the only Halphen pencils on $X$.
\end{proposition}
\begin{proof}
Due to  Lemmas~\ref{lemma:smooth-points} and
Corollary~\ref{corollary:Ryder-a1}, we may assume  that
$$
\mathbb{CS}\Big(X, \frac{1}{n}\mathcal{M}\Big)\subset\Big\{O_{1}, O_{2}, O_{3}, O_{4}, O_{5}, O_{6}, O_7,  P\Big\}.%
$$
If it contains the point $O_i$, we consider the Kawamata blow up
$\pi_i:Y_i\to X$  at the point $O_i$ with weights $(1,1,1)$. The
proof of  Lemma~\ref{lemma:n-18-O1-O6} then shows
$\mathcal{M}=\mathcal{P}_i$.

From now, we suppose that the set $\mathbb{CS}(X,
\frac{1}{n}\mathcal{M})$ consists of the point $P$. Then, the set
$\mathbb{CS}(Y, \frac{1}{n}\mathcal{M}_Y)$ must contain the
singular point that is contained in the exceptional divisor of the
birational morphism $\pi$.

Let $\sigma:W\to Y$ be the Kawamata blow up at this point. Then,
the base locus of the pencil $\mathcal{P}_{W}$ consists of the
irreducible curve $\bar{C}_W$. Moreover,
$\mathcal{P}_{W}\sim_{\mathbb{Q}} -2K_W$, $-K_W\cdot \bar{C}_W<0$,
and $\mathcal{M}_{W}\sim_{\mathbb{Q}} -nK_W$. Therefore,
Theorem~\ref{theorem:main-tool} gives us the identity
$\mathcal{M}=\mathcal{P}$.
\end{proof}
\begin{K3proposition}\label{K3proposiotn:n-22}
A general surface in each of the pencils $\mathcal{P}$,
$\mathcal{P}_{1}$, $\mathcal{P}_{2}$, $\mathcal{P}_{3}$,
$\mathcal{P}_{4}$, $\mathcal{P}_{5}$, $\mathcal{P}_{6}$,  and
$\mathcal{P}_{7}$ is birational to a smooth K3 surface.
\end{K3proposition}
\begin{proof}
 A general surface in the pencil
$\mathcal{P}_i$ is a compactification of a double cover of
$\mathbb{C}^2$ ramified along a sextic curve. By
Theorem~\ref{theorem:CPR}, it must be birational to a smooth K3
surface. For the pencil $\mathcal{P}$, use the same argument with
the exceptional divisor of $\sigma$ as in K3-Proposition~\ref{K3proposition:n-18-i}.
\end{proof}

\section{Case $\gimel=28$, hypersurface of degree $15$ in
$\mathbb{P}(1,3,3,4,5)$.}\index{$\gimel=28$} \label{section:n-28}

 The threefold $X$ is a general
hypersurface of degree $15$ in $\mathbb{P}(1,3,3,4,5)$ with
$-K_{X}^{3}=\frac{1}{12}$. The singularities of  $X$ consist of
five points $O_{1}$, $O_{2}$, $O_{3}$, $O_{4}$, and $O_{5}$ that
are quotient singularities of type $\frac{1}{3}(1,1,2)$ and one
point $P$ that is a quotient singularity of type
$\frac{1}{4}(1,3,1)$.

There is a commutative diagram
$$
\xymatrix{
&&Y\ar@{->}[dl]_{\pi}\ar@{->}[dr]^{\eta}&\\
&X\ar@{-->}[rr]_{\psi}&&\mathbb{P}(1,3,3),&&}
$$
where \begin{itemize}

\item $\psi$ is the natural  projection,

\item $\pi$ is the Kawamata blow up at the point $P$ with weights
$(1,3,1)$,

\item $\eta$ is an elliptic fibration.
\end{itemize}

The hypersurface $X$ can be given by the equation
$$
t^{3}z+t^2f_{7}(x,y,z,w)+tf_{11}(x,y,z,w)+f_{15}(x,y,z,w)=0,
$$
where $f_{i}(x,y,z,t)$ is a general quasihomogeneous polynomial of
degree $i$. Let $\mathcal{P}$ be the pencil of surfaces that are
cut out on the hypersurface $X$ by the equations $\lambda
x^{3}+\mu z=0$, where $(\lambda :\mu)\in\mathbb{P}^{1}$.
As in the previous cases, the  linear system $\mathcal{P}$ is a
Halphen pencil.

The hypersurface $X$ can also be given by the equation
$$
xg_{14}(x,y,z,t)+tg_{11}(x,y,z,t)+wg_{10}(x,y,z,t)+yg_{12}(y,z)=0
$$
such that the point $O_{1}$ is given by the equations $x=y=t=w=0$,
where $g_{i}$ is a general quasihomogeneous polynomial of degree
$i$. Let $\mathcal{P}_{1}$ be the pencil of surfaces that are cut
out on the hypersurface $X$ by the pencil $\lambda x^{3}+\mu y=0$,
where $(\lambda :\mu)\in\mathbb{P}^{1}$. Then, the linear system-
$\mathcal{P}_{1}$ is a Halphen pencil. The base locus of
$\mathcal{P}_{1}$ does not contain the points $O_{2}$, $O_{3}$,
$O_{4}$, and $O_{5}$. Similarly, we can construct a Halphen pencil
$\mathcal{P}_{i}$ such that $\mathcal{P}_{i}\subset|-3K_{X}|$ and
the base locus of the pencil $\mathcal{P}_{i}$ contains the point
$O_{i}$.

\begin{proposition}
\label{proposition:n-28} The linear systems
 $\mathcal{P}$, $\mathcal{P}_{1}$, $\mathcal{P}_{2}$,
$\mathcal{P}_{3}$, $\mathcal{P}_{4}$, and $\mathcal{P}_{5}$ are
the only Halphen pencils on $X$.
\end{proposition}
\begin{proof}
The proof is the same as that of
Proposition~\ref{proposition:n-22}.
\end{proof}
\begin{K3proposition}\label{K3proposiotn:n-28}
A general surface in each of the pencils $\mathcal{P}$,
$\mathcal{P}_{1}$, $\mathcal{P}_{2}$, $\mathcal{P}_{3}$,
$\mathcal{P}_{4}$,  and $\mathcal{P}_{5}$ is birational to a
smooth K3 surface.
\end{K3proposition}
\begin{proof}
Essentially, the proof is the same as the proof of
K3-Proposition~\ref{K3proposition:n-18-i}. Instead of a conic, we
however consider a cubic curve on $\mathbb{P}(1,1,2)$ or
$\mathbb{P}(1,1,3)$, which has a rational irreducible component.
\end{proof}

\newpage

\part{Fano threefold hypersurfaces with infinitely many  Halphen
pencils.}\label{section:infinite-Halphens}

\section{Case $\gimel=1$, hypersurface of degree $4$ in
$\mathbb{P}^{4}$.}\index{$\gimel=01$}\label{section:n-1}

Let $X$ be a general quartic hypersurface in $\mathbb{P}^{4}$. It
is smooth and the log pair $(X, \frac{1}{n}\mathcal{M})$ is
canonical (Theorem~3.6 in \cite{Co00}).

\begin{proposition}
\label{proposition:n-1} Every Halphen pencil is contained in
$|-K_{X}|$.
\end{proposition}

Let us prove Proposition~\ref{proposition:n-1}. Suppose that the
set $\mathbb{CS}(X, \frac{1}{n}\mathcal{M})$ contains a curve $Z$
and does not contain any point. Then, we have
$$\mathrm{mult}_{Z}(\mathcal{M})=n.$$
It follows from Lemma~\ref{lemma:curves} that $\mathrm{deg}(Z)\leq
4$.

\begin{lemma}
\label{lemma:n-1-plane-curve} The curve $Z$ is contained in a
two-dimensional linear subspace of $\mathbb{P}^{4}$.
\end{lemma}

\begin{proof}
Suppose that the curve $Z$ is not contained in any plane in
$\mathbb{P}^{4}$. Then, the degree of the curve $Z$ is either $3$
or $4$. If the degree is $3$, then the curve is smooth. If the
degree is $4$, then the curve  can be singular but the
singularities consist of only one double point.

Suppose that $Z$ is smooth. Let $\alpha:U\to X$ be the blow up
along the curve $Z$ and  $F$ be its exceptional divisor. Then, the
base locus of the linear system
$|\alpha^{*}(-\mathrm{deg}(Z)K_{X})-F|$ does not contain any curve
but
$$
\Big(\alpha^{*}(-\mathrm{deg}(Z)K_{X})-F\Big)\cdot D_{1}\cdot D_{2}<0,%
$$
where $D_{1}$ and $D_{2}$ are general surfaces of the linear
system $\mathcal{M}_U$, which is a contradiction.

Suppose that the curve $Z$ is a quartic curve with a double point
$P$. Let $\beta:W\to X$ be the composition of the blow up at the
point $P$ with the blow up along the proper transform of the curve
$Z$. Let $G$ and $E$ be the exceptional divisors of $\beta$ such
that $\beta(E)=Z$ and $\beta(G)=P$. Then, the base locus of the
linear system $|\beta^{*}(-4K_{X})-E-2G|$ does not contain any
curve  but
$$
\Big(\beta^{*}(-4K_{X})-E-2G\Big)\cdot D_{1}\cdot D_{2}<0
$$
where $D_{1}$ and $D_{2}$ are general surfaces of the linear
system $\mathcal{M}_W$, which is a contradiction.
\end{proof}

\begin{lemma}
\label{lemma:n-1-curve-is-line} If the curve $Z$ is a line, then
the pencil  $\mathcal{M}$ is contained in $|-K_{X}|$.
\end{lemma}

\begin{proof}
Let $\pi:V\to X$ be the blow up along the line $Z$. Then, the
linear system  $|-K_{V}|$ is base-point-free and induces an
elliptic fibration $\eta:V\to\mathbb{P}^{2}$. Therefore,
$\mathcal{M}_V$ is contained in  fibers of $\eta$. In particular,
the base locus of the pencil $\mathcal{M}_V$ does not contain
curves not contracted by the morphism $\eta$.

The set $\mathbb{CS}(V, \frac{1}{n}\mathcal{M}_V)$ is not empty by
the Theorem~\ref{theorem:Noether-Fano}. However, it does not
contain any point  because we assume that the set  $\mathbb{CS}(X,
\frac{1}{n}\mathcal{M})$ does not contain points. Hence, there is
an irreducible curve $L\subset V$ such that
$\mathrm{mult}_{L}(\mathcal{M}_V)=n$ and $\eta(L)$ is a point.

The pencil $\mathcal{M}_V$ is the pull-back of a pencil
$\mathcal{P}$ on $\mathbb{P}^{2}$ via the morphism $\eta$ such
that
$\mathcal{P}\sim_{\mathbb{Q}}\mathcal{O}_{\mathbb{P}^{2}}(n)$.
Hence, the equality $\mathrm{mult}_{L}(\mathcal{M}_V)=n$ implies
that the multiplicity of the pencil $\mathcal{P}$ at the point
$\eta(L)$ is $n$, which implies that $n=1$.
\end{proof}

Thus, we may assume that the set $\mathbb{CS}(X,
\frac{1}{n}\mathcal{M})$ does not contain lines. Moreover,
 the pencil $\mathcal{M}$ is contained
in $|-K_{X}|$ if  $Z$ is a plane quartic curve by
Theorem~\ref{theorem:main-tool}. Thus, we may assume that $Z$ is
either a plane cubic curve or a conic.

\begin{lemma}
\label{lemma:n-1-curve-is-cubic} If the curve $Z$ is a cubic, then
$\mathcal{M}$ is a pencil in $|-K_{X}|$.
\end{lemma}

\begin{proof}
Let $\mathcal{P}$ be the pencil in $|-K_{X}|$ that contains all
surfaces passing through the cubic curve $Z$ and $D$ be a general
surface in $\mathcal{P}$. Then, $D$ is a smooth K3 surface but the
base locus of the pencil $\mathcal{P}$ consists of the curve $Z$
and some line $L\subset X$. We have
$$
\mathcal{M}\Big\vert_{D}=nZ+\mathrm{mult}_{L}(\mathcal{M})L+\mathcal{B}\equiv
nZ+nL,
$$
where $\mathcal{B}$ is a pencil on $D$ without fixed components.
On the other hand, we have $L^{2}=-2$ on the surface $D$, which
implies that $\mathrm{mult}_{L}(\mathcal{M})=n$ and
$\mathcal{B}=\varnothing$. Hence, we have
$\mathcal{M}=\mathcal{P}$ by Theorem~\ref{theorem:main-tool}.
\end{proof}

Therefore, we may assume that the curve $Z$ is a conic. Let $\Pi$
be the plane in $\mathbb{P}^{4}$ that contains the conic $Z$.

\begin{lemma}
\label{lemma:n-1-curve-is-conic} If $\Pi\cap X=Z$, then
$\mathcal{M}$ is a pencil in $|-K_{X}|$.
\end{lemma}

\begin{proof}
Let $\alpha:U\to X$ be the blow up along the curve $Z$ and $D$ be
a general surface of the pencil $|-K_{U}|$. Then, $D$ is a smooth
K3 surface but the base locus of the pencil $|-K_{U}|$ consists of
the irreducible curve $L$ such that $\alpha(L)=Z$ and
$-K_{U}\vert_{D}\equiv L$. Therefore, we have
$$
\mathcal{M}_U\Big\vert_{D}\equiv nL,
$$
but $L^{2}=-2$ on the surface $D$. Hence, we have
$\mathcal{M}_U=|-K_{U}|$ by Theorem~\ref{theorem:main-tool}.
\end{proof}

In the case when the set-theoretic intersection $\Pi\cap X$
contains a curve different from a conic $Z$, the arguments of the
proof of Lemma~\ref{lemma:n-1-curve-is-cubic} easily imply that
$\mathcal{M}$ is a pencil in $|-K_{X}|$. Therefore, we may assume
that the set $\mathbb{CS}(X, \frac{1}{n}\mathcal{M})$ contains a
point $P$ of the quartic $X$.

Let $M_{1}$ and $M_{2}$ be two general surfaces in $\mathcal{M}$.
Then, the inequality $\mathrm{mult}_{P}(M_{1}\cdot M_{2})\geq
4n^{2}$ holds (\cite{Co00} and \cite{Pu98a}). On the other hand,
the degree of the cycle $M_{1}\cdot M_{2}$ is $4n^{2}$, which
implies that $\mathrm{mult}_{P}(M_{1}\cdot M_{2})=4n^{2}$. In
particular, the support of the cycle $M_{1}\cdot M_{2}$ consists
of the union of all lines passing through the point $P$, which
implies that there are at most finitely many lines on the quartic
$X$ passing through the point $P$. Moreover, the equality
$\mathrm{mult}_{P}(\mathcal{M})=2n$ holds (\cite{Ch00a},
\cite{Ch03a}, Corollary~C.14 in \cite{Ch05}).

\begin{lemma}
\label{lemma:n-1-lines} For a line $L$ on $X$ passing through $P$,
$\frac{n}{3}\leq\mathrm{mult}_{L}(\mathcal{M})\leq \frac{n}{2}$.
\end{lemma}

\begin{proof}
Let $D$ be a general hyperplane section of $X$ that passes through
the line $L$ and $M$ be a general surface in $\mathcal{M}$. Then,
$D$ is a smooth K3 surface and
$$
M\Big\vert_{D}=\mathrm{mult}_{L}(\mathcal{M})L+\Delta,
$$
where $\Delta$ is an effective divisor such that
$\mathrm{mult}_{P}(\Delta)\geq 2n-\mathrm{mult}_{L}(\mathcal{M})$.
On the other hand, we have $L^{2}=-2$ on the surface $D$. Hence,
we have
$$
n+2\mathrm{mult}_{L}(\mathcal{M})=L\cdot\Delta\geq\mathrm{mult}_{P}(\Delta)\geq 2n-\mathrm{mult}_{L}(\mathcal{M}),%
$$
which implies $\frac{n}{3}\leq\mathrm{mult}_{L}(\mathcal{M})$.

Let $T$ be the hyperplane section  tangent to the quartic $X$ at
the point $P$. Then, $T$ has isolated singularities and the point
$P$ is an isolated double point of the surface $T$ because
$\mathrm{mult}_{P}(M)=2n$. The cycle $T\cdot D$ is reduced and
consists of the line $L$ and possibly reducible cubic curve
$Z\subset D$ that passes through the point $P$. Thus, we have
$$
3n=\big(\mathrm{mult}_{L}(\mathcal{M})L+\Delta\big)\cdot
Z=3\mathrm{mult}_{L}(\mathcal{M})+\Delta\cdot Z\geq
3\mathrm{mult}_{L}(\mathcal{M})+2n-\mathrm{mult}_{L}(\mathcal{M}),
$$
which implies $\mathrm{mult}_{L}(\mathcal{M})\leq \frac{n}{2}$.
\end{proof}

Therefore,  any curve containing the point $P$ cannot belong to
the set $\mathbb{CS}(X, \frac{1}{n}\mathcal{M})$. Let $\pi:V\to X$
be the blow up at the point $P$ and $E$ be the exceptional divisor
of the blow up $\pi$. In addition, let $B_{i}$ be the proper
transform of the divisor $M_{i}$ by $\pi$. Then, the equalities
$\mathrm{mult}_{P}(M_{1}\cdot M_{2})=4n^{2}$ and
$\mathrm{mult}_{P}(\mathcal{M})=2n$ imply that
$$
B_{1}\cdot
B_{2}=\sum_{i=1}^{k}\mathrm{mult}_{\bar{L}_{i}}\big(B_{1}\cdot B_{2})\bar{L}_{i},%
$$
where $k$ is a number of lines on $X$ that passes through the
point $P$ and $\bar{L}_{i}$ is an irreducible curve such that
$\pi(\bar{L}_{i})$ is a line on $X$ that passes through the point
$P$.

\begin{lemma}
\label{lemma:n-1-curves} Let $Z$ be an irreducible curve on $X$
that is not a line passing through the point $P$. Then,
$$\mathrm{deg}(Z)\geq 2\mathrm{mult}_{P}(Z),$$
where the equality  holds only if the proper transform $Z_V$ does
not intersect the curve $\bar{L}_{i}$ for any $i$.
\end{lemma}

\begin{proof}
The proper transform $Z_V$ is not contained in $B_{i}$ because the
base locus of the pencil $\mathcal{M}_V$ consists of the curves
$\bar{L}_{1},\cdots, \bar{L}_{k}$. Hence, we have
$$
0\leq B_{i}\cdot Z_V\leq
n\big(\mathrm{deg}(Z)-2\mathrm{mult}_{P}(Z)\big),
$$
which concludes the proof.
\end{proof}

Note that so far we never use the generality of the quartic $X$
besides its smoothness. In the following we assume that there are
at most $3$ lines on $X$ passing though a given point of $X$ and
every line on $X$ has normal bundle
$\mathcal{O}_{\mathbb{P}^{1}}(-1)\oplus\mathcal{O}_{\mathbb{P}^{1}}$.
It follows from the proof of Proposition~1 in \cite{Pu98a} that
the former condition is satisfied on a general quartic threefold.
The latter  condition is also satisfied on a general quartic
threefold by  \cite{Col79}. Moreover, the article \cite{Col79}
shows that the latter  condition is equivalent to the following:
no two-dimensional linear subspace of $\mathbb{P}^{4}$ is tangent
to the quartic $X$ along a line. In particular, we see that no
hyperplane section of $X$ can be singular at three points that are
contained in a single line.

\begin{lemma}
\label{lemma:n-1-mult}
 For a line $L$ in $X$ passing through $P$,
$\mathrm{mult}_{L}(\mathcal{M})= \frac{n}{2}$.
\end{lemma}

\begin{proof} By Lemma~\ref{lemma:n-1-lines}, it is enough to show
$\mathrm{mult}_{L}(\mathcal{M})\geq \frac{n}{2}$.

Let $\alpha:W\to X$ be the blow up along the line $L$ and $F$ be
the exceptional divisor of the blow up $\alpha$. Then, the surface
$F$ is the rational ruled surface $\mathbb{F}_{1}$.

Let $\Delta$ be the irreducible curve on the surface $F$ such that
$\Delta^{2}=-1$ and $Z$ be the fiber of the restricted morphism
$\pi\vert_{F}:F\to L$ over the point $P$. Then, $F\vert_{F}\equiv
-(\Delta +Z)$, which implies that
$$
\mathcal{M}_W\Big\vert_{F}\equiv
nZ+\mathrm{mult}_{L}(\mathcal{M})(\Delta+Z).
$$

Let $\beta:U\to W$ be the blow up along the curve $Z$ and $G$ be
the exceptional divisor of $\beta$. Then, the exceptional divisor
$E$ of $\pi$  is the proper transform of the divisor $G$ on the
threefold $V$. Hence, we have
$$
\mathrm{mult}_{Z}(\mathcal{M}_W)=2n-\mathrm{mult}_{L}(\mathcal{M}),
$$
which implies that $\mathrm{mult}_{Z}(\mathcal{M}_W\vert_{F})\geq
2n-\mathrm{mult}_{L}(\mathcal{M})$. Therefore, we have
$$
n+\mathrm{mult}_{L}(\mathcal{M})\geq
2n-\mathrm{mult}_{L}(\mathcal{M}),
$$
which gives $\mathrm{mult}_{L}(\mathcal{M})\geq \frac{n}{2}$.
\end{proof}

Let ${T}$ be a hyperplane section of $X$ that is singular at the
point $P$. Then, ${T}$ has only isolated singularities. Moreover,
we have $\mathrm{mult}_{P}({T}\cdot M_{i})=4n$, which implies that
the point $P$ is an isolated double point of the surface $T$. Put
${L}_{i}=\pi(\bar{L}_{i})$.

\begin{lemma}
\label{lemma:n-1-ODP} The point $P$ is not an ordinary double
point of the surface $T$.
\end{lemma}

\begin{proof}
Suppose that the point $P$ is  an ordinary double point of the
surface ${T}$. Let us show that this assumption leads us to a
contradiction.

Let $H_{i}$ be a general hyperplane section of the quartic $X$
that passes through the line ${L}_{i}$. Then,
$$
H_{i}\cdot T={L}_{i}+Z_{i},
$$
where $Z_{i}$ is a cubic curve. The cubic curve $Z_{i}$ intersect
the line ${L}_{i}$ at the point $P$ and at some smooth point of
the surface ${T}$ because ${L}_{i}$ does not contain three
singular points of the surface ${T}$. Hence, we have
$$
{L}_{i}^{2}=H_{i}\cdot L_{i}-Z_{i}\cdot L_{i}<-\frac{1}{2}.
$$

The proper transform $T_V$ has isolated singularities and normal.
Moreover, the inequality ${L}_{i}^{2}<-\frac{1}{2}$ implies that
$\bar{L}_{i}^{2}<-1$.

Let $M$ be a general surface in $\mathcal{M}$. The support of the
cycle ${T}\cdot M$ consists of the union of all lines on $X$
passing through the point $P$ because $\mathrm{mult}_{P}({T}\cdot
M)=4n$. Thus, the equalities $\mathrm{mult}_{P}({T})=2n$ and
$\mathrm{mult}_{P}(M)=2n$ implies that the support of the cycle
$T_V\cdot M_V$ consists of the union of the curves
$\bar{L}_{1},\cdots,\bar{L}_{k}$. Hence, we have
$$
M_{V}\Big\vert_{T_V}=\sum_{i=1}^{k}m_{i}\bar{L}_{i},
$$
but $M_{V}\cdot \bar{L}_{l}=-n$ and $\bar{L}_{i}\cdot
\bar{L}_{j}=0$ for $i\ne j$. Hence, we have
$$
-n=M_{V}\cdot \bar{L}_{j}=\sum_{i=1}^{k}m_{i}\bar{L}_{i}\cdot
\bar{L}_{j}=m_{j}\bar{L}_{j}^{2},
$$
which implies that $m_{j}<n$.

Let $H$ be the proper transform of a general hyperplane section of
$X$ on the threefold $V$. Then,
$$
4n=M_{V}\cdot T_V\cdot H=\sum_{i=1}^{k}m_{i}\bar{L}_{i}\cdot
H=\sum_{i=1}^{k}m_{i}<kn,
$$
which implies that $k>4$. Thus, the threefold $X$ has at least
five lines that pass through the point $P$, which is a
contradiction.
\end{proof}

Thus, the point $P$ is not an ordinary double point on the surface
${T}$. Therefore, there is a hyperplane section $Z$ of the quartic
surface ${T}$ with $\mathrm{mult}_{P}(Z)\geq 3$. Hence, the curve
$Z$ is reducible by Lemma~\ref{lemma:n-1-curves}. Moreover, the
curve $Z$ is reduced and $\mathrm{mult}_{P}(Z)=3$ by our
assumption of generality of the quartic $X$.

\begin{lemma}
\label{lemma:n-1-3-lines} The curve $Z$ is not a union of four
lines.
\end{lemma}

\begin{proof}
Suppose that the curve $Z$ is a union of four lines. Then, one
component of $Z$ is a line $L$ that does not pass through the
point $P$. Then, $L$ intersects $M_{i}$  at least three points
that are contained in the union of the lines
${L}_{1},\cdots,L_{k}$. On the other hand, we have $M_{i}\cdot
L=n$, which implies that $L$ is contained in $M_{i}$ by
Lemma~\ref{lemma:n-1-mult}, which is impossible because the base
locus of $\mathcal{M}$ is the union of the lines
$L_{1},\cdots,L_{k}$.
\end{proof}

The curve $Z$ is not a union of an irreducible cubic curve and a
line due to Lemma~\ref{lemma:n-1-curves}. Hence, the curve $Z$ is
a union of two different lines passing through the point $P$ and a
conic that also passes through the point $P$, which is impossible
by Lemma~\ref{lemma:n-1-curves}. Hence, we have completed the
proof of Proposition~\ref{proposition:n-1}.

\section{Case $\gimel=2$, hypersurface of degree $5$ in
$\mathbb{P}(1,1,1,1,2)$.}\label{section:n-2}\index{$\gimel=02$}

The threefold $X$ is a general hypersurface of degree $5$ in
$\mathbb{P}(1,1,1,1,2)$ with $-K_X^3=\frac{5}{2}$. It has only one
singular point $O$ at $(0:0:0:0:1)$ which is a quotient
singularity of type $\frac{1}{2}(1,1,1)$. The hypersurface $X$ can
be given by the equation
$$
w^{2}f_{1}\big(x,y,z,t\big)+wf_{3}\big(x,y,z,t\big)+f_{5}\big(x,y,z,t\big)=0,
$$
where  $f_{i}$ is a homogeneous polynomial of degree $i$.

There is a commutative diagram
$$
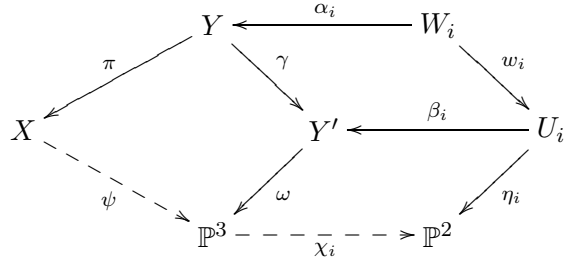
\xymatrix{
&&&Y\ar@{->}[dll]_{\pi}\ar@{->}[dr]^{\gamma}&&W_{i}\ar@{->}[ll]_{\alpha_{i}}\ar@{->}[dr]^{w_{i}}&&\\%
&X\ar@{-->}[drr]_{\psi}&&&Y'\ar@{->}[dl]^{\omega}&&U_{i}\ar@{->}[ll]_{\beta_{i}}\ar@{->}[ld]^{\eta_{i}}&\\
&&&\mathbb{P}^{3}\ar@{-->}[rr]_{\chi_{i}}&&\mathbb{P}^{2}&&}
$$
where \begin{itemize} \item $\psi$ is the natural projection,

\item $\pi$ is the Kawamata blow up at the point $O$ with weights
$(1,1,1)$,

\item $\gamma$ is the birational morphism that contracts $15$
smooth rational curves $L_{1},\cdots , L_{15}$ to $15$ isolated
ordinary double points $P_{1},\cdots , P_{15}$ of the variety
$Y'$, respectively,

\item $\alpha_{i}$ is the blow up  along the curve $L_{i}$,

\item  $\beta_{i}$ is the  blow up at the point $P_{i}$,

\item  $w_{i}$ is a birational morphism,

\item $\omega$ is a double cover of $\mathbb{P}^3$ branched over a
sextic surface $R\subset\mathbb{P}^{3}$,

\item $\chi_{i}$ is the projection from  the point
$\omega(P_{i})$,

\item  $\eta_{i}$ is an elliptic fibration.

\end{itemize}
The surface $R$ is given by the equation
$$
f_{3}\big(x,y,z,t\big)^{2}-4f_{1}\big(x,y,z,t\big)f_{5}\big(x,y,z,t\big)=0\subset\mathbb{P}^{3}\cong\mathrm{Proj}\Big(\mathbb{C}[x,y,z,t]\Big).
$$
It has $15$ ordinary double points $\omega(P_{1}),\cdots ,
\omega(P_{15})$ that are given by the equations
$$
f_{3}\big(x,y,z,t\big)=f_{1}\big(x,y,z,t\big)=f_{5}\big(x,y,z,t\big)=0\subset\mathbb{P}^{3}.
$$
We may assume that the curves on $\mathbb{P}^3$ defined by
$f_{3}=f_{1}=0$ and $f_{5}=f_{1}=0$ are irreducible.

For the convenience, let $M$ and $M'$ be general surfaces in the
pencil $\mathcal{M}$.
\begin{lemma}
\label{lemma:n-2-smooth-points} The set $\mathbb{CS}(X,
\frac{1}{n}\mathcal{M})$ does not contain any smooth point of $X$.
\end{lemma}

\begin{proof}
Suppose that the set $\mathbb{CS}(X, \frac{1}{n}\mathcal{M})$
contains a smooth point $P$ of $X$. Let $D$ be a general surface
of the linear system $|-K_{X}|$ that passes through the point $P$.
The surface $D$ does not contain an irreducible component of the
cycle $M\cdot M^{\prime}$  if none of $\pi(L_i)$ passes through
the point $P$.  In particular, in such a case, we see
$$
\mathrm{mult}_{P}\Big(M\cdot M^{\prime}\Big)\leq M\cdot M^{\prime}\cdot D=-n^{2}K_{X}^{3}=\frac{5}{2}n^{2},%
$$
which is impossible by \cite{Pu98a}. Thus, we may assume that  the
curve $\pi(L_{1})$ passes through the point $P$.

Let us use the arguments of the article \cite{CPR}. Put
$L=\pi(L_{1})$ and
$$
\mathcal{M}\Big\vert_{D}=\mathcal{L}+\mathrm{mult}_{L}\big(\mathcal{M}\big)L,
$$
where $\mathcal{L}$ is a pencil on the surface $D$ without fixed
curves. Then, the point $P$ is a center of log canonical
singularities of the log pair $(D,
\frac{1}{n}\mathcal{M}\vert_{D})$ by the Shokurov connectedness
principle (\cite{Co00} and \cite{Sho93}). It implies that
$$
\mathrm{mult}_{P}\Big(\Lambda_{1}\cdot \Lambda_{2}\Big)
\geq 4n\Big(n-\mathrm{mult}_{L}\big(\mathcal{M}\big)\Big)%
$$
by Theorem~3.1 in \cite{Co00}, where $\Lambda_{1}$ and
$\Lambda_{2}$ are general curves in $\mathcal{L}$. The equality
$$
\Lambda_{1}\cdot \Lambda_{2}=\frac{5}{2}n^{2}-\mathrm{mult}_{L}\big(\mathcal{M}\big)n-\frac{3}{2}\mathrm{mult}^{2}_{L}\big(\mathcal{M}\big)%
$$
holds on the surface $D$ because $L^{2}=-\frac{3}{2}$ on the
surface $D$. Hence, we have
$$
\frac{5}{2}n^{2}-\mathrm{mult}_{L}\big(\mathcal{M}\big)n-\frac{3}{2}\mathrm{mult}^{2}_{L}\big(\mathcal{M}\big)\geq 4n\Big(n-\mathrm{mult}_{L}\big(\mathcal{M}\big)\Big),%
$$
which gives $\mathrm{mult}_{L}(\mathcal{M})=n$. Thus, the set
$\mathbb{CS}(X, \frac{1}{n}\mathcal{M})$ contains the curve
$\pi(L_{1})$.

The set $\mathbb{CS}(X, \frac{1}{n}\mathcal{M})$ contains the
point $O$ by Lemma~\ref{lemma:Kawamata}. Then,
$\mathcal{M}_{W_{1}}\sim_{\mathbb{Q}}-nK_{W_{1}}$ because
$\mathrm{mult}_{L}(\mathcal{M})=n$, which implies that each
surface of $\mathcal{M}_{W_{1}}$ is contracted to a curve by the
elliptic fibration $\eta_{1}\circ w_{1}$. On the other hand, the
set $\mathbb{CS}(W_{1}, \frac{1}{n}\mathcal{M}_{W_{1}})$ contains
a subvariety of the threefold $W_{1}$ that dominates the point
$P$.

Let $E_{1}$ be the exceptional divisor of $\alpha_{1}$. Then,
$E_{1}\cong\mathbb{P}^{1}\times\mathbb{P}^{1}$ and the pencil
$\mathcal{M}_{W_{1}}\vert_{E_{1}}$ does not have fixed components
because $E_{1}$ is a section of the elliptic fibration
$\eta_{1}\circ\omega_{1}$ and the base locus of the pencil
$\mathcal{M}_{W_{1}}$ does not contain curves not contracted by
the elliptic fibration $\eta_{1}\circ w_{1}$. Thus, the set
$\mathbb{CS}(W_{1}, \frac{1}{n}\mathcal{M}_{W_{1}})$ contains a
point $Q$ of the surface $E_{1}$ such that
$\pi\circ\alpha_{1}(Q)=P$.

The point $Q$ is a center of log canonical singularities of the
log pair $(E_{1}, \frac{1}{n}\mathcal{M}_{W_{1}}\vert_{E_{1}})$ by
the Shokurov connectedness principle (\cite{Co00} and
\cite{Sho93}). Let $\Delta_{1}$ and $\Delta_{2}$ be  general
curves in $\mathcal{M}_{W_{1}}\vert_{E_{1}}$. Then, the inequality
$$
2n^{2}=\mathrm{mult}_{Q}\Big(\Delta_{1}\cdot \Delta_{2}\Big)\geq 4n^{2}%
$$
holds by Theorem~3.1 in \cite{Co00}, which is a contradiction.
\end{proof}

\begin{lemma}
\label{lemma:n-2-curves-in-nonsingular-locus} If the set
$\mathbb{CS}(X, \frac{1}{n}\mathcal{M})$ contains a curve
$\Lambda$ not passing through the singular point $O$, then the
pencil $\mathcal{M}$ is contained in $|-K_{X}|$.
\end{lemma}

\begin{proof}
We have $\mathrm{mult}_{\Lambda}(\mathcal{M})=n$ and $-K_{X}\cdot
\Lambda\leq 2$  by Corollary~\ref{corollary:curves}.

Suppose that $-K_{X}\cdot \Lambda=2$ and $\psi(\Lambda)$ is a
line. Then, the line $\psi(\Lambda)$ passes through a unique
singular point of $R$. Hence, we may assume that the curve
$\Lambda$ intersects $\pi(L_{i})$ only for $i=1$.

Let $\mathcal{D}$ be the pencil in the linear system $|-K_{X}|$
consisting of surfaces that pass through the curve $\Lambda$ and
$D$ be a general surface of the pencil $\mathcal{D}$. Then, the
surface $D$ is smooth in the outside of the singular point $O$,
the point $O$ is an ordinary double point of the surface $D$, and
the base locus of the pencil $\mathcal{D}$ consists of the curve
$\Lambda$ and the curve $\pi(L_{1})$. Put $L=\pi(L_{1})$. Then,
$$
\mathcal{M}\Big\vert_{D}=\mathrm{mult}_{\Lambda}\big(\mathcal{M}\big)\Lambda+
\mathrm{mult}_{L}\big(\mathcal{M}\big)L+\mathcal{L}\equiv
n\big(\Lambda+L\big),
$$
where $\mathcal{L}$ is a pencil with no fixed curves. It gives
$\mathcal{M}=\mathcal{D}$ by Theorem~\ref{theorem:main-tool}
because the inequality $L^{2}<0$ holds on the surface $D$.

We may assume that either the equality $-K_{X}\cdot \Lambda=1$
holds or $\psi(\Lambda)$ is a conic, which implies that $\Lambda$
is smooth. Let $\sigma:\breve{X}\to X$ be the blow up along the
curve $\Lambda$ and $G$ be its exceptional divisor.

Suppose that $-K_{X}\cdot \Lambda=2$. Then, $\Lambda$ is cut, in
the set-theoretic sense, by the surfaces of the linear system
$|-2K_{X}|$ that pass through the curve $\Lambda$. Moreover, the
scheme-the\-o\-re\-tic intersection of two  general surfaces of
the linear system $|-2K_{X}|$ passing through the curve $\Lambda$
is reduced at a generic point of the curve $\Lambda$, which
implies that the divisor $\sigma^*(-2K_{X})-G$ is nef by
Lemma~5.2.5 in \cite{CPR}. However, we obtain an absurd inequality
$$
-3n^{2}=\Big(\sigma^*(-2K_{X})-G\Big)\cdot M_{\breve{X}}\cdot M^{\prime}_{\breve{X}}\geq 0.%
$$
 Therefore, the equality $-K_{X}\cdot \Lambda=1$
holds, which implies that $|-K_{\breve{X}}|$ is a pencil.

Suppose that $\psi(\Lambda)$ is not contained in the plane
$f_{1}(x,y,z,t)=0$. Then, $\psi(\Lambda)$ contains a unique
singular point of the surface $R\subset\mathbb{P}^{3}$. Hence, we
may assume that the curve $\Lambda$ intersects $\pi(L_{i})$ only
for $i=1$. It implies that the base locus of the linear system
$|-K_{\breve{X}}|$ consists of the irreducible curves
$\breve{\Lambda}$ and $\breve{L}_{1}$ such that
$(\psi\circ\sigma)(\breve{\Lambda})=\psi(\Lambda)$ and
$\sigma(\breve{L}_{1})=\pi(L_{1})$. Let $\breve{D}$ be a general
surface in $|-K_{\breve{X}}|$. Then, we can consider the curves
$\breve{\Lambda}$ and $\breve{L}_{1}$  as divisors on $\breve{D}$.
We have

$$
\breve{\Lambda}^{2}=-2,\ \breve{L}_{1}^{2}=-\frac{3}{2},\ \breve{\Lambda}\cdot\breve{L}_{1}=1,%
$$
which implies the negative-definiteness of the intersection form
of $\breve{\Lambda}$ and $\breve{L}_{1}$. Because
$$
\mathcal{M}_{\breve{X}}\Big\vert_{\breve{D}}\equiv
-nK_{\breve{X}}\Big\vert_{\breve{D}}\equiv
n\big(\breve{\Lambda}+\breve{L}_{1}\big),%
$$
it follows from Theorem~\ref{theorem:main-tool} that
$\mathcal{M}_{\breve{X}}=|-K_{\breve{X}}|$.

Finally, we suppose that the line $\psi(\Lambda)$ is contained in
the plane $f_{1}(x,y,z,t)=0$. In particular, the line
$\psi(\Lambda)$ is not contained in the surface $R$ because the
curve $f_{3}=f_{1}=0$ is irreducible. Moreover, the line
$\psi(\Lambda)$ contains exactly three singular points of the
ramification surface\footnote{In fact, we may assume that no three
points of the set $\mathrm{Sing}(R)$ are collinear.}; otherwise
the point $O$ would belong to the curve $\Lambda$. Thus, the curve
$\Lambda$ intersects exactly three curves among the curves
$L_{1},\cdots, L_{15}$; otherwise $\Lambda$ would contain the
point $O$.

We may assume that $\Lambda$ intersects the curves $\pi(L_{1})$,
$\pi(L_{2})$, and $\pi(L_{3})$, which means that the points
$\omega(P_{1})$, $\omega(P_{2})$, $\omega(P_{3})$ are contained in
$\psi(\Lambda)$. The base locus of $|-K_{\breve{X}}|$ consists of
the curves $\breve{L}_{1}$, $\breve{L}_{2}$, $\breve{L}_{3}$ such
that 
$\sigma(\breve{L}_{i})=\pi(L_{i})$. The curves $\breve{L}_{1}$,
$\breve{L}_{2}$, $\breve{L}_{3}$  can be contracted on the surface
$\breve{D}$ to a singular point of type $\mathbb{D}_{4}$, which
implies that their intersection form is negative-definite. Hence,
we have $\mathcal{M}_{\breve{X}}=|-K_{\breve{X}}|$ by
Theorem~\ref{theorem:main-tool} .
\end{proof}

The equivalence $\mathcal{M}_{Y}\sim_{\mathbb{Q}}-nK_{Y}$ holds by
Lemma~\ref{lemma:Kawamata}. It implies that the set
$\mathbb{CS}(Y, \frac{1}{n}\mathcal{M}_{Y})$ contains no point of
$Y$ due to Lemmas~\ref{lemma:Cheltsov-Kawamata} and
\ref{lemma:n-2-smooth-points}. Let $\mathcal{M}_{Y'}$ be the
push-forward of the pencil $\mathcal{M}_Y$ by the birational
morphism $\gamma$. Then,
$\mathcal{M}_{Y'}\sim_{\mathbb{Q}}-nK_{Y'}$, the log pair $(Y',
\frac{1}{n}\mathcal{M}_{Y'})$ has canonical singularities but it
follows from Theorem~\ref{theorem:Noether-Fano} that the
singularities of the log pair $(Y', \frac{1}{n}\mathcal{M}_{Y'})$
are not terminal.

\begin{lemma}
\label{lemma:n-2-curves-downstairs} If the set $\mathbb{CS}(Y',
\frac{1}{n}\mathcal{M}_{Y'})$ contains an irreducible curve
$\Gamma$ with $-K_{Y'}\cdot\Gamma\ne 1$, then the pencil
$\mathcal{M}$ is contained in $|-K_{X}|$.
\end{lemma}

\begin{proof}
Let $D$ be a general divisor in $|-K_{Y'}|$. In addition, let
$M_{Y'}=\gamma(M_Y)$ and $M_{Y'}'=\gamma(M_Y')$. Then,
$$
2n^{2}=D\cdot M_{Y'}\cdot M^{\prime}_{Y'}\geq
\mathrm{mult}_{\Gamma}\Big(M_{Y'}\cdot M^{\prime}_{Y'}\Big)D\cdot
\Gamma
\geq -n^{2}K_{Y'}\cdot\Gamma%
$$
because $\mathrm{mult}_{\Gamma}(\mathcal{M}_{Y'})=n$. Therefore,
the inequality $-K_{Y'}\cdot\Gamma\leq 2$ holds.

Suppose that $-K_{Y'}\cdot\Gamma=2$ but the curve $\omega(\Gamma)$
is a line. Let $\mathcal{T}$ be the linear subsystem of the linear
system $|-K_{Y'}|$ consisting of surfaces passing through the
curve $\Gamma$ and $T$ be a general surface in the pencil
$\mathcal{T}$. Then, the base locus of the pencil $\mathcal{T}$
consists of the curve $\Gamma$ and the rational map induced by the
pencil $\mathcal{T}$ is the composition of the double cover
$\omega$ with the projection from the line $\omega(\Gamma)$. On
the other hand, we have
$$
2n=D\cdot T\cdot M_{Y'}\geq \mathrm{mult}_{\Gamma}\Big(T\cdot
M_{Y'}\Big)D\cdot \Gamma\geq
 -nK_{Y'}\cdot\Gamma,%
$$
which implies that the support of the cycle $T\cdot M_{Y'}$ is
contained in $\Gamma$. Thus, we have
$\mathcal{M}_{Y'}=\mathcal{T}$ by Theorem~\ref{theorem:main-tool}.

For now, we suppose that $-K_{Y'}\cdot\Gamma=2$ but the curve
$\omega(\Gamma)$ is a conic. Then, $\Gamma$ is smooth and
$\omega\vert_{\Gamma}$ is an isomorphism. Moreover, the curve
$\Gamma$ contains at most $2$ singular points of the threefold
$Y'$ if the curve $\omega(\Gamma)$ is not contained in the plane
$f_{1}(x,y,z,t)=0$, and the curve $\Gamma$ contains at most $6$
singular points of the threefold $Y'$ otherwise. We may assume
that $\Gamma$ passes through $P_{1},\cdots ,P_{k}$, where $0\leq
k\leq 6$. The equality $k=0$ means that $\Gamma$ lies in the
smooth locus of the threefold $Y'$.

Let $\beta:V\to Y'$ be the blow up at the points $P_{1},\cdots
,P_{k}$, and $E_{i}$ be the exceptional divisor of the blow up
$\beta$ with $\beta(E_{i})=P_{i}$. The exceptional divisor
$E_{i}$ is isomorphic to $\mathbb{P}^{1}\times\mathbb{P}^{1}$. The
proper transform $\Gamma_V$ intersects the surface $E_{i}$
transversally at a single point, which we denote by $Q_{i}$.

Let $\upsilon:W\to V$ be the blow up along the curve ${\Gamma}_V$
and $G$ be the exceptional divisor of the birational morphism
$\upsilon$. In addition, let $A_{i}$ and $B_{i}$ be the fibers of
the natural projections of the surface $E_{i}$ that pass through
the point $Q_{i}$, and $\bar{A}_{i}$ and $\bar{B}_{i}$ be the
proper transforms of the curves $A_{i}$ and $B_{i}$ on the
threefold $W$, respectively. Then, we can flop the curves
$\bar{A}_{i}$ and $\bar{B}_{i}$.

Let $\upsilon_1:U\to W$ be the blow up along the curves
$\bar{A}_{1}, \bar{B}_{1}, \cdots , \bar{A}_{k}, \bar{B}_{k}$.
Also, let $F_{i}$ and $H_{i}$ be the exceptional divisors of
$\upsilon_1$ such that $\upsilon_1(F_{i})=\bar{A}_{i}$ and
$\upsilon_1(H_{i})=\bar{B}_{i}$. Then, all the exceptional
divisors are isomorphic to  $\mathbb{P}^{1}\times\mathbb{P}^{1}$.
There is a birational morphism $\upsilon_1^{\prime}:U\to
W^{\prime}$ such that $\upsilon_1^{\prime}(F_{i})$ and
$\upsilon_1^{\prime}(H_{i})$ are rational curves but
$\upsilon_1^{\prime}\circ\upsilon_1^{-1}$ is not biregular in a
neighborhood of $\bar{A}_{i}$ and $\bar{B}_{i}$. Let
$E_{i}^{\prime}$ be the proper transform of $E_{i}$ on the
threefold $W^{\prime}$. Then, we can contract the surface
$E_{i}^{\prime}$ to a singular point of type $\frac{1}{2}(1,1,1)$.

Let $\upsilon^{\prime}:W^{\prime}\to V^{\prime}$ be the
con\-trac\-tion of $E_{1}^{\prime},\cdots , E_{k}^{\prime}$ and
$G^{\prime}$ be the proper trans\-form of the surface $G$ on the
threefold $V^{\prime}$. Then, there is a birational morphism
$\beta^{\prime}:V^{\prime}\to Y'$ that contracts the divisor
$G^{\prime}$ to the curve $\Gamma$. Hence, we constructed the
commutative diagram
$$
\xymatrix{
&&W\ar@{->}[dl]_{\upsilon}&&U\ar@{->}[ll]_{\upsilon_1}\ar@{->}[dr]^{\upsilon_1^{\prime}}\\%
&V\ar@{->}[dr]_{\beta}&&&&W^{\prime}\ar@{->}[dl]^{\upsilon^{\prime}}\\%
&&Y'&&V^{\prime}.\ar@{->}[ll]_{\beta^{\prime}}&}
$$
The threefold $V^{\prime}$ is projective. Its singularities
consist of $15-k$ ordinary double points and $k$ singular points
of type $\frac{1}{2}(1,1,1)$. However, it is not
$\mathbb{Q}$-factorial because the threefold $Y'$ is not
$\mathbb{Q}$-factorial.

The construction of the birational morphism $\beta^{\prime}$
implies that
$$
\mathcal{M}_{V^{\prime}}\sim_{\mathbb{Q}} -n\beta'^*(K_{Y'})-nG^{\prime}\sim_{\mathbb{Q}}-nK_{V^{\prime}}.%
$$
Let $D^{\prime}$ be a general surface of the linear system
$|{\beta^{\prime}}^{*}(-4K_{Y'})-G^{\prime}|$. Then, the divisor
$D^{\prime}$ is nef by Lemma~5.2.5 in \cite{CPR}. The construction
of the birational morphism $\beta^{\prime}$ implies that
$$
0>\Big(-4+\frac{k}{2}\Big)n^{2}=\Big({\beta^{\prime}}^{*}(-4K_{Y'})-G^{\prime}\Big)
\cdot\Big({\beta^{\prime}}^{*}(-nK_{Y'})-nG^{\prime}\Big)^{2}
=D^{\prime}\cdot M_{V^{\prime}}\cdot M_{V^{\prime}}^{\prime}\geq 0,%
$$
where $M_{V'}$ and $M_{V'}'$ are the proper transforms of $M_{Y'}$
and $M_{Y'}'$ by the birational morphism $\beta'$. We have
obtained a contradiction.
\end{proof}

\begin{lemma}
\label{lemma:n-2-lines-downstairs} If the set $\mathbb{CS}(Y',
\frac{1}{n}\mathcal{M}_{Y'})$ contains  a curve $\Gamma$ with
$-K_{Z}\cdot\Gamma=1$, then the pencil $\mathcal{M}$ is contained
in $|-K_{X}|$.
\end{lemma}

\begin{proof}
The curve $\omega(\Gamma)$ is a line in $\mathbb{P}^3$. The
restricted morphism $\omega\vert_{\Gamma}:\Gamma\to\omega(\Gamma)$
is an isomorphism. The curve $\Gamma$ contains at most one
singular point of $Y'$ if $\omega(\Gamma)$ is not contained in the
plane $f_{1}(x,y,z,t)=0$, and the curve $\omega(\Gamma)$ contains
at most three singular points of the threefold $Y'$ otherwise. We
may assume that $\Gamma$ contains $P_{1},\cdots ,P_{k}$, where
$0\leq k\leq 3$. Here, the equality $k=0$ means that $\Gamma$ lies
in the smooth locus of the threefold $Y'$.

Suppose that the line $\omega(\Gamma)$ is not contained in $R$.
Let $D$ be a general surface in $|-K_{Y'}|$ that passes through
the curve $\Gamma$. Then,
$$
\mathcal{M}_{Y'}\Big\vert_{D}=\mathrm{mult}_{\Gamma}\big(\mathcal{M}_{Y'}\big)\Gamma+\mathrm{mult}_{\Omega}\big(\mathcal{M}_{Y'}\big)\Omega+\mathcal{L},
$$
where $\mathcal{L}$ is a pencil without fixed curves and $\Omega$
is a smooth rational curve different from $\Gamma$ such that
$\omega(\Omega)=\omega(\Gamma)$. Moreover, the surface $D$ is
smooth in the outside the points $P_{1},\cdots ,P_{k}$ but the
points $P_{1},\cdots ,P_{k}$ are isolated ordinary double points
of the surface $D$. We have $\Omega^{2}=-2+\frac{k}{2}$ on the
surface $D$ and
$$
\Big(n-\mathrm{mult}_{\Omega}(\mathcal{M}_{Y'})\Big)\Omega^{2}=\Big(\mathrm{mult}_{\Gamma}(\mathcal{M}_{Y'})-n\Big)\Gamma\cdot \Omega+L\cdot\Omega=L\cdot\Omega\geq 0,%
$$
where $L$ is a general curve in $\mathcal{L}$. Therefore, the
equality $\mathrm{mult}_{\Omega}(\mathcal{M}_{Y'})=n$ holds, which
easily implies that $\mathcal{M}_{Y'}$ is a pencil in $|-K_{Y'}|$
due to Theorem~\ref{theorem:main-tool}.

Finally, we suppose that $\omega(\Gamma)$ is contained in the
ramification surface of $\omega$. It implies that $\omega(\Gamma)$
is not contained in the plane $f_{1}(x,y,z,y)=0$. The proof of
Lemma~\ref{lemma:n-2-curves-downstairs} shows the existence of a
birational morphism $\upsilon^{\prime}:W^{\prime}\to V^{\prime}$
that  contracts a single irreducible divisor $G^{\prime}$ to the
curve $\Gamma$, the surface $G^{\prime}$ contains $k$ singular
points of the threefold $V^{\prime}$ of type $\frac{1}{2}(1,1,1)$,
and $\upsilon^{\prime}$ is the blow up of $\Gamma$ at a generic
point of $\Gamma$.

Let $D^{\prime}$ be a general surface in $|-K_{V^{\prime}}|$.
Then, $\omega\circ\beta^{\prime}(D^{\prime})$ is a plane that
passes through $\omega(\Gamma)$. It implies that the base locus of
the pencil $|-K_{V^{\prime}}|$ consists of an irreducible curve
$\Gamma^{\prime}$ such that
$\beta^{\prime}(\Gamma^{\prime})=\Gamma$ and
$$
D^{\prime}\cdot\Gamma^{\prime}=-K_{V^{\prime}}^{3}=-2+\frac{k}{2}.
$$
Then, one can easily see that
$\mathcal{M}_{V^{\prime}}=|-K_{V^{\prime}}|$ by
Theorem~\ref{theorem:main-tool}. Hence, the linear system
$\mathcal{M}$ is a pencil in $|-K_{X}|$.
\end{proof}

\begin{proposition}
\label{proposition:n-2} Every Halphen pencil is contained  in
$|-K_{X}|$.
\end{proposition}
\begin{proof}
Let $\mathcal{M}_{U_i}$ be the push-forward of the pencil
$\mathcal{M}_{W_i}$ by the morphism $w_i$. Due to the previous
arguments, we may assume that
$$
P_{1}\in\mathbb{CS}\Big(Y', \frac{1}{n}\mathcal{M}_{Y'}\Big)\subseteq\Big\{P_{1},\cdots , P_{15}\Big\},%
$$
which implies that
$\mathcal{M}_{U_{1}}\sim_{\mathbb{Q}}-nK_{U_{1}}$ by Theorem~3.10
in \cite{Co00}. Therefore, each member in the pencil
$\mathcal{M}_{U_{1}}$ is contracted to a curve by the elliptic
fibration $\eta_{1}$. Therefore, the base locus of the pencil
$\mathcal{M}_{U_{1}}$ does not contain curves that are not
contracted by $\eta_{1}$. On the other hand, the singularities of
the log pair $(U_{1}, \frac{1}{n}\mathcal{M}_{U_{1}})$ are not
terminal by Theorem~\ref{theorem:Noether-Fano}.

The proof of Lemma~\ref{lemma:n-2-smooth-points} implies that the
set $\mathbb{CS}(U_{1}, \frac{1}{n}\mathcal{M}_{U_{1}})$ does not
contain a smooth point of the exceptional divisor of $\beta_1$.
Therefore, the set $\mathbb{CS}(U_{1},
\frac{1}{n}\mathcal{M}_{U_{1}})$ contains a singular point of the
threefold $U_{1}$, which implies that
$$
\Big\{P_{1},P_{i}\Big\}\subseteq\mathbb{CS}\Big(Y', \frac{1}{n}\mathcal{M}_{Y'}\Big)\subseteq\Big\{P_{1},\cdots , P_{15}\Big\},%
$$
for some $i\ne 1$. Thus, each member in the pencil
$\mathcal{M}_{U_{i}}$ is contracted to a curve by  the elliptic
fibration $\eta_{i}$, which implies that $\mathcal{M}$ is a pencil
in $|-K_{X}|$.
\end{proof}

\section{Case $\gimel=3$, hypersurface of degree $6$ in
$\mathbb{P}(1,1,1,1,3)$.} \label{section:n-3}\index{$\gimel=03$}

Let $X$ be a general hypersurface of degree $6$ in
$\mathbb{P}(1,1,1,1,3)$ with $-K_X^3=2$. It is smooth. It cannot
be birationally transformed into an elliptic fibration
(\cite{ChPa05}).

\begin{proposition}
\label{proposition:n-3} Every Halphen pencil is  contained in
$|-K_{X}|$.
\end{proposition}

\begin{proof}
It follows from Lemma~\ref{lemma:smooth-points} that the set
$\mathbb{CS}(X, \frac{1}{n}\mathcal{M})$ does not contain any
point of $X$. Hence, it must contain a curve $Z$. Then, the
in\-equa\-li\-ty $$\mathrm{mult}_{Z}(\mathcal{M})\geq n$$ holds.

For general surfaces $M_1$ and $M_2$ in $\mathcal{M}$ and a
general surface $D$ in $|-K_{X}|$, we have
$$
2n^2=M_1\cdot M_2\cdot D\geq \mathrm{mult}^2_{Z}(\mathcal{M})(-K_{X}\cdot Z)\geq n^2,%
$$
which implies that there are different surfaces $D_1$ and $D_2$ in
the linear system $|-K_{X}|$ such that the intersection $D_1\cap
D_2$ contains the curve $Z$.

Let $\mathcal{P}$ be the pencil in $|-K_{X}|$ consisting of
surfaces passing through the curve $Z$.

Suppose that $-K_X\cdot Z=2$.  For a general surface $D'$ in the
pencil $\mathcal{P}$, the inequality
\[2n=M_1\cdot D'\cdot D\geq 2n\]
implies that
$\mathrm{Supp}(M_1)\cap\mathrm{Supp}(D')\subset\mathrm{Supp}(Z)$.
It follows from Theorem~\ref{theorem:main-tool} that the linear
system $\mathcal{M}$ is the pencil in $|-K_{X}|$ consisting of
surfaces that pass through $Z$.

Now, we suppose that $-K_{X}\cdot Z=1$. The generality of $X$
implies that the general surface $D$ in $|-K_X|$ is smooth and
that the intersection $D_{1}\cap D_{2}$ consists of the curve $Z$
and an irreducible curve $Z'$ such that $Z\ne Z'$. Hence, we have
$Z^{2}=-2$ on the surface $D$ and $\mathcal{M}\vert_{D}\equiv
nZ+nZ'$. Therefore, the inequality
$\mathrm{mult}_Z(\mathcal{M})\geq n$ implies the identity
$\mathcal{M}=\mathcal{P}$ by Theorem~\ref{theorem:main-tool}.
\end{proof}

\section{Case $\gimel=4$, hypersurface of degree $6$ in
$\mathbb{P}(1,1,1,2,2)$.} \label{section:n-4}\index{$\gimel=04$}

The weighted hypersurface $X$ is defined by a general
quasihomogeneous polynomial of degree $6$ in
$\mathbb{P}(1,1,1,2,2)$ with $-K_X^3=\frac{3}{2}$. The
singularities of the hypersurface $X$ consist of points $P_{1}$,
$P_{2}$, $P_{3}$ that are quotient singularities of types
$\frac{1}{2}(1,1,1)$. The hypersurface $X$ can be given by the
equation
$$
w^{2}t+\Big(t^{2}+tf_{2}(x,y,z)+f_{4}(x,y,z)\Big)w+f_{6}(x,y,z,t)=0
$$
such that $P_{1}$ is given by the equations $x=y=z=t=0$, where
$f_{i}$ is a general quasihomogeneous polynomial of degree $i$.

There is a commutative diagram
$$
\xymatrix{
&&&U\ar@{->}[dll]_{\beta}\ar@{->}[rr]^{\alpha}&&X\ar@{-->}[ddrr]^{\psi}\ar@{-->}[ddll]_{\chi}&&Y\ar@{->}[dd]^{\eta}\ar@{->}[ll]_{\pi}\\
&W\ar@{->}[rrd]_{\omega}&&&&&\\
&&&\mathbb{P}(1,1,1,2)\ar@{-->}[rrrr]_{\xi}&&&&\mathbb{P}^{2},&}
$$
where \begin{itemize} \item $\psi$ is the natural projection,

\item $\pi$ is the composition of the Kawamata blow ups at the
points $P_{1}$, $P_{2}$, and $P_{3}$,

 \item $\eta$ is an elliptic fibration

\item $\alpha$ is the Kawamata blow up of the point $P_{1}$,

\item $\xi$ and $\chi$ are the natural projections,

\item $\beta$ is a birational morphism,

\item $\omega$ is a double cover ramified along an  octic surface
$R\subset\mathbb{P}(1,1,1,2)$.
\end{itemize}
The surface $R$ is given by the equation
$$
\Big(t^{2}+tf_{2}(x,y,z)+f_{4}(x,y,z)\Big)^{2}-4tf_{6}(x,y,z,t)=0\subset\mathbb{P}(1,1,1,2)\cong\mathrm{Proj}\big(\mathbb{C}[x,y,z,t]\big),
$$
which implies that the surface $R$ has exactly $24$ isolated
ordinary double points given by the equations
$$
t=t^{2}+tf_{2}(x,y,z)+f_{4}(x,y,z)=f_{6}(x,y,z,t)=0.
$$
The birational morphism $\beta$ contracts $24$ smooth rational
curves $C_{1},\cdots , C_{24}$ to isolated ordinary double points
of the variety $W$ that dominate the singular points of $R$.

 It easily follows from
Theorems~\ref{theorem:Noether-Fano}, \ref{theorem:Ryder},
Lemmas~\ref{lemma:curves}, \ref{lemma:Cheltsov-Kawamata}, and
\ref{lemma:smooth-points} that either the set $\mathbb{CS}(X,
\frac{1}{n}\mathcal{M})$ contains an irreducible curve passing
through a singular point of $X$ or the set $\mathbb{CS}(X,
\frac{1}{n}\mathcal{M})$ consists of a single singular point of
$X$. In particular, we may assume that the set $\mathbb{CS}(X,
\frac{1}{n}\mathcal{M})$ contains the point $P_1$.

\begin{proposition}
\label{proposition:n-4} Every Halphen pencil on $X$ is contained
in $|-K_{X}|$.
\end{proposition}

\begin{proof}

Suppose that the set $\mathbb{CS}(X, \frac{1}{n}\mathcal{M})$
contains an irreducible curve $Z$ that passes through $P_{1}$.
Then, it follows from Theorem~\ref{theorem:Ryder} that the linear
system $\mathcal{M}$ is a pencil in $|-K_{X}|$ in the case when
$-K_{X}\cdot Z=\frac{3}{2}$. Therefore, we may assume that the
curve $Z$ is contracted by the rational map $\psi$ to a point.
Also, we may assume that  either $-K_{X}\cdot Z=\frac{1}{2}$ or
$-K_{X}\cdot Z=1$.

Let $\mathcal{B}$ be the pencil in $|-K_{X}|$ consisting of
surfaces passing through $Z$. In addition, let $B$ and
$B^{\prime}$ be general surfaces in $\mathcal{B}$. Then, the cycle
$B\cdot B^{\prime}$ is reduced and contains the curve $Z$. Put
$\tilde{Z}=B\cdot B^{\prime}$ and let $\tilde{Z}_W$ be the image
of the curve $\tilde{Z}_U$ by the birational morphism $\beta$.
Then, $\omega(\tilde{Z}_{W})$ is a ruling of the cone
$\mathbb{P}(1,1,1,2)$. In particular, the curve
$\omega(\tilde{Z}_{W})$ contains at most one singular point of the
surface $R$.

There are exactly $24$ rulings of the cone $\mathbb{P}(1,1,1,2)$
that pass through the singular points of the surface $R$. Thus, we
may assume that the curve $\tilde{Z}_{W}$ is irreducible in the
case when the curve $\omega(\tilde{Z}_{W})$ passes through a
singular point of the surface $R$. Moreover, the surface $B_{W}$
that is the image of the surface $B_U$ by $\beta$
 has an isolated ordinary double point at the point
$\beta(C_{i})$ in the case when
$\omega\circ\beta(C_{i})\in\omega(\tilde{Z}_{W})$. Therefore, the
cycle $\tilde{Z}$ consists of two irreducible components.

Let $\bar{Z}$ be the irreducible component of $\tilde{Z}$ that is
different from $Z$. Then, the generality of the hypersurface $X$
implies that $\bar{Z}^{2}<0$ on the surface $B$, but
$M\vert_{B}\equiv nZ+n\bar{Z}$. On the other hand, we have
$$
M\Big\vert_{B}=m_{1}Z+m_{2}\bar{Z}+F,
$$
where and $m_{1}$ and $m_{2}$ are natural numbers and $F$ is an
effective divisor on $B$ whose support contains neither the curve
$Z$ nor the curve $\bar{Z}$. We have
$$m_{1}\geq\mathrm{mult}_{Z}(\mathcal{M})\geq n$$ and
$$
(n-m_{2})\bar{Z}\equiv F+(m_{1}-n)Z,
$$
They imply that $m_{2}=m_{1}=n$ and the support of the cycle
$M\cdot B$ is contained in $Z\cup\bar{Z}$. Therefore, the identity
$\mathcal{M}=\mathcal{B}$ follows from
Theorem~\ref{theorem:main-tool}.

For now, we suppose that the set $\mathbb{CS}(X,
\frac{1}{n}\mathcal{M})$ consists of the point $P_{1}$. It follows
from Lemma~\ref{lemma:Kawamata} that
$\mathcal{M}_{U}\sim_{\mathbb{Q}} -nK_{U}$. Therefore, the set
$\mathbb{CS}(U, \frac{1}{n}\mathcal{M}_{U})$ is not empty by
Theorem~\ref{theorem:Noether-Fano}. Let $E$ be the exceptional
divisor of $\alpha$. Then, $E\cong\mathbb{P}^{2}$ and the set
$\mathbb{CS}(U, \frac{1}{n}\mathcal{M}_{U})$ contains a line $L$
on the surface $E$ by Lemma~\ref{lemma:Cheltsov-Kawamata}.

Let $Z$ be the curve $S^t_{U}\cap E$. Then, $Z$ does not contain
the curve $L$, the surface $S^t_{U}$ contains every curve $C_{i}$,
and  the curve $Z$ is a smooth plane quartic curve. The
hypersurface $X$ is general by assumption. In particular, the
surface $S^t_{U}$ is smooth along the curve $C_{i}$, the morphism
$\beta\vert_{S^t_{U}}$ contracts the curve $C_{i}$ to a smooth
point of the surface $S^t_{W}$ which is the image of $S^t_U$ by
$\beta$. Moreover, we may assume that the intersection $L\cap Z$
contains at least one point of the curve $Z$ that is not contained
in $\cup_{i=1}^{24}C_{i}$. Indeed, it is enough to assume that the
set $\cup_{i=1}^{24}(C_{i}\cap Z)$ does not contain bi-tangent
points of the plane quartic curve $Z$.

Let $M^{\prime}$ be a general surface in $\mathcal{M}$ and $D$ be
a general surface in $|-2K_{U}|$. Then,
$$
2n^{2}=D\cdot M_{U}\cdot M_{U}^{\prime}\geq2\mathrm{mult}_{L}(M_{U}\cdot M_{U}^{\prime})\geq 2\mathrm{mult}_{L}(M_{U})\mathrm{mult}_{L}(M_{U}^{\prime})\geq 2n^{2},%
$$
which implies that the support of the cycle $M_{U}\cdot
M_{U}^{\prime}$ is contained in the union of the curve $L$ and
$\cup_{i=1}^{24}C_{i}$. Hence, we have
$$
\mathcal{M}_{U}\Big\vert_{S^t_{U}}=\mathcal{D}+\sum_{i=1}^{24}m_{i}C_{i},
$$
where $m_{i}$ is a natural number and $\mathcal{D}$ is a pencil
without fixed components. Let $P$ be a point of $L\cap Z$ that is
not contained in $\cup_{i=1}^{24}C_{i}$. For general curves
$D_{1}$ and $D_{2}$ in $\mathcal{D}$,
$$
n^{2}-\sum_{i=1}^{24}m_{i}^{2}=D_{1}\cdot D_{2}\geq\mathrm{mult}_{P}(D_{1})\mathrm{mult}_{P}(D_{2})\geq n^{2},%
$$
which implies that $m_{1}=m_2=\cdots=m_{24}=0$. Therefore, we have
$M_{U}\cdot M_{U}^{\prime}=n^{2}L$, which is impossible because
the suppose of the cycle $M\cdot M^{\prime}$ must contain a curve
on $X$.
\end{proof}

\section{Case $\gimel=5$, hypersurface of degree $7$ in
$\mathbb{P}(1,1,1,2,3)$.} \label{section:n-5}\index{$\gimel=05$}

The threefold $X$ is a general hypersurface of degree $7$ in
$\mathbb{P}(1,1,1,2,3)$ with $-K_{X}^{3}=\frac{7}{6}$.  The
singularities of the hypersurface $X$ consist of  two points $P$
and $Q$ that the are quotient singularities of types
$\frac{1}{2}(1,1,1)$ and  $\frac{1}{3}(1,1,2)$, respectively. The
hypersurface $X$ can be given by the equation
$$
w^{2}z+f_{4}(x,y,z,t)w+f_{7}(x,y,z,t)=0,
$$
where  $f_{i}$ is a quasihomogeneous polynomial of degree $i$.
Hence, the point $P$ is given by the equations $x=y=z=w=0$ and the
point $Q$ is given by the equations $x=y=z=t=0$.

There is a commutative diagram
$$
\xymatrix{
&&W\ar@{->}[dll]_{\pi}\ar@{->}[d]^{\alpha}&&&\\%
Y\ar@{->}[dr]_{\omega}&&X\ar@{-->}[dl]_{\xi}\ar@{-->}[dr]^{\psi}&&U\ar@{->}[llu]_{\beta}&\\
&\mathbb{P}(1,1,1,2)\ar@{-->}[rr]_{\chi}&&\mathbb{P}^{2},&&V\ar@{->}[ll]^{\eta}\ar@{->}[lu]_{\gamma}}
$$
where \begin{itemize}

\item $\psi$ is the natural projection,

\item $\alpha$ is the Kawamata blow up at the point $Q$ with
weights $(1,1,2)$,

\item $\beta$ is the Kawamata blow up with weights $(1,1,1)$ at
the point of $W$ whose image to $X$ is the point $P$,

\item $\gamma$ is the Kawamata blow up with weights $(1,1,1)$ at
the singular point of $U$ whose image to $X$ is  the point $Q$,

\item $\eta$ is an elliptic fibration,

\item $\pi$ is a birational morphism,

\item $\chi$ and $\xi$ are the natural projections,

\item $\omega$ is a double cover of $\mathbb{P}(1,1,1,2)$ ramified
along an octic surface $R$.

\end{itemize}
The generality of $X$ implies that the birational morphism $\pi$
contracts $14$ smooth irreducible rational curves $C_{1},\cdots ,
C_{14}$ into $14$ isolated ordinary double points $P_{1},\cdots ,
P_{14}$ of the variety $Y$, respectively. The double cover
$\omega$ is branched over the octic surface
$R\subset\mathbb{P}(1,1,1,2)$ that is given by the equation
$$
f_{4}(x,y,z,t)^{2}-4zf_{7}(x,y,z,t)=0\subset\mathbb{P}(1,1,1,2)\cong\mathrm{Proj}\big(\mathbb{C}[x,y,z,t]\big),
$$
which has $14$ isolated ordinary double points
$\omega(P_{1}),\cdots,\omega(P_{14})$.

Now let us prove the following result, which is due to
\cite{Ry02}.

\begin{proposition}
\label{proposition:n-5} Every Halphen pencil is contained in
$|-K_{X}|$.
\end{proposition}

\begin{proof}
It follows from Theorem~\ref{theorem:Ryder} and the generality of
$X$ that the linear system $\mathcal{M}$ is a pencil in $|-K_{X}|$
if the set $\mathbb{CS}(X, \frac{1}{n}\mathcal{M})$ contain a
curve (the proof of Proposition~\ref{proposition:n-4}). Hence, we
may assume that the set $\mathbb{CS}(X, \frac{1}{n}\mathcal{M})$
consists of singular points of $X$ by
Lemma~\ref{lemma:smooth-points}.

It easily follows from Theorem~\ref{theorem:Noether-Fano},
Lemmas~\ref{lemma:curves}, and \ref{lemma:Cheltsov-Kawamata} that
$\mathbb{CS}(X, \frac{1}{n}\mathcal{M})\ne \{P\}$.

Suppose that $\mathbb{CS}(X, \frac{1}{n}\mathcal{M})=\{P, Q\}$.
Let $E$ be the exceptional divisor of $\alpha$. Then, the surface
$E$ is a quadric cone. It follows from
Theorem~\ref{theorem:Noether-Fano}, Lemmas~ \ref{lemma:curves},
\ref{lemma:Kawamata}, and \ref{lemma:Cheltsov-Kawamata} that
$\mathcal{M}_{U}\sim_{\mathbb{Q}} -nK_{U}$ and the set
$\mathbb{CS}(U, \frac{1}{n}\mathcal{M}_{U})$ contains a ruling $L$
of the cone $E_U$. Hence, it follows from
Theorem~\ref{theorem:main-tool} and the proof of
Lemma~\ref{lemma:curves} imply that $\mathcal{M}_{U}$ is the
pencil consisting of surfaces in $|-K_{U}|$ that contain the curve
$L$ because $-K_{U}\cdot L=\frac{1}{2}$ and
$-K_{U}^{3}=\frac{1}{2}$.

Now, we suppose that $\mathbb{CS}(X,
\frac{1}{n}\mathcal{M})=\{Q\}$. Let $O$ be the singular point of
the variety $W$ whose image to $X$ is the point $Q$ and $\bar{P}$
be the singular point of $W$ whose image to $X$ is the point $P$.
It follows from Theorem~\ref{theorem:Noether-Fano},
Lemmas~\ref{lemma:curves}, \ref{lemma:Kawamata}, and
\ref{lemma:Cheltsov-Kawamata} that
$\mathcal{M}_{Y}\sim_{\mathbb{Q}} -nK_{Y}$ and that the set
$\mathbb{CS}(Y, \frac{1}{n}\mathcal{M}_{Y})$ contains an
irreducible curve $Z$ such that $Z$ is the image of a ruling of
the cone $E$. Hence, the equality $-K_{Y}\cdot Z=\frac{1}{2}$
holds and $\chi\circ\omega(Z)$ is a point.

The variety $\mathbb{P}(1,1,1,2)$ is a cone over the Veronese
surface. Hence, the curve $\omega(Z)$ is a ruling of the cone
$\mathbb{P}(1,1,1,2)$ and the point
$(\omega\circ\pi)(\bar{P})=(\omega\circ\pi)(O)$ is the vertex of
$\mathbb{P}(1,1,1,2)$. The generality of the hypersurface $X$
implies the existence of the irreducible curve $\bar{Z}$ on the
variety $Y$ such that $\bar{Z}\ne Z$, $\omega(\bar{Z})=\omega(Z)$,
and $\pi(\bar{P})\in \bar{Z}$.

Let $D$ be a general surface in $|-K_{Y}|$ that contains $Z$.
Then, the inequality $\bar{Z}^{2}<0$ holds on $D$, but
$\mathcal{M}_{Y}\vert_{D}\equiv n\bar{Z}+nZ$. Now the proof of
Proposition~\ref{proposition:n-4} implies that $\mathcal{M}_{Y}$
is the pencil consisting of surfaces in $|-K_{Y}|$ that contain
the curve $Z$.
\end{proof}

\section{Case $\gimel=6$, hypersurface of degree $8$ in
$\mathbb{P}(1,1,1,2,4)$.} \label{section:n-6}\index{$\gimel=06$}

The threefold $X$ is a general hypersurface  of degree $8$ in
$\mathbb{P}(1,1,1,2,4)$ with $-K_X^3=1$.  Its singularities
consist of points $P$ and $Q$ that are quotient singularities of
type $\frac{1}{2}(1,1,1)$.

We have a commutative diagram
$$
\xymatrix{
&&U\ar@{->}[dl]_{\pi}\ar@{->}[dr]^{\eta}&&\\%
&X\ar@{-->}[rr]_{\psi}&&\mathbb{P}^{2},&}
$$
where \begin{itemize} \item $\psi$ is the natural projection,

\item $\pi$ is the composition of the Kawamata blow ups at the
singular points $P$ and $Q$,

\item $\eta$ is an elliptic fibration.
\end{itemize}

\begin{proposition}
\label{proposition:n-6}  Every Halphen pencil is contained  in
$|-K_{X}|$.
\end{proposition}
\begin{proof}
The log pair $(X, \frac{1}{n}\mathcal{M})$ has terminal
singularities at a smooth point of the hypersurface $X$ by
Lemma~\ref{lemma:smooth-points}. Moreover, it easily follows from
Theorem~\ref{theorem:Noether-Fano}, Lemmas~\ref{lemma:curves}, and
\ref{lemma:Cheltsov-Kawamata} that the set $\mathbb{CS}(X,
\frac{1}{n}\mathcal{M})$ contains an irreducible curve $Z$. Then,
it follows from Theorem~\ref{theorem:Ryder} that the curve $Z$ is
a fiber of the projection $\psi$, which easily implies that the
linear system $\mathcal{M}$ is a pencil in the linear system
$|-K_{X}|$ by Theorem~\ref{theorem:main-tool} in the case when the
equality $-K_{X}\cdot Z=1$ holds.

To conclude the proof, we assume that $-K_{X}\cdot Z=\frac{1}{2}$.
Let $D$ be a general surface in $|-K_{X}|$ that contains $Z$.
Then, the surface $D$ is smooth in the outside of the points $P$
and $Q$ which are isolated ordinary double points on $D$. Let $F$
be a fiber of the rational map $\psi$ over the point $\psi(Z)$.
Then, $F$ consists of two irreducible components. Let $\bar{Z}$ be
the component of $F$ different from $Z$. Then, the generality of
$X$ implies that $\bar{Z}^{2}<0$ and the proof of
Proposition~\ref{proposition:n-4} implies that the pencil
$\mathcal{M}$ consists  of surfaces in $|-K_{X}|$ that contain the
curve $Z$.
\end{proof}

\section{Case $\gimel=8$, hypersurface of degree $9$ in
$\mathbb{P}(1,1,1,3,4)$.} \label{section:n-8}\index{$\gimel=08$}

 The threefold  $X$ is a
 general hypersurface of degree $9$ in $\mathbb{P}(1,1,1,3,4)$
with $-K_{X}^{3}=\frac{3}{4}$. The singularities of the
hypersurface $X$ consist of the singular point $O$ that is a
quotient singularity of type $\frac{1}{4}(1,1,3)$. The
hypersurface $X$ can be given by the equation
$$
w^{2}z+\Big(f_{2}\big(x,y,z\big)t+f_{5}\big(x,y,z\big)\Big)w+f_{9}\big(x,y,z,t\big)=0,
$$
where $f_{i}$ is a quasihomogeneous polynomial of degree $i$.
Thus, the point $O$ is given by the equations $x=y=z=t=0$.
Furthermore, we may assume that the polynomials $f_{2}(x,y,0)$ and
$f_{5}(x,y,0)$ are co-prime.

There is a commutative diagram
$$
\xymatrix{
&&&W\ar@{->}[dll]_{\sigma}\ar@{->}[d]^{\alpha}&&U\ar@{->}[ll]_{\beta}&&&\\%
&Y\ar@{->}[dr]_{\omega}&&X\ar@{-->}[dl]_{\xi}\ar@{-->}[drrr]^{\psi}&&&&V\ar@{->}[ld]^{\eta}\ar@{->}[llu]_{\gamma}&\\
&&\mathbb{P}(1,1,1,3)\ar@{-->}[rrrr]_{\chi}&&&&\mathbb{P}^{2},&&}
$$
where \begin{itemize}

\item $\xi$, $\psi$ and $\chi$ are the natural projections,

\item $\alpha$ is the Kawamata blow up at the point $O$ with
weights $(1,1,3)$,

\item $\beta$ is the Kawamata blow up at the singular point of the
variety $W$ that is a quotient singularity of type
$\frac{1}{3}(1,1,2)$,

\item $\gamma$ is the Kawamata blow up at the singular point of
the variety $U$ that is a quotient singularity of type
$\frac{1}{2}(1,1,1)$,

\item $\eta$ is an elliptic fibration,

\item $\sigma$ is a birational morphism that contracts $15$ smooth
rational curves to $15$ isolated ordinary double points
$P_{1},\cdots , P_{15}$ of the variety $Y$,

\item $\omega$ is a double cover of $\mathbb{P}(1,1,1,3)$ branched
over a surface $R$ of degree $10$.
\end{itemize}

The surface $R$ is given by the equation
$$
\Big(f_{2}\big(x,y,z\big)t+f_{5}\big(x,y,z\big)\Big)^{2}-4zf_{9}\big(x,y,z,t\big)=0\subset\mathbb{P}\big(1,1,1,3\big)\cong\mathrm{Proj}\Big(\mathbb{C}[x,y,z,t]\Big).
$$
It has $15$ ordinary double points
 given by
$z=tf_{2}+f_{5}=f_{9}=0$. Let $P_{1},\cdots,P_{15}$ be the points
of $Y$ whose image via the double cover $\omega$ are the $15$
ordinary double points of the surface $R$.

Let $E$ be the exceptional divisor of $\alpha$ and $F$ be the
exceptional divisor of $\beta$. In addition, let $P$ be the
singular point of $W$ and  $Q$ be the singular point of $U$. Then,
the surface $\omega\circ\sigma(E)$ is given by the equation $z=0$
and $\mathbb{P}(1,1,1,3)$ is a cone whose vertex is the point
$\omega\circ\sigma(P)$. The generality of the polynomials $f_{5}$
and $f_{2}$ implies that the surface $R$ does not contain the
rulings of $\mathbb{P}(1,1,1,3)$ that are contained in the surface
$\omega\circ\sigma(E)$.

It follows from Lemma~\ref{lemma:smooth-points} that the set
$\mathbb{CS}(X, \frac{1}{n}\mathcal{M})$ does not contain smooth
points of the hypersurface $X$. Therefore,  by
Theorem~\ref{theorem:Ryder} and Lemma~\ref{lemma:Kawamata}, it
must contain the point $O$.

\begin{lemma}
\label{lemma:n-8-curves} If the set $\mathbb{CS}(X,
\frac{1}{n}\mathcal{M})$ contains a curve, then $\mathcal{M}$ is a
pencil in $|-K_{X}|$.
\end{lemma}

\begin{proof}
Let $L$ be a curve in $\mathbb{CS}(X, \frac{1}{n}\mathcal{M})$.
Then, it follows from Theorem~\ref{theorem:Ryder} that there are
two different surfaces $D$ and $D^{\prime}$ in $|-K_{X}|$ such
that $L$ is a component of the cycle $D\cdot D^{\prime}$.
Moreover, the generality of $X$ implies that the cycle $D\cdot
D^{\prime}$ is reduced and contains at most two components.

Let $\mathcal{P}$ be the pencil in $|-K_{X}|$  generated by
surfaces $D$ and $D^{\prime}$. From
Theorem~\ref{theorem:main-tool} and the proof of
Lemma~\ref{lemma:curves}, we obtain $\mathcal{M}=\mathcal{P}$ if
$-K_{X}\cdot L=\frac{3}{4}$. Hence, we may assume that either
$-K_{X}\cdot L=\frac{1}{4}$ or $-K_{X}\cdot L=\frac{1}{2}$. Thus,
the cycle $D\cdot D^{\prime}$ contains a component $L^{\prime}$
such that
$$
-K_{X}\cdot L^{\prime}=\frac{3}{4}+K_{X}\cdot L
$$
and $L^{\prime}\ne L$. We consider only the case when $-K_{X}\cdot
L=\frac{1}{4}$ because the case $-K_{X}\cdot L=\frac{1}{2}$ is
simpler and similar.

The proper transform $S^z_{W}$ contains the curve $L_{W}$ because
$$
S^z_{W}\sim_{\mathbb{Q}} \alpha^*(-K_{X})-\frac{5}{4}E, \ \ \
E\cdot L_{W}\geq \frac{1}{3}.$$ Thus, either the curve $L_{W}$ is
contracted by $\sigma$ or the curve $\omega(L_{Y})$ is a ruling of
the cone $\mathbb{P}(1,1,1,3)$ contained in the surface
$\omega\circ\sigma(E)$, where $L_Y$ is the image of $L_W$ by
$\sigma$.

Suppose that the curve $L_{W}$ is not contracted by $\sigma$.
Then, the curve $\omega(L_{Y})$ is not contained in the surface
$R$, which implies that $\omega(L_{Y})$ contains at most one
singular point of the surface $R$  different from the point
$\omega\circ\sigma(P)$. Moreover, the curve $\omega(L_{Y})$ must
contain a singular point of $R$ different from
$\omega\circ\sigma(P)$ because $-K_{X}\cdot L=\frac{3}{4}$
otherwise. Thus, we may assume that the curve  $\omega(L_{Y})$
contains the point $\omega(P_{1})$.

Let $D_Y=\sigma(D_W)$ and $D_Y'=\sigma(D_W')$. Then, the point
$P_{1}$ is an isolated ordinary double point of the surface
$D_{Y}$. Thus, wee see that the curve $L^{\prime}_{W}$ is
contracted to the point $P_{1}$ by $\sigma$ and
$$
D_{Y}\cdot D^{\prime}_{Y}=L_{Y}+\bar{L}_{Y},
$$
where $\bar{L}_{Y}$ is a ruling of $E\cong\mathbb{P}(1,1,3)$. In
particular, we have $-K_{X}\cdot L^{\prime}=\frac{1}{4}$, which
contradicts the equality $-K_{X}\cdot L=\frac{1}{4}$. Hence, the
curve $L_{W}$ is  contracted by $\sigma$, which implies that the
curve $L_{W}^{\prime}$ is not contracted by $\sigma$ and the curve
$\omega(L^{\prime}_{Y})$ is a ruling of the cone
$\mathbb{P}(1,1,1,3)$ that is contained in the surface
$\omega\circ\sigma(E)$, where $L_Y'$ is the image of the curve
$L'_W$ by $\sigma$. The curve $\omega(L^{\prime}_{Y})$ is not
contained in the surface $R$. It implies that
$\omega(L^{\prime}_{Y})$ contains at most one singular point of
the surface $R$ different from the point $\omega\circ\sigma(P)$.
The curve $\omega(L^{\prime}_{Y})$ must contain a singular point
of $R$ different from $\omega\circ\sigma(P)$ because $-K_{X}\cdot
L^{\prime}\ne \frac{3}{4}$. Thus, we may assume that the curve
$\omega(L^{\prime}_{Y})$ contains the point $\omega(P_{1})$.

The point $P_{1}$ is an isolated ordinary double point of the
surface $D_{Y}$ and the curve $L_{W}$ is contracted to the point
$P_{1}$ by $\sigma$. Hence, we have
$$
D_{Y}\cdot D^{\prime}_{Y}=L^{\prime}_{Y}+\bar{L}_{Y}',
$$
where $\bar{L}_{Y}'$ is a ruling of $E\cong\mathbb{P}(1,1,3)$.
Therefore, the intersection $L_{W}\cap L_{W}^{\prime}$ consists of
a point $O^{\prime}$ such that $O^{\prime}\not\in E$, and hence
the intersection $L\cap L'$ contains the point
$\alpha(O^{\prime})$ that is different from $O$.

The surface $D$ is normal and it is smooth at the point
$\alpha(O^{\prime})$. On the other hand, the equality
$(L+L^{\prime})\cdot L^{\prime}=\frac{1}{2}$ holds on the surface
$D$, which implies that the inequality $L^{\prime}\cdot
L^{\prime}<0$ holds on the surface $D$. Therefore, we have
$$
\mathcal{M}\Big\vert_{D}=m_{1}L+m_{2}L^{\prime}+\mathcal{L}\equiv  nL+nL^{\prime},%
$$
where $\mathcal{L}$ is a pencil on $D$ that does not have fixed
components, and $m_1$ and $m_{2}$ are natural numbers such that
$m_{1}\geq n$. In particular, we have
$$
0\leq(m_{1}-n)L^{\prime}\cdot L+\mathcal{L}\cdot L^{\prime}= (n-m_{2})L^{\prime}\cdot L^{\prime},%
$$
which implies that $m_{2}=m_{1}=n$ and
$\mathcal{M}\vert_{D}=nL+nL^{\prime}$ because $L^{\prime}\cdot
L^{\prime}<0$. It follows from Theorem~\ref{theorem:main-tool}
that $\mathcal{M}=\mathcal{P}$.
\end{proof}

Therefore, we may assume that the set $\mathbb{CS}(X,
\frac{1}{n}\mathcal{M})$ consists of the singular point $O$, which
implies that the set $\mathbb{CS}(W, \frac{1}{n}\mathcal{M}_{W})$
contains the point $P$ by Theorem~\ref{theorem:Noether-Fano} and
Lemma~\ref{lemma:Cheltsov-Kawamata}.

\begin{lemma}
\label{lemma:n-8-second-floor} If the set  $\mathbb{CS}(W,
\frac{1}{n}\mathcal{M}_{W})$ consists of the point $P$, then
$\mathcal{M}$ is a pencil in $|-K_{X}|$.
\end{lemma}

\begin{proof}
Our assumption implies that the set $\mathbb{CS}(U,
\frac{1}{n}\mathcal{M}_{U})$ contains the point $Q$ by
Theorem~\ref{theorem:Noether-Fano} and
Lemma~\ref{lemma:Cheltsov-Kawamata}, and hence the set
$\mathbb{CS}(V, \frac{1}{n}\mathcal{M}_{V})$ is not empty by
Theorem~\ref{theorem:Noether-Fano}. However, the set
$\mathbb{CS}(V, \frac{1}{n}\mathcal{M}_{V})$ does not contain any
subvariety of the exceptional divisor of $\gamma$ by
Lemmas~\ref{lemma:curves} and \ref{lemma:Cheltsov-Kawamata}. Thus,
the set $\mathbb{CS}(U \frac{1}{n}\mathcal{M}_{U})$ contains an
element different from the point $Q$.

The surface $F$ is a quadric cone and it follows from
Lemma~\ref{lemma:Cheltsov-Kawamata} that the set $\mathbb{CS}(U,
\frac{1}{n}\mathcal{M}_{U})$ contains a ruling $L$ of the cone
$F$. Let $\mathcal{B}$ be the linear system consisting of surfaces
in $|-K_{U}|$ that contain the curve $L$. Then, $\mathcal{B}$ is a
pencil because the curve $L$ is contracted by the map
$\eta\circ\gamma^{-1}$ to a point.

Let $D$ be a general surface in $|-rK_{U}|$ for $r\gg 0$ and
$M_{U}$ and $B$ be general surfaces in $\mathcal{M}_{W}$ and
$\mathcal{B}$, respectively. Then,
$$
\frac{rn}{2}=D\cdot M_{U}\cdot B\geq \mathrm{mult}_{L}(M_{U}\cdot
B)(D\cdot L)\geq \mathrm{mult}_{L}(M_{U})
\mathrm{mult}_{L}(B)(D\cdot L)\geq\frac{rn}{2},
$$
which implies that the support of the effective cycle $M_{U}\cdot
B$ consists of the curve $L$ and a cycle $Z$ such that $D\cdot
Z=0$. On the other hand, the divisor $-K_{U}$ is big and big.
Hence, it follows from Theorem~\ref{theorem:main-tool} that
$\mathcal{M}_{U}=\mathcal{B}$, which implies that  $\mathcal{M}$
is a pencil in $|-K_{X}|$.
\end{proof}

\begin{proposition}
\label{proposition:n-8} Every Halphen pencil is contained in
$|-K_{X}|$
\end{proposition}
\begin{proof}
Due to the previous arguments, we may assume that the set
$\mathbb{CS}(W, \frac{1}{n}\mathcal{M}_{W})$ contains a subvariety
$Z$ different from the point $P$. Then, $Z$ is contained in the
surface $E$ that is a cone over the smooth rational curve of
degree $3$. Moreover, it follows from
Lemma~\ref{lemma:Cheltsov-Kawamata} that $Z$ is a ruling of $E$.
Put $\bar{Z}=\sigma(Z)$. Then, $\bar{Z}\in\mathbb{CS}(Y,
\frac{1}{n}\mathcal{M}_{Y})$ and $-K_{Y}\cdot
\bar{Z}=\frac{1}{3}$.

The curve $\omega(\bar{Z})$ is a ruling of the cone
$\mathbb{P}(1,1,1,3)$. Let $\mathcal{B}$ be the linear system
consisting of surfaces in $|-K_{Y}|$ that contain the curve
$\bar{Z}$. Then, $\mathcal{B}$ is a pencil whose base locus
consists of the curve $\bar{Z}$ and an irreducible smooth rational
curve $L$ on the variety $Y$ such that $L$ is different from the
curve $\bar{Z}$ and $\omega(L)=\omega(\bar{Z})$.

Let $B$ be a general surface in $\mathcal{B}$. Then, $B$ is smooth
in the outside of the points
$$
\sigma(P)\cup\Big(\{P_{1},\cdots ,P_{15}\}\cap (\bar{Z}\cup L)\Big),%
$$
and every singular point of $B$ different from $\sigma(P)$ is an
isolated ordinary double point.

The generality of $X$ implies that the curve $\omega(L)$ does not
contain more than one singular point of $R$ different from the
point $\omega\circ\sigma(P)$. Thus, arguing as in the proof of
Lemma~\ref{lemma:n-8-curves}, we see that the inequality $L^{2}<0$
holds on the surface $B$ if the intersection $L\cap \bar{Z}$
contains a point different from $\sigma(P)$. On the other hand,
the curve $\omega(L)$ does not contain singular points of $R$
different from the point $\omega\circ\sigma(P)$ if $L\cap
\bar{Z}=\{\sigma(P)\}$. Thus,  the inequality $L^{2}<0$ holds on
the surface $B$ as well if the intersection $L\cap \bar{Z}$
consists of the point $\sigma(P)$ because the curve $L$ is an
image of the curve $\bar{Z}$ via the biregular involution of the
surface $B$ and the curve $\bar{Z}$ is contracted on the surface
$B$.

The equivalence $\mathcal{M}_{Y}\vert_{B}\equiv n\bar{Z}+nL$ holds
on the surface $B$, which implies that the support of the cycle
$M_{Y}\cdot B$ is the union of the curves $\bar{Z}$ and $L$
because $\mathrm{mult}_{\bar{Z}}(\mathcal{M}_{Y})\geq n$. Hence,
it follows from Theorem~\ref{theorem:main-tool} that
$\mathcal{M}_{Y}=\mathcal{B}$. Thus, the linear system
$\mathcal{M}$ is a pencil in $|-K_{X}|$. We have completed the
proof.
\end{proof}

\section{Case $\gimel=10$, hypersurface of degree $10$ in
$\mathbb{P}(1,1,1,3,5)$.} \label{section:n-10}\index{$\gimel=10$}

The threefold $X$ is a general hypersurface of degree $10$ in
$\mathbb{P}(1,1,1,3,5)$ with $-K_{X}^{3}=\frac{2}{3}$. The
singularities of the hypersurface $X$ consist of one point $O$
that is a quotient singularity of type $\frac{1}{3}(1,1,2)$.

The hypersurface $X$ can be given by the equation
$$
w^{2}=t^{3}z+t^{2}f_{4}(x,y,z)+tf_{7}(x,y,z)+f_{10}(x,y,z),
$$
where  $f_{i}$ is a general quasihomogeneous polynomial of degree
$i$. In particular, we may assume that the polynomials
$f_{4}(x,y,0)$ and $f_{7}(x,y,0)$ are co-prime  and the polynomial
$f_{7}^{2}(x,y,0)-f_{4}(x,y,0)f_{10}(x,y,0)$ is reduced, i.e., it
has $14$ distinct linear factors.

There is a commutative diagram
$$
\xymatrix{
&&&&W\ar@{->}[dll]_{\alpha}&&&\\%
&&X\ar@{->}[dl]_{\xi}\ar@{-->}[drrr]^{\psi}&&&&Y\ar@{->}[ld]^{\eta}\ar@{->}[llu]_{\beta}&\\
&\mathbb{P}(1,1,1,3)\ar@{-->}[rrrr]_{\chi}&&&&\mathbb{P}^{2},&&}
$$
where \begin{itemize}

\item $\xi$, $\psi$, and $\chi$ are the natural projections,

\item  $\alpha$ is the Kawamata blow up at the point  $O$ with
weights $(1,1,2)$,

\item $\beta$ is the Kawamata blow up at the singular point of the
variety $W$,

\item $\eta$ is an elliptic fibration.

\end{itemize}

\begin{lemma}
\label{lemma:n-10-curves-exclusion} If the set $\mathbb{CS}(X,
\frac{1}{n}\mathcal{M})$  contains a curve, then $\mathcal{M}$ is
a pencil in $|-K_{X}|$.
\end{lemma}

\begin{proof}
Let $L$ be a curve on $X$ that is contained in $\mathbb{CS}(X,
\frac{1}{n}\mathcal{M})$. Then, $-K_{X}\cdot L\leq \frac{2}{3}$ by
Lemma~\ref{lemma:curves}. Moreover, the proof of
Lemma~\ref{lemma:curves} together with
Theorem~\ref{theorem:main-tool} and Lemma~\ref{theorem:Ryder}
implies  that $\mathcal{M}$ is a pencil in $|-K_{X}|$ in the case
when $-K_{X}\cdot L=\frac{2}{3}$. Hence, we may assume that
$-K_{X}\cdot L=\frac{1}{3}$, which implies that the curve $L$ is
contracted by the rational map $\psi$ to a point.

The variety $\mathbb{P}(1,1,1,3)$ is a cone whose vertex is the
point $\xi(O)$. The curve $\xi(L)$ is a ruling of the cone
$\mathbb{P}(1,1,1,3)$. The generality of the hypersurface $X$
implies that $\xi(L)$ is not contained in the ramification divisor
of $\xi$. Thus, there is an irreducible curve $\bar{L}$ on the
variety $X$ such that $\bar{L}$ is different from $L$ but
$\xi(L)=\xi(\bar{L})$.

Because $-K_{W}\cdot L_{W}=0$,  the curve $L$ is one of $28$
curves that are cut out on $X$ by the equations
$z=f_{7}^{2}(x,y,z)-4f_{4}(x,y,z)f_{10}(x,y,z)=0$. Therefore, the
generality of the hypersurface $X$ implies that the intersection
$L\cap \bar{L}$ consists of the point $O$ and another distinct
point $P$.

Let $\mathcal{B}$ be the pencil consisting of surfaces in
$|-K_{X}|$ that contain both $L$ and $\bar{L}$ and $B$ be a
general surface in $\mathcal{B}$. Then, $B$ is smooth at the point
$P$. Thus, the equality $(L+\bar{L})\cdot\bar{L}=\frac{1}{3}$
holds on the surface $B$, which implies that $\bar{L}^{2}<0$. On
the other hand, we have
$$
\mathcal{M}\Big\vert_{B}= m_{1}L+m_{2}\bar{L}+\mathcal{L}\equiv  nL+n\bar{L},%
$$
where $\mathcal{L}$ is a pencil on $B$ without fixed components,
and $m_1$ and $m_{2}$ are natural numbers such that $m_{1}\geq n$.
Now the inequalities $\bar{L}^{2}<0$ and $m_{1}\geq n$ imply that
$m_{2}=m_{1}=n$ and $\mathcal{M}\vert_{B}= nL+n\bar{L}$.
Therefore, it follows from Theorem~\ref{theorem:main-tool} that
$\mathcal{M}=\mathcal{B}$.
\end{proof}

\begin{proposition}
\label{proposition:n-10} Every Halphen pencil is contained in
$|-K_{X}|$
\end{proposition}
\begin{proof}
Due to Lemmas~\ref{lemma:smooth-points} and
\ref{lemma:n-10-curves-exclusion}, we may assume that
$\mathbb{CS}(X, \frac{1}{n}\mathcal{M})=\{O\}$. Let $Q$ be the
unique singular point of the variety $W$. Then, the set
$\mathbb{CS}(W, \frac{1}{n}\mathcal{M}_{W})$ contains the point
$Q$ by Theorem~\ref{theorem:Noether-Fano} and
Lemma~\ref{lemma:Cheltsov-Kawamata}.

Each member in the linear system $\mathcal{M}_{Y}$ is contracted
to a curve by the elliptic fibration $\eta$ and the set
$\mathbb{CS}(Y, \frac{1}{n}\mathcal{M}_{Y})$ is not empty by
Theorem~\ref{theorem:Noether-Fano}. Moreover, the set
$\mathbb{CS}(Y, \frac{1}{n}\mathcal{M}_{Y})$ does not contain  any
subvariety of the exceptional divisor of $\beta$ by
Lemmas~\ref{lemma:curves} and \ref{lemma:Cheltsov-Kawamata}. Thus,
the set $\mathbb{CS}(W, \frac{1}{n}\mathcal{M}_{W})$ must contain
an element other than the point $Q$.

Let $E$ be the exceptional divisor of $\alpha$. Then, $E$ is a
quadric cone  and it follows from
Lemma~\ref{lemma:Cheltsov-Kawamata} that the set $\mathbb{CS}(W,
\frac{1}{n}\mathcal{M}_{W})$ contains a ruling $Z$ of the cone
$E$. Then, the proper transform $Z_Y$ is contracted by $\eta$ to a
point.

Let $\mathcal{T}$ be the pencil consisting of surfaces in
$|-K_{W}|$ that contain the curve $Z$ and $\gamma:W\dasharrow V$
be the dominant rational map induced by the linear system
$|-rK_{X}|$ for $r\gg 0$. The pencil $\mathcal{T}$ is the proper
transform of a pencil contained in $|-K_{X}|$, the map $\gamma$ is
a birational morphism, and $V$ is a hypersurface  of degree $12$
in $\mathbb{P}(1,1,1,4,6)$.

Let $D$ be a general surface in $|-rK_{X}|$, and $M_{W}$ and $T$
be general surfaces in $\mathcal{M}_{W}$ and $\mathcal{T}$,
respectively. Then,
$$
\frac{rn}{2}=D\cdot M_{W}\cdot T\geq \mathrm{mult}_{Z}(M_{W}\cdot
T)(D\cdot L)\geq \mathrm{mult}_{Z}(M_{W})
\mathrm{mult}_{Z}(T)(D\cdot Z)\geq\frac{rn}{2},
$$
which implies that the support of the effective cycle $M_{W}\cdot
T$ is contained in the union of the curve $Z$ and a finite union
of curves contracted by the morphism $\gamma$. Now it follows from
Theorem~\ref{theorem:main-tool} that
$\mathcal{M}_{W}=\mathcal{T}$, which completes the proof.
\end{proof}

\section{Case $\gimel=14$,  hypersurface of degree $12$ in $\mathbb{P}(1,1,1,4,6)$.}%
\label{section:n-14}\index{$\gimel=14$}

Let $X$ be a general hypersurface of degree $12$ in
$\mathbb{P}(1,1,1,4,6)$ with $-K_X^3=\frac{1}{2}$. It has only one
singular point $P$ that is a quotient singularity of type
$\frac{1}{2}(1,1,1)$.

We have an elliptic fibration as follows:
$$ 
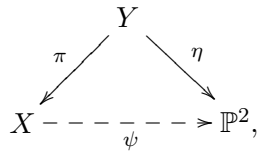
\xymatrix{
&&Y\ar@{->}[dl]_{\pi}\ar@{->}[dr]^{\eta}&&\\
&X\ar@{-->}[rr]_{\psi}&&\mathbb{P}^2,&}
$$
where \begin{itemize}

\item $\psi$ is the natural  projection,

\item

$\pi$ is the Kawamata blow up at the point $P$ with weights
$(1,1,1)$,

\item $\eta$ is an elliptic fibration. \end{itemize}

\begin{proposition}
\label{proposition:n-14} Every Halphen pencil on $X$ is contained
in $|-K_{X}|$.
\end{proposition}

\begin{proof}
The log pair $(X, \frac{1}{n}\mathcal{M})$ is not terminal by
Theorem~\ref{theorem:Noether-Fano}. However, it is terminal at a
smooth point by Lemma~\ref{lemma:smooth-points}.

Suppose that the log pair $(X, \frac{1}{n}\mathcal{M})$ is not
terminal along a curve $Z\subset X$. Then, the in\-equa\-li\-ty
$$\mathrm{mult}_{Z}(\mathcal{M})\geq n$$ holds.

For general surfaces $M_1$ and $M_2$ in $\mathcal{M}$ and a
general surface $D$ in $|-K_{X}|$, we have
$$
\frac{n^2}{2}=M_1\cdot M_2\cdot D\geq \mathrm{mult}_{Z}^2(\mathcal{M})(-K_{X}\cdot Z)\geq \frac{n^2}{2},%
$$
which implies that the curve $Z$ is a fiber of the rational map
$\psi$. For a  general surface $D'$ in $|-K_{X}|$ that contains
the curve $Z$,
$$
\frac{n}{2}=M_1\cdot D\cdot D'\geq \frac{n}{2},%
$$
which implies that
$\mathrm{Supp}(M_1)\cap\mathrm{Supp}(D')\subset\mathrm{Supp}(Z)$.
It follows from Theorem~\ref{theorem:main-tool} that the linear
system $\mathcal{M}$ is the pencil in $|-K_{X}|$ consisting of
surfaces that pass through $Z$.

Suppose that the log pair $(X, \frac{1}{n}\mathcal{M})$ is not
terminal at the point $P$. Because $-K_Y^3=0$ and
$\mathcal{M}_Y\sim_{\mathbb{Q}}-nK_Y$, every surface in the pencil
$\mathcal{M}_Y$ is contracted to a curve by the morphism $\eta$.
The log pair $(Y, \frac{1}{n}\mathcal{M}_Y)$ is not terminal along
a curve contained in the exceptional divisor of $\pi$ by
Theorem~\ref{theorem:Noether-Fano} and
Lemma~\ref{lemma:Cheltsov-Kawamata}. However, it contradicts
Lemma~\ref{lemma:curves}.
\end{proof}

\newpage

\part{The Table.}
\label{section:weighted-Fanos}\index{The Table}

Before we explain the table, we should mention that all the
contents, except the numbers of Halphen pencils, are obtained from
\cite{CPR}.

 We tabulate the singular points of
the hypersurface
\[X=X_d\subset\mathbb{P}(1, a_1, a_2, a_3, a_4),\]
and the
number of Halphen pencils on $X$, \emph{i.e.}, the number of ways
in which the hypersurface $X$ is birationally transformed into a
 fibration by surfaces of Kodaira dimension zero.

The contents in the entries on the first row and the second column
is the number of Halphen pencils. These pencils define rational
maps  a generic fiber of which is birational to a smooth K3 surface.

The contents in the entries from the second rows explain the
singular points on $X$.  The first column tabulates the types of
singularities. The second column shows the numbers $b$ and $c$ in
Proposition~\ref{proposition:special-singular-points} when we take
the Kawamata blow up at a given point. The divisors $B$ and $E$
are the anticanonical divisor  and the exceptional divisor ,
respectively, on the Kawamata blow up at a given singular point.
For simplicity, we keep the divisors $bB+cE$ only when $B^3<0$. The blank entries simply mean $B^3\geq 0$.
In such cases, we do not need the divisors $bB+cE$ for the present article.

\begin{center}

\begin{entry}
\no 1; \reln 4; \gens 1,1,1,1,1; \degree 4; \inv \textbf{F}_0;
\ell \infty; \kkk  \infty;\type \text{smooth}; \excl \text{N/A};
\untw none;
\end{entry}

\begin{entry}
\no 2; \reln 5; \gens 1,1,1,1,2; \degree 5/2; \inv
\textbf{F}_1;\ell 15; \kkk  \infty;\type P_4=\recip2(1,1,1); \untw
Q.I. $xw^2,3,5$;
\end{entry}

\begin{entry}
\no 3; \reln 6; \gens 1,1,1,1,3; \degree 2;\inv \textbf{F}_0;\ell
0 ; \kkk \infty; \type \text{smooth}; \excl \text{N/A}; \untw
none;
\end{entry}

\begin{entry}
\no 4; \reln 6; \gens 1,1,1,2,2; \degree 3/2; \inv
\hat{\textbf{F}}_3;\ell 1 ; \kkk \infty; \type
P_3P_4=3\times\recip2(1,1,1); \untw Q.I. $tw^2,4,6$;
\end{entry}

\begin{entry}
\no 5; \reln 7; \gens 1,1,1,2,3; \degree 7/6; \inv
\textbf{F}_2;\ell 1 ; \kkk \infty; \type P_4=\recip3(1,1,2); \untw
Q.I. $xw^2,4,7$; \typea P_3=\recip2(1,1,1); \untwa Q.I.
$*wt^2,5,7$;
\end{entry}

\begin{entry}
\no 6; \reln 8; \gens 1,1,1,2,4; \degree 1; \inv \textbf{F}_2;\ell
1 ;\kkk \infty;\type P_3P_4=2\times\recip2(1,1,1); \untw Q.I.
$wt^2,6,8$;
\end{entry}

\begin{entry}
\no 7; \reln 8; \gens 1,1,2,2,3; \degree 2/3;\inv
\textbf{F}_5;\ell 5 ; \kkk 1; \type P_4=\recip3(1,1,2); \untw Q.I.
$tw^2,6,8$; \typea P_2P_3=4\times\recip2(1,1,1); \untwa E.I.
$*tw^2-zt^3,5,8$;
\end{entry}

\begin{entry}
\no 8; \reln 9; \gens 1,1,1,3,4; \degree 3/4;\inv
\textbf{F}_1;\ell 1 ; \kkk \infty; \type P_4=\recip4(1,1,3); \untw
Q.I. $xw^2,5,9$;
\end{entry}

\begin{entry}
\no 9; \reln 9; \gens 1,1,2,3,3; \degree 1/2;\inv
\hat{\textbf{F}}_3; \ell 2 ;\kkk 1; \type P_2=\recip2(1,1,1);
\excl ; \typea P_3P_4=3\times\recip3(1,1,2); \untwa Q.I.
$tw^2,6,9$;
\end{entry}

\begin{entry}
\no 10; \reln 10; \gens 1,1,1,3,5; \degree 2/3;\inv
\textbf{F}_0;\ell 1 ; \kkk \infty; \type P_3=\recip3(1,1,2); \excl
; \untw none;
\end{entry}

\begin{entry}
\no 11; \reln 10; \gens 1,1,2,2,5; \degree 1/2;\inv
\textbf{F}_0;\ell 5 ;\kkk 1; \type P_2P_3=5\times\recip2(1,1,1);
\excl ; \untw none;
\end{entry}

\begin{entry}
\no 12; \reln 10; \gens 1,1,2,3,4; \degree 5/12;\inv
\textbf{F}_2;\ell 1 ; \kkk 1; \type P_4=\recip4(1,1,3); \untw Q.I.
$zw^2,6,10$; \typea P_3=\recip3(1,1,2); \untwa Q.I. $*t^2w,7,10$;
\typeb P_2P_4=2\times\recip2(1,1,1); \exclb   B;
\end{entry}

\begin{entry}
\no 13; \reln 11; \gens 1,1,2,3,5; \degree 11/30;\inv
\textbf{F}_2; \ell 1; \kkk 1 ;\type P_4=\recip5(1,2,3); \untw Q.I.
$xw^2,6,11$; \typea P_3=\recip3(1,1,2); \untwa Q.I. $*t^2w,8,11$;
\typeb P_2=\recip2(1,1,1); \exclb   B;
\end{entry}

\begin{entry}
\no 14; \reln 12; \gens 1,1,1,4,6; \degree 1/2;\inv
\textbf{F}_0;\ell 1 ; \kkk \infty; \type
P_3P_4=1\times\recip2(1,1,1); \excl ; \untw none;
\end{entry}

\begin{entry}
\no 15; \reln 12; \gens 1,1,2,3,6; \degree 1/3;\inv
\textbf{F}_2;\ell 1 ; \kkk 1; \type P_3P_4=2\times\recip3(1,1,2);
\untw Q.I. $wt^2,9,12$; \typea P_2P_4=2\times\recip2(1,1,1);
\excla   B;
\end{entry}

\begin{entry}
\no 16; \reln 12; \gens 1,1,2,4,5; \degree 3/10;\inv
\textbf{F}_1;\ell 1 ; \kkk 1; \type P_4=\recip5(1,1,4); \untw Q.I.
$zw^2,8,12$; \typea P_2P_3=3\times\recip2(1,1,1); \excla   B;
\end{entry}

\begin{entry}
\no 17; \reln 12; \gens 1,1,3,4,4; \degree 1/4;\inv
\hat{\textbf{F}}_3;\ell 4 ; \kkk 1; \type
P_3P_4=3\times\recip4(1,1,3); \untw Q.I. $tw^2,8,12$;
\end{entry}

\begin{entry}
\no 18; \reln 12; \gens 1,2,2,3,5; \degree 1/5;\inv
\textbf{F}_1;\ell 1 ;\kkk 7; \type P_4=\recip5(1,2,3); \untw Q.I.
$yw^2,7,12$; \typea P_1P_2=6\times\recip2(1,1,1); \excla
2B;
\end{entry}

\begin{entry}
\no 19; \reln 12; \gens 1,2,3,3,4; \degree 1/6;\inv
\textbf{F}_0;\ell 4 ;  \kkk 1;\type P_1P_4=3\times\recip2(1,1,1);
\excl   6B+E ; \untw none; \typea
P_2P_3=4\times\recip3(1,2,1); \excla ;
\end{entry}

\begin{entry}
\no 20; \reln 13; \gens 1,1,3,4,5; \degree 13/60;\inv
\textbf{F}_3; \ell 2 ;\kkk 1; \type P_4=\recip5(1,1,4); \untw Q.I.
$zw^2,8,13$; \typea P_3=\recip4(1,1,3); \untwa Q.I. $*t^2w,9,13$;
\typeb P_2=\recip3(1,1,2); \untwb E.I. $zw^2-*tz^3,8,13$;
\end{entry}

\begin{entry}
\no 21; \reln 14; \gens 1,1,2,4,7; \degree 1/4;\inv
\textbf{F}_0;\ell 1 ;\kkk 1; \type P_3=\recip4(1,1,3); \excl
; \untw none; \typea P_2P_3=3\times\recip2(1,1,1);
\excla   B;
\end{entry}

\begin{entry}
\no 22; \reln 14; \gens 1,2,2,3,7; \degree 1/6;\inv
\textbf{F}_0;\ell 1 ; \kkk 8; \type P_3=\recip3(1,2,1); \excl
; \untw none; \typea P_1P_2=7\times\recip2(1,1,1); \excla
  2B;
\end{entry}

\begin{entry}
\no 23; \reln 14; \gens 1,2,3,4,5; \degree 7/60;\inv
\textbf{F}_2;\ell 1; \kkk 1; \type P_4=\recip5(1,2,3); \untw Q.I.
$tw^2,9,14$; \typea P_3=\recip4(1,3,1); \untwa E.I.
$tw^2-yt^3,9,14$; \typeb P_2=\recip3(1,2,1); \exclb   2B;
\typec P_1P_3=3\times\recip2(1,1,1); \exclc   4B+E;
\end{entry}

\begin{entry}
\no 24; \reln 15; \gens 1,1,2,5,7; \degree 3/14;\inv
\textbf{F}_1;\ell 1 ; \kkk 1; \type P_4=\recip7(1,2,5); \untw Q.I.
$xw^2,8,15$; \typea P_2=\recip2(1,1,1); \excla   B;
\end{entry}

\begin{entry}
\no 25; \reln 15; \gens 1,1,3,4,7; \degree 5/28; \inv
\textbf{F}_2;\ell 1;\kkk 1; \type P_4=\recip7(1,3,4); \untw Q.I.
$xw^2,8,15$; \typea P_3=\recip4(1,1,3); \untwa Q.I. $*wt^2,11,15$;
\end{entry}

\begin{entry}
\no 26; \reln 15; \gens 1,1,3,5,6; \degree 1/6;\inv
\textbf{F}_1;\ell 2 ;\kkk 1; \type P_4=\recip6(1,1,5); \untw Q.I.
$zw^2,9,15$; \typea P_3P_4=2\times\recip3(1,1,2); \excla ;
\end{entry}

\begin{entry}
\no 27; \reln 15; \gens 1,2,3,5,5; \degree 1/10;\inv
\hat{\textbf{F}}_3;\ell 1 ; \kkk 1; \type
P_3P_4=3\times\recip5(1,2,3); \untw Q.I. $tw^2,10,15$; \typea
P_1=\recip2(1,1,1); \excla   5B+E;
\end{entry}

\begin{entry}
\no 28; \reln 15; \gens 1,3,3,4,5; \degree 1/12;\inv
\textbf{F}_0;\ell 1 ; \kkk 6; \type P_3=\recip4(1,3,1); \excl
; \untw none; \typea P_1P_2=5\times\recip3(1,1,2); \excla
  3B;
\end{entry}

\begin{entry}
\no 29; \reln 16; \gens 1,1,2,5,8; \degree 1/5;\inv
\textbf{F}_0;\ell 1 ;\kkk 1; \type P_3=\recip5(1,2,3); \excl
; \untw none; \typea P_2P_4=2\times\recip2(1,1,1);
\excla   B;
\end{entry}

\begin{entry}
\no 30; \reln 16; \gens 1,1,3,4,8; \degree 1/6;\inv
\textbf{F}_2;\ell 2 ; \kkk 1; \type P_3P_4=2\times\recip4(1,1,3);
\untw Q.I. $wt^2,12,16$; \typea P_2=\recip3(1,1,2); \excla ;
\end{entry}

\begin{entry}
\no 31; \reln 16; \gens 1,1,4,5,6; \degree 2/15;\inv
\textbf{F}_2;\ell 2 ; \kkk 1; \type P_4=\recip6(1,1,5); \untw Q.I.
$zw^2,10,16$; \typea P_3=\recip5(1,1,4); \untwa Q.I.
$*wt^2,11,16$; \typeb P_2P_4=1\times\recip2(1,1,1); \exclb
B;
\end{entry}

\begin{entry}
\no 32; \reln 16; \gens 1,2,3,4,7; \degree 2/21;\inv
\textbf{F}_1;\ell 1 ;\kkk 1; \type P_4=\recip7(1,3,4); \untw Q.I.
$yw^2,9,16$; \typea P_2=\recip3(1,2,1); \excla   2B; \typeb
P_1P_3=4\times\recip2(1,1,1); \exclb   4B+E;
\end{entry}

\begin{entry}
\no 33; \reln 17; \gens 1,2,3,5,7; \degree 17/210;\inv
\textbf{F}_2; \ell 1 ;\kkk 1; \type P_4=\recip7(1,2,5); \untw Q.I.
$zw^2,10,17$; \typea P_3=\recip5(1,2,3); \untwa Q.I.
$*wt^2,12,17$; \typeb P_2=\recip3(1,2,1); \exclb   2B; \typec
P_1=\recip2(1,1,1); \exclc   7B+2E;
\end{entry}

\begin{entry}
\no 34; \reln 18; \gens 1,1,2,6,9; \degree 1/6;\inv
\textbf{F}_0;\ell 1 ; \kkk 1; \type P_3P_4=1\times\recip3(1,1,2);
 \untw none; \typea P_2P_3=3\times\recip2(1,1,1);
\excla   B;
\end{entry}

\begin{entry}
\no 35; \reln 18; \gens 1,1,3,5,9; \degree 2/15; \inv
\textbf{F}_0;\ell 1 ; \kkk 1; \type P_3=\recip5(1,1,4); \excl
; \untw none; \typea P_2P_4=2\times\recip3(1,1,2);
\excla   B;
\end{entry}

\begin{entry}
\no 36; \reln 18; \gens 1,1,4,6,7; \degree 3/28; \inv
\textbf{F}_2;\ell 2 ;\kkk 1; \type P_4=\recip7(1,1,6); \untw Q.I.
$zw^2,11,18$; \typea P_2=\recip4(1,1,3); \untwa E.I.
$zw^2-tz^3,11,18$; \typeb P_2P_3=1\times\recip2(1,1,1); \exclb
  B;
\end{entry}

\begin{entry}
\no 37; \reln 18; \gens 1,2,3,4,9; \degree 1/12;\inv
\textbf{F}_0;\ell 1;\kkk 1 ; \type P_3=\recip4(1,3,1); \excl
; \untw none; \typea P_2P_4=2\times\recip3(1,2,1); \excla
  2B; \typeb P_1P_3=4\times\recip2(1,1,1); \exclb
4B+E;
\end{entry}

\begin{entry}
\no 38; \reln 18; \gens 1,2,3,5,8; \degree 3/40; \inv
\textbf{F}_2;\ell 1 ;\kkk 1; \type P_4=\recip8(1,3,5); \untw Q.I.
$yw^2,10,18$; \typea P_3=\recip5(1,2,3); \untwa Q.I.
$*wt^2,13,18$; \typeb P_1P_4=2\times\recip2(1,1,1); \exclb
10B+3E ;
\end{entry}

\begin{entry}
\no 39; \reln 18; \gens 1,3,4,5,6; \degree 1/20;\inv
\textbf{F}_0;\ell 1;\kkk 1 ; \type P_3=\recip5(1,4,1); \excl
; \untw none; \typea P_2=\recip4(1,3,1); \excla
3B; \typeb P_2P_4=\recip2(1,1,1); \exclb   3B+E; \typec
P_1P_4=3\times\recip3(1,1,2); \exclc   20B+3E ;
\end{entry}

\begin{entry}
\no 40; \reln 19; \gens 1,3,4,5,7; \degree 19/420;\inv
\textbf{F}_2;\ell 1 ;\kkk 1; \type P_4=\recip7(1,3,4); \untw Q.I.
$tw^2,12,19$; \typea P_3=\recip5(1,3,2); \untwa E.I.
$tw^2-zt^3,12,19$; \typeb P_2=\recip4(1,3,1); \exclb   3B;
\typec P_1=\recip3(1,1,2); \exclc   7B+E;
\end{entry}

\begin{entry}
\no 41; \reln 20; \gens 1,1,4,5,10; \degree 1/10;\inv
\textbf{F}_2; \ell 1 ;\kkk 1; \type P_3P_4=2\times\recip5(1,1,4);
\untw Q.I. $wt^2,15,20$; \typea P_2P_4=\recip2(1,1,1); \excla
  B;
\end{entry}

\begin{entry}
\no 42; \reln 20; \gens 1,2,3,5,10; \degree 1/15;\inv
\textbf{F}_2; \ell 1 ;\kkk 1; \type P_2=\recip3(1,2,1); \excl
  2B; \typea P_3P_4=2\times\recip5(1,2,3); \untwa Q.I.
$tz^2,15,20$; \typeb P_1P_4=2\times\recip2(1,1,1); \exclb
10B+3E;
\end{entry}

\begin{entry}
\no 43; \reln 20; \gens 1,2,4,5,9; \degree 1/18;\inv
\textbf{F}_1;\ell 1;\kkk 1; \type P_4=\recip9(1,4,5); \untw Q.I.
$yw^2,11,20$; \typea P_1P_2=5\times\recip2(1,1,1); \excla
4B+E;
\end{entry}

\begin{entry}
\no 44; \reln 20; \gens 1,2,5,6,7; \degree 1/21;\inv
\textbf{F}_2;\ell 2 ;\kkk 1; \type P_4=\recip7(1,2,5); \untw Q.I.
$tw^2,13,20$; \typea P_3=\recip6(1,5,1); \untwa E.I.
$tw^2-yt^3,13,20$; \typeb P_1P_3=3\times\recip2(1,1,1); \exclb
  6B+2E;
\end{entry}

\begin{entry}
\no 45; \reln 20; \gens 1,3,4,5,8; \degree 1/24; \inv
\textbf{F}_1;\ell 1 ;\kkk 2; \type P_4=\recip8(1,3,5); \untw Q.I.
$zw^2,12,20$; \typea P_1=\recip3(1,1,2); \excla   8B+E;
\typeb P_2P_4=2\times\recip4(1,3,1); \exclb   3B;
\end{entry}

\begin{entry}
\no 46; \reln 21; \gens 1,1,3,7,10; \degree 1/10;\inv
\textbf{F}_1; \ell 1 ;\kkk 1; \type P_4=\recip{10}(1,3,7); \untw
Q.I. $xw^2,11,21$;
\end{entry}

\begin{entry}
\no 47; \reln 21; \gens 1,1,5,7,8; \degree 3/40;\inv
\textbf{F}_1;\ell 1 ;\kkk 1; \type P_4=\recip8(1,1,7); \untw Q.I.
$zw^2,13,21$; \typea P_2=\recip5(1,2,3); \excla ;
\end{entry}

\begin{entry}
\no 48; \reln 21; \gens 1,2,3,7,9; \degree 1/18;\inv
\textbf{F}_1;\ell 1 ;\kkk 2; \type P_4=\recip9(1,2,7); \untw Q.I.
$zw^2,12,21$; \typea P_1=\recip2(1,1,1); \excla   3B+E;
\typeb P_2P_4=2\times\recip3(1,2,1); \exclb   2B;
\end{entry}

\begin{entry}
\no 49; \reln 21; \gens 1,3,5,6,7; \degree 1/30;\inv
\textbf{F}_0;\ell 2 ; \kkk 1; \type P_3=\recip6(1,5,1); \excl
; \untw none; \typea P_2=\recip5(1,3,2); \excla ; \typeb P_1P_3=3\times\recip3(1,2,1); \exclb   7B+E;
\end{entry}

\begin{entry}
\no 50; \reln 22; \gens 1,1,3,7,11; \degree 2/21;\inv
\textbf{F}_0; \ell 1 ;\kkk 1; \type P_3=\recip7(1,3,4); \excl
; \untw none; \typea P_2=\recip3(1,1,2); \excla
B;
\end{entry}

\begin{entry}
\no 51; \reln 22; \gens 1,1,4,6,11; \degree 1/12;\inv
\textbf{F}_0;\ell 2 ; \kkk 1; \type P_3=\recip6(1,1,5); \excl
; \untw none; \typea P_2=\recip4(1,1,3); \excla ; \typeb P_2P_3=\recip2(1,1,1); \exclb   B;
\end{entry}

\begin{entry}
\no 52; \reln 22; \gens 1,2,4,5,11; \degree 1/20;\inv
\textbf{F}_0; \ell 1;\kkk 1 ;\type P_3=\recip5(1,4,1); \excl
; \untw none; \typea P_2=\recip4(1,1,3); \excla
2B; \typeb P_1P_2=5\times\recip2(1,1,1); \exclb   4B+E;
\end{entry}

\begin{entry}
\no 53; \reln 24; \gens 1,1,3,8,12; \degree 1/12;\inv
\textbf{F}_0; \ell 1 ;\kkk 1; \type P_3P_4=1\times\recip4(1,1,3);
 \untw none; \typea P_2P_4=2\times\recip3(1,1,2);
\excla   B;
\end{entry}

\begin{entry}
\no 54; \reln 24; \gens 1,1,6,8,9; \degree 1/18;\inv
\textbf{F}_1;\ell 1 ; \kkk 1; \type P_4=\recip9(1,1,8); \untw Q.I.
$zw^2,15,24$; \typea P_2P_4=1\times\recip3(1,1,2); \excla
B; \typeb P_2P_3=\recip2(1,1,1); \exclb   B;
\end{entry}

\begin{entry}
\no 55; \reln 24; \gens 1,2,3,7,12; \degree 1/21;\inv
\textbf{F}_0; \ell 1 ;\kkk 2; \type P_3=\recip7(1,2,5); \excl
; \untw none; \typea P_2P_4=2\times\recip3(1,2,1);
\excla   2B; \typeb P_1P_4=2\times\recip2(1,1,1); \exclb
  3B+E;
\end{entry}

\begin{entry}
\no 56; \reln 24; \gens 1,2,3,8,11; \degree 1/22;\inv
\textbf{F}_1; \ell 1 ;\kkk 1; \type P_4=\recip{11}(1,3,8); \untw
Q.I. $yw^2,13,24$; \typea P_1P_3=3\times\recip2(1,1,1); \excla
  3B+E;
\end{entry}

\begin{entry}
\no 57; \reln 24; \gens 1,3,4,5,12; \degree 1/30;\inv
\textbf{F}_0; \ell 1;\kkk 2 ;\type P_3=\recip5(1,3,2); \excl
; \untw none; \typea P_2P_4=2\times\recip4(1,3,1); \excla
  3B; \typeb P_1P_4=2\times\recip3(1,1,2); \exclb
12B+2E;
\end{entry}

\begin{entry}
\no 58; \reln 24; \gens 1,3,4,7,10; \degree 1/35;\inv
\textbf{F}_2; \ell 1 ;\kkk 2; \type P_4=\recip{10}(1,3,7); \untw
Q.I. $zw^2,14,24$; \typea P_3=\recip7(1,3,4); \untwa Q.I.
$*wt^2,17,24$; \typeb P_2P_4=\recip2(1,1,1); \exclb   3B+E;
\end{entry}

\begin{entry}
\no 59; \reln 24; \gens 1,3,6,7,8; \degree 1/42;\inv
\textbf{F}_0;\ell 1 ;\kkk 1; \type P_3=\recip7(1,6,1); \excl
; \untw none; \typea P_2P_4=\recip2(1,1,1); \excla
3B+E; \typeb P_1P_2=4\times\recip3(1,1,2); \exclb   6B+E;
\end{entry}

\begin{entry}
\no 60; \reln 24; \gens 1,4,5,6,9; \degree 1/45;\inv
\textbf{F}_1;\ell 0; \kkk 2; \type P_4=\recip9(1,4,5); \untw Q.I.
$tw^2,15,24$; \typea P_2=\recip5(1,4,1); \excla   4B; \typeb
P_3P_4=1\times\recip3(1,1,2); \exclb   5B+E; \typec
P_1P_3=2\times\recip2(1,1,1); \exclc   18B+7E ;
\end{entry}

\begin{entry}
\no 61; \reln 25; \gens 1,4,5,7,9; \degree 5/252;\inv
\textbf{F}_2; \ell 1 ;\kkk 1; \type P_4=\recip9(1,4,5); \untw Q.I.
$tw^2,16,25$; \typea P_3=\recip7(1,5,2); \untwa E.I.
$tw^2-yt^3,16,25$; \typeb P_1=\recip4(1,3,1); \exclb   9B+E;
\end{entry}

\begin{entry}
\no 62; \reln 26; \gens 1,1,5,7,13; \degree 2/35;\inv
\textbf{F}_0; \ell 1 ;\kkk 1; \type P_3=\recip7(1,1,6); \excl
; \untw none; \typea P_2=\recip5(1,2,3); \excla ;
\end{entry}

\begin{entry}
\no 63; \reln 26; \gens 1,2,3,8,13; \degree 1/24;\inv
\textbf{F}_0;\ell 1 ;\kkk 1; \type P_3=\recip8(1,3,5); \excl
; \untw none; \typea P_2=\recip3(1,2,1); \excla
2B; \typeb P_1P_3=3\times\recip2(1,1,1); \exclb   3B+E;
\end{entry}

\begin{entry}
\no 64; \reln 26; \gens 1,2,5,6,13; \degree 1/30;\inv
\textbf{F}_0; \ell 2 ; \kkk 1;\type P_3=\recip6(1,5,1); \excl
; \untw none; \typea P_2=\recip5(1,2,3); \excla ; \typeb P_1P_3=4\times\recip2(1,1,1); \exclb   6B+2E;
\end{entry}

\begin{entry}
\no 65; \reln 27; \gens 1,2,5,9,11; \degree 3/110;\inv
\textbf{F}_1;\ell 1 ; \kkk 1; \type P_4=\recip{11}(1,2,9); \untw
Q.I. $zw^2,16,27$; \typea P_2=\recip5(1,4,1); \excla   2B;
\typeb P_1=\recip2(1,1,1); \exclb   11B+4E;
\end{entry}

\begin{entry}
\no 66; \reln 27; \gens 1,5,6,7,9; \degree 1/70;\inv
\textbf{F}_0;\ell 1;\kkk 2 ; \type P_3=\recip7(1,5,2); \excl
; \untw none; \typea P_2=\recip6(1,5,1); \excla
5B; \typeb P_1=\recip5(1,1,4); \exclb   12B+E ; \typec
P_2P_4=1\times\recip3(1,2,1); \exclc   5B+E;
\end{entry}

\begin{entry}
\no 67; \reln 28; \gens 1,1,4,9,14; \degree 1/18;\inv
\textbf{F}_0; \ell 1 ;\kkk 1; \type P_3=\recip9(1,4,5); \excl
; \untw none; \typea P_2P_4=\recip2(1,1,1); \excla
  B;
\end{entry}

\begin{entry}
\no 68; \reln 28; \gens 1,3,4,7,14; \degree 1/42;\inv
\textbf{F}_2; \ell 1; \kkk 1; \type P_1=\recip3(1,1,2); \excl
  7B+E; \typea P_3P_4=2\times\recip7(1,3,4); \untwa Q.I.
$wt^2,21,28$; \typeb P_2P_4=\recip2(1,1,1); \exclb   3B+E;
\end{entry}

\begin{entry}
\no 69; \reln 28; \gens 1,4,6,7,11; \degree 1/66;\inv
\textbf{F}_1; \ell 1;\kkk 2 ;\type P_4=\recip{11}(1,4,7); \untw
Q.I. $zw^2,17,28$; \typea P_2=\recip6(1,1,5); \excla   4B;
\typeb P_1P_2=2\times\recip2(1,1,1); \exclb   22B+9E ;
\end{entry}

\begin{entry}
\no 70; \reln 30; \gens 1,1,4,10,15; \degree 1/20;\inv
\textbf{F}_0;\ell 1 ; \kkk 1; \type P_2=\recip4(1,1,3); \excl
  B; \untw none; \typea P_3P_4=\recip5(1,1,4); \excla ; \typeb P_2P_3=1\times\recip2(1,1,1); \exclb   B;
\end{entry}

\begin{entry}
\no 71; \reln 30; \gens 1,1,6,8,15; \degree 1/24;\inv
\textbf{F}_0; \ell 1 ;\kkk 1;\type P_3=\recip8(1,1,7); \excl
 ; \untw none; \typea P_2P_4=1\times\recip3(1,1,2);
\excla   B; \typeb P_2P_3=1\times\recip2(1,1,1); \exclb
  B;
\end{entry}

\begin{entry}
\no 72; \reln 30; \gens 1,2,3,10,15; \degree 1/30;\inv
\textbf{F}_0;\ell 1;\kkk 1 ; \type P_3P_4=1\times\recip5(1,2,3);
 \untw none; \typea P_2P_4=2\times\recip3(1,2,1);
\excla   2B; \typeb P_1P_3=3\times\recip2(1,1,1); \exclb
  3B+E;
\end{entry}

\begin{entry}
\no 73; \reln 30; \gens 1,2,6,7,15; \degree 1/42;\inv
\textbf{F}_0; \ell 1;\kkk 1 ;\type P_3=\recip7(1,6,1); \excl
; \untw none; \typea P_2P_4=1\times\recip3(1,2,1); \excla
  2B; \typeb P_1P_2=5\times\recip2(1,1,1); \exclb
6B+2E;
\end{entry}

\begin{entry}
\no 74; \reln 30; \gens 1,3,4,10,13; \degree 1/52;\inv
\textbf{F}_1;\ell 1 ; \kkk 2; \type P_4=\recip{13}(1,3,10); \untw
Q.I. $zw^2,17,30$; \typea P_2=\recip4(1,3,1); \excla   3B;
\typeb P_2P_3=\recip2(1,1,1); \exclb   3B+E;
\end{entry}

\begin{entry}
\no 75; \reln 30; \gens 1,4,5,6,15; \degree 1/60;\inv
\textbf{F}_0; \ell 0;\kkk 1;\type P_1=\recip4(1,1,3); \excl
10B+E ; \untw none; \typea P_3P_4=1\times\recip3(1,1,2); \excla
  5B+E; \typeb P_2P_4=2\times\recip5(1,4,1); \exclb
4B; \typec P_1P_3=2\times\recip2(1,1,1); \exclc   5B+2E;
\end{entry}

\begin{entry}
\no 76; \reln 30; \gens 1,5,6,8,11; \degree 1/88;\inv
\textbf{F}_2; \ell 1;\kkk 2 ;\type P_4=\recip{11}(1,5,6); \untw
Q.I. $tw^2,19,30$; \typea P_3=\recip8(1,5,3); \untwa E.I.
$tw^2-zt^3,22,30$; \typeb P_2P_3=1\times\recip2(1,1,1); \exclb
  5B+2E;
\end{entry}

\begin{entry}
\no 77; \reln 32; \gens 1,2,5,9,16; \degree 1/45;\inv
\textbf{F}_0; \ell 1 ;\kkk 1; \type P_3=\recip9(1,2,7); \excl
; \untw none; \typea P_2=\recip5(1,4,1); \excla
2B; \typeb P_1P_4=2\times\recip2(1,1,1); \exclb   18B+7E;
\end{entry}

\begin{entry}
\no 78; \reln 32; \gens 1,4,5,7,16; \degree 1/70;\inv
\textbf{F}_0; \ell 1;\kkk 1 ;\type P_3=\recip7(1,5,2); \excl
; \untw none; \typea P_2=\recip5(1,4,1); \excla
4B; \typeb P_1P_4=2\times\recip4(1,1,3); \exclb   7B+E;
\end{entry}

\begin{entry}
\no 79; \reln 33; \gens 1,3,5,11,14; \degree 1/70;\inv
\textbf{F}_1;\ell 1 ; \kkk 2; \type P_4=\recip{14}(1,3,11); \untw
Q.I. $zw^2,19,33$; \typea P_2=\recip5(1,1,4); \excla   3B;
\end{entry}

\begin{entry}
\no 80; \reln 34; \gens 1,3,4,10,17; \degree 1/60;\inv
\textbf{F}_0; \ell 1 ;\kkk 2; \type P_3=\recip{10}(1,3,7); \excl
; \untw none; \typea P_2=\recip4(1,3,1); \excla
3B; \typeb P_1=\recip3(1,1,2); \exclb   10B+2E; \typec
P_2P_3=1\times\recip2(1,1,1); \exclc   3B+E;
\end{entry}

\begin{entry}
\no 81; \reln 34; \gens 1,4,6,7,17; \degree 1/84;\inv
\textbf{F}_0; \ell 1; \kkk 2 ;\type P_3=\recip7(1,4,3); \excl
; \untw none; \typea P_2=\recip6(1,1,5); \excla
4B; \typeb P_1=\recip4(1,3,1); \exclb   7B+E; \typec
P_1P_2=2\times\recip2(1,1,1); \exclc   12B+5E ;
\end{entry}

\begin{entry}
\no 82; \reln 36; \gens 1,1,5,12,18; \degree 1/30;\inv
\textbf{F}_0; \ell 1; \kkk 1; \type P_2=\recip5(1,2,3); \excl
; \untw none; \typea P_3P_4=1\times\recip6(1,1,5); \end{entry}

\begin{entry}
\no 83; \reln 36; \gens 1,3,4,11,18; \degree 1/66;\inv
\textbf{F}_0;\ell 1 ;  \kkk 1; \type P_3=\recip{11}(1,4,7); \excl
; \untw none; \typea P_2P_4=1\times\recip2(1,1,1);
\excla   3B+E; \typeb P_1P_4=2\times\recip3(1,1,2); \exclb
  6B+18E;
\end{entry}

\begin{entry}
\no 84; \reln 36; \gens 1,7,8,9,12; \degree 1/168;\inv
\textbf{F}_0; \ell 0 ;\kkk 2;\type P_2=\recip8(1,7,1); \excl
  7B; \untw none; \typea P_1=\recip7(1,2,5); \excla
8B; \typeb P_3P_4=1\times\recip3(1,1,2); \exclb   8B+2E;
\typec P_2P_4=1\times\recip4(1,3,1); \exclc   7B+E;
\end{entry}

\begin{entry}
\no 85; \reln 38; \gens 1,3,5,11,19; \degree 2/165;\inv
\textbf{F}_0; \ell 1 ;\kkk 1; \type P_3=\recip{11}(1,3,8); \excl
; \untw none; \typea P_2=\recip5(1,1,4); \excla
3B; \typeb P_1=\recip3(1,2,1); \exclb   5B+E;
\end{entry}

\begin{entry}
\no 86; \reln 38; \gens 1,5,6,8,19; \degree 1/120;\inv
\textbf{F}_0; \ell 1;\kkk 2 ; \type P_3=\recip8(1,5,3); \excl
; \untw none; \typea P_2=\recip6(1,5,1); \excla
5B; \typeb P_1=\recip5(1,1,4); \exclb   18B+2E ; \typec
P_2P_3=1\times\recip2(1,1,1); \exclc   5B+2E;
\end{entry}

\begin{entry}
\no 87; \reln 40; \gens 1,5,7,8,20; \degree 1/140;\inv
\textbf{F}_0; \ell 0 ;\kkk 1;\type P_2=\recip7(1,1,6); \excl
  5B; \untw none; \typea P_3P_4=1\times\recip4(1,1,3); \excla
  7B+E; \typeb P_1P_4=2\times\recip5(1,2,3); \exclb
20B+3E;
\end{entry}

\begin{entry}
\no 88; \reln 42; \gens 1,1,6,14,21; \degree 1/42;\inv
\textbf{F}_0; \ell 1 ;\kkk 1; \type P_3P_4=1\times\recip7(1,1,6);
 \untw none; \typea P_2P_4=1\times\recip3(1,1,2);
\excla   B; \typeb P_2P_3=1\times\recip2(1,1,1); \exclb
  B;
\end{entry}

\begin{entry}
\no 89; \reln 42; \gens 1,2,5,14,21; \degree 1/70;\inv
\textbf{F}_0;\ell 1;\kkk 1 ; \type P_2=\recip5(1,4,1); \excl
  2B; \untw none; \typea P_3P_4=1\times\recip7(1,2,5); \excla
; \typeb P_1P_3=3\times\recip2(1,1,1); \exclb
5B+2E;
\end{entry}

\begin{entry}
\no 90; \reln 42; \gens 1,3,4,14,21; \degree 1/84;\inv
\textbf{F}_0; \ell 1; \kkk 1;\type P_2=\recip4(1,3,1); \excl
  3B; \untw none; \typea P_3P_4=\recip7(1,3,4); \excla ; \typeb P_1P_4=2\times\recip3(1,1,2); \exclb   21B+5E;
\typec P_2P_3=\recip2(1,1,1); \exclc   3B+E;
\end{entry}

\begin{entry}
\no 91; \reln 44; \gens 1,4,5,13,22; \degree 1/130;\inv
\textbf{F}_0; \ell 1 ;\kkk 2; \type P_3=\recip{13}(1,4,9); \excl
; \untw none; \typea P_2=\recip5(1,3,2); \excla
4B; \typeb P_1P_4=1\times\recip2(1,1,1); \exclb   5B+2E;
\end{entry}

\begin{entry}
\no 92; \reln 48; \gens 1,3,5,16,24; \degree 1/120;\inv
\textbf{F}_0; \ell 1; \kkk 1;\type P_2=\recip5(1,1,4); \excl
  3B; \untw none; \typea P_3P_4=1\times\recip8(1,3,5); \excla
; \typeb P_1P_4=2\times\recip3(1,2,1); \exclb
5B+E;
\end{entry}

\begin{entry}
\no 93; \reln 50; \gens 1,7,8,10,25; \degree 1/280;\inv
\textbf{F}_0; \ell 0 ;\kkk 2;\type P_2=\recip8(1,7,1); \excl
  7B; \untw none; \typea P_1=\recip7(1,3,4); \excla
8B; \typeb P_3P_4=1\times\recip5(1,2,3); \exclb   8B+E;
\typec P_2P_3=1\times\recip2(1,1,1); \exclc   7B+3E;
\end{entry}

\begin{entry}
\no 94; \reln 54; \gens 1,4,5,18,27; \degree 1/180;\inv
\textbf{F}_0; \ell 1; \kkk 1; \type P_2=\recip5(1,3,2); \excl
  4B; \untw none; \typea P_1=\recip4(1,1,3); \excla
18B+3E; \typeb P_3P_4=1\times\recip9(1,4,5); \exclb ;
\typec P_1P_3=1\times\recip2(1,1,1); \exclc   5B+2E;
\end{entry}

\begin{entry}
\no 95; \reln 66; \gens 1,5,6,22,33; \degree 1/330;\inv
\textbf{F}_0; \ell 1; \kkk 2; \type P_1=\recip5(1,2,3); \excl
  6B; \untw none; \typea P_3P_4=1\times\recip{11}(1,5,6);
\excla ; \typeb P_2P_4=\recip3(1,2,1); \exclb
5B+E; \typec P_2P_3=1\times\recip2(1,1,1); \exclc   5B+2E;
\end{entry}

\end{center}


\end{document}